\documentclass[11pt, reqno]{amsbook}

\usepackage{amsmath}
\usepackage{amssymb}
\usepackage{mathrsfs}
\usepackage{bussproofs}

\usepackage[colorlinks]{hyperref}
\newcommand{\Long}  [2]                {#2}
\newcommand{\PaperState}               {\Long}
\newcommand{\IfPaperState}  [2]        {\PaperState{#1}{#2}}

\newcommand{\EM}                       { {\mathsf{EM}} }
\newcommand{\HA}                       { {\mathsf{HA}} }
\newcommand{\PA}                       { {\mathsf{PA}} }
\newcommand{\EA}                          {{\mathsf{EA}}}
\newcommand{\Nat}                      { {\tt N} }
\newcommand{\Bool}                     { {\tt Bool} }
\newcommand{\State}                    { {\tt S} }
\newcommand{\NatSet}                   {\mathbb{N}}
\newcommand{\BoolSet}                  {\mathbb{B}}
\newcommand{\StateSet}                 {\mathbb{S}}
\newcommand{\SystemT}                  {\mathcal{T}}
\newcommand{\SystemTG}                  {\mathsf{T}}
\newcommand{\SystemB}                {{\mathsf{B}}}
\newcommand{\True}                     { {\tt{True}} }
\newcommand{\False}                    { {\tt{False}} }
\newcommand{\Realizer}                 { {\mbox{\tiny Real}} }
\newcommand{\StructureN}               { {\mathcal N} }
\newcommand{\Class}                    {\mbox{\tiny Class}}
\newcommand{\Learn}                    {\mbox{\tiny Learn}}
\newcommand{\SystemTState}             {\SystemT_\State}
\newcommand{\SystemTClass}             {{\SystemT_{\Class}}}
\newcommand{\SystemTLearn}             {{\SystemT_{\Learn}}}

\newcommand{\Language}                 {\mathcal{L}}
\newcommand{\LanguageClass}            {\Language_{\Class}}
\newcommand{\LanguageLearn}            {\Language_{\Learn}}
\newcommand{\comment}[1]{}
\newcommand{\proj}                     { {p} }
\newcommand{\CupSem}                   { {\mathcal U} }
\newcommand{\Add}                      { {\mathsf{Add}} }
\newcommand{\add}                      { {\mathsf{add}} }

\newcommand{\makenum}        [1]       { {\overline{#1}} }
\newcommand{\makestate}      [1]       { {|{#1}|^{-1}} }
\newcommand{\PRclass}                    {{\mathcal{PCF}_{\Class}}}
\newcommand{\PRlearn}                    {{\mathcal{PCF}_{\Learn}}} 
\newcommand{\PCFclass}                    {{\mathcal{PCF}_{\Class}}}
 
\newcommand{\BR}                       {{\mathsf{BR}}}

\newcommand{\fun}                   {{\NatSet\rightarrow\NatSet}}
\newcommand{\seqfun}                   {{\Nat\rightarrow (\Nat\rightarrow\Nat)}}
\newcommand{\funnat}                   { (\Nat\rightarrow\Nat)\rightarrow \Nat}
\newcommand{\seqin}[1]                   { \Nat\rightarrow {#1}}
\newcommand{\funn}                   { {\Nat\rightarrow\Nat}}
\newcommand{\ifthen} [3]                        { {\mathsf{if}\ {#1}\ \mathsf{then}\ {#2}\ \mathsf{else}\ {#3} } }
\newcommand{\id}             {\mathsf{id}}
\newcommand{\gmc}                      {{\ \mathsf{gmc}_s\ }}
\newcommand{\M}                      {\mathsf{M}}
\newcommand{\update}                 {{\mathcal U}}
\newcommand{\rec}                          {{\mathsf{R}}}
\newcommand{\ifn}                          {{\mathsf{if}}}
\newcommand{\suc}                      {{\mathsf{S}}}

\newcommand{\model}[1]                  {{[\![#1]\!]_s}}
\newcommand{\Phic}                        {{\mathsf{\Phi}}}

\newcommand{\den}[1]                 {{|{#1}|}}
\newcommand{\Chi}                    {{\mathsf{X}}}
\newcommand{\E}[1]                   {{\mathsf{E}_{#1}}}
\newcommand{\oneback}                          {{\mathsf{1back}}}
\newcommand{\substitution} [1]         { {\overline{#1}} }
\newcommand{\zerof}                 { {\mathsf{0}}}
\newcommand{\rank}                    { {\mathsf{rk}} }

\newtheoremstyle{examplestyle}% name of the style to be used
  {6mm}% measure of space to leave above the theorem. topsep:=default
  {6mm}% measure of space to leave below the theorem. topsep:=default
  {}% name of font to use in the body of the theorem
  {}% measure of space to indent
  {\sc}% name of head font
  {.}% punctuation between head and body
  {.5em}% space after theorem head
  {\thmname{#1} \thmnumber{#2} \thmnote{(#3)}}% Manually specify head
  
  \newtheoremstyle{theorems}% name of the style to be used
  {6mm}% measure of space to leave above the theorem. topsep:=default
  {6mm}% measure of space to leave below the theorem. topsep:=default
  {\itshape}% name of font to use in the body of the theorem
  {}% measure of space to indent
  {\sc}% name of head font
  {.}% punctuation between head and body
  {.5em}% space after theorem head
  {\thmname{#1} \thmnumber{#2} \thmnote{(#3)}}% Manually specify head

\theoremstyle{theorems}
\newtheorem{theorem}{Theorem}[section]

\newtheorem{proposition}[theorem]{Proposition}

\theoremstyle{examplestyle}
\newtheorem{definition}{Definition}[section]
\newtheorem{remark}[theorem]{Remark}
\newtheorem{example}[theorem]{Example}

\theoremstyle{theorems}
\newtheorem{thm}{Theorem}[section]
\newtheorem{cor}[thm]{Corollary}
\newtheorem{lem}[thm]{Lemma}
\newtheorem{prop}[thm]{Proposition}

\theoremstyle{examplestyle}

\newtheorem{rem}[thm]{Remark}

\newtheorem{exa}[thm]{Example}

\newtheorem{defi}[thm]{Definition}

\newtheorem{fact}{Fact}

\numberwithin{section}{chapter}
\numberwithin{equation}{chapter}

\makeindex

\begin{document}

\frontmatter

\title{Universit\`a degli Studi di Torino\\
Dipartimento di Informatica\\
and\\
Queen Mary, University of London\\
School of Electronic Engineering and Computer Science \\
%(joint with: Queen Mary, University of London\\
%School of Electronic Engineering and Computer Science)
LEARNING, REALIZABILITY AND GAMES IN CLASSICAL ARITHMETIC}

%    Remove any unused author tags.

%    author one information
\author{Federico Aschieri\\ PhD Thesis}
\comment{\address{Dipartimento di Informatica\\ }
\address{Universit\`a di Torino}
\address{School of Electronic Engineering and Computer Science\\}
\address{Queen Mary, University of London}
\address{Supervisors: prof. Stefano Berardi and dr. Paulo Oliva}}

\address{Supervisors: prof. Stefano Berardi and dr. Paulo Oliva}
\address{Reviewers: prof. Ulrich Berger and prof. Pierre-Louis Curien}
\address{Defended on March, 25, 2011 (London) and April, 11, 2011 (Turin)}
%\subjclass[2000]{Primary }
%    For books to be published after 1 January 2010, you may use
%    the following version:


\date{30 November, 2010}

\begin{abstract}
In this dissertation we provide mathematical evidence that the concept of \emph{learning} can be used to give a new and intuitive computational semantics of classical proofs in various fragments of Predicative Arithmetic. 

 First, we extend Kreisel modified realizability to a classical fragment of first order Arithmetic, Heyting Arithmetic plus $\EM_1$ (Excluded middle axiom restricted to $\Sigma^0_1$ formulas). We introduce a new realizability semantics we call ``Interactive Learning-Based Realizability''. Our realizers are \emph{self-correcting} programs, which learn from their errors and evolve through time, thanks to their ability of perpetually questioning, testing and extending their knowledge. Remarkably, that capability is entirely due to classical principles when they are applied on top of intuitionistic logic.
 
Secondly, we extend the class of learning based realizers to a classical version $\PCFclass$ of $\mathcal{PCF}$ and, then,  compare the resulting notion of realizability with Coquand game semantics and prove a full soundness and completeness result. In particular, we show there is a one-to-one correspondence between realizers and recursive winning strategies in the 1-Backtracking version of Tarski games. 
 
Third,  we provide a complete and fully detailed constructive analysis of learning as it arises in learning based realizability for $\HA+\EM_1$, Avigad's update procedures and epsilon substitution method for Peano Arithmetic $\PA$. We present new constructive techniques to bound the length of learning processes and we apply them to reprove - by means of our theory - the classic result of G\"odel that provably total functions of $\PA$ can be represented in G\"odel's system $\SystemTG$.

Last, we give an axiomatization of the kind of learning that is needed to computationally interpret Predicative classical second order Arithmetic. Our work is an extension of Avigad's and generalizes the concept of update procedure to the transfinite case. Transfinite update procedures have to learn values of transfinite sequences of non computable functions in order to extract witnesses from classical proofs. 

 \end{abstract}

\maketitle

%   Dedication.  If the dedication is longer than a line or two,
%    remove the centering instructions and the line break.
\cleardoublepage
\thispagestyle{empty}
\vspace*{13.5pc}
\begin{center} A Margherita\\[2pt] 
.\\[2pt]
\comment{``Quand on s'est connus,\\[2pt]
Quand on s'est reconnus,\\[2pt]
Pourquoi se perdre de vue,\\[2pt]
Se reperdre de vue?\\[2pt]
Quand on s'est retrouv\'es,\\[2pt]
Quand on s'est r\'echauff\'es,\\[2pt]
Pourquoi se s\'eparer?\\[2pt]}
``Alors tous deux on est repartis\\[2pt]
Dans le tourbillon de la vie\\[2pt]
On \`a continu\'e \`a tourner\\[2pt]
Tous les deux enlac\'es\\[2pt]
Tous les deux enlac\'es\\[2pt]
Tous les deux enlac\'es."% Dedication text (use \\[2pt] for line break if necessary)
\end{center}
%\cleardoublepage

%    Change page number to 6 if a dedication is present.
\setcounter{page}{4}

\textbf{Acknowledgments}

This dissertation would simply not exist without the contribution of my supervisors, Stefano Berardi and Paulo Oliva. 

Stefano, I have to thank you so much. You have deeply influenced my work with your profound ideas on proof theory: you know, my accomplishment would not have been possible if I hadn't started from your great insights. You also  have invested an extraordinary amount of time in checking my work, teaching me how to write mathematics, supporting and helping me in many situations, well beyond your duties. 

Paulo, thanks for your support,  your brilliant guidance and suggestions: you have led me to exciting research in my year at Queen Mary.

I wish to thank also Gianluigi Bellin:  without you too this thesis would not have been produced. You have supported and  helped me greatly in the early stages of my studies. You have introduced me to my future topics of research and had a determinant influence on the direction of my career. Thanks also for our dinners in London.\\

\comment{Finally, I thank my family, for the constant care, support and love. }

Ringrazio anche la mia famiglia, per la presenza costante, per tutto l'affetto e il supporto che non mi hanno mai fatto mancare. 

Infine, grazie Margherita, per essere speciale come sei,  per avermi voluto bene sempre, perch\'e quando si tratta di te niente passa, niente scorre, niente cambia: questa tesi \`e dedicata a te, come il mio cuore lo \`e da quando ci siamo conosciuti.

\tableofcontents

%    Include unnumbered chapters (preface, acknowledgments, etc.) here.
\mainmatter
%    Include main chapters here.

\chapter{Introduction}\label{chapter-introduction}

\section{A Computational Semantics of Classical Proofs}
In this dissertation we provide mathematical evidence that the concept of \emph{learning} can be used to give a new and intuitive computational semantics of classical proofs in various fragments of Predicative Arithmetic. The main definite result in this sense is a new realizability semantics for $\HA+\EM_1$, which we call ``Interactive Learning Based Realizability" (shortly, \emph{learning based realizability}). $\HA+\EM_1$ is first order intuitionistic Heyting Arithmetic with the principle of excluded middle over $\Sigma^0_1$ formulas, a classical axiom that eluded computational semantics for a long time and, probably,  never had an intuitive one.

 $\EM_1$ has already been extensively studied in terms of Curry-Howard correspondence (see \cite{Sorensen}), that is, in terms of computational constructs that can be associated to it in order to extract computational information from classical proofs. Thanks to this approach, the constructive content of classical logic can be interestingly attained in terms of \emph{proof transformations}, or \emph{reduction rules} applied to the corresponding computational constructs. While these results are satisfying from the computational point of view, we feel they are not satisfying in terms of human \emph{understanding}. Though this issue may at first appear marginal, it is not: classical proofs can now be ``executed" thanks to Curry-Howard correspondence, but extracted programs are still very difficult to understand. 
 Nowadays, there are research programs whose aim is precisely to \emph{understand}  programs extracted from classical proofs, since without this understanding programs cannot be analyzed, neither optimized nor improved, to begin with. Without high level grasp of a program, it is like having a black box, that in some mysterious way gives always correct answers. This phenomenon probably happens because classical principles are well understood only in terms of the computational devices they are associated to. Therefore, when one has to interpret a classical proof, he is forced to formalize it and then to extract the corresponding program. While this approach is feasible with small proofs, it is unmanageable with complex ones without great technical effort. This seems to explain why computational interpretations of classical logic are still not universally used by mathematicians and computer scientists.  
 
As for ourselves, we think that proof theory should  offer a \emph{proof semantics}: that is, not only a way of extracting the computational content of classical proofs, but \emph{also} a high level explanation of what are the general ideas used by the programs extracted from proofs and what is the constructive \emph{meaning} of quantifiers, logical connectives and axioms in the framework of classical logic. Only in this way, the general mathematician or computer scientist could extract intuitively the computational content hidden in some proof he wishes to analyze. Such a semantics was put forward a long time ago for intuitionistic logic - starting from Heyting semantics, passing through Kleene realizability to arrive to Kreisel modified realizability - but we think that in the case of classical logic there are still no ultimate results.
 
 In this dissertation, we put a major effort to lay the ground for a \emph{proof semantics} of Predicative Arithmetic based on the concept of \emph{learning}. A learning based realizability interpretation of $\HA+\EM_1$ will be completely explained and formalized, compared with game semantics and thoroughly analyzed by constructive means. We will define a computational model of ``intelligent" \emph{self-correcting} programs, which learn from their errors and evolve through time, thanks to their ability of perpetually questioning, testing and extending their knowledge. Remarkably, that capability is entirely due to classical principles when they are applied on top of intuitionistic logic. We shall thus conclude that the computational content of classical fragment $\HA+\EM_1$ can be described in terms of learning. Moreover, we will introduce a more general concept of learning \emph{by levels}, that generalizes Avigad's one \cite{Avigad} and will serve as a foundation for possible extensions of our learning based realizability to full first order Peano Arithmetic and even Predicative Analysis. The learning based computation interpretation of classical logic  is a new and exciting field of research, but just started: in this thesis we give substantial contributions, but the path to follow is still long.\\
 
 In this chapter we give an overview of what are proof transformations, what are proofs semantics and what are the general ideas of learning in classical logic. We start with a short review of known results about intuitionistic logic, in order to explain how the issue ``proof transformation vs. proof semantics" was solved in that case. We then explain the ideas behind epsilon substitution method, Coquand's game semantics and Avigad's update procedures, the three major sources of inspiration for our work and where the concept of learning first appeared, implicitly in the first case, explicitly and much more elegantly in the second and the third. These pioneering contributions are also the underpinnings of the more recent work of Berardi and de' Liguoro \cite{BerardiLiguoro}, which will be the starting point of this thesis. We conclude with a synopsis of the contribution and the structure of our dissertation.

\section{Proof Transformations and Proof Semantics in Intuitionistic Arithmetic}

Intuitionistic logic was the main method of reasoning used in mathematics until the eighteenth century. Existence of objects with special properties was mainly established by explicitly constructing objects with those properties. With Dedekind, Cantor, Hilbert, Brouwer, Frege and many others, non constructive reasoning became central, thanks to its great conceptual and simplifying power. Non constructive methods were immediately questioned by mathematicians such as Kronecker and others, but soon became generally employed and accepted. After the famous days of paradoxes in set theory and mathematics, however, a new season of concerns and foundational efforts began in mathematical logic. Brouwer, in particular, advocated the need for \emph{intuitionistic logic}, which rejected some classical principles such as the excluded middle and impredicative definitions.

Intuitionistic logic was born as a \emph{constructive logic}, but what does that mean? Syntactically, it is presented just as a set of axioms and inference rules in a formal language. And the reason why, for example, excluded middle is left out remains inevitably obscure, if a semantics for assertions in intuitionistic logic is not provided; the excluded middle is a very natural and intuitive principle, after all, known in logic since Aristotle. 

The issue is informally solved through the so called \emph{Heyting semantics} (see for example Girard \cite{Girard}). According to Brouwer, in order to assert a statement, one has to provide some sort of \emph{construction}. What  is a construction? The following is Heyting's answer:
\begin{enumerate}
\item A construction of an atomic formula is a \emph{computation} showing its truth.\\
\item A construction of $A\land B$ is a pair formed by a construction of $A$ and a construction of $B$.\\
\item A construction of $A\lor B$ is a pair whose first element is a boolean $i$ such that: if $i=\True$, then the second element of the pair is a construction of $A$, if $i=\False$, then the second element of the pair is a construction of $B$.\\
\item A construction of $A\rightarrow B$ is an algorithm taking as input a construction of $A$ and returning as output a construction of $B$.\\
\item A construction of $\forall x^\Nat A(x)$ is an algorithm taking as input a natural number $n$ and returning as output a construction of $A(n)$.\\
\item A  construction of $\exists x^\Nat A(x)$ is a pair whose first element is a number $n$ such that the the second element of the pair is a construction of $A(n)$.\\ 

\end{enumerate}   

Thanks to Heyting semantics, one immediately recognizes why the excluded middle $A\lor \lnot A$ is rejected by intuitionistic logic: a construction of $A\lor \lnot A$ would require to decide which one among $A, \lnot A$ is true, which is not generally possible in algorithmic way. In Heyting semantics, each logical connective and quantifier is given a computational meaning, which is enough clear as to be intuitively used by humans with the aim of devising constructive proofs. 

\subsection{Proof Transformations in $\HA$}

Let us consider the formal system $\HA$, Heyting Arithmetic, which is the intuitionistic version of the usual first order Peano Arithmetic $\PA$. Since $\HA$ is born as a constructive logic, one expects it to be sound with respect to Heyting semantics. But since as a formal system $\HA$ does not mention at all the concept of Heyting construction, it is not obvious, given any provable formula, how to extract a construction from any of its proofs. 

The first method for implicitly extracting such a computational content is due to Gentzen \cite{Gentzen1}, \cite{Gentzen2} and falls under the category of \emph{proof transformations}. The idea is the following. One defines a set of transformation rules mapping proofs of a formula into proofs of the same formula. He then iterates these rules from the initial proof until he finds a proof with a special syntactical structure, and from this proof he obtains the constructive content ``implicit" in the initial proof. We will not consider Gentzen technique (called ``cut elimination") but instead a later one due to Prawitz, known as \emph{normalization}, which applies to natural deduction proofs. We start from recalling the inference rules of $\HA$.\\

\begin{enumerate}

\item %AXIOM
$\begin{array}{c}   \hline   A
\end{array}\ \ \ \ $ where $A$ is a Peano axiom.\\

\item %CONJUNCTION
$\begin{array}{c}  A\ \ \  B\\ \hline

A\wedge B
\end{array}\ \ \ \ $
$\begin{array}{c}   A\wedge B\\ \hline  A
\end{array}\ \ \ \ $
$\begin{array}{c}   A\wedge B\\ \hline  B
\end{array}$\\\\

\item %IMPLICATION
$\begin{array}{c}  A\rightarrow B\ \ \  A \\ \hline
 B
\end{array}\ \ \ \ $
$\begin{array}{c}  B\\ \hline 
A\rightarrow B
\end{array}$\\\\

\item %DISJUNCTION
$\begin{array}{c}   A\\ \hline 
 A\vee B
\end{array}\ \ \ \ $
$\begin{array}{c}  B\\ \hline 
 A\vee B
\end{array}$\\
$\begin{array}{c}  A\vee B\ \ \  C\ \ \ 
C\\ \hline  C
\end{array}$\\

\item %FOR ALL
$\begin{array}{c} \forall \alpha^\Nat A\\ \hline   A[t/\alpha]
\end{array} $
$\begin{array}{c}   A\\ \hline
\forall \alpha^\Nat A
\end{array}$\\
where $t$ is a term in the language of $\HA$\\
\item %EXISTS
$\begin{array}{c}  A[t/\alpha^\Nat]\\ \hline  
\exists
\alpha^\Nat. A
\end{array}$ \ \ \ \
$\begin{array}{c}  \exists \alpha^\Nat. A\ \ \  C\\
\hline
C
\end{array} $\\

\item %INDUCTION
$\begin{array}{c}   A(0)\ \ \  \forall \alpha^\Nat.
A(\alpha)\rightarrow A(\suc(\alpha))\\ \hline 
\forall
\alpha^\Nat A
\end{array}\ \ \ \ $\\

\end{enumerate}
A natural deduction proof is a tree whose nodes are formulas obtained from their children by inference rules. Leaves of a deduction tree are either hypotheses or discharged hypotheses and all the standard management of discharged ones is assumed here (see for example van Dalen \cite{van Dalen}). As usual, quantifier rules must be applied under the usual restrictions on the quantified variables.

 We are interested in a set of rules that allow to transform any proof of any formula $A$ in a special \emph{normal form} proof of $A$. We limit ourselves to the rules for implication and induction and we refer to \cite{van Dalen} for a complete treatment.
The proof transformations defined are the following. 
 In any deduction, a sub-deduction 
 \begin{prooftree}
\AxiomC{$\mathcal{D}_1$}
\noLine
\UnaryInfC{$B$}
\UnaryInfC{$A\rightarrow B$}

\AxiomC{$\mathcal{D}_2$}
\noLine
\UnaryInfC{$A$}
\BinaryInfC{$B$}
\end{prooftree}
is replaced by
\begin{prooftree}
\AxiomC{$\mathcal{D}_1[\mathcal{D}_2/A]$}
\noLine
\UnaryInfC{$B$}
\end{prooftree}
where $\mathcal{D}_1[\mathcal{D}_2/A]$ is the deduction obtained from $D_1$ by replacing all the hypotheses $A$ with the deduction $\mathcal{D}_2$ of $A$. Moreover, in any deduction, a sub-deduction
\begin{prooftree}
\AxiomC{$\mathcal{D}_1$}
\noLine
\UnaryInfC{$A(0)$}

\AxiomC{$\mathcal{D}_2$}
\noLine
\UnaryInfC{$\forall\alpha^\Nat. A(\alpha)\rightarrow A(\suc(\alpha))$}

\BinaryInfC{$\forall \alpha^\Nat A $}
\UnaryInfC{$A(n)$}
\end{prooftree}
where $n$ is a numeral, is replaced by a deduction $\mathcal{D}^n$ of $A(n)$, where by induction $\mathcal{D}^0$ is defined as
\begin{prooftree}
\AxiomC{$\mathcal{D}_1$}
\noLine
\UnaryInfC{$A(0)$}
\end{prooftree}
and  for any numeral $m$, $\mathcal{D}^{\suc(m)}$ is defined as

\begin{prooftree}

\AxiomC{$\mathcal{D}_2$}

\noLine
\UnaryInfC{$\forall\alpha^\Nat. A(\alpha)\rightarrow A(\suc(\alpha))$}
\UnaryInfC{$A(m)\rightarrow A(\suc(m))$}

\AxiomC{$\mathcal{D}^m$}
\noLine
\UnaryInfC{$A(m)$}
\BinaryInfC{$A(\suc(m))$}
\end{prooftree}

It is possible to show that any natural deduction proof without hypotheses of a closed formula $A$ can be transformed into a closed normal proof of $A$. Just by inspection and some straigthforward reasoning on the resulting proof shape, one can show that if $A=\exists x^\Nat B$, then the last inference rule must be of the form
\[\begin{array}{c}  B[t/\alpha^\Nat]\\ \hline  
\exists
\alpha^\Nat. B
\end{array}\] 
Hence one automatically finds a witness $t$ for $\exists x^\Nat B$. If $A=B_1\lor B_2$, then the last inference rule must be of the form
\[\begin{array}{c}   B_i\\ \hline 
 B_0\vee B_1
\end{array}\]
with $i\in{0,1}$: again a witness a required by Heyting semantics. The constructive information of the original proof of $A$ can thus be found by normalization. 

In general, by normalization, it is possible to show that every formula provable in $\HA$ has a construction, in the sense of Heyting. There is a remark to be done: without Heyting semantics, the process of normalization would be perfectly non intelligible. One would only see a sequence of transformations performed on a initial proof until a normal form proof is found and, magically, a witness pops out of nothing. Moreover, normalization is not a technique that can be used in an intuitive way, because it requires full formalization of  proofs and then to understand how the normalization process proceeds.
\subsection{Proof Semantics for $\HA$}

We now explain how a formal proof semantics for $\HA$ can be formulated in terms of Kreisel modified realizability \cite{Kreisel}. Modified realizability is a formalization of Heyting semantics, carefully carved and designed for $\HA^\omega$ (but we shall consider its restriction to $\HA$). The idea is to restrict the class of algorithms used in Heyting's definition: only algorithm representable in G\"odel's system $\SystemTG$ (see chapter \ref{chapter-technicalpreliminaries}) are allowed. This is an important restriction: only bounded iteration is explicitly used by such algorithms and this property rules out the possibility that the computational content of proofs may be found by blind search and other constructively unjustified techniques (as it might happen with trivial Kleene-style realizers \cite{Kleene}). 

We start by associating types to formulas as to mirror the structure of the programs used  in Heyting semantics.
 
 \begin{defi}
\label{definition-TypesForRealizers} (Types for realizers) For each
arithmetical formula $A$ we define a type $|A|$ of $\SystemTG$ by
induction on $A$:
\begin{enumerate}
\item
$|P(t_1,\ldots,t_n)|=\Nat$,\\
\item
$|A\wedge B|=|A|\times |B|$,\\
\item
$|A\vee B|= \Bool\times (|A|\times |B|)$,\\
\item
$|A\rightarrow B|=|A|\rightarrow |B|$,\\
\item
$|\forall x A|=\Nat\rightarrow |A|$,\\
\item
$|\exists x A|= \Nat\times |A|$\\
\end{enumerate}
\end{defi}

We now define the Kreisel modified realizability relation $t\Vdash C$. It is clear that a modified realizer represents a construction in the sense of Heyting.
\begin{defi}[Modified Realizability]
\label{lemma-IndexedRealizabilityAndRealizability}
Assume $t$ is a closed term of G\"odel's system $\SystemTG$ (see chapter \ref{chapter-technicalpreliminaries}), $C$ is a closed formula, and $t:|C|$. Let $\vec{t} = t_1, \ldots, t_n : \Nat$. We define the relation $t\Vdash C$ by induction and by cases according to the form of $C$:

\begin{enumerate}
\item
$t\Vdash P(\vec{t})$ if and only if 
$P(\vec{t})={\True}$\\

\item
$t\Vdash {A\wedge B}$ if and only if $\pi_0t \Vdash {A}$ and $\pi_1t\Vdash {B}$\\

\item
$t\Vdash {A\vee B}$  if and only if either $\proj_0t={\True}$  and $\proj_1t\Vdash A$, or $\proj_0t={\False}$  and $\proj_2t\Vdash B$\\

\item
$t\Vvdash {A\rightarrow B}$ if and only if for all $u$, if $u\Vdash {A}$,
then $tu\Vdash_s{B}$\\

\item
$t\Vdash {\forall x A}$ if and only if for all numerals $n$,
$t{n}\Vdash A[{n}/x]$\\
\item

$t\Vdash \exists x A$ if and only for some numeral $n$, $\pi_0t= {n}$  and $\pi_1t \Vdash A[{n}/x]$\\

\end{enumerate}

\end{defi}
 
Thanks to Kreisel modified realizability and to the correspondence between typed lambda terms and natural deduction proofs, one can give a new intuitive meaning to the normalization process for intuitionistic proofs. One does that by decorating inference rules with terms of system $\SystemTG$ in such a way that each of the natural deduction transformation rules we have previously seen correspond to a normalization step in the associated term. We shall present in chapter \ref{chapter-learningbasedrealizability} the full decoration, here we consider only the case for $\rightarrow$-rules:
\[
\begin{array}{c}  u\vdash A\rightarrow B\ \ \ t\vdash A \\ \hline
ut\vdash B
\end{array}\ \ \ \ 
\begin{array}{c}  u\vdash B\\ \hline \lambda x^{|A|}u\vdash
A\rightarrow B
\end{array}\]
We see that, finally, the idea of Heyting construction, in the form of modified realizability, is explicitly associated to the proof. Each inference rule either uses a realizer to compute something or defines a realizer. If we consider a decorated deduction $\mathcal{D}$
 \begin{prooftree}
\AxiomC{$\mathcal{D}_1$}
\noLine
\UnaryInfC{$u\vdash B$}
\UnaryInfC{$\lambda x^{|A|} u\vdash A\rightarrow B$}

\AxiomC{$\mathcal{D}_2$}
\noLine
\UnaryInfC{$t\vdash A$}
\BinaryInfC{$(\lambda x^{|A|} u)t\vdash B$}
\end{prooftree}
the previously described proof transformation rule of $\mathcal{D}$ in $\mathcal{D}_1[\mathcal{D}_2/A]$ corresponds to the normalization step $(\lambda x^{|A|} u) t=u[t/x]$, which is simply the evaluation of a function at its argument.

The point to be made here is the following. First, one has a \emph{local} interpretation of each inference rule in terms of realizability, that helps to understand what the rule itself means and what it is going to accomplish from the computational point of view. Moreover, every axiom and every formula in general is given a computational meaning through realizability. In this way, one understands what arithmetical assertions mean from the constructive standpoint and what one has to do in order to constructively prove them. Without such a semantics, it would not be possible for a human to understand in a simple way what is the intuitive computational content that a proof offers. The goal of this dissertation is to extend this result to $\HA+\EM_1$. 

\section{Learning in Classical Arithmetic: A Brief History}

In this dissertation we provide a realizability semantics of $\HA+\EM_1$ classical proofs in terms of learning. The first question we have to answer is therefore the following: \emph{what is to be learned by realizers of classical proofs}? Surprisingly, the answer was anticipated a long time ago by Hilbert and his epsilon substitution method.

\subsection{Learning in the Epsilon Substitution Method} 

Kreisel's modified realizers have not enough computational power for deciding truth of formulas and thus they are prevented  to realize classical principles such as the excluded middle; this limitation in turn prevents them to decorate classical proofs with the aim of finding witnesses for classically provable existential statements. However, as proven by Kreisel himself \cite{Kreisel0}, one can extract from a classical proof in $\PA$ of any $\Pi^0_2$ formula $\forall x^\Nat \exists y^\Nat Pxy$, with $P$ decidable, a non trivial algorithm that given any number $n$ finds a witness $m$ such that $Pnm$ is true. This result implies that even non constructive proofs of existence hide constructive information.   

Hilbert's idea (for a modern account, see \cite{Avigad})  to circumvent the apparent inability of constructive methods  to computationally interpret classical proofs, at least for $\Pi^0_2$ provable formulas, was to introduce non computable functions and to define through them classical witnesses. The role of those non computable functions, called epsilon substitutions, is to assign values to some non effectively evaluable \emph{epsilon terms}.  For each arithmetical formula $A(x)$ one introduces a term $\epsilon x A(x)$, whose intended denotation can be any number $n$ such that $A(n)$ is true. Thus,  $\epsilon x A(x)$ reads as ``an $x$ such that $A(x)$". For every term $t$, one has therefore an axiom
\[ A(t)\rightarrow A(\epsilon x A(x))\]
which captures completely the intended meaning of an epsilon term. Such an axiom is said to be a \emph{critical formula}. Epsilon terms are not effectively computable, so one starts with an epsilon substitution $S$ that gives to them dummy values. In order to extract the computational content of a classical proof, one has to satisfy all the critical formulas appearing in  the proof. Fortunately, there is only a finite number of them in any proof and so nothing in principle makes the goal impossible. It turns out that one can \emph{learn} values of epsilon terms by counterexamples. Suppose, for instance, that a critical formula 
\[ A(t)\rightarrow A(\epsilon x A(x))\]
is false under some substitution $S$. Then, if we denote with $S(\epsilon x A(x))$ the value associated to $\epsilon x A(x)$ by $S$, we have that $A(S(\epsilon x A(x)))$ is false. However, if $t$ evaluates to $n$ under $S$, $A(n)$ must be true! So, one updates the substitution $S$ as to give to $\epsilon x A(x)$ the value $n$, because he has \emph{learned} a witness through a counterexample. Things, however, are not so easy: if there are many critical formulas, this attempt of making true one of them, may make false another. Hilbert's \emph{Ansatz} (approach) was to show that a specially defined series of this learning steps must terminate.

Hilbert's idea is brilliant, but it has never been used to give a semantics of classical proofs. As to ourselves, we shall use it and we can now anticipate that the goal of our learning based realizers will be to \emph{learn values of non computable functions}. 

\subsection{Learning in Coquand Game Semantics}

The concept of learning in classical logic has been beautifully improved and reframed in terms of special Tarski games, in which the participating players are allowed to correct their moves: this is Coquand games semantics \cite{Coquand}. The interest, for our purposes, of this semantics, lies in the fact that the concept of recursive winning strategy in a standard Tarski game is nothing but a rephrasing of realizability: an adaptation of Tarski games to classical logic offers an opportunity to translate back the new intuitions in a realizability semantics, which is more suitable to proof interpretations. 

We first review what is a  \emph{Tarski game} (actually, the concept was explicitly defined by Hintikka (see \cite{Hintikka}) and it is a sort of folklore). In a Tarski game there are two players, Eloise and Abelard, and a negation-and-implication-free formula $B$ on the board (we assume that for every atomic formula its negation can also be expressed as an atomic formula). Intuitively, Eloise tries to show that the formula is true, while Abelard tries to show that it is false. Turns, sets of possible moves and winners are defined accordingly to the form of $B$:
\begin{enumerate}
\item If $B=\exists x^\Nat A(x)$, then Eloise has to choose a numeral $n$ and $B$ is replaced by the formula $A(n)$.\\
\item If $B=\forall x^\Nat A(x)$, then Abelard has to choose a numeral $n$ and $B$ is replaced by the formula $A(n)$.\\
\item If $B=A_1\lor A_2$, then Eloise has to choose a numeral $i\in\{0,1\}$ and $B$ is replaced by the formula $A_i$.\\
\item If $B=A_1\land A_2$, then Abelard has to choose a numeral $i\in\{0,1\}$ and $B$ is replaced by the formula $A_i$.\\
\item If $B$ is atomic and true, Eloise wins the game, otherwise Abelard wins.\\
\end{enumerate}

Informally, Eloise has a \emph{recursive winning strategy} in this Tarski game, if she has an algorithm for choosing her next moves that enables her to win every play. It is possible to show that a formula has a construction in the sense of Heyting if and only if Eloise has a recursive winning strategy in its associated Tarski game. In this sense, Tarski games  rephrase Heyting semantics. 

It is not surprising, then, that Eloise does not have a recursive winning strategy for every instance of the excluded middle $\EM_1$. Coquand solved this impasse by allowing Eloise to backtrack, i.e. to erase her moves and return to a previous state of the game (actually, when interpreting the cut rule, Coquand allowed also Abelard to backtrack, but we will not consider this case). Again, \emph{learning} by counterexamples enters the scene. It is possible to prove that Eloise has a recursive winning strategy in the backtracking version of the Tarski game associated to $\EM_1$. Suppose, for example, that 
\[\EM_1:=\forall x^\Nat.\ \exists y^\Nat Pxy \lor \forall y^\Nat \lnot Pxy\]   
and consider any play: we show how Eloise can win. Abelard has to move and, for some $n$, chooses the formula
\[\exists y^\Nat Pny \lor \forall y^\Nat \lnot Pny\] 
Then it is the turn of Eloise, who believes that no witness for $\exists y^\Nat Pny$ can be found. So she chooses 
\[ \forall y^\Nat \lnot Pny\] 
Again, Abelard for some $m$ chooses
\[\lnot Pnm\]
If $Pnm$ is false, then Eloise wins. But if $Pnm$ is true, she would have lost the standard Tarski game. However, according to the new rules, she can now backtrack to a previous position. Observe that  Abelard has falsified Eloise's belief in the non existence of witnesses for the formula $\exists y^\Nat Pny$, \emph{actually providing} one witness. So Eloise backtracks to the position 
\[\exists y^\Nat Pny \lor \forall y^\Nat \lnot Pny\] 
and this time chooses 
\[\exists y^\Nat Pny\]
followed by
\[Pnm\]
and so she wins. 

We remark how close is this kind of learning to the one in epsilon substitution method. In both cases, some formula wished to be true is false. But  from its falsehood one can always learn a new positive fact: a witness for an existential statement.

\subsection{Learning in Avigad's Update Procedures}

In \cite{Avigad}, Avigad has formulated an abstract axiomatization of learning as it implicitly appears in the epsilon substitution method for first order Peano Arithmetic. He has explicitly introduced the non computable functions needed by the epsilon method to elicit the computational content of classical proofs and formulated in a clear way the notion of \emph{update procedure} that formalizes what we call ``learning by levels".

We give a definition slightly different from Avigad's, but analogous. Intuitively, a $k$-ary update procedure, with $k\in \NatSet^+$, is a functional which takes as input a finite sequence $f=f_1, \ldots, f_k$ of functions approximating some oracles $\Phic_1, \ldots, \Phic_k$, such that each one of those functions is defined in terms of the previous ones. Then, it uses those functions to compute some witnesses for some provable $\Sigma_1^0$ formula of $\PA$. Afterwards, it checks whether the result of its computation is sound. If it is not, it identifies some wrong value $f_i(n)$ used in the computation and  corrects it with a new one.

\begin{definition}[Update Procedures]\label{definition-finiteupdateprocintro}

A \emph{$k$-ary  update procedure}, $k\in\NatSet^+$,  is a continuous function $\update: (\NatSet\rightarrow \NatSet)^k\rightarrow \NatSet^3\cup\{\emptyset\}$ (i.e., its output is determined by a finite number of values of the input functions) such that the following holds: 
\begin{enumerate}
\item for all function sequences ${f}={f}_1,\ldots, f_k$
\[\update {f}=(i,n,m) \implies 1\leq i\leq k\]
\item for all function sequences ${f}={f}_1,\ldots, {f}_k$ and ${g}={g}_1,\ldots, {g}_k$, for all $1\leq i<k$, if \\

i) for all $j<i$, ${f}_j={g}_j$; \\

ii) $\update {f}=(i,n,m)$, ${g}_i(n)=m$ and $\update {g}=(i,h,l)$\\\\
then $h\neq n$. 
\end{enumerate}
If $\update$ is a $k$-ary  update procedure, a \emph{zero} for $\update$ is a sequence ${f}={f}_1,\ldots, f_k$ of functions such that $\update f=\emptyset$.

\end{definition}

Condition (2) of definition \ref{definition-finiteupdateprocintro} means that the values of the $i$-th function depend on the values of some of the functions $f_j$, with $j<i$, and learning on level $i$ is possible only if all the lower levels $j$ have ``stabilized". 
In particular, if $\update$ is a $k$-ary update procedure and $f:(\NatSet\rightarrow \NatSet)^k$ is a sequence of functions approximating the oracles $\Phic_1, \ldots, \Phic_k$, there are two possibilities: either $f$ is a fine approximation and then $\update f=\emptyset$; or $f$ is not  and then $\update f=(i,n,m)$, for some numerals $n,m$: $\update$ says the function $f_i$ should be updated as to output $m$ on input $n$. Moreover, if $\update f = (i, n, m)$, one in a sense has \emph{learned} that $\Phic_i(n)=m$: by definition of update procedure, if $g$ is a function sequence agreeing with $f$ in its first $i-1$ elements, $g_i$ is another candidate approximation of $\Phic_i$ and $g_i(n)=m$, then $\update g$ does not represent a request to modify the value of $g_i$ at point $n$, for $\update g=(i, h,l)$ implies $h\neq n$.

The main theorem about update procedures is that they always have zeros and these latter can be computed through  learning processes \emph{guided} by the former. Intuitively a zero of an update procedure represents a good approximation of the oracles used in a computation, and in particular a good enough one to yield some sought classical witness.

\section{Learning in Classical Arithmetic: Contributions and Structure of This Dissertation}

The aim of this dissertation is to study and \emph{describe} the computational content of classical proofs in terms of \emph{learning}. In particular, the contributions and the structure of this dissertation are the following.

\subsection{Chapter \ref{chapter-learningbasedrealizability}} Our first contribution is to put together in a novel way the ideas contained in the epsilon substitution method, Coquand game semantics and Avigad's update procedures in order to define a realizability semantics of proofs which extends in a simple way Kreisel modified realizability to the classical system $\HA+\EM_1$: we shall call it (interactive) \emph{learning based realizability}. 

In a few words, learning based realizability describes a way of making oracle computations more effective, through the use of approximations of oracle values and learning of new values by counterexamples. A learning based realizer is in the first place a term of system $\SystemTClass$, which is a simple extension G\"odel's system $\SystemTG$ plus some oracles of the same Turing degree of an oracle for the Halting problem. Of course, if a realizer  was only this, it would be ineffective and hence useless. Therefore, learning based realizers are computed with respect to \emph{approximations} of the oracles of $\SystemTClass$ and thus effectiveness is recovered. Since approximations may be sometimes inadequate, results of computations may be wrong. But a learning based realizer is also a \emph{self-correcting} program, able to spot incorrect oracle values used during computations and to correct them with \emph{right} values. The new values are \emph{learned}, surprisingly, by  realizers of  $\EM_1$ and all the oracle values needed during each particular computation are acquired through learning processes. Here is the fundamental insight: classical principle may be computationally interpreted as learning devices.

 Our realizability semantics allows not only to extract realizers as usual by decorating classical proofs, but also to understand the intuitive meaning,  behaviour and goals of the extracted realizers. 
 
\subsection{Chapter \ref{chapter-realizabilityandgames}} Our second contribution is, first, to extend the class of learning based realizers to a classical version $\PCFclass$ of $\mathcal{PCF}$ and, then, to compare the resulting notion of realizability with Coquand game semantics and prove a full soundness and completeness result. In particular, we show there is a one-to-one correspondence between realizers and recursive winning strategies in the 1-Backtracking version of Tarski games.  

The soundness theorem should be useful to understand the significance and see possible uses of learning based realizability.  The idea is that playing games represents a way of challenging realizers and of  seeing how they react to the challenge by learning from failure and counterexamples. 

The proof of the completeness theorem in our view, moreover, has an interesting feature. In a sense, it is  the first application of the ideas of learning based realizability to a concrete non trivial classical proof, which is our version of the one given by Berardi et al. \cite{BerCoq}. That proof classically shows that if Eloise has \emph{recursive} winning strategy in the 1-Backtracking Tarski game associated to a formula $A$, then she also has a winning strategy in the Tarski game associated to $A$ (but a non computable strategy, only recursive in an oracle for the Halting problem). We manage to associate a constructive content to this seemingly ineffective proof and find out that it hides a learning mechanism to gain correct oracle values from failures and counterexamples. We then transform this learning mechanism into a learning based realizer.   

\subsection{Chapter \ref{chapter-constructiveanalysislearning}} Our third contribution is a complete and fully detailed constructive analysis of learning as it arises in learning based realizability for $\HA+\EM_1$, Avigad's update procedures and epsilon substitution method for Peano Arithmetic $\PA$. We present new constructive techniques to bound the length of learning processes and we apply them to reprove - by means of our theory - the classic result of Kreisel  that provably total functions of $\PA$ can be represented in G\"odel's system $\SystemTG$. An interesting novelty is that we develop type-theoretic techniques to reason and prove in a new way theorems about update procedures and epsilon substitution method. A notable byproduct of our work is the introduction of a ``constructive" non standard model of G\"odel's system $\SystemTG$. Our analysis is also a first step toward an extension of learning based realizability to full $\PA$. 

\subsection{Chapter \ref{chapter-learninganalysis}} Our last contribution is an axiomatization of the kind of learning that is needed to computationally interpret Predicative classical second order Arithmetic. Our work is an extension of Avigad's and generalizes the concept of update procedure to the transfinite case. Transfinite update procedures have to learn values of transfinite sequences of non computable functions in order to extract witnesses from classical proofs. We shall present several proofs of the fact that transfinite update procedures have zeros. The last one uses methods of type theory and in particular bar recursion: the algorithms presented are powerful and yet quite simple. The interest of our results is twofold. First, we extend Avigad's intuitive description of the learning content of the epsilon substitution method to the second order case. Secondly, we take a first step toward the extension of learning based realizability to Predicative second order Arithmetic, since we have isolated the concept of learning that will have to be employed.

\chapter{Technical Preliminaries}\label{chapter-technicalpreliminaries}

\section{Background}
\label{section-Background}

This dissertation is almost self-contained from the technical point of view. We cover the needed background here by reviewing what will be our main technical tool: G\"odel's System $\SystemTG$.

\subsection{G\"odel's system $\SystemTG$}
G\"odel's system $\SystemTG$ (see \cite{Girard}, for example) is simply typed $\lambda$-calculus, enriched with natural numbers, booleans, conditional $
 \ifn_T$ and primitive recursion $\rec_T$ in all types together with their associated reduction rules. We start by defining the types of system $\SystemTG$.

\begin{definition}[Types of System $\SystemTG$]
\label{definition-Types} The set of \emph{types} of system $\SystemTG$ is defined inductively as follows:

\begin{enumerate} 

\item $\Nat$ and $\Bool$ are types.\\

\item If $U, V$ are types, $U\times V$, $U\rightarrow V$ are types.

\end{enumerate}
\end{definition}

The type $\Nat$ represents the set $\NatSet$ of natural numbers and $\Bool$ represents the set $\BoolSet = \{\mbox{True},\mbox{False}\}$ of booleans, while product types $T \times U$ and arrows types $T \rightarrow U$ represent respectively cartesian products and function spaces. We assume $\rightarrow$ associates to the left:
\[T\rightarrow  U\rightarrow V=T\rightarrow  (U\rightarrow V)\]

\begin{definition}[Terms of System $\SystemTG$] 
\label{definition-Terms} We define the terms of system $\SystemTG$ inductively as follows:
\begin{enumerate}

\item For all types $U$, the variables $x_0^U,\ldots, x_n^U,\ldots$ are terms of type $U$.\\

\item $0$ is a term of type $\Nat$. $\True$ and $\False$ are terms of type $\Bool$. \\

\item For every type $T$, $\ifn_U$ is a term (constant) of type $U:=\Bool\rightarrow T \rightarrow T\rightarrow T$. Terms of the form $\ifn_Tt_1t_2t_3$ will be written in the more legible form $\ifthen{t_1}{t_2}{t_3}$.\\

\item For every type $T$, $\rec_U$ is a term (recursion constant) of type $U:=T\rightarrow (\Nat \rightarrow (T \rightarrow T))\rightarrow \Nat\rightarrow T$. The type $T$ in $\rec_T$ will be omitted whenever inferable from the context.\\

\item If $t$ is of type $\Nat$, then $\suc(t)$ is a term of type $\Nat$.\\

\item If $t$ and $u$ are terms of types respectively $U\rightarrow V$  and $U$, then $tu$ is a term of type $V$.\\

\item  If $t$ is a term of type $U\times V$, then  $\pi_0 t$ and $\pi_1 t$ are terms of types respectively $U$ and $V$.\\

\item If $u$ and $v$ are terms of types respectively $U$  and $V$, then $\langle u, v\rangle $ is a term of type $U\times V$.\\

\item  If $v$ is a term of type $V$ and $x^U$ a variable, then $\lambda x^U v$ is a term of type $U\rightarrow V$.

\end{enumerate}
\end{definition}

As usual, given two terms $u$ and $t$ and a variable $x$, with $u[t/x]$ we shall denote the term resulting from  $u$ by replacing all free occurrences of $x$ in $u$ with $t$, avoiding capture of variables.

We now give the reduction rules that explain the computational meaning of the syntax of system $\SystemTG$.

\begin{definition}[One Step $\beta_0$ Reduction and $\beta$ Reduction]
\label{definition-Reduction} We define a binary relation $\beta_0$ between terms of system $\SystemTG$ as the least relation satisfying the following properties:

\begin{enumerate}

\item If $u\;\beta_0\; u'$ and $v\;\beta_0\; v'$, then for all terms $u, v$, it holds that $uv\;\beta_0\; u'v$,  $uv\;\beta_0\; uv'$,  $\langle u, v\rangle\;\beta_0\; \langle u',v\rangle$ and $\langle u, v\rangle\;\beta_0\; \langle u,v'\rangle$.\\

\item If $u\;\beta_0\; u'$, then $\lambda x^T u\;\beta_0\;\lambda x^T u'$ and $\pi_i u\;\beta_0\;\pi_i u'$ for $i=0,1$.\\

\item $(\lambda x^T u)t\;\beta_0\; u[t/x]$.\\

\item  $\pi_i \langle u_0, u_1\rangle\;\beta_0\; u_i$ for $i=0,1$.\\

\item $\rec uv0\;\beta_0\; u$.\\

\item $Ruv\suc(t)\;\beta_0\; vt(\rec uvt)$.\\

\item $\ifthen{\True}{u_1}{u_2}\;\beta_0\; u_1$ and $\ifthen{\False}{u_1}{u_2}\;\beta_0\; u_2$.\\

\item $t\ \beta\ t'$ if $t=t'$ or $t\ \beta_0\ u_1\ \beta_0\ u_2\ \beta_0\ldots\beta_0\ u_n\ \beta_0\ t'$ for some terms $u_1,u_2\ldots,u_n$.\\

\item We say that $t$ is in normal form if $t\ \beta_0\ t'$ does not hold, for every $t'$.

\end{enumerate}
\end{definition}

We now define the equality rules for terms of system $\SystemTG$. Throughout  the dissertation we will write $u=t$ if the equation is provable by means of the following rules.  

\begin{definition}[Equational Theory of $\SystemT$]

\label{definition-EquationalTheory} We list the axioms of equality for system $\SystemTG$:
\begin{enumerate}

\item $t=t$.\\

\item If $t=u$, then $u=t$.\\

\item If $t=u$ and $u=v$, then $t=v$.\\

\item If $u=u'$ and $v=v'$, then $uv=u'v'$ and $\langle u, v\rangle=\langle u',v'\rangle$.\\

\item If $u=u'$, then $\lambda x^V u=\lambda x^V u'$ and $\pi_i u=\pi_i u'$ for $i=0,1$.\\

\item $(\lambda x^T u)t=u[t/x]$.\\

\item  $\pi_i \langle u_0, u_1\rangle=u_i$ for $i=0,1$.\\

\item $\rec uv0=u$.\\

\item $\rec uv\suc (t)=vt(\rec uvt)$.\\

\item $\ifthen{\True}{u_1}{u_2}= u_1$ and $\ifthen{\False}{u_1}{u_2}= u_2$.\\

\end{enumerate}
\end{definition}

Every term of G\"odel's system $\SystemTG$ has a unique normal form (see \cite{Girard}).

\begin{theorem}[Normalization and Church-Rosser Property]
\label{StrongNormChurch}
For every term $t$ of system $\SystemTG$ there exists $n\in\NatSet$ such that  $t\ \beta_0\ t_1\ \beta_0\ldots \beta_0\ t_m$ implies $m\leq n$. Hence, $t$ has a normal form.\\ Moreover, if $t\ \beta\ t_1$ and $t\ \beta\ t_2$, then there exists $t'$ such that $t_1\ \beta\ t'$ and $t_2\ \beta\ t'$. As consequence, $t$ has a unique normal form.
\end{theorem}

A term is \emph{closed} if it has no free variables; a \emph{numeral} is a term of the form $\suc^n(0)$, with $n\in\NatSet$, having inductively defined $\suc^0(0):=0$ and $\suc^{n+1}(0):=\suc(\suc^n(0))$. We will constantly use the following characterization of normal forms (see again \cite{Girard}). 

\begin{theorem}[Normal Form Characterization for System $\SystemTG$]\label{theorem-normalformpropertyT}
\label{lemma-normalform} Assume $A$ is an atomic type. Then any closed normal term $t$ of  $\SystemTG$ of type $A$ is either a numeral ${n}:\Nat$ or a boolean $\True,\False:\Bool$.
\end{theorem}

\chapter{Interactive Learning-Based Realizability
for Heyting Arithmetic with
$\EM_1$}\label{chapter-learningbasedrealizability}

\begin{abstract}
In this chapter, we extend Kreisel modified realizability to a classical fragment of first order Arithmetic, Heyting Arithmetic plus $\EM_1$ (Excluded middle axiom restricted to $\Sigma^0_1$ formulas). In particular, we introduce a new realizability semantics we call ``Interactive Learning-Based Realizability''. Our realizers are \emph{self-correcting} programs, which are able to learn something from every failure and use the new knowledge to correct themselves. We build our semantics over Avigad's fixed point result
\cite{Avigad},
but the same semantics may be defined over different constructive
interpretations of classical arithmetic (in \cite{Berardi},
continuations are used). Our notion of realizability extends
intuitionistic realizability and differs from it only in the atomic case: we interpret atomic realizers as ``learning agents''.
\end{abstract}

\section{Introduction}
From now on, we will call $\HA$ Heyting Intuitionistic Arithmetic, with a language including one symbol for each primitive recursive predicate or function (see \cite{van Dalen} or section \ref{section-ALearningBasedRealizability}). We call $\Sigma^0_1$-formulas the set of all formulas $\exists x.P(x,y)$ for some primitive recursive predicate $P$, and $\EM_1$ the {\em Excluded middle axiom restricted to $\Sigma^0_1$-formulas}. For a detailed study of the intuitionistic consequences of the sub-classical axiom $\EM_1$ we refer to \cite{BerCoqKohlLICS}.

This chapter is based on Aschieri and Berardi \cite{AB}. We extend Berardi and de' Liguoro (\cite{Ber2005}, \cite{Berardi}) notion of atomic realizability - originally conceived for quantifier free primitive recursive Arithmetic plus $\EM_1$ - to full predicate logic, namely Heyting Arithmetic with $\EM_1$ ($\HA + \EM_1$). Our idea is to interpret classical proofs as constructive proofs on a
suitable structure $\StructureN$ for natural numbers and maps of G\"odel's system $\SystemTG$, by applying to the semantics of Arithmetic the idea of ``finite approximation'' used to interpret Herbrand's Theorem. We extend intuitionistic realizability to a new notion of realizability, which we call ``Interactive learning-based Realizability". We provide a term assignment for the standard natural deduction system of $\HA + \EM_1$, which is surprisingly equal in all respects to that of $\HA$, but for the fact that we have non-trivial realizers for
atomic formulas and a new realizer for $\EM_1$.

Our semantics is ``local'': we do not introduce a global variable representing an ``external" goal, different for each particular proof one wants to interpret, as in continuation interpretation, in Friedman's $A$-translation and in Krivine's Classical Realizability. The goal of realizers is  fixed, ``internal", and is either to provide right constructions or to learn new information about excluded middle. In this way, we interpret $\EM_1$ and thus classical proofs locally and step-by-step, in order to solve a major problem of all computational interpretations: global illegibility, which means that, even for simple classical proofs, it is extremely difficult to understand how each step of the extracted program is related to the ideas of the proof, and what it is the particular task performed by each subprogram of the extracted program. The main sources of inspiration of this chapter are works of
Kleene, Hilbert, Coquand, Hayashi, Berardi and de' Liguoro and Avigad.

\textit{Intuitionistic Realizability revisited.} Recall chapter \ref{chapter-introduction}. In \cite{Kleene}, Kleene
introduced the notion of realizability, a formal semantics for
intuitionistic arithmetic. Later, Kreisel \cite{Kreisel} defined modified realizability, the same notion but with respect to G\"odel's system $\SystemTG$ instead of Kleene's formalism of partial recursive functions.  Realizability is nothing but a formal
version of Heyting semantics for intuitionistic logic, translated
into the language of arithmetic.

Intuitively, realizing a closed arithmetical formula $A$  means
exhibiting a computer program - called realizer - able to calculate
all  relevant information about the truth of $A$.  Hence, realizing
a formula $A\lor B$ means realizing $A$ or  realizing $B$, after
calculating \textit{which one} of the two is actually realized;
realizing a formula $\exists x A(x)$ means computing a numeral $n$ - called a witness -  and realizing $A(n)$.

These two cases are indeed the only ones in which we have relevant
information to calculate about the truth of the corresponding
formula, and there is a  decision to be made: realizing a formula
$\forall x A$ means exhibiting an algorithm which takes as input a numeral $n$ and gives as output realizers of $A(n)$;
realizing a formula $A\wedge B$ means realizing $A$ and realizing $B$; realizing $A\rightarrow B$ means providing an algorithm which takes as input realizers of $A$ and gives realizers of $B$; in these cases we provide no information about the formula we realize and we only take the inputs we will use for realizing existential or disjunctive formulas. Finally, realizing an atomic formula means that the formula is true: in this case, the realizer does nothing at all.

\IfPaperState{}{Hence, intuitionistic realizability closely follows
Tarski's definition of truth - the only difference being
effectiveness: for instance, while Tarski, to assert that $\exists
x
A$ is true, contented himself to know that there exists some $n$
such that  $A(n)$ is true, Kleene asked for a program that
calculates an $n$ such that $A(n)$ is
true.}

Intuitionistic natural deduction rules are perfectly suited to
preserve realizability. In order to actually build realizers from
intuitionistic natural deductions, it suffices to give realizers for
the axioms. Since our goal is to interpret classical connectives
using Heyting and Kleene interpretation of intuitionistic
connectives, then a first, quite naive idea would be the following:
if we devised realizers for Excluded Middle, we would be able to extend
realizability to all classical arithmetic.

Unfortunately, from the work of Turing it is
well known that not every instance of Excluded Middle is
realizable.
If $Txyz$ is Kleene's predicate, realizing $\forall x \forall y.
\exists z Txyz \lor \forall z \neg Txyz$ implies exhibiting an
algorithm which for every $n,m$ calculates whether or not the
$n$-th
Turing machine halts on input $m$: the halting problem would be
decidable. Hence, there is no hope of computing with effective
programs all the information about the truth  of Excluded
Middle.

However, not all is lost. A key observation is
the
following. Suppose we had a realizer
$O$ of the Excluded Middle and we made a natural deduction of a
formula $\exists x A$ actually using  Excluded Middle; then, we
would be able to extract from the proof a program $u$, containing
$O$ as subprogram, able to compute the witness for $\exists x A$.
Given the effectiveness of $u$, after a finite number of steps -
and
more importantly, after a finite number of calls to $O$ - $u$ would
yield the required witness. It is thus clear that $u$, to perform
the calculation, would use only a \textit{finite piece of
information about the Excluded Middle}. This fundamental fact gives us hope: maybe there is not always necessity of fully realizing Excluded Middle, since in finite computations only a finite amount of information is used. If we were able to gain that finite information during the computation, as it is the case in the proof of Herbrand's Theorem, we could adapt intuitionistic realizability to Classical Logic.\\

\textit{Herbrand's Theorem and the idea of ``finite approximation".} (A corollary of) Herbrand's Theorem says that if a universal first order theory $T$, in a suitable language supporting definition by cases, proves a statement $\exists x P(x)$, then one can extract from any proof a term $t$ and closed instances $A_1,\ldots, A_n$ of some universal formulas of $T$ such that $A_1\land \ldots \land A_n\rightarrow P(t)$ is a propositional tautology. So, even using classical logic, one can define witnesses. The problem is that the functions occurring in $t$ may not be computable, because the language of $T$ is allowed to contain arbitrary functions. However, given the finiteness of the information needed about any function used during any finite computation of $t$, in order to carry out actual calculations one would only have to find finite approximations of the non-computable functions involved, thus recovering effectiveness. We choose to follow this intuition:  we will add non-computable functions to our language for realizers and exploit the existence of these ideal objects in order to find concrete computational solutions.

This general idea dates back to Hilbert's $\epsilon$-substitution method (for a neat reformulation of the $\epsilon$-method see for example Avigad \cite{Avigad}). As noted by Ackermann \cite{Ackermann}, the $\epsilon$- substitution method may be used to compute witnesses of provable existential statements of first order Peano Arithmetic. The procedure is simple: introduce Skolem functions (equivalently, $\epsilon$-terms) and correspondent quantifier free Skolem axioms in order to reduce any axiom to a quantifier free form; take a $PA$-proof of a sentence $\exists x P(x)$ and translate it into a proof using as axioms only universal formulas; then apply Herbrand's theorem to the resulting proof, obtaining a quantifier free proof of $P(t)$, for some term $t$ of the extended language; finally, calculate a suitable finite approximation of the Skolem functions occurring in $t$ and calculate from $t$ an $n$ such that $P(n)$ holds.

However, while proofs in quantifier free style are very simple combinatorial objects, they lose the intuitive appeal, the general concepts, the structure of high level proofs. Hence, it may be an impossible task to understand extracted programs. Moreover we have a computational syntactic method but no semantics of proofs and logical operators based on the idea of ``finite approximation'', as the realizability interpretations are based on the idea of ``construction''. However, in the $\epsilon$-method, albeit only for quantifier free formulas, we see in action the method of intelligent \textit{learning}, driven by the Skolem axioms used in the proofs. One of the aims of this chapter is to extend this ``semantics of learning" from atomic propositions to individuals, maps, logical connectives and quantifiers of full natural deduction proofs. An important contribution comes from Coquand \cite{Coquand}.\\

\textit{Coquand's Game Semantics for Classical Arithmetic.} Computing all relevant information about the truth of a
given formula $A$ is not always possible. In \cite{Coquand} and in the context of game semantics, Coquand introduced a new key idea around this problem: the correspondence between backtracking (in game theory, retracting a move) and ``learning'', a refinement of the idea of ``finite approximation''. If we cannot
compute all the right information about the truth of a formula, maybe we could do this if we were allowed to make finitely many mistakes and to learn from them.

Suppose, for instance, we have the formula $\forall x. \exists y
Pxy\lor \forall y \neg Pxy$, but we have no algorithm which, for all numeral $n$ given as input, outputs false if $\forall y \neg Pny$ holds and outputs true if $\exists y Pny$ holds. Then we may describe a learning algorithm $r$ as follows.  Initially, for all $n$ given as input, $r$ outputs false. Intuitively, $r$ is initially persuaded - following the principle ``if I don't see, I do not believe" - that for all numeral $n$ there is no numeral $m$ such that $Pnm$ holds. Hence, when asked for his opinion about the formula $\exists y Pny\lor \forall y \neg Pny$, $r$ always says: $\exists y Pny$ is false. However, if someone - an opponent of $r$ - to show that $r$ is wrong, comes out with an $m$ such that $Pnm$ holds, $r$ realizes indeed to be mistaken, and stores the information ``$Pnm$ is true". Then, the next time being asked for an opinion about $\exists y Pny\lor \forall y \neg Pny$, $r$ will say: true. In other words, such $r$, after at most one ``mind changing'', would be able to learn the correct answer to any question of the form: ``which one among $\exists y Pny$, $\forall y \neg Pny$ does hold?". This is actually learning by counterexamples and is the key idea behind
Coquand's semantics.

Our question is now: can we formulate a realizability notion based on learning by counterexamples in order to extend Kreisel's interpretation to all individuals, maps and connectives of the sub-classical Arithmetic $\HA +  \EM_1$? Following Hayashi \cite{Hayashi}, in our solution we modify the notion of individual, in such a way that individuals change with time, and realizers ``interact'' with them.\\

\textit{Hayashi's Proof Animation and Realizability.} In \cite{Hayashi}, Hayashi explains a notion of realizability for a sub-classical arithmetic, called limit computable mathematics.
Basing his analysis on ideas of Gold \cite{Gold}, he defines a Kleene's style
notion of realizability equal to the original one but for the fact
that the notion of individual changes: the
witnesses of existential and disjunctive formulas are calculated by
a stream of guesses and ``learned in the limit'' (in the sense that
the limit of the stream is a correct witness). An individual $a$
is therefore a computable map $a:\mathbb{N} \rightarrow \mathbb{N}$,
with $a(t)$ representing the value of the individual at time $t$.
%The technical device that makes his interpretation working is the
%class of limiting recursive functions introduced in Gold
%\cite{Gold}.

For instance, how would Hayashi realize the formula $\forall x. \exists y Pxy\lor \forall y \neg Pxy$? He would define an algorithm $H$ as follows. Given any numeral $n$, $H$ would calculate the truth value of $\forall y\leq n Pny$. Then the correct answer to the question: ``which one among $\exists y Pny$, $\forall y \neg Pny$ does hold?" is learned in the limit by computing $P(n,0)$, $P(n,1)$, $P(n,2)$,\ldots, $P(n,k)$,\ldots and thus  producing a stream of guesses either of the form false, false, false,\ldots, true, true,\ldots, true,\ldots or of the form false, false, false, \ldots, false, \ldots, the first stabilizing in the limit to true, the second to false. Hayashi's idea is to perform a completely
blind
and exhaustive search: in such a way, the correct answer is
guaranteed to be eventually learned (classically). Hayashi's realizers do not
learn in an efficient way: in Hayashi's notion of realizability the
only learning device is to look through all possible cases.
Instead,
we want to combine the idea of individual as limit, taken from Hayashi, with notion of learning in which the stream of guesses is
{\em driven by the proof itself}, as in Coquand's game semantics. For the quantifier-free fragment, this was done by Berardi \cite{Ber2005} and Berardi-de' Liguoro \cite{Berardi}.\\

\textit{Realizability Based on Learning: Berardi-de'
Liguoro interpretation.}
%As long as one investigates learning and
%the process of correcting hypotheses by means of counterexamples, it is natural to
We explain \cite{Berardi} using Popper's ideas \cite{Popper} as a metaphor.
According to Popper, a scientific theory relies on a set of
unproved - and unprovable - hypotheses and, through logic, makes
predictions susceptible of being falsified by experiments. If a prediction is falsified,
some hypothesis is
incorrect. In front of a counterexample to a theory's prediction,
one must modify the set of hypotheses and  build a better theory,
which will be tested by experiments, and so on. Laws of Nature are
universal statements, that cannot be verified, but are suitable to
falsification.
%\IfPaperState{}
%%%%%%%%
%% LONG VERSION
%%%%%%%%
%{Schematically: \[\mbox{HYPOTHESES + LOGIC }\longmapsto
%\mbox{PREDICTIONS}\] \[\mbox{WRONG PREDICTIONS}\longmapsto
%\mbox{WRONG HYPOTHESES}\] because we assume soundness of Logic.\\
%}
%%%%%%%%
%% END LONG VERSION
%%%%%%%%
We may explain the link between falsifiable hypotheses and
$\EM_1$. For every $n$, given an instance $\exists y.Pny  \vee
\forall
y.\neg Pny$  of $\EM_1$ (with $P$ atomic), we may formulate an
hypothesis about which side of the disjunction is true. If we know
that $Pnm$ is true for some $m$, we know that $\exists y.Pny$ is
true. Otherwise we may assume $\forall y. \neg Pny$ as hypothesis,
because it is a falsifiable hypothesis.

%Here \EM_1 comes into the scene. If $Pxy$ is atomic, for every
%$n$,
%we know that  $\forall y \neg Pny$ is true or false; sure, we
%cannot mechanically decide which one of the two cases holds, but
%we
%may formulate hypotheses.

We formalize the process of making hypotheses about
$\EM_1$ by a finite state of knowledge, called $S$,
collecting the instances $Pnm$ which we know to hold, e.g. by direct calculation.  If we have evidence that $Pnm$ holds for some $m$ (that is, $Pnm\in S$) we know that $\exists y Pny$ is true; in the other case, we assume that $\forall y \neg Pny$ is true. So $S$ defines a set of hypotheses on $\EM_1$, of the form $\forall y \neg Pny$: universal falsifiable statements. Using $S$ a realizer $r$ may effectively decide which side of a given instance of $\EM_1$ is true,
at the price of making mistakes: to decide if $\forall y \neg Pny$ is true, $r$ looks for any $Pnm$ in the finite state $S$ and outputs ``false'' if the research is successful, ``true'' otherwise. If and when from an hypothesis $\forall y \neg Pny$ we obtain some false conclusion $\neg Pnm$, the realizer $r$ returns the additional knowledge: ``$Pnm$ is true'', to be added to $S$.\\

\textit{Extending Berardi-de'
Liguoro interpretation to $\HA+\EM_1$.} In our chapter, we interpret each classical proof $p$ of $A$ in $\HA + \EM_1$ by a ``learning realizer'' $r$. $r$ returns a ``prediction'' of the truth of this formula, based on the information in $S$, and some additional knowledge in the case the prediction is effectively falsified.  For example, in front of
a formula $\exists x. A \wedge B$, a realizer $r$ predicts that $A(n)\wedge B(n)$ is true for some numeral $n$ (and since $n$ depends on $S$, in our model we change the notion of individual, interpreting ``numbers'' as  computable maps from the set of bases of knowledge to $\NatSet$).
Then $r$ predicts, say, that $B(n)$ is true, and so on, until $r$ arrives at some atomic formula, say $\neg Pnm$. Either $Pnm$ is actually true, or $r$ is able to effectively find one or more flawed hypothesis $\forall x.\neg Q_1n_1x, \ldots, \forall x.\neg Q_kn_kx$ among the hypotheses used to predict that $Pnm$ is true, and for each flawed hypothesis one counterexample $ Q_1n_1m_1,\ldots, Q_kn_km_k$. In this case, $r$ requires to enlarge our state of knowledge $S$ by including the information ``$Q_1n_1m_1$ is true", \ldots, ``$Q_kn_km_k$ is true".

Our Interactive Realizability differs from Intuitionistic
Realizability in the notion of individual (the value of an individual may depend on our knowledge state), and in the
realizability relation for the atomic case. In our interpretation, to realize an atomic formula does not mean that the formula is true, but that the realizer requires to extend our state of knowledge $S$ if the formula is not true. The realizer is thought as a learning device. Each extension of $S$ may change the value of the individuals which are parameters of the atomic formula, and therefore may make the atomic formula false again. Then the realizer requires to extend $S$ again, and so forth. The convergence of this ``interaction'' between a realizer and a group of individuals follows by Avigad's fixed point thm. \cite{Avigad} (a constructive proof may be found in \cite{Ber2005}), and it is the analogue of the termination of Hilbert's $\epsilon$-substitution method.

%The idea of using finite bases of knowledge and to use them to ``decide'' Excluded Middle comes from Berardi and de' Liguoro \cite{Ber2005}, \cite{Berardi}.

%In \cite{Berardi} there are only
%realizers for atomic formulas, which have the task of extending the
%current knowledge and hence are not trivial: they embody learning
%strategies.
%Actually, the case of atomic formulas is the crucial
%step of our interpretation.\\

\textit{Why the Arithmetic $\HA + \EM_1$} instead of considering the full Peano Arithmetic? We have two main reasons. First, we observe that $\EM_1$ enjoys a very good property: the information about its truth can be computed in the limit, in the sense of Gold \cite{Gold}, as we saw en passant when discussing Hayashi's realizability. This implies that witnesses for existential and disjunctive statements too can be  learned in the limit, as shown in Hayashi \cite{Hayashi}. In  chapter \ref{chapter-realizabilityandgames} we show that realizers which we will be able to extract from proofs have a straightforward interpretation as winning strategies in 1-Backtracking games \cite{BerCoq}, which are the most natural and simple instances of Coquand's style games. Secondly, a great deal of mathematical theorems are proved by using $\EM_1$ alone (\cite{BerCoqKohlLICS}, \cite{BerAPAL}). Third, as shown in chapter \ref{chapter-constructiveanalysislearning}, already $\HA+\EM_1$ - plus G\"odel's double negation translation  - suffices to interpret all provably total recursive functions of $\PA$, with the advantage of eliminating the extra step of Friedman or Dialectica translation (see Kohlenbach for \cite{Kohl} these latter)

%\IfPaperState{}{For
%instance, a lot of theorems can be derived as consequences of the
%minimum principle for functions over $\NatSet$. Remarkably,
%Hilbert's famous basis theorem can be proved with $\EM_1$, and
%interpreted with a learning in the limit program, as remarked by
%Hayashi by reading Hilbert's original notebooks.}

\textit{Plan of the Chapter}. The chapter is organized as follows. In \S \ref{section-TheTermCalculus} we define the term calculus in which our realizers will be written: a version of G\"odel's system $\SystemTG$, extended with some syntactic sugar, in order to represent bases of knowledge (which we shall call states) and to manipulate them. Then we prove a convergence property for this calculus (as in Avigad \cite{Avigad} or in \cite{Ber2005}). In \S \ref{section-ALearningBasedRealizability}, we introduce the notion of realizability and prove our Main Theorem, the Adequacy Theorem: ``if a closed arithmetical formula is provable in $\HA + \EM_1$, then it is realizable''.
In \S \ref{section-conclusion} we conclude the discussion about our notion of realizability by comparing it with other notions of realizability for classical logic, then we consider some possible future work.

%In \S \ref{section-Examples} we conclude the chapter with examples of program extraction, which are the motivation of this chapter: we show how our semantics yields legible, intuitive, effective programs from classical proofs.

\section{The Term Calculus}\label{section-TheTermCalculus}

In this section we formalize the intuition of ``learning realizer'' we discussed in the introduction.

We associate to any instance $\exists y Pxy\lor \forall y \neg Pxy$ of $\EM_1$ (Excluded Middle restricted to $\Sigma^0_1$-formulas) two functions $\chi_P$ and $\varphi_P$. The function $\chi_P$ takes a knowledge state $S$, a numeral $n$, and  returns a guess for the truth value of $\exists y. Pny$. When this guess is ``true'' the function $\varphi_P$ returns a witness $m$ of $\exists y.Pny$. The guess for the truth value of $\exists y. Pny$ is computed w.r.t. the knowledge state $ S $, and it may be wrong. For each constant $s$ denoting some knowledge state $ S $, the function $\lambda x^\Nat\ \chi_P( s ,x)$ is some ``approximation'' of an ideal map $\lambda x^\Nat\ \Chi_P(x)$, the \emph{oracle} returning the truth value of $\exists y. Pxy$. In the same way, the function $\lambda x^\Nat\ \phi_P( s ,x)$ is some ``approximation'' of an ideal map $\lambda x^\Nat\ \Phic_P(x)$, the \emph{Skolem map} for $\exists y. Pxy$, returning some $y$ such that $Pxy$ if any, and $0$ otherwise. The Skolem axioms effectively used by a given proof take the place of a set of experiments testing the correctness of the predictions made by $\varphi_P(s,x), \chi_P(s,x)$ about $\Chi_P(x),\Phic_P(x)$ (we do not check the correctness of $\varphi_P, \chi_P$ in an exhaustive way, but only on the values required by the Skolem axioms used by a proof).

Our Term Calculus is an extension of G\"odel's system $\SystemTG$ (see chapter \ref{chapter-technicalpreliminaries} or \cite{Girard}).  From now on, if $t, u$ are terms of $\SystemTG$ with $t=u$ we denote provable equality in $\SystemTG$. If $k \in \NatSet$, the numeral denoting $k$ is the closed normal term $\makenum{k} = \suc^k(0)$ of type $\Nat$.  All closed normal terms of type $\Nat$ are numerals (see chapter \ref{chapter-technicalpreliminaries}).  We recall that any closed normal term of type $\Bool$ in $\SystemTG$ is ${\True}$ or ${\False}$.

We introduce a notation for ternary projections: if $T = A \times (B \times C)$, with $p_0, p_1, p_2$ we respectively denote the terms $\pi_0$, $\lambda x:T.\pi_0(\pi_1(x))$, $\lambda x:T.\pi_1(\pi_1(x))$.
If $u = \langle u_0,\langle u_1,u_2\rangle \rangle : T$, then $\proj_iu=u_i$ in $\SystemTG$ for $i=0,1,2$. We abbreviate $\langle u_0,\langle u_1,u_2\rangle \rangle :T $ with $\langle u_0,u_1,u_2\rangle : T$. We formalize the idea of ``finite information about $\EM_1$'' by the notion of {\em state of knowledge}.

\begin{defi}[States of
Knowledge and Consistent Union]\label{definition-StateOfKnowledge}We define: \begin{enumerate}
\item
A $k$-ary {\em predicate} of $\SystemTG$ is any closed normal term $P:\Nat^{k}\rightarrow \Bool$ of $\SystemTG$.\\

\item
An atom is any triple $\langle P,\vec{n},{m}\rangle $, where $P$ is a $(k+1)$-ary predicate, $\vec{n},m$ are $k+1$ numerals, and $P\vec{n}m = \True$ in $\SystemTG$.\\

\item
Two atoms $\langle P,\vec{n},{m}\rangle $, $\langle P',\vec{n}',{m'}\rangle $ are {\em consistent} if $P = P'$ and $\vec{n} = \vec{n}'$ imply $m = m'$.\\

\item
A state of knowledge, shortly a {\em state}, is any finite set $S$ of pairwise consistent atoms.\\

\item
Two states $S_1, S_2$ are consistent if $S_1 \cup S_2$ is a state.\\

\item
$\StateSet$ is the set of all states of knowledge.\\
\item
The {\em consistent union} $S_1 \CupSem S_2$ of $S_1, S_2 \in \StateSet$ is $S_1 \cup S_2 \in \StateSet$ minus all atoms of $S_2$ which are inconsistent with some atom of $S_1$.\\
\end{enumerate}
\end{defi}

We think of an atom $\langle P,\vec{n},{m}\rangle $ as the code of a witness for $\exists y.P(\vec{n},y)$. Consistency condition allows at most one witness for each $\exists y.P(\vec{n},y)$ in each knowledge state $S$. Two states $S_1, S_2$ are consistent if and only if each atom of $S_1$ is consistent with each atom of $S_2$.

$S_1 \CupSem S_2$ is an non-commutative operation: whenever an atom of $S_1$ and an atom of $S_2$ are inconsistent, we arbitrarily keep the atom of $S_1$ and we reject the atom of $S_2$, therefore for some $S_1, S_2$ we have $S_1 \CupSem S_2 \not = S_2 \CupSem S_1$. $\CupSem$ is a ``learning strategy'', a way of selecting a consistent subset of $S_1 \cup S_2$. It is immediate to show that $\CupSem$ is an associative operation on the set of consistent states, with neutral element $\emptyset$, with upper bound $S_1 \cup S_2$, and returning a non-empty state whenever $S_1 \cup S_2$ is non-empty.

\begin{lem} \label{lemma-Cup}
Assume $i \in \NatSet$ and $S_1, \ldots, S_i \in \StateSet$.
\begin{enumerate}
\item
$S_1 \CupSem \ldots \CupSem S_i \subseteq S_1 \cup \ldots \cup S_i$

\item
$S_1 \CupSem \ldots \CupSem S_i = \emptyset$ implies $S_1 = \ldots = S_i = \emptyset$.
\end{enumerate}
\end{lem}

%\proof
%\begin{enumerate}
%\item
%By definition of $\CupSem$.
%\item
%It is enough to prove the result for $i=2$. Assume $S_1 \CupSem S_2 = \emptyset$. By definition of $\CupSem$ we have $S_1 \subseteq S_1 \CupSem S_2 = \emptyset$, therefore $S_1 = \emptyset$. Thus, no atom is inconsistent with some atom of $S_1$. We conclude that $\emptyset = S_1 \CupSem S_2 = S_1 \cup (S_2 \setminus \mbox{\{all atoms of $S_2$ inconsistent with some atom of $S_1$\}}) = \emptyset \cup (S_2 \setminus\emptyset) = S_2$.
%\end{enumerate}
%\qed

In fact, the whole realizability Semantics is a Monad \cite{BerardiLiguoroMonadi}. In \cite{BerardiLiguoroMonadi}, it is proved that our realizability Semantics is parametric with respect to the definition we choose for $\CupSem$. Any associative operation $\CupSem$, with neutral element $\emptyset$ and satisfying the two properties of Lemma \ref{lemma-Cup}, defines a different but sound realizability Semantics, corresponding to a different ``learning strategy''. An immediate consequence of Lemma \ref{lemma-Cup} is:

\begin{lem} \label{lemma-consistency}
Assume $S, S_1, S_2 \in \StateSet$.
\begin{enumerate}
\item
If $S$ is consistent with $S_1, S_2$, then $S$ is consistent with $S_1 \CupSem S_2$.
\item
If $S$ is disjoint with $S_1, S_2$, then $S$ is disjoint with $S_1 \CupSem S_2$.
\end{enumerate}
\end{lem}

For each state of knowledge $S$ we assume having a unique constant $s$ denoting it. We denote the state denoted by a constant $s$ with $\den{s}$ and as usual with $|\_|^{-1}$ the inverse of $|\_|$; that is, $||s||^{-1}=s$. We assume $\varnothing$ is the state constant denoting the empty state $\emptyset$; that is, $\den{\varnothing}=\emptyset$. We define with \[\SystemTState  = \SystemTG + \State + \{s\ |\ \den{s} \in \StateSet\}\] the extension of $\SystemTG$ with one atomic type $\State$ denoting $\StateSet$, and a constant $s$  for each state $S \in \StateSet$, and {\em no} new reduction rule. We denote states by $S, S', \ldots$ and state constants by $s, s', \ldots$. Any closed normal form of type $\Nat, \Bool, \State$ in $\SystemTState$ is, respectively, some numeral $n$, some boolean $\True, \False$, some state constant $s$. Computation on states will be defined by some suitable set of algebraic reduction rules we call ``functional''.

\begin{defi}[Functional set of rules]\label{definition-functional0}
Let $C$ be any set of constants, each one of some type $A_1\rightarrow \ldots \rightarrow A_n\rightarrow A$, for some $A_1,\ldots,A_n, A \in\{ \Bool, \Nat, \State\}$. We say that $\mathcal{R}$ is a {\em functional set of reduction rules} for $C$ if $\mathcal{R}$ consists, for all $c\in C$ and all ${a_1}:A_1,\ldots, {a_n}:A_n$ closed normal terms of $\SystemT_\State$, of exactly one rule $c {a_1}\ldots {a_n}\mapsto {a}$, for some closed normal term ${a}:A$ of $\SystemT_\State$.
\end{defi}

\begin{thm} \label{theorem-ExtStrongNormalization}
Assume that $\mathcal{R}$ is a functional set of reduction rules for $C$ (def. \ref{definition-functional0}). Then $\SystemTState + C + \mathcal{R}$ enjoys: i) strong normalization; ii) weak-Church-Rosser (uniqueness of normal forms) for all closed terms of atomic types.

\end{thm}

\proof \emph{(Sketch)} For strong normalization, see Berger \cite{Berger} (the constants ${s}:\State$ and $c \in C$ are trivially strongly computable). For weak Church-Rosser property, we start from the fact that there is the canonical set-theoretical model $\mathcal{M}$ of $\SystemTState + C + \mathcal{R}$. The interpretation of $\Bool, \Nat,\State$ in $\mathcal{M}$ consists of all closed normal form of these types. Arrows and pairs are interpreted set-theoretically. Each constant $c \in C$ is interpreted by some map $f_c$, defined by $f_c(a_1,\ldots,a_n)=a$ for all reduction rules $(c a_1\ldots a_n\mapsto a) \in \mathcal{R}$. Assume $u,v:A$ are closed normal term, $A = \Bool,\Nat$, or $\State$ is an atomic type, and $u,v$ are equal in $\SystemTState + C + \mathcal{R}$, in order to prove that $u, v$ are the same term. $u,v$ are equal in $\mathcal{M}$ because $\mathcal{M}$ is a model of $\SystemTState + C + \mathcal{R}$. By induction on $w$ we prove that if $w$ is a closed normal form of atomic type $\SystemT + C + \mathcal{R}$, then $w$ is a numeral, or $\True,\False$, or a state constant, and therefore $w$ is interpreted by itself in $\mathcal{M}$. From $u,v$ equal in $\mathcal{M}$ we conclude that $u,v$ are the same term of $\SystemTState + C + \mathcal{R}$.\\
\qed

We define two extensions of $\SystemT_\State$: an extension $\SystemTClass$ with symbols denoting the non-computable maps $\Chi_P, \Phic_P$ and no computable reduction rules, another extension $\SystemTLearn$, with the computable approximations $\chi_P,\phi_P$ of $\Chi_P, \Phic_P$, and a computable set of reduction rules. We use the elements of $\SystemTClass$ to represent non-computable realizers, and the elements of $\SystemTLearn$ to represent a computable ``approximation'' of a realizer. In the next definition, we denote terms of type $\State$ by $\rho, \rho', \ldots$.

\begin{defi} \label{definition-TermLanguageL1}
Assume $P:\Nat^{k+1}\rightarrow \Bool$ is a $k+1$-ary predicate of $\SystemTG$. We introduce the following constants:
\begin{enumerate}
\item
$\Chi_P:\Nat^k\rightarrow \Bool$ and $\Phic_P: \Nat^k \rightarrow \Nat$.\\

\item
$\chi_P:\State \rightarrow \Nat^k\rightarrow \Bool$
and
$\varphi_P:\State \rightarrow \Nat^k\rightarrow \Nat$.\\

\item
$\Cup:\State\rightarrow \State \rightarrow \State$.\\

\item
$\Add_P:\Nat^{k+1} \rightarrow \State$ and $\add_P:\State \rightarrow \Nat^{k+1} \rightarrow \State$.\\

\end{enumerate}

%
%With the expression $\Add_P$ we denote the term \[\lambda \vec{n},m\ {\tt if} \ P\vec{n}m\ {\tt then}\ ({\tt if}\ \Chi_P\vec{n}\ {\tt then}\ \emptyset\ {\tt else}\ \{\langle P, \vec{n}, m\rangle\})\ {\tt else}\ \emptyset\]
%and, if $s$ is a state, with $\add_P{s}$ we denote the term
%
%\[\lambda \vec{n},m\ {\tt if} \ P\vec{n}m\ {\tt then}\ ({\tt if}\ \chi_P{s}\vec{n}\ {\tt then}\ \emptyset\ {\tt else}\ \{\langle P, \vec{n}, m\rangle\})\ {\tt else}\ \emptyset\]

We denote $\Cup\rho_1\rho_2$ with $\rho_1\Cup\rho_2$. \\\begin{enumerate}
\item
$\Xi_\State$ is the set of all constants $\chi_P,\varphi_P, \Cup, \add_P$.\\

\item
$\Xi$ is the set of all constants $\Chi_P,\Phic_P, \Cup, \Add_P$.\\

\item
$\SystemTClass = \SystemT_\State + \Xi$.\\

\item
A term $t \in \SystemTClass$ has state $\varnothing$ if it has no state constant different from $\varnothing$.\\
\end{enumerate}
\end{defi}

Let $\vec{t} = t_1\ldots t_k$. We interpret $\chi_P{s} \vec{t}$ and $\varphi_P{s}\vec{t} $ respectively as a ``guess'' for the values of the oracle and the Skolem map $\Chi_P$ and $\Phic_P$ for $\exists y.P\vec{t}y$, guess computed w.r.t. the knowledge state denoted by the constant $s$.  There is no set of computable reduction rules for the constants $\Phic_P, \Chi_P \in \Xi$, and therefore no set of computable reduction rules for $\SystemTClass$.

 If $s_1, s_2$ are state constants, we  interpret $s_1 \Cup s_2$ as denoting the consistent union $\den{s_1} \CupSem \den{s_2}$. $\Add_P$ denotes the map constantly equal to the empty state $\emptyset$. $\add_Ps \vec{n}m $ denotes the empty state $\emptyset$ if we cannot add the atom $\langle P, \vec{n},m\rangle$ to $\den{s}$, either because $\langle P,\vec{n},l\rangle \in \den{s}$ for some numeral $l$, or because $P\vec{n}m={\False}$; $\add_P{s} \vec{n}m $ denotes the state $\{\langle P, \vec{n},m \rangle\}$ otherwise. We define a system $\SystemTLearn$ with reduction rules over $\Xi_\State$ by a functional reduction set $\mathcal{R}_\State$.

\begin{defi}[The System $\SystemTLearn$] \label{definition-EquationalTheoryL1}
Let  $s, s_1, s_2$ be state constants. Let $\langle P, \vec{n},m \rangle$ be an atom. $\mathcal{R}_\State$ is the following functional set of reduction rules for $\Xi_\State$:

\[\begin{aligned}\chi_P{s}\vec{n} &\mapsto
\begin{cases}\True &\mathsf{if}\ \exists m.\ \langle P,\vec{n},{m}\rangle \in \den{s}\\ 
\False &\mathsf{otherwise}
\end{cases}\\\\
\varphi_P{s}\vec{n} &\mapsto
\begin{cases}m &\mathsf{if}\ \exists m.\ \langle P,\vec{n},{m}\rangle \in \den{s}\\ 
0 &\mathsf{otherwise}
\end{cases}\\\\
\add_P{s}\vec{n}{m} &\mapsto
\begin{cases} \varnothing &\mathsf{if}\ \exists l.\ \langle P,\vec{n},{l}\rangle \in \den{s} \lor P\vec{n}m=\False\\
\makestate{\{\langle P,\vec{n},{m} \rangle\}} &\mathsf{otherwise}
\end{cases}\\\\
{s_1}\Cup{s_2} &\mapsto s_3, \text{ where $s_3$ is the state constant such that $\den{s_3}=\den{s_1} \CupSem \den{s_2}$}
\end{aligned}
\]
We define $\SystemTLearn = \SystemT_\State + \Xi_\State + \mathcal{R}_\State$.
\end{defi}

\textbf{Remark.} $\SystemTLearn$ is nothing but $\SystemTState$ with some ``syntactic sugar''. By Theorem \ref{theorem-ExtStrongNormalization}, $\SystemTLearn$ is strongly normalizing and has the weak Church-Rosser property for closed term of atomic types.  $\SystemTLearn$ satisfies a Normal Form Property.

\begin{lem}[Normal Form Property for $\SystemTLearn$] \label{lemma-normalform} Assume $A$ is either an atomic type or a product type. Then any closed normal term $t \in \SystemTLearn$ of type $A$ is: a numeral ${n}:\Nat$, or a boolean $\True,\False:\Bool$, or a state constant $s:\State$, or a pair $\langle u,v \rangle: B \times C$.
\end{lem}

\proof \emph{(Sketch)} By induction over $t$. For some $\vec{v}$, either $t$ is $(\lambda \vec{x}.u)(\vec{v})$, or $t$ is $\langle u,w\rangle(\vec{v})$, or $t$ is $x(\vec{v})$ for some variable $x$, or $t$ is $c(\vec{v})$ for some constant $c$, and either $c={0}, \mbox{S}, {\True}, {\False}, {s}, \linebreak R_T, {\tt if}_T, \pi_i$ is some constant of $\SystemTState$, or $c \in \Xi_\State$. If $t=(\lambda \vec{x}.u)(\vec{v})$, then $t$ has an arrow type if $\vec{v} = \emptyset$, while $t$ is not normal if $\vec{v} \not = \emptyset$. If $t=\langle u,w\rangle(\vec{v})$, then $\vec{v}=\emptyset$ and we are done.
If $t = x(\vec{v})$ then $t$ is not closed. The only case left is $t=c(\vec{u}):A$. $A$ is not an arrow type, therefore all arguments of $c$ are in $\vec{u}$. If $t=0$ we are done, if $t=\mbox{S}(u)$ we apply the induction hypothesis, if $t={\True},{\False}:\Bool$ or $t={s}: \State$ or $t = \langle u,v \rangle$ we are done. Otherwise either $t=R_T(n,f,a)\vec{t}, {\tt if}_T(b,a_1,a_2)\vec{t}, \pi_i(v)\vec{t}$, or $t = \chi_P(u,\vec{w}):\Nat$, or $t= \varphi_P(u,\vec{w}) :\Nat$, or $t=\Cup(u_1,u_2) :\State$, or $t=\add_P(u,\vec{w}): \State$. The proper subterms $n, w_1, \ldots, w_k:\Nat$, $b:\Bool$, $v:A \times B$, $u,u_1,u_2:\State$ of $t$ have atomic or product type and are closed normal. By induction hypothesis they are, respectively, a numeral, a boolean, a pair, a state constant. In all cases, $t$ is not normal.

 \qed

Let $t_1 t_2 \in \SystemTLearn$ be two closed terms of type $\State$. We abbreviate ``the normal forms of $t_1, t_2$ denote two states which are consistent and disjoint'' by: $t_1, t_2$ are consistent and disjoint. $\varnothing, s$ are consistent and disjoint for every state constant $s$. The maps denoted by $\Cup, \add_P$ preserve the relation ``to be consistent and disjoint''.

\begin{lem}\label{lemma-disjoint}
Assume $s, s_1, s_2$ are state constants and $\langle P, \vec{n},m\rangle$ is an atom.
\begin{enumerate}
\item
$s, (\add_Ps\vec{n}m)$ are consistent and disjoint.\\
\item
Assume $s,s_1$ are consistent and disjoint, and $s,s_2$ are consistent and disjoint. Then $s, s_1 \Cup s_2$ are consistent and disjoint.
\end{enumerate}
\end{lem}

\proof
\begin{enumerate}
\item
 If $\add_Ps\vec{n}m$ denotes the empty state the thesis is immediate. Otherwise $\add_Ps\vec{n}m$ denotes $\{\langle P, \vec{n},m\rangle\}$ and $\langle P, \vec{n},l\rangle \not \in \den{s}$ for all numerals $l$. Then $\{\langle P, \vec{n},m\rangle\}$ is consistent and disjoint with $\den{s}$.
\item
By Lemma \ref{lemma-consistency}.
\end{enumerate}
\qed

Each (in general, non-computable) term $t \in \SystemTClass$ is associated to a set $\{t[{s}]\ | $s$ \mbox{ is a } \linebreak \mbox{state constant}\} \subseteq \SystemTLearn$ of computable terms we call its ``approximations'', one for each state constant $s$.

\begin{defi}[Approximation at state $s$] Assume $t \in \SystemTClass$ and  $s$ is a state constant. We call ``approximation of $t$ at state $s$'' the term $t[{s}]$ of $\SystemTLearn$ obtained from $t$ by replacing each constant $\Chi_P$ with $\chi_P{s}$, each constant $\Phic_P$ with $\varphi_P{s}$, each constant $\Add_P$ with $\add_P{s}$.
\end{defi}

We interpret any $t[{s}] \in \SystemTLearn$ as a learning process evaluated w.r.t. the information taken from a state constant $s$ (the same $s$ for the whole term).
\\

Assume $t \in \SystemTClass$ is closed, $t:\State$ and $s$ is a state constant. Then $t[{s}]$ is a closed term of $\SystemTLearn$, and its normal form, by the Normal Form Property \ref{lemma-normalform}, is some state constant ${s}'$. We conclude $t[{s}]={s}'$ in $\SystemTLearn$. We prove that $s, s'$ are consistent and disjoint.

\begin{lem}\label{lemma-consistentdisjoint}
Assume $s$ is a state constant, $t \in \SystemTClass$,
$t:\State$ is closed, and all state constants in $t$ are consistent and disjoint with $s$.
\begin{enumerate}
\item
If $t[s]$ reduces to $t'[s]$, then all state constants in $t'$ are consistent and disjoint with $s$.\\
\item
$s,t[s]$ are consistent and disjoint.\\
\item
If $u\in \SystemTClass$, $u:\State$ and all state constants in $u$ are $\varnothing$, then $s, u[s]$ are consistent and disjoint.
\end{enumerate}
\end{lem}
%We define a notion of validity over $t \in \SystemTClass$. A closed term $t$ of atomic type is valid if $t[s]$ reduces to a numeral, $\True, \False$, or a state constant $s'$ such that $s, s'$ denote two states which are disjoint and consistent. A closed term of product type is valid if its two projections are, a closed term of arrow type is valid if it maps closed valid terms into closed valid terms. A term is valid if all its closed substitution with closed valid terms are valid. By induction over $t \in \SystemTClass$ we may prove that $t$ is valid. The only non-trivial cases are $t = \Add_P \vec{n} m$ or $t = t_1 \Cup t_2$. In this case $t[s] = \add_Ps \vec{n} m$ or $t[s] = t_1[s] \Cup t_2[s]$. By induction hypothesis $s, t_1[s]$ and $s, t_2[s]$ are disjoint and consistent. We apply Lemma \ref{lemma-disjoint}.

\proof
\begin{enumerate}
\item
It is enough to consider a one-step reduction. Suppose that $t[s]$ reduces to $t'[s]$ by contraction of a redex $r$ of $t[s]$. If $r$ is $(\lambda x u)t$ or $\rec_Tuv\suc(w)$ or $ {\tt if}_T(b,a_1,a_2)$ or $\pi_i\langle v_1,v_2\rangle$ or $ \chi_Ps\vec{n}$, or $ \varphi_Ps\vec{n}$, then its contractum $r'$ does not contain any new state constant; hence, all state constants in $t'$ are consistent and disjoint with $s$. If $r$ is $s_1\Cup s_2$ or $\add_Ps\vec{n}m$, then both $s, s_1$ and $s,s_2$ are consistent and disjoint state constants by hypothesis on $t$; therefore, by Lemma \ref{lemma-disjoint}, in both cases $s$ and the contraction of $r$ are consistent and disjoint;  so all state constants in $t'$ are consistent and disjoint with $s$.\\
\item
Every reduct of $t[s]$ is $t'[s]$ for some $t' \in \SystemTClass$. If $t[s]$ reduces to a normal form $t'[s] \equiv s'$, then the only possibility is $t' \equiv s'$. By the previous point $1$, we conclude that $s'$ is consistent and disjoint with $s$.\\
\item
By the previous point $2$, and the fact that the only state constant $\varnothing$ in $u$ is consistent and disjoint with any $s$.
\end{enumerate}
\qed

We introduce now a notion of convergence for families of terms  $\{t[{s_i}]\}_{i \in \NatSet} \subseteq \SystemTLearn$, defined by some $t \in \SystemTClass$ and indexed over a set of state constants $\{s_i\}_{i\in\NatSet}$. Informally, ``$t$ convergent" means that the normal form of $t[{s}]$ eventually stops changing when the knowledge state $s$ increases. If $s_1, s_2$ are state constants, we write $s_1 \le s_2$ for $\den{s_1}\subseteq \den{s_2}$. We say that a sequence $\{s_i\}_{i\in\NatSet}$ of state constants is a weakly increasing chain of states (is w.i. for short), if $s_i\le s_{i+1}$ for all $i\in\NatSet$.

\begin{defi}[Convergence]
\label{definition-Convergence} Assume
that $\{s_i\}_{i\in\NatSet} $ is a w.i. sequence of state constants,
and $u \in \SystemTClass$.
\begin{enumerate}

\item
  $u$ converges in $\{s_i\}_{i\in\NatSet}$ if $\exists i\in\NatSet.
\forall j\geq i.u[s_j]=u[s_{i}]$ in $\SystemTLearn$.\\

\item
$u$ converges if $u$ converges in every w.i. sequence of state constants.
\end{enumerate}
\end{defi}

We remark that if  $u$ is convergent, we do not ask that $u$
is convergent to the {\em same} value on {\em all} w.i. chain of states. The value learned by $u$ may depend on the information contained in the particular chain of state constants by which $u$ gets the knowledge. The chain of states, in turn, is selected by the particular definition we use for the ``learning strategy'' $\CupSem$. Different ``learning strategies'' may learn different values.

\begin{thm}[Stability Theorem] \label{theorem-StabilityTheorem0}
Assume $t \in \SystemTClass$ is a closed term of atomic type $A$ ($A\in\{\Bool,\Nat,\State\}$). Then $t$ is convergent.
\end{thm}

\proof \emph{(Classical)}. Assume $S$ is any consistent and possibly infinite set of atoms. We define some (in general, \emph{not} computable) functional reduction set $\mathcal{R}(S)$ for the set $\Xi$ of constants and for $\SystemTClass$. The reductions for $\Chi_P, \Phic_P, \Add_P$ are those for $\chi_P, \phi_P, \add_P$ in $\SystemTLearn$:

\begin{enumerate}
\item
If $\langle P,\vec{n},{m}\rangle \in S$, then
$(\Chi_P\vec{n}\mapsto {\True}), (\Phic_P\vec{n}\mapsto{m}) \in \mathcal{R}(S)$, else \linebreak
$(\Chi_P\vec{n}\mapsto {\False}), (\Phic_P\vec{n}\mapsto {0}) \in \mathcal{R}(S)$.\\

%\item
%${s_1}\Cup{s_2} \mapsto {s}_3$, where $s_3$ is defined as in  denotes the state $\{\langle P_1, \vec{n}_1, {m_1}\rangle, \ldots, \langle P_j,\vec{n}_j,{m_j}\rangle\}$ such that $\langle P_i,\vec{n}_i,{m_i}\rangle \in s_3$ if and only if either $\langle P_i,\vec{n}_i,{m_i}\rangle \in s_1$ or  $\langle P_i,\vec{n}_i,{m_i}\rangle \in s_2$ and $\langle P_i,\vec{n}_i,{m}'\rangle \notin s_1$ for every numeral $m'$.

\item
$\Add_P\vec{n}{m} \mapsto \varnothing\in \mathcal{R}(S)$ if either $\langle P,\vec{n},l \rangle \in S$ for some numeral $l$ or $P\vec{n}{m} = {\False}$; otherwise,  $\Add_P\vec{n}{m} \mapsto \makestate{\{\langle P,\vec{n},{m} \rangle\}}\in \mathcal{R}(S)$.
\end{enumerate}
and the reduction for $\Cup$ in $\mathcal{R}(S)$ is the reduction for $\Cup$ in $\SystemTLearn$. By theorem \ref{theorem-ExtStrongNormalization}, $\SystemTClass + \mathcal{R}(S)$ is strongly normalizing and weak-CR for all closed terms of atomic type, for any consistent set of atoms $S$. For the rest of the proof, let $\{s_i\}_{i\in\NatSet}$ be a w.i. chain of state constants. Assume $t \in \SystemTClass$ is a closed term of atomic type $A$.\\ {\em Claim}. For any state constant $s$, the map $u \mapsto u[{s}]$ is a bijection from the reduction tree of $t$ in $\SystemTClass + \mathcal{R}(\den{s})$ to the reduction tree of $t[{s}]$ in $\SystemTLearn$.\\ {\em Proof of the Claim.} By induction over the reduction tree of $t[{s}]$. Every reduction $\beta,\pi,{\tt if}_T, R_T, \Cup$ over $t[{s}]$ may be obtained from the same reduction over $t$. All occurrences of $\chi_P, \varphi_P, \add_P$ in the reduction tree of $t[{s}]$ are of the form  $\chi_P{s}, \varphi_P{s}, \add_P{s}$, therefore every reduction over $\chi_P, \varphi_P, \add_P$ may be obtained from the corresponding reduction over $\Chi_P,\Phic_P,\Add_P$.

Assume now $a$ is the (unique, by weak-CR) normal form of $t$ in $\SystemTClass + \mathcal{R}(\den{s})$. By the \emph{Claim}, $a[{s}]$ is the normal form of $t[{s}]$ in $\SystemTLearn$. Since $a$ is normal in $\SystemTClass + \mathcal{R}(\den{s})$, there is no $\Chi_P,\Phic_P,\Add_P$ in $a$. Thus $a$ and $a[{s}]$ are the same term: $t$ and $t[{s}]$ have the same normal form respectively in $\SystemTClass + \mathcal{R}(\den{s})$ and in $\SystemTLearn$. Let $\{s_i\}_{i \in \NatSet}$ be a given sequence of state constants. Define $S_\omega = \cup_{i \in \NatSet}\den{s_i}$. By strong normalization, the reduction tree of $t$ in $\SystemTClass + \mathcal{R}(S_\omega)$ is finite. Therefore in this reduction tree are used only \textit{finitely many} reduction rules from $\mathcal{R}(S_\omega)$, and for some numeral $n$ it is equal to the reduction tree of $t$ in $\SystemTClass + \mathcal{R}(\den{s_n})$, and in $\SystemTClass + \mathcal{R}(\den{s_m})$ for all $m \ge n$. We deduce that for all $m \ge n$ the normal forms of $t$ in $\SystemTClass + \mathcal{R}(\den{s_m})$ are the same. Thus, the normal form in $\SystemTLearn$ of all $t[s_m]$ with $m \ge n$ are the same, as we wished to show.
\qed

%%%%%%%%%%%%%%%%%

\begin{rem} The idea of the proof of theorem \ref{theorem-StabilityTheorem0} corresponds exactly to the intuition of the introduction.  During any computation, the oracles $\Chi_P$ and $\Phic_P$ are consulted a finite number of times and hence asked for a finite number of values. When our state of knowledge is great enough, we can substitute the oracles with their approximation $\chi_P{s}$ and $\varphi_P{s}$ for some state constant $s$, and we will obtain the same oracle values and hence the same results.

The proof, though non constructive, is short and explains well why the result is true. However, provided we replace the notion of convergence used in this chapter with the intuitionistic notion introduced in \cite{Ber2005},
we are able to reformulate and prove theorem \ref{theorem-StabilityTheorem0} in a purely intuitionistic way, achieving thus a constructive description of learning in $\HA +\EM_1$. Being the intuitionistic proof much more elaborated and less intuitive than the present one and connected with other foundationally interesting results, it will be the subject of chapter \ref{chapter-constructiveanalysislearning}.

Our proof of convergence follows the pattern
of Avigad's one in \cite{Avigad}.
A closed term $t\in \SystemTClass$ of atomic
type and in the constants $c_1,\ldots, c_n\in \Xi$,
may be seen as a functional $F_t$ which
maps functions $f_1,\ldots, f_n$ of the same type
of $c_1,\ldots c_n$  into an object of atomic type:
$F_t(f_1,\ldots, f_n)$ is defined as the normal
form of $t$ in $\SystemTClass + \mathcal{R}$,
where $\mathcal{R}
=\{ c_ia_1\ldots a_n\mapsto a\ |\ f_i(a_1,\ldots,a_n)=
a \mbox{ and $i\in\{1,\ldots, n\}$}\}$. $F_t$
is continuous in the sense of Avigad. Moreover, since $\Chi_P$ and $\Add_P$ have a set-theoretical definition in terms of $\Phic_P$, we may assume $F_t$ depends only on the functions which define in $\mathcal{R}$ the reduction rules for $\Phic_{P_1},\ldots\Phic_{P_n}$. Then, if $t$ is
of type $\State$, it is not difficult to see that
$F_t$ represents an update procedure with respect to any of its argument.
The fact that $F_t$ is an update procedure implies convergence
for $t$ and the zero theorem
\ref{FixedPointProperty}.
\end{rem}

As last result of this section, we prove that if we start from any state constant $s$ and we repeatedly apply any closed term $t: \State$ of $\SystemTClass$ of state $\varnothing$ (see definition \ref{definition-TermLanguageL1}), we obtain a ``zero" of $t$, that is a state constant $s_n$ such that $t[s_n]=\varnothing$. We interpret this by saying that any term $t$ represents a terminating learning process.

\begin{thm}[Zero Theorem]\label{FixedPointProperty}
Let $t:\State$ be a closed term of $\SystemTClass$ of state $\varnothing$ and $s$ any state constant. Define, by induction on $n$, a sequence $\{s_n\}_{n\in\NatSet}$ of state constants such that: $s_0=s$  and $s_{n+1}=s_n\Cup t[s_n]$. Then, there exists an $n$ such that $t[s_n]=\varnothing$.
\end{thm}

\proof ${s_0}, {s_1}, {s_2}, \ldots$ is a weakly increasing chain of state constants by construction. By theorem \ref{theorem-StabilityTheorem0}, $t$ converges over this chain: there exists $k\in\NatSet$ such that for every $j\geq k$, $t[s_j]=t[s_k]$. By choice of $k$
\[\begin{aligned}s_{k+2}&=s_{k+1} \Cup t[s_{k+1}]\\
&=(s_k\Cup t[s_k])\Cup t[s_{k+1}]\\
&=(s_k\Cup t[s_k])\Cup t[s_k]\\
&=s_k\Cup t[s_k]\\
&=s_{k+1}
\end{aligned}\]
Since $s_{k+2}=s_{k+1}$ and $s_{k+1}$, $t[s_{k+1}]$ are consistent and disjoint by lemma \ref{lemma-consistentdisjoint}, we conclude $t[s_{k+1}]=\varnothing$.

\qed

\section{An Interactive Learning-Based Notion of
Realizability}\label{section-ALearningBasedRealizability}
In this section we introduce the notion of realizability for $\HA +
\EM_1$, Heyting Arithmetic plus Excluded Middle on
$\Sigma^0_1$-formulas, then we prove our main Theorem, the Adequacy
Theorem: {\em ``if a closed arithmetical formula is provable in
$\HA
+ \EM_1$, then it is realizable''}. \IfPaperState{For proofs we
refer to
\cite{ExtendedVersion}.}{}

We first define the formal system $\HA + \EM_1$, from now on ``Extended $\EM_1$ Arithmetic''. We represent atomic predicates of $\HA + \EM_1$ with (in general, non-computable) closed terms of $\SystemTClass$ of type $\Bool$. Terms of $\HA + \EM_1$ may include function symbols $\Chi_P$, $\Phic_P$ denoting non-computable functions: oracles and Skolem maps for $\Sigma^0_1$-formulas $\exists x.Px\vec{{n}}$, with $P$ predicate of $\SystemTG$. We remark that our realizability can be formulated already for the standard language of Arithmetic: we add non computable functions to the language for greater generality.  We assume having in $\SystemTG$ some terms $\Rightarrow_\Bool: \Bool,\Bool\rightarrow\Bool, \neg_\Bool: \Bool \rightarrow \Bool, \ldots$, implementing boolean connectives. If $t_1, \ldots, t_n, t \in \SystemTG$ have type $\Bool$ and are made from free variables all of type $\Bool$, using boolean connectives, we say that $t$ is a tautological consequence of $t_1, \ldots, t_n$ in $\SystemTG$ (a tautology if $n=0$) if all boolean assignments making $t_1, \ldots, t_n$ equal to ${\True}$ in $\SystemTG$ also make $t$ equal to ${\True}$ in
$\SystemTG$.

\begin{defi}[Extended $\EM_1$ Intuitionistic Arithmetic: $\HA + \EM_1$]\label{definition-extendedarithmetic}
The language $\LanguageClass$ of $\HA + \EM_1$ is defined as follows.
\begin{enumerate}

\item
The terms of $\LanguageClass$ are all $t \in \SystemTClass$ with state $\varnothing$, such that $t:\Nat$ and $FV(t) \subseteq \{x_1^\Nat, \ldots, x_n^\Nat\}$ for some $x_1, \ldots, x_n$.\\

\item
The atomic formulas of $\LanguageClass$ are all $Qt_1\ldots t_n \in \SystemTClass$, for some $Q:\Nat^{n}\rightarrow \Bool$ {\em closed term of $\SystemTClass$} of state $\varnothing$, and some terms $t_1,\ldots,t_n$ of $\LanguageClass$.\\

\item
The formulas of $\LanguageClass$ are built from atomic formulas of $\LanguageClass$ by the connectives $\lor,\land,\rightarrow \forall,\exists$ as usual.\\
\end{enumerate}

%$\LanguageClass^\sigma$ is the subset of the terms and formulas
%$\LanguageClass$ without constants of type $\State$ and whose only
%free
%variable is $\sigma:\State$.

A formula of $\HA$ is a formula of $\HA + \EM_1$ in which all predicates and terms are terms of $\SystemTG$.

Deduction rules for $\HA + \EM_1$
are as in van Dalen \cite{van Dalen}, with: %\begin{itemize}
%\item
{\em (i)} an axiom schema for $\EM_1$;
%\item
{\em (ii)} the induction rule;
%\item
{\em (iii
)} as Post rules:
all axioms of equality and ordering on $\Nat$, all
equational axioms of $\SystemTG$, and one schema for each
tautological consequences of $\SystemTG$. {\em (iv)} the axiom schemas for oracles: $P(\vec{t},t) \Rightarrow_\Bool \Chi_P \vec{t}$ and for Skolem maps: $\Chi_P \vec{t} \Rightarrow_\Bool P( \vec{t},(\Phic_P \vec{t}))$, for any predicate $P$ of $\SystemTG$.

%\end{itemize}
\end{defi}

We denote with $\bot$ the atomic formula ${\False}$ and will sometimes write  a generic atomic formula as $P(t_1,\ldots, t_n)$ rather than in the form $Pt_1\ldots t_n$. Finally, since any arithmetical formula has only variables of type $\Nat$, we shall freely omit their types, writing for instance $\forall x. A$ in place of $\forall x^\Nat. A$. Post rules cover many rules with atomic assumptions and conclusion as we find useful, for example, the rule: ``if $f(z)\leq 0$ then $f(z)=0$''.

We defined $\Rightarrow_\Bool:\Bool,\Bool\rightarrow \Bool$ as a term implementing implication, therefore, to be accurate, the axiom $P(t_1,\ldots, t_n,t) \Rightarrow_\Bool \Chi_P t_1\ldots t_n$ is not an implication between two atomic formulas, but it is equal to the single atomic formula $Qt_1\ldots t_nt$, where \[ Q = \lambda x_1^{\Nat}\ldots \lambda
x_{n+1}^{\Nat}. \Rightarrow_\Bool (Px_1\ldots x_nx_{n+1})(\Chi_P x_1\ldots x_{n+1})\]

Similarly, $\lnot_{\Bool}P(\vec{t},t)$ will denote a single atomic formula. Any atomic formula $A$ of $\LanguageClass$ is a boolean term of $\SystemTClass$, therefore for any state constant $s$ we may form the ``finite approximation'' $A[{s}]:\Bool, A[{s}] \in \SystemTLearn$ of $A$. In $A[{s}]$ we replace all oracles $\Chi_P$ and all Skolem maps $\Phic_P$ we have in $A$ by their finite approximation $\chi_P{s}, \phi_Ps$, computed with respect to the state constant $s$. We denote with $\LanguageLearn$ the set of all expressions $A[s]$ with $A \in \LanguageClass$ and $s$ a state constant. All $A[s] \in \LanguageLearn$ may be interpreted by first order arithmetical formulas having all closed atomic subformulas decidable.

Using the metaphor explained in the introduction, we use a set of falsifiable hypotheses determined by $s$ to predict a computable truth value $A[{s}]:\Bool$ in $\SystemTLearn$ for an atomic formula $A \in \LanguageClass$ that we cannot effectively evaluate. Our definition of realizability  provides a formal semantics  for the Extended
Intuitionistic Arithmetic $\HA+\EM_1$, and therefore also for the more usual language of Arithmetic $\HA$, in which all functions represent recursive maps.

\begin{defi}[Types for realizers]
\label{definition-TypesForRealizers} For each
arithmetical formula $A$ we define a type $|A|$ of $\SystemTG$ by
induction on $A$:
\begin{enumerate}
\item
$|P(t_1,\ldots,t_n)|=\State$,\\
\item
$|A\wedge B|=|A|\times |B|$,\\
\item
$|A\vee B|= \Bool\times (|A|\times |B|)$,\\
\item
$|A\rightarrow B|=|A|\rightarrow |B|$,\\
\item
$|\forall x A|=\Nat\rightarrow |A|$,\\
\item
$|\exists x A|= \Nat\times |A|$\\
\end{enumerate}
\end{defi}

We now define the realizability relation $t\Vvdash A$, where $t \in \SystemTClass$, $A \in \LanguageClass$, $t$ has state $\varnothing$ and $t:|A|$.

\begin{defi}[Realizability]
\label{lemma-IndexedRealizabilityAndRealizability}
Assume $s$ is a state constant, $t\in \SystemTClass$ is a closed term of state $\varnothing$, $C \in \LanguageClass $ is a closed formula, and $t:|C|$. Let $\vec{t} = t_1, \ldots, t_n : \Nat$. We define first the relation $t\Vvdash_s C$ by induction and by cases according to the form of $C$:

\begin{enumerate}
\item
$t\Vvdash_s P(\vec{t})$ if and only if $t[s]  = \varnothing$ in $\SystemTLearn$ implies
$P(\vec{t})[{s}]={\True}$\\

\item
$t\Vvdash_s{A\wedge B}$ if and only if $\pi_0t \Vvdash_s{A}$ and $\pi_1t\Vvdash_s{B}$\\

\item
$t\Vvdash_s {A\vee B}$  if and only if either $\proj_0t[{s}]={\True}$ in $\SystemTLearn$ and $\proj_1t\Vvdash_s A$, or $\proj_0t[{s}]={\False}$ in $\SystemTLearn$ and $\proj_2t\Vvdash_s B$\\

\item
$t\Vvdash_s {A\rightarrow B}$ if and only if for all $u$, if $u\Vvdash_s{A}$,
then $tu\Vvdash_s{B}$\\

\item
$t\Vvdash_s {\forall x A}$ if and only if for all numerals $n$,
$t{n}\Vvdash_s A[{n}/x]$\\
\item

$t\Vvdash_s \exists x A$ if and only for some numeral $n$, $\pi_0t[{s}]= {n}$ in $\SystemTLearn$ and $\pi_1t \Vvdash_s A[{n}/x]$\\

\end{enumerate}
We define $t\Vvdash A$ if and only if for all state constants $s$, $t\Vvdash_sA $.

\end{defi}

The definition of $\Vvdash$ formalizes all the idea we sketched in the introduction. A realizer is a term $t$ of $\SystemTClass$, possibly containing the non-computable functions $\Chi_P, \Phic_P$; if such functions were computable, $t$ would be an intuitionistic realizer. Since in general $t$ is not computable, we calculate its approximation $t[s]$ at state $s$, which is a term of $\SystemTLearn$, and we require it to satisfy the indexed-by-state realizability clauses. Realizers of disjunctions and existential statements provide a witness, which is an
individual depending on an actual state of knowledge, representing all the hypotheses
used to approximate the non-computable. The actual behavior
of a realizer depends upon the current state of knowledge. The state is used only when there is relevant information about the truth of a given formula to be computed: the truth value $P(t_1,\ldots,t_n)[{s}]$ of an atomic formula and the disjunctive  witness $\proj_0t[{s}]$ and the existential witness $\pi_0u [{s}]$ are computed w.r.t. the constant state $s$. A realizer $t$ of $A\lor B$ uses $s$ to predict which one between $A$ and $B$ is realizable (if $\proj_0t[{s}]={\True}$ then $A$ is realizable, and if $\proj_0t[{s}]={\False}$ then $B$ is realizable). A realizer $u$ of $\exists  x A$ uses $s$ to predict that $\pi_0u[{s}]$ equals an ${n}$, some witness for $\exists x A$ (i.e. that $A[{n}/x]$ is realizable). These predictions need not  be always correct; hence, it is possible that a realized atomic formula is actually false; we may have $t \Vvdash_s P$ and $P[s]={\False}$ in $\SystemTLearn$. If an atomic formula, although predicted to be true, is indeed false, then we have encountered a counterexample and so our theory is wrong, our approximation  still inadequate; in this case, $t[{s}] \not = \varnothing$ by definition of $t\Vvdash_s P$, and the atomic realizer $t$ takes $s$ and extends it to a larger state $s'$, union of $s$ and $t[{s}]$. That is to say: if something goes wrong, we must learn from our
mistakes. The point is that after every learning, the actual state of knowledge grows, and if we ask to the same realizer new predictions, we will obtain ``better'' answers.

Indeed, we can say more about this last point. Suppose for instance that $t\Vvdash A\lor B$ and let $\{s_i\}_{i\in\NatSet}$ be a w.i. sequence. Then, since $t: \Bool\times |A|\times |B|$, then $\proj_0t:\Bool$ is a closed term of $\SystemTClass$, converging in $\{s_i\}_{i\in\NatSet}$ to a boolean; thus the sequence of predictions  $\{\proj_0t[s_i]\}_{i\in \NatSet}$ eventually stabilizes, and hence a witness is eventually learned in the limit.

In the atomic case, in order to have $t\Vvdash_s P(t_1,\ldots, t_n)$, we require that if $t[{s}] = \varnothing$, then $P(t_1,\ldots,t_n)[{s}] = {\True}$ in $\SystemTLearn$. That is to say: if $t$ has no new information to add to $s$,  then $t$ must assure the truth of $P(t_1,\ldots,t_n)$ w.r.t. $s$. By the zero theorem \ref{FixedPointProperty}, when $t:\State$ is closed, there is plenty of state constants $s$ such that $t[s]=\varnothing$; hence search for truth will be for us {\em computation of a zero}, driven by the excluded-middle instances and the Skolem axioms used by the proof, rather than exhaustive search for counterexamples. In chapter \ref{chapter-constructiveanalysislearning} we will prove that, actually, zeros for terms of $\SystemTClass$ can be computed by learning processes whose length can be bounded through constructive reasoning.\\

It is useful to give a slightly different definition of indexed realizability, which in some situations is \emph{slightly} easier to reason with. The difference with definition \ref{lemma-IndexedRealizabilityAndRealizability} is only that the relation we are going to define now is between terms of $\SystemTLearn$ and formulas of $\LanguageLearn$, which are from the beginning approximations at some state $s$ respectively of terms of $\SystemTClass$ and formulas of $\LanguageClass$.
\begin{defi}[Variant of Indexed Realizability]
\label{definition-AuxiliaryRealizability} Let $s$ be a state constant. Assume $t \in \SystemTLearn$ and $A \in \LanguageLearn$ are of the form $t = t'[s], A = A'[s]$ for some closed $t' \in \SystemTClass$ of state $\varnothing$ and some closed $A' \in \LanguageClass$. We define $t \Vdash_s A$ for any state constant $s$ by induction on $A$.
\begin{enumerate}
\item $t\Vdash_s P(t_1,\ldots, t_n)$ if and only if  $t = \varnothing$ in $\SystemTLearn$ implies
$P(t_1,\ldots,t_n)={\True}$\\

\item
$t\Vdash_s {A\wedge B}$ if and only if $\pi_0t\Vdash_s{A}$ and
$\pi_1t\Vdash_s{B}$\\

\item
$t\Vdash_s {A\vee B}$ if and only if: either $p_0t={\True}$ in $\SystemTLearn$ and
$p_1t\Vdash_s{A}$, or $p_0t={\False}$ and $p_2t\Vdash_s{B}$\\

\item
$t\Vdash_s {A\rightarrow B}$ if and only if for all $u$, if $u\Vdash_s{A}$,
then $tu\Vdash_s{B}$\\

\item
$t\Vdash_s {\forall x A}$ if and only if for all numerals $n$,
$t{n}\Vdash_s A[{n}/x]$\\

\item
$t\Vdash_s {\exists x A}$ if and only if  for some numeral $n$ $\pi_0t= {n}$ in $\SystemTLearn$ and $\pi_1t\Vdash_s A[{n}/x]$

\end{enumerate}

\end{defi}

The realizability relation is compatible with equality in $\SystemTLearn$:

\begin{lem}\label{RealizabilityEquality}

If $t_1\Vdash_s A[u_1/x]$, $t_1=t_2$ and $u_1=u_2$ in $\SystemTLearn$, then $t_2\Vdash_s A[u_2/x]$

\end{lem}

\proof By straightforward induction on $A$.

\qed

We can now characterize $\Vvdash$ in the following way.   
%
%%%%%%%%%%%%%%%%%%%%%%%%%%%%%%%%%%%%%%%%%%%%%%%%
% THIS IS THE COMPLETE CHARACTERIZATION WE HAVE%
%%%%%%%%%%%%%%%%%%%%%%%%%%%%%%%%%%%%%%%%%%%%%%%%

\begin{lem}[Alternative Characterization  of Realizability]
\label{lemma-IndexedRealizabilityAndRealizability}
Assume $t\in \SystemTClass$ is a closed term, $A \in \LanguageClass $ is a closed formula, and $t:|A|$. Then 
\[t\Vvdash A\text{ if and only if  for all state constants } s,\ t[s]\Vdash_s A[s]\]
\end{lem}

\proof By definition unfolding and by induction on $A$, one shows that $t\Vvdash_sA$ if and only if $t[s]\Vdash_s A[s]$.

\qed

\begin{exa} The most remarkable feature of our Realizability Semantics is the existence of a  realizer $\E{P}$ for $\EM_1$. Assume that $P$ is a predicate of $\SystemTG$ and define  \[\E{P}:=\lambda \vec{\alpha}^{\Nat}
\langle  \Chi_P\vec{\alpha},\ \langle
\Phic_P\vec{\alpha},\
\varnothing \rangle ,\ \lambda n^{\Nat}\
\Add_P \vec{\alpha}n\rangle   \]
Indeed $\E{P}$ realizes its associated instance of $\EM_1$.
\begin{prop}[Realizer $\E{P}$ of $\EM_1$]\label{RealizerOfEM1}
 \[\E{P}\Vvdash \forall \vec{x}.\ \exists y\ P(\vec{x}, y)\vee \forall y \neg_\Bool P(\vec{x},y)\]
\end{prop}

\proof Let $\vec{m}$ be a vector of numerals. $\E{P}\vec{m} [{s}]$ is equal to \[ \langle  \chi_P{s}\vec{m},\ \langle \varphi_P{s}\vec{m},\ \varnothing \rangle ,\ \lambda n^{\Nat}\ \add_P{s}\vec{m}n \rangle  \] and we want to prove that \[\E{P}\vec{m}[s]\Vdash_s  \exists y\ P(\vec{m}, y)\vee \forall y \neg_\Bool P(\vec{m},y)\]
We have $\proj_0\E{P}\vec{m}[{s}]=\chi_P{s}\vec{m}$ in $\SystemTLearn$. There are two cases.
\begin{enumerate}
\item $\chi_P{s} \vec{m}={\True}$. Then $\langle P, \vec{m}, {n} \rangle \in |s|$ for some numeral $n$ such that $P(\vec{m},{n}) = {\True}$, and we have to prove \[\proj_1 \E{P}{m}[{s}] \Vdash_s  \exists y\ P(\vec{m}, y)\] By definition of $\varphi_Ps\vec{m}$
\[\proj_1\E{P}{m}[{s}] = \langle \varphi_P{s} \vec{m},\varnothing\rangle = \langle {n}, \varnothing \rangle\] Thus
\[\pi_0(\proj_1\E{P}{m})[{s}] =  \pi_0\langle {n}, \varnothing \rangle = {n}\] and 
\[\pi_1(\proj_1\E{P}{m})[s] \Vdash_s P(\vec{m}, {n})\] because $P(\vec{m}, {n})={\True}$. We conclude 
\[\proj_1\E{P}{m}[{s}] \Vdash_s  \exists y\ P(\vec{m}, y)\] 

\item $ \chi_P{s} \vec{m}={\False}$. Then $\langle  P, \vec{m}, l\rangle  \not
\in |s|$ for all numerals $l$.
We have to prove
\[\proj_2\E{P}\vec{m}[{s}]=\lambda n\ \add_P{s}\vec{m}n\
\Vdash_s
\forall y \neg_\Bool P({m},y)\] i.e. that given any numeral $n$ \[ \add_ P{s}\vec{m}{n}  \Vdash_s \neg_\Bool
P({m},{n})\] By the definition of realizer in this case, we have to assume that $\add_P{s} \vec{m}{n} = \varnothing$, and prove that $\neg_\Bool
P(\vec{m}, {n}) [{s}] = {\True}$. The substitution $(.)[{s}]$ has an empty effect over $P(\vec{m}, {n})$, therefore we have to prove that $\neg_\Bool P(\vec{m}, {n}) = {\True}$, that is, that $P(\vec{m},{n})= {\False}$. Assume for contradiction that $P(\vec{m},{n}) ={\True}$. We already proved that $\langle P, \vec{m}, l \rangle \not \in |s|$, for all numerals $l$: from this and $P(\vec{m},{n}) ={\True}$ we deduce that by definition $\add_P {s} \vec{m}{n} = \makestate{\{\langle P, \vec{m}, {n} \rangle \}}$, contradiction.\\
\end{enumerate}
\qed
%
% Finally, for the Stability
%Theorem  $\E{P}\in \| | \forall \vec{x}.\ \exists y\ P(\vec{x},
%y)\vee
%\forall y \neg_\Bool P(\vec{x},y)|\|$.\\
%

$\E{P}$ works according to the ideas we sketched in the introduction.
It uses $\chi_P$ to make predictions about which one between
$\exists y\ P(\vec{m}, y)$ and $\forall y \neg_\Bool P(\vec{m},y)$ is true. $\chi_P$, in turn, relies on  the constant $s$ denoting the actual state to make its own prediction. If $\chi_P{s}{m}={\False}$, given any $n$, $ \neg_\Bool P({m},{n})$ is predicted to be true; if it is not
the case, we have a counterexample and $\add_P$ requires to extend the state with $\langle P,\vec{m},{n}\rangle $. On the contrary, if $\chi_P{s}{m}={\True}$, there is unquestionable evidence that $\exists y P(\vec{m},y)$ holds; namely, there is some numeral $n$ such that $\langle P,\vec{m},{n}\rangle$ is in $s$; then $\varphi_P$ is called, and it returns $\varphi_P{s} \vec{m}={n}$.

This is the basic mechanism by which we implement learning: every state extension is linked with an assumption about an
instance of $\EM_1$ which we used and turned out to be wrong (this is the only way to come across a counterexample); in next computations, the actual state will be bigger, the realizer will not do the same error,  and hence will be ``wiser".
\end{exa}
As usual for a realizability interpretation, we may extract from any
realizer $t\Vvdash {\forall x. \exists y. P(x,y)}$, with $P \in \SystemTG$, some recursive map $f$ from the set of numerals to the set of numerals, such that $P(n,f(n))$ for all numerals $n$. 

\begin{exa}[Program Extraction via Learning Based Realizability]\label{theorem-extraction0} Let $t$ be a term of $\SystemTClass$ and suppose that $t\Vvdash \forall x^\Nat \exists y^\Nat Pxy$, with $P$ atomic. Then, from $t$ one can effectively define a recursive function $f$ from the set of numerals to the set of numerals such that for every numeral $n$, $Pn(f(n))=\True$.
\end{exa}

\proof Let \[v:=\lambda m^\Nat\ \pi_1(tm)\]
$v$ is of type $\Nat\rightarrow \State$. 
By zero theorem \ref{FixedPointProperty}, there exists a recursive function $\mathsf{zero}$ from the set of numerals to the set of state constants such that $vn[\mathsf{zero}(n)]=\varnothing$ for every numeral $n$.  Define $f$ as the function
\[m\mapsto \pi_0(tm)[\mathsf{zero}(m)]\]
and fix a numeral $n$. By unfolding the definition of realizability with respect to the state $\mathsf{zero}(n)$, we have that 

\[tn\Vvdash_{\mathsf{zero}(n)} \exists y^\Nat Pny\] 
and hence
\[\pi_1(tn)\Vvdash_{\mathsf{zero}(n)}\ Pn(f(n))\]
that is to say
\[vn[\mathsf{zero}(n)]=\varnothing\implies Pn(f(n))=\True\]
and therefore \[Pn(f(n))=\True\] 
 which is the thesis. 

\qed

\begin{rem} In chapter \ref{chapter-constructiveanalysislearning} we shall prove that the map $f$ constructed in example \ref{theorem-extraction0} is even definable in G\"odel's $\SystemTG$. This result formally proves that our realizability interpretation is a constructive semantics and that $f$ is not a brute force search algorithm. More precisely, we can argue as follows. The numeral $f(n)$ is computed by finding a zero of $vn=\pi_1(tn)$, i.e. a state $s$ such that $\pi_1(tn)[s]=\varnothing$. By the Zero theorem \ref{FixedPointProperty}, this zero is computed step by step by constructing the sequence $s_0=\varnothing$, $s_{n+1}=s_n\Cup t[s_n]$ and stopping at the first $m$ such that $\pi_1(tn)[s_m]=\varnothing$. \\First, we observe that each portion of $s_m$ is efficiently constructed: for each $n$, the state $s_n$ is efficiently extended to $s_{n+1}$, through the addition of new oracle values learned by $t$ by counterexamples, i.e. by the falsification of some excluded middle or Skolem axiom instances. No brute force search whatsoever, thus, for new oracle values: they are all efficiently produced by $t[s_n]$, which, modulo some trivial coding, \emph{is}  a term of G\"odel's $\SystemTG$ and just cannot search blindly for oracle values, since it is a primitive recursive functional of finite type.\\ Secondly, an upper bound to $m$ can be computed in G\"odel's $\SystemTG$, as proven in chapter \ref{chapter-constructiveanalysislearning}, theorem \ref{theorem-zerotclass}. Moreover, this upper bound results from a constructive proof of the Zero theorem. Hence, $s_m$ and thus $f(n)$ can be defined by a primitive recursive functional of finite type, which, again, by construction cannot explore blindly the infinite search space of the knowledge states in order to find a zero.

Moreover, from the low level computational point of view and in the language of $\epsilon$-substitution method, our realizers represent convergent procedures to find out a ``solving substitution", i.e. a state representing an approximation of Skolem functions (i.e., $\epsilon$-terms) which makes true the Skolem axioms instances used in a proof of an existential statement. The advantage of our semantics is the possibility of defining such procedures directly from high level proofs, by means of Curry-Howard correspondence, hence avoiding the roundabout route which forces to use a quantifier free deduction system. In the case of a provable formula in the language of Peano Arithmetic (that is, one not containing the symbols $\Chi_P$ or $\Phic_P$) we do not need at all to modify the language of its proof and to use the Skolem axioms $\chi, \varphi$.
\end{rem}

Now we explain how to turn each proof $\mathcal{D}$ of a formula $A
\in \LanguageClass$ in $\HA + \EM_1$ into a realizers $\mathcal{D}^*$
of the same $A$. By induction on $\mathcal{D}$, we define a
``decoration with
realizers'' $\mathcal{D}^\Realizer$ of $\mathcal{D}$, in which each
formula $B$ of $\mathcal{D}$ is replaced by a new statement $u
\vdash B$, for some $u \in \SystemTClass$ of state $\varnothing$. If $t \vdash A$ is the
conclusion of $\mathcal{D}^\Realizer$, we set $\mathcal{D}^* = t$.
Then we will prove that if $\mathcal{D}$ is closed and without
assumptions, then $\mathcal{D}^* \in \SystemTClass$ and
$\mathcal{D}^* \Vvdash A$. The decoration  $\mathcal{D}^\Realizer$
of
$\mathcal{D}$ with realizers is completely standard: we have new
realizers only for $\EM_1$ and for atomic formulas. For
notation simplicity, if $x_i$ is the label for the set of
occurrences of some assumption $A_i$ of $\mathcal{D}$, we use $x_i$
also as a name of one free variable in $\mathcal{D}^*$ of type
$|A_i|$. If $T$ is any type of $\SystemTState$, we denote with $d^T$ a dummy term of type $T$, defined by $d^\Nat = 0$, $d^\Bool = \False$, $d^\State = \varnothing$, $d^{A \rightarrow B} = \lambda \_^A.d^B$ (with $\_^A$ any variable of type $A$), $d^{A \times B} = \langle d^A, d^B\rangle$.

\begin{defi}[Term Assignment Rules for $\HA + \EM_1$]
\label{definition-NaturalDeductionandTermAssignmentRulesforExtendedArithmetic}Assume $\mathcal{D}$ is a proof of $A \in \LanguageClass$ in $\HA +
\EM_1$, with free assumptions $A_1, \ldots, A_n$ denoted by proof variables
$x_1^{A_1},
\ldots, x_n^{A_n}$ and free integer variables $\alpha_1^\Nat, \ldots,
\alpha_m^\Nat$. By induction on $\mathcal{D}$, we define a
decorated
proof-tree $\mathcal{D}^\Realizer$, in which each formula $B$
is replaced by $u \vdash B$ for some $u \in \SystemTClass$, and the
conclusion $A$ with some $t \vdash A$, with $FV(t) \subseteq
\{x_1^{|A_1|}, \ldots, x_1^{|A_1|},\alpha_1^\Nat, \ldots,
\alpha_m^\Nat \}$. Eventually we set $\mathcal{D}^*=t$.\\
%
%Let $A_1,A_2,\ldots, A_n,\ldots$ be an effective enumeration of
%formulas of $\Language$. For each formula $A_i$ and for every
%$n\in\mathbb{N}$ we define the symbol $x_{n}^{A_i}$ as a shorthand
%for $x_{(n,i)}^{|A_i|}$, where $(\cdot, \cdot
%):\mathbb{N}^2\rightarrow \mathbb{N}$ is a recursive bijection.
%his way,  $x_{n}^{A_i}$ - notwithstanding its being a variable of
%$\SystemT_1$ and of type $|A_i|$ - does not lose information about
%$A_i$, since from the index $(n,i)$ we can effectively recover $i$
%and hence $A_i$. If we are not interested in the index $n$, we
%simply write $x^{A_i}$ rather than $x_n^{A_i}$.
%%For the sake of faithful proof representation, we assume that,

%every formula $A$,  there is a special variable $x_A^{|A|}$ of

%$|A|$; we will denote $x_A^{|A|}$ by $x_A$.
%
%As usual, we define a system of inference rules and we shall write
%$u\vdash A$ if that expression can be obtained by means of those
%inference rules. We call $D$ a natural deduction of $u\vdash A$ if
%it is a proof tree ending with $u\vdash A$. The string $u\vdash A$
%should be read as: ``the term $u: |A|$ codifies a proof of $A$
%under %the assumptions $A_1,\ldots, A_n$ given by its free
%variables of the %form $x^{A_1},\ldots, x^{A_n}$ ''.
%%We assign to each assumption of the proof and to each
%%conclusion of  each inference rule some term $u$ of
%$\SystemT_1^{-}$ such that
%%if the premises of the deduction are among $A_1,\ldots, A_n$,
%then
%%$FV(u)\subseteq \{x_1^{|A_1|},\ldots, x_n^{|A_n|}, \alpha_1^\Nat,
%\ldots, \alpha_k^\Nat, \sigma\}$, for some $x_i$ and
%$\alpha_i$.

\begin{enumerate}

\item %AXIOM
$\begin{array}{c}   \hline  x^{|A|}\vdash A
\end{array}\ \ \ \ $ if $\mathcal{D}$ consists of a single free
assumption $A\in \LanguageClass$ labeled $x^A$.\\

\item %CONJUNCTION
$\begin{array}{c}  u\vdash A\ \ \ t\vdash B\\ \hline \langle
u,t\rangle \vdash
A\wedge B
\end{array}\ \ \ \ $
$\begin{array}{c}  u\vdash A\wedge B\\ \hline \pi_0 u\vdash A
\end{array}\ \ \ \ $
$\begin{array}{c}  u\vdash A\wedge B\\ \hline \pi_1 u\vdash B
\end{array}$\\\\

\item %IMPLICATION
$\begin{array}{c}  u\vdash A\rightarrow B\ \ \ t\vdash A \\ \hline
ut\vdash B
\end{array}\ \ \ \ $
$\begin{array}{c}  u\vdash B\\ \hline \lambda x^{|A|}u\vdash
A\rightarrow B
\end{array}$\\\\

\item %DISJUNCTION
$\begin{array}{c}  u\vdash A\\ \hline \langle {\True},u,d^{B}\rangle
\vdash A\vee B
\end{array}\ \ \ \ $
$\begin{array}{c}  u\vdash B\\ \hline \langle  {\False},d^{A}, u\rangle
\vdash A\vee B
\end{array}$\\
$\begin{array}{c}  u\vdash A\vee B\ \ \ w_1\vdash C\ \ \ w_2\vdash
C\\ \hline  \ifthen{\proj_0 u}{(\lambda x^{|A|} w_1)(\proj_1u)}{
(\lambda x^{|B|} w_2)(\proj_2 u)}\vdash C
\end{array}$\\

where $d^A$ and $d^B$ are dummy closed terms of $\SystemTClass$ of type $|A|$ and
$|B|$\\ %and $x^{|A|}$, $y^{|B|}$ are variables respectively
%associated to premises $A$ and $B$ of the deductions of $w_1\vdash
%C$ and $w_2\vdash C$.

\item %FOR ALL
$\begin{array}{c}  u\vdash \forall \alpha A\\ \hline  ut\vdash A[t/\alpha]
\end{array} $
$\begin{array}{c}  u\vdash A\\ \hline \lambda \alpha^{\Nat} u\vdash
\forall \alpha A
\end{array}$

where $t$ is a term of $\LanguageClass$ and $\alpha^{\Nat}$ does not occur
free in any free assumption $B$ of the subproof of $\mathcal{D}$ of
conclusion $A$.\\

\item %EXISTS
$\begin{array}{c}  u\vdash A[t/\alpha^\Nat]\\ \hline  \langle
t,u\rangle \vdash
\exists
\alpha^\Nat. A
\end{array}$ \ \ \ \
$\begin{array}{c}  u\vdash \exists \alpha^\Nat. A\ \ \ t\vdash C\\
\hline
(\lambda \alpha^{\Nat}\lambda x^{|A|}\ t)(\pi_0 u)(\pi_1 u)\vdash C
\end{array} $\\

where $\alpha^{\Nat}$ is not free in $C$
nor in any free assumption $B$ different from $A$ in the
subproof of $\mathcal{D}$ of conclusion $C$.\\

\item %INDUCTION
$\begin{array}{c}  u\vdash A(0)\ \ \ v\vdash \forall \alpha.
A(\alpha)\rightarrow A(\suc(\alpha))\\ \hline \lambda \alpha^{\Nat} \rec uv\alpha\vdash
\forall
\alpha A
\end{array}\ \ \ \ $\\

\item %ATOMIC RULES
$\begin{array}{c}  u_1\vdash A_1\ u_2\vdash A_2\ \cdots \ u_n\vdash
A_n\\ \hline u_1\Cup u_2\Cup\cdots\Cup u_n\vdash A
\end{array}$

where $n > 0 $ and $A_1,A_2,\ldots,A_n,A$ are atomic
formulas of $\LanguageClass$, and the rule is a Post rule for equality
or
ordering, or
a tautological
consequence. \\

\item %AXIOMS OF HA
$\begin{array}{c}  \hline \varnothing \vdash A
\end{array}$

where $A$ is an atomic axiom of $\HA + \EM_1$ (an axiom
of equality or of ordering or a tautology or an equation of
$\SystemTG$) \\

\item %1-EXCLUDED MIDLE
$\begin{array}{c}\hline \E{P} \vdash  \forall \vec{x}.\ \exists y\
P(\vec{x}, y)\vee \forall y \neg_\Bool P(\vec{x},y)
\end{array}$

where $P$ is a predicate of $\SystemTG$.\\ 
\item
$\begin{array}{c} \hline \Add_P \vec{t},t \vdash
P(\vec{t},t) \Rightarrow_\Bool \Chi_P \vec{t}
\end{array}$, $\chi$-Axiom\\

\item
$\begin{array}{c} \hline \varnothing \vdash \Chi_P \vec{t}
\Rightarrow_\Bool P(\vec{t},(\Phic_P \vec{t}))
\end{array}$, $\varphi$-Axiom

\end{enumerate}
\end{defi}

The term decorating the conclusion of a Post rule is of the form
$u_1\Cup \cdots \Cup u_n$. In this case, we have $n$ different
realizers, whose learning capabilities are put together through a
sort of union.
By Lemma \ref{lemma-Cup}.2, if $u_1\Cup \cdots \Cup u_n[{s}] = \varnothing$, then $u_1[{s}] = \ldots = u_n[{s}] = \varnothing$,
i.e. all $u_i$ ``have nothing to learn''. In that case, each $u_i$ must guarantee $A_i$ to be true, and therefore the conclusion of the Post rule is true, because true premises $A_1, \ldots, A_n$ spell a true conclusion
$A$.

We now prove our main theorem, that every theorem of $\HA +
\EM_1$ is realizable.

\begin{thm}[Adequacy Theorem]\label{Adeguacy Theorem}
Suppose that $\mathcal{D}$ is a proof of $A$ in
the system $\HA + \EM_1$ with free assumptions
$x_1^{A_1},\ldots,x_n^{A_n}$ and free variables
$\alpha_1:{\Nat},\ldots,\alpha_k:{\Nat}$. Let $w =
\mathcal{D}^*$.
For all state constants $s$ and for all numerals $n_1,\ldots,n_k$, if
\[t_1[s]\Vdash_s A_1[{n}_1/\alpha_1\cdots
{n}_k/\alpha_k][{s}], \ldots, t_n[s]\Vdash_s
A_n[{n}_1/\alpha_1\cdots
{n}_k/\alpha_k][{s}]\] then \[w[t_1/x_1^{|A_1|}\cdots
t_n/x_n^{|A_n|}\  {n}_1/\alpha_1\cdots
{n}_k/\alpha_k][{s}]\Vdash_s A[{n}_1/\alpha_1\cdots
{n}_k/\alpha_k][{s}]\]
\end{thm}

\IfPaperState
%%%%%%%%%%%%%%%%%%%%
% SHORT VERSION
%%%%%%%%%%%%%%%%%%%
{
{\bf Proof} By induction on $w$ (see \cite{ExtendedVersion}).
}
%%%%%%%%%%%%%
% END OF SHORT VERSION
%%%%%%%%%%%%%
%%%%%%%%%%%%%
% LONG VERSION
%%%%%%%%%%%%%%%
{
\proof Notation: for any term $v$ and formula $B$, we
denote
\[v[t_1/x_1^{|A_1|}\cdots t_n/x_n^{|A_n|}\
{n}_1/\alpha_1\cdots
{n}_k/\alpha_k][{s}]\]
with $\substitution{v}$ and $B
[{n}_1/\alpha_1\cdots {n}_k/\alpha_k][{s}]$ with
$\substitution{B}$. We have $|\substitution{B}| = |B|$ for all formulas
$B$. We denote with $=$ the provable equality in $\SystemTLearn$.
We proceed by induction on $w$. Consider the last rule in the
derivation $\mathcal{D}$:

\begin{enumerate}
\item
If it is the rule for variables, then
$w=x_i^{|A_i|}=x^{|\substitution{A_i}|}$ and
$A=A_i$. So $\substitution{w}=t_i\Vdash_s
\substitution{A_i}=\substitution{A}$.\\

\item
If it is the $\wedge I$ rule, then $w=\langle u,t\rangle $,
$A=B\wedge
C$, $u\vdash B$ and $t\vdash C$. Therefore, $\substitution{w}=
\langle \substitution{u},\substitution{t}\rangle $. By induction
hypothesis,
$\pi_0\substitution{w}=\substitution{u}\Vdash_s \substitution{B}$ and
$\pi_1\substitution{w}=\substitution{t}\Vdash_s \substitution{C}$; so, by
definition, $\substitution{w}\Vdash_s
\substitution{B}\wedge\substitution{C}=\substitution{A}$.\\

  \item If it is a $\wedge E$ rule, say left, then $w=\pi_0 u$ and
$u\vdash A\wedge B$. So $\substitution{w}=\pi_0 \substitution{u}\Vdash_s
\substitution{A}$, because $\substitution{u}\Vdash_s \substitution{A}\wedge
\substitution{B}$ by induction hypothesis.\\

  \item If it is the $\rightarrow E$ rule, then $w=ut$, $u\vdash
B\rightarrow A$ and $t\vdash B$. So
$\substitution{w}=\substitution{u}\substitution{t}\Vdash_s \substitution{A}$, for
$\substitution{u}\Vdash_s \substitution{B}\rightarrow \substitution{A}$ and
$\substitution{t}\Vdash_s \substitution{B}$ by induction hypothesis.\\

  \item If it is the $\rightarrow I$ rule, then $w=\lambda x^{|B|}
u$,
$A=B\rightarrow C$ and $u\vdash C$. Thus, $\substitution{w}=\lambda x^{|B|} \substitution{u}$. Suppose now that $t\Vdash_s
\substitution{B}$;
%we have to prove that $\substitution{w}t\Vdash_s
%\substitution{C}$.
by induction hypothesis on $u$,
$\substitution{w}t=\substitution{u}[t/x^{|B|}] \Vdash_s \substitution{C}$.\\

%
%
%
%\[\substitution{w}t=\substitution{u}[t/x^{B}]=u[{n}_1/\alpha_1\cdots
%{n}_k/\alpha_k\ x^{\substitution{B}}/x^B][t_1/x_1^{A_1}\cdots
%t_n/x_n^{A_n}\ t/x^{\substitution{B}}]\Vdash_s \substitution{C}\]
%
%since $u[{n}_1/\alpha_1\cdots
%{n}_k/\alpha_k\ x^{\substitution{B}}/x^B]\vdash \substitution{C}$
%by Proposition \ref{fact on
%deductions} ($A_1,\ldots, A_n$ are closed).
%
%
%

\item If it is a $\vee I$ rule, say left, then $w=\langle
{\True},u,d^C\rangle $,
$A=B\vee C$ and $u\vdash B$. So,
$\substitution{w}=\langle {\True},\substitution{u},d^{C}\rangle $ and hence
$\proj_0\substitution{w}={\True}$. We indeed verify that
$\proj_1\substitution{w}=\substitution{u}\Vdash_s\substitution{B}$ with the help
of induction hypothesis.\\

\item If it is a $\vee E$ rule, then \[w= \ifthen{\proj_0 u}{
(\lambda x^{|B|} w_1)\proj_1u}{(\lambda y^{|C|} w_2)\proj_2u} \]
 and  $u\vdash B\vee C, w_1\vdash D,w_2\vdash D, A=D$. So,
\[\substitution{w}= \ifthen{\proj_0 \substitution{u}}{(\lambda
x^{|B|}\substitution{w_1})\proj_1\substitution{u}}{
(\lambda
y^{|C|}\substitution{w_2})\proj_2 \substitution{u}} \]
Assume $\proj_0\substitution{u}={\True}$. Then by inductive hypothesis
$\proj_1 \substitution{u} \Vdash_s \substitution{B}$, and
%Moreover, by Proposition
%$\ref{fact on deductions}$, $w_1[{n}_1/\alpha_1\cdots
%{n}_k/\alpha_k\ x^{\substitution{B}}/x^B]\vdash \substitution D$
%and hence
again by induction hypothesis,
$\substitution{w}=\substitution{w}_1[\proj_1\substitution{u}/x^{|\substitution{B}|}]\Vdash_s
\substitution{D}$.
Symmetrically, if $\proj_0\substitution{u}={\False}$, then
$\substitution{w}\Vdash_s \substitution{D}$.\\

\item If it is the $\forall E$ rule, then $w=ut$, $A=B[t/\alpha]$
and $u\vdash \forall \alpha B$. So,
$\substitution{w}=\substitution{u}\substitution{t}$. For some numeral $n$ we have
${n}=\substitution{t}$. By inductive hypothesis  $\substitution{u}\Vdash_s
\forall\alpha \substitution{B}$, therefore
$\substitution{u}\substitution{t}=\substitution{u}{n}\Vdash_s
\substitution{B}[{n}/\alpha] = \substitution{B}[\substitution{t}/\alpha]=\substitution{A}$.\\

  \item If it is the $\forall I$ rule, then $w=\lambda
\alpha^{\Nat}u$, $A=\forall \alpha B$ and $u\vdash B$. So,
$\substitution{w}=\lambda \alpha^{\Nat} \substitution{u}$. Let $n$ be a numeral; we have to prove that
$\substitution{w}{n}=\substitution{u}[{n}/\alpha]\Vdash_s
\substitution{B}[{n}/\alpha]$, which is true, indeed, by
induction hypothesis.\\

  \item If it is the $\exists E$ rule, then $w=(\lambda
\alpha^{\Nat}\lambda x^{|B|} t)(\pi_0u)( \pi_1u)$, $t\vdash A$ and
$u\vdash \exists \alpha^{\Nat}. B$. Assume ${n} = \pi_0
{u}$, for some numeral $n$. Then
\[\substitution{t}[{n}/\alpha^{\Nat},\pi_1\substitution{u}/
x^{|\substitution{B}[{n}/\alpha^{\Nat}]|}]\Vdash_s
\substitution{A}[{n}/\alpha]=A\] by
inductive hypothesis, whose application being justified by the
fact, also by induction, that $\substitution{u}\Vdash_s \exists
\alpha^{\Nat}.
\substitution{B}$ and hence $\pi_1\substitution{u}\Vdash_s
\substitution{B}[{n}/\alpha^{\Nat}]$. We thus obtain
\[\substitution{w}=\substitution{t}[\pi_0
\substitution{u}/\alpha^{\Nat}\
\pi_1\substitution{u}/x^{|B|}]
\Vdash_s
\substitution{A}[{n}/\alpha]=A\]

  \item If it is the $\exists I$ rule, then $w=\langle  t,u\rangle
$, $A=\exists
\alpha
B$, $u\vdash B[t/\alpha]$. So, $\substitution{w}=\langle
\substitution{t},\substitution{u}\rangle $; and, indeed, $\pi_1
\substitution{w}=\substitution{u}\Vdash_s \substitution{B}[\pi_0\substitution{w}/\alpha]=\substitution{B}[\substitution{t}/\alpha]$ since by induction hypothesis
$\substitution{u}\Vdash_s \substitution{B}[\substitution{t}/\alpha]$.\\

  \item If it is the induction rule, then $w=\lambda \alpha^{\Nat}\
\rec uv\alpha$, $A=\forall \alpha B$, $u\vdash B(0)$ and $v\vdash
\forall \alpha. B(\alpha)\rightarrow B(\suc(\alpha))$. So,
$\substitution{w}=\lambda \alpha^{\Nat}
\rec\substitution{u}\substitution{v}\alpha$. Now let $n$ be a numeral. A
plain induction on $n$ shows that
$\substitution{w}{n}=\rec\substitution{u}\substitution{v}{n}\Vdash_s
\substitution{B}[{n}/\alpha]$, for $\substitution{u}\Vdash_s
\substitution{B}(0)$ and $\substitution{v}{i}\Vdash_s
\substitution{B}({i})\rightarrow
\substitution{B}({\suc(i)})$ for all numerals $i$ by
induction hypothesis.\\

  \item If it is a Post rule, then $w=u_1\Cup u_2\Cup \cdots\Cup
u_n$ and     $u_i\vdash  A_i$. So,
$\substitution{w}=\substitution{u}_1\Cup \substitution{u}_2\Cup \cdots\Cup
\substitution{u}_n$. Suppose now that $\substitution{w}[{s}] = \varnothing$; then we
have to prove that $\substitution{A}={\True}$. It suffices to prove that $\substitution{A}_1 =\substitution{A}_2 =\cdots= \substitution{A}_n ={\True}$. By Lemma \ref{lemma-Cup} we have $\substitution{u}_1 =\cdots= \substitution{u}_n =\varnothing$ and by induction hypothesis $\substitution{A}_1 =\cdots=\substitution{A}_n={\True}$, since $\substitution{u}_i\Vdash_s \substitution{A}_i$, for $i=1,\ldots,n$.\\

\item If it is a $\chi$-axiom rule, then  $w=\Add_Pt_1\ldots t_nt $ and \[A= P(t_1,\ldots, t_n,t)
\Rightarrow
\Chi_P
 t_1\ldots t_n\] Let $\vec{t} = \substitution{t}_1,\ldots,\substitution{t}_n$. For some numeral $m$ we have $m = \substitution{t}$. Suppose by contradiction that $\substitution{w} = \varnothing$ and
$P(\vec{t},
\substitution{t}) = P(\vec{t},
{m}) ={\True}$ and $\chi_P s \vec{t}={\False}$. From $\chi_P s \vec{t}={\False}$ we get $\langle P,\vec{t}, {m}'\rangle \not \in s$ for all numerals $m'$. We deduce $\substitution{w} =\add_Ps\vec{t}
{m} = \makestate{\{\langle P, \vec{t},{m}\rangle\}}$, contradiction.\\

\item
$w$ realizes an $\EM_1$ axiom: this is Proposition \ref{RealizerOfEM1}.\\

\item If it is a $\varphi$-axiom rule, then $w=\varnothing$ and
\[A=\Chi_P t_1\ldots t_n \Rightarrow P( t_1,\ldots,
t_n,(\Phic_P t_1\ldots t_n))\] We have
$\substitution{w}=\varnothing$. Let us denote $\vec{t} = {\substitution{t}_1}\ldots \substitution{t}_n$. Suppose that $\chi_P{s}\vec{t} ={\True}$. Then for some numeral $m$ we have $\langle P,\vec{t},{m}\rangle
\in s$ and $P\vec{t}{m} = {\True}$ and $\varphi_P{s}\vec{t} = {m}$. By definition of $\varphi_P$ we have  \[P( \vec{t},(\varphi_P{s} \vec{t}))={\True}\] We conclude that $\substitution{A} = {\True}$.
\end{enumerate}
}
\qed

%%%%%%%%%%%%%%%%%%
% END OF LONG VERSION
%%%%%%%%%%%%%%%%%

\begin{cor}\label{Realizability Theorem}
If $A$ is a closed formula provable in $\HA + \EM_1$, then there
exists $t\in \SystemTClass$ such that $t\Vvdash A$.
\end{cor}

\section{Conclusion and further works}
\label{section-conclusion}

Many notions of realizability for Classical Logic already exists. A notion similar to our one in spirit and motivations is Goodman's notion of Relative realizability \cite{Goodman}. However, there is an intrinsic difference between our solution and Goodman's solution. Goodman uses forcing to obtain a ``static'' description of learning. His ``possible worlds'' are learning states, but {\em there is no explicit operation updating a world to a larger word}. The dynamic aspect of learning (which is represented by a winning strategy in Game Semantics) is therefore lost. Using our realizability model, a realizer of an atomic formula, instead of being a trivial map, is {\em a map extending worlds}, whose fixed points are the worlds in which the atomic formula is true. Extending a world represents, in our realizability Semantics, the idea of ``learning by trial-and-error'' that we have in game semantics, while fixed points represent the final state of the game.

A second notion related to our realizability Semantics is Avigad's idea of ``update procedure'' \cite{Avigad}. A state $s$ in our chapter corresponds to a finite model of skolem maps in Avigad. An ``update procedure'' is a construction ``steering'' the future evolution of a finite partial model $s$ of skolem maps, to which our individuals belong, in a wanted direction. The main difference with our work is that we express this idea formally, by interpreting an ``update procedure'' as a realizer (in the sense of Kreisel) for a Skolem axiom. Another important difference is that our realizability relation is defined for all first-order formulas with Skolem maps, while the theory of ``update procedures'' is defined only for quantifier-free formulas with Skolem maps.

Another difference with the other realizability or Kripke models for Classical Logic is in the notion of individual and in the equality between individuals. Assume that $m$ is the output of a skolem map for $\exists y.P(n,y)$, with $P$ decidable, and $m = \{m[s] | s \in \State\}$ a family of values depending on the finite partial model $s$. Then our realizer for Skolem axioms ``steers'' the evolution of $s$ towards some universe in which the axiom $\exists y.P(n,y) \Rightarrow P(n,m[s])$ is true. Modifying the evolution of $s$ may modify the value of $m[s]$. In our realizability Semantics we introduce a notion of individuality which is ``dynamical'' (depending on a state $s$) and ``interactive'' (the value of the individual depends on what a realizer does). This second aspect is new. A realizer may ``try'' to equate an individual $a = \{a[s] | s \in \State\}$ with another individual $b = \{b[s] | s \in \State\}$. Whenever this is possible, the realizer defines a construction over the evolution of the universe $s$ producing such an effect, while a random evolution of $s$ (without an ``interaction'' with the realizer) does not guarantee that eventually we have $a[s] = b[s]$. This is why, in our realizability model, even equality among concrete objects  is not a ``statical'' fact, but it is the effect of applying a realizer (which is a construction over the evolution of the state or ``world'' $s$). In the other models either equality is ``static'', or, even when it is ``dynamical'', and it changes with time, it is not ``interactive'': the final truth value of an equality is not the effect of the application of the realizer, but it is eventually the same in all future evolutions of the current world.
\\

Many aspects of this chapter will require some further work. A challenging idea is to iterate the construction we had for $\EM_1$, in order to provide a learning model for the entire classical Arithmetic. In this case the leading concepts would be the game-theoretical notion of \emph{``level of backtracking''}, introduced in \cite{BerCoq} and \cite{BerardiLiguoro}, a notion related to the more informal notion of \emph{non-monotonic learning}.

Another aspect deserving further work is comparing the programs extracted from classical proofs with our method and with other methods, say, with Friedman $A$-translation. Our interpretation, explaining in term of learning how the extracted program work, should allow us to modify and improve the extracted program in a way impossible for the more formal (but very elegant) $A$-translation.

We remarked that our interpretation is implicitly parametric with respect to the operation $\CupSem$ merging the realizers of two atomic formulas. As explained in \cite{BerardiLiguoroMonadi}, by choosing different variant of this operation we may study different evaluation strategies for the extracted programs: sequential and parallel, left-to-right and right-to-left, confluent and non-confluent. We would like to study whether by choosing a particular evaluation strategy we may extract a more efficient program.

\chapter{Learning Based Realizability and 1-Backtracking Games}\label{chapter-realizabilityandgames}

\begin{abstract}  We prove a soundness and completeness result for learning based realizability with respect to 1-Backtracking Coquand game semantics. First, we prove that interactive learning based classical realizability  is sound with respect to Coquand game semantics. In particular, any realizer of an implication-and-negation-free arithmetical formula embodies a winning recursive strategy for the 1-Backtracking version of Tarski games. We also give examples of realizer and winning strategy extraction for some classical proofs. Secondly, we extend our notion of realizability to a total recursive learning based realizability and show that the notion is complete with respect to Coquand semantics, when it is restricted to 1-Backtracking games.
 \end{abstract}

\section{Introduction}

In this chapter we show that learning based realizability (see chapter \ref{chapter-learningbasedrealizability}) relates to 1-Backtracking Tarski games as intuitionistic realizability (see Kleene \cite{Kleene}) relates to Tarski games, when one considers implication-and-negation-free formulas. The relationship we refer to is between realizability on one hand, and existence of winning strategies on the other.
In particular, it is known that a negation-and-implication-free arithmetical formula is Kleene realizable if and only if Eloise has a recursive winning strategy in the associated Tarski game. We show as well that an implication-and-negation-free arithmetical formula is ``learning realizable" if and only if Eloise has recursive winning strategy in the associated 1-Backtracking Tarski game.

It is well known that Tarski games (which were actually  introduced by Hintikka, see \cite{Hintikka} and definition \ref{definition-TarskiGame}) are just a simple way of rephrasing the concept of classical truth in terms of a game between two players - the first one, Eloise, trying to show the truth of a formula, the second, Abelard, its falsehood - and that a Kleene realizer gives a recursive winning strategy to the first player. The result is quite expected: since a realizer gives a way of computing all the information about the truth of a formula, the player trying to prove the truth of that formula has a recursive winning strategy. However, not at all \emph{any} classically provable arithmetical formula allows a winning recursive strategy for that player; otherwise, the decidability of the Halting problem would follow. 

In \cite{Coquand}, Coquand introduced a new game semantics for Peano Arithmetic, centered on the concept of ``Backtracking Tarski game": a special Tarski game in which players have the additional possibility of correcting their moves and backtracking to a previous position of the game anytime they wish. Coquand then showed that for any provable negation-and-implication-free arithmetical formula $A$, Eloise has a \emph{recursive} winning strategy in the Backtracking Tarski game associated to $A$. Remarkably, a proof in Peano Arithmetic thus hides a non trivial computational content that can be described as a recursive strategy that produces witnesses in classical Arithmetic by interaction and learning.

In the first part of this chapter, we show that  learning based realizers have direct interpretation as recursive winning strategies in 1-Backtracking Tarski games (which are a particular case of Coquand games: see Berardi et al \cite{BerCoq} and definition \ref{definition-1BacktrackingGames} below). The result was wished, because interactive learning based realizers, by design, are similar to strategies in games with backtracking: they improve their computational ability by learning from interaction and counterexamples in a convergent way; eventually, they gather enough information about the truth of a formula to win its associated game. 

An interesting but incomplete step towards our result was the Hayashi realizability \cite{Hayashi1}. Indeed, a realizer in the sense of Hayashi  represents a recursive winning strategy in 1-Backtracking games. However, from the computational point of view, Hayashi  realizers do not relate to 1-Backtracking games in a significant way: Hayashi winning strategies work by exhaustive search and, actually, do not learn from the game and from the \emph{interaction} with the other player. As a result of this issue, constructive upper bounds on the length of games cannot be obtained, whereas using our realizability it is possible. For example, in the case of the 1-Backtracking Tarski game for the formula $\exists x \forall y f(x)\leq f(y)$, the Hayashi realizer checks all the natural numbers to be sure that an $n$ such that $\forall y f(n)\leq f(y)$ is eventually found. On the contrary, our realizer yields a strategy for Eloise which bounds the number of backtrackings by $f(0)$, as shown in this paper; moreover, what the strategy learns is uniquely determined by interaction with the other player. In this case, the Hayashi strategy is the same one suggested by the classical \emph{truth} of the formula, whereas ours is the constructive strategy suggested by its classical \emph{proof}.

Since learning based realizers are extracted from proofs in $\HA+ \EM_1$ (Heyting Arithmetic with excluded middle over existential sentences, see chapter \ref{chapter-learningbasedrealizability}), one also has an interpretation of classical proofs as strategies with 1-Backtracking. Moreover, studying  learning based realizers in terms of 1-Backtracking games also sheds light on their behaviour and offers an interesting case study in program extraction and interpretation in classical arithmetic. 

In the second part of the chapter, we extend the class of learning based realizers from a classical version of G\"odel's system $\SystemTG$ to a classical version of $\mathcal{PCF}$ and define a more general ``total recursive learning based realizability". This step is analogous to the (conceptual, rather than chronological) step leading from Kreisel realizability to Kleene realizability: one extends the computational power of realizers. We then prove a completeness theorem: for every implication-and-negation-free arithmetical formula $A$, if Eloise has recursive winning strategy in the 1-Backtracking Tarski game associated to $A$, then $A$ is also realizable. \\

The \emph{plan of the chapter} is the following. In section \S \ref{calculusandrealizability}, we recall the definitions and results from chapter \ref{chapter-learningbasedrealizability} that we shall need in the present one. In section \S \ref{gamesrealizability}, we prove our first main theorem: a realizer of an arithmetical formula embodies a winning strategy in its associated 1-Backtracking Tarski game. In section \S \ref{examples}, we extract realizers from two classical proofs and study their behavior as learning strategies. In section \S \ref{completeness}, we define an extension of the learning based realizability of chapter \ref{chapter-learningbasedrealizability} and in section \S \ref{section-completenessproof} prove its completeness with respect to 1-Backtracking Tarski games.

\section{Learning-Based Realizability for the Standard Language of Arithmetic}\label{calculusandrealizability}

 In this chapter, we will use a standard language of Arithmetic: the symbols $\Chi_P, \Phic_P$ will not occur in the language of formulas, but only in realizers. We recall the definition and results we need here.

\begin{definition}[Convergence]
\label{definition-Convergence} Assume
that $\{s_i\}_{i\in\NatSet} $ is a w.i. sequence of state constants,
and $u, v \in \SystemTClass$.
\begin{enumerate}

\item
  $u$ converges in $\{s_i\}_{i\in\NatSet}$ if $\exists i\in\NatSet.
\forall j\geq i.u[s_j]=u[s_{i}]$ in $\SystemTLearn$.\\

\item
$u$ converges if $u$ converges in every w.i. sequence of state constants.
\end{enumerate}
\end{definition}

We will make use of the following two theorems of chapter \ref{chapter-learningbasedrealizability}.

\begin{theorem}[Stability Theorem] \label{theorem-StabilityTheorem}
Assume $t \in \SystemTClass$ is a closed term of atomic type $A$ ($A\in\{\Bool,\Nat,\State\}$). Then $t$ is convergent.
\end{theorem}

\begin{theorem}[Zero Theorem]\label{Fixed Point Property}
Let $t:\State$ be a closed term of $\SystemTClass$ of state $\varnothing$ and $s$ any state constant. Define, by induction on $n$, a sequence $\{s_n\}_{n\in\NatSet}$ of state constants such that: $s_0=s$  and $s_{n+1}=s_n\Cup t[s_n]$. Then, there exists an $n$ such that $t[s_n]=\varnothing$.

\end{theorem}
We now define a language for Peano Arithmetic and then formulate a realizability relation between terms of $\SystemTClass$ and formulas of the language.

\begin{definition}[The language $\mathcal{L}$ of Peano Arithmetic] \label{definition-extendedarithmetic}  We define:

\begin{enumerate}

\item
The terms of $\mathcal{L}$ are all terms $t$ of G\"odel's system $\SystemTG$, such that $t:\Nat$ and $FV(t) \subseteq \{x_1^\Nat, \ldots, x_n^\Nat\}$ for some $x_1, \ldots, x_n$.\\

\item
The atomic formulas of $\mathcal{L}$ are all terms  $Qt_1\ldots t_n$ of G\"odel's system $\SystemTG$, for some $Q:\Nat^{n}\rightarrow \Bool$ {\em closed term of $\SystemTG$}, and some terms $t_1,\ldots,t_n$ of $\mathcal{L}$.\\

\item
The formulas of $\mathcal{L}$ are built from atomic formulas of $\mathcal{L}$ by the connectives $\lor,\land,\rightarrow \forall,\exists$ as usual.
\end{enumerate}

\end{definition}

We now define the types realizers as in chapter \ref{chapter-learningbasedrealizability} (we only use a different notation in order to avoid confusion in the rest of the chapter).

\begin{definition}[Types for realizers]
\label{definition-TypesForRealizers} For each
arithmetical formula $A$ we define a type $[ A]$ of $\SystemTG$ by
induction on $A$:
%\begin{enumerate}
%\item
$[ P(t_1,\ldots,t_n)]=\State$,
%\item
$[ A\wedge B]=[ A] \times [ B]$,
%\item
$[ A\vee B] = \Bool\times ([ A] \times [ B])$,
%\item
$[ A\rightarrow B] =[ A] \rightarrow [ B] $,
%\item
$[ \forall x A] =\Nat\rightarrow [ A] $,
%\item
$[ \exists x A] = \Nat\times [ A]$
%\end{enumerate}
\end{definition}We give the simplified notion of learning-based realizability we shall use in the following.
\begin{definition}[Learning-Based Realizability]
\label{lemma-IndexedRealizabilityAndRealizability}
Assume $s$ is a state constant, $t\in \SystemTClass$ is a closed term of state $\varnothing$, $A \in \mathcal{L} $ is a closed formula, and $t: [A]$. Let $\vec{t} = t_1, \ldots, t_n : \Nat$.

\begin{enumerate}
\item
$t\Vvdash_s P(\vec{t})$ if and only if $t[s]  = \varnothing$ in $\SystemTLearn$ implies
$P(\vec{t})={\True}$\\

\item
$t\Vvdash_s{A\wedge B}$ if and only if $\pi_0t \Vvdash_s{A}$ and $\pi_1t\Vvdash_s{B}$\\

\item
$t\Vvdash_s {A\vee B}$  if and only if either $\proj_0t[{s}]={\True}$ in $\SystemTLearn$ and $\proj_1t\Vvdash_s A$, or $\proj_0t[{s}]={\False}$ in $\SystemTLearn$ and $\proj_2t\Vvdash_s B$\\

\item
$t\Vvdash_s {A\rightarrow B}$ if and only if for all $u$, if $u\Vvdash_s{A}$,
then $tu\Vvdash_s{B}$\\

\item
$t\Vvdash_s {\forall x A}$ if and only if for all numerals $n$,
$t{n}\Vvdash_s A[{n}/x]$\\
\item

$t\Vvdash_s \exists x A$ if and only for some numeral $n$, $\pi_0t[{s}]= {n}$ in $\SystemTLearn$ and $\pi_1t \Vvdash_s A[{n}/x]$\\
\end{enumerate}
We define $t \Vvdash A$ if and only if $t\Vvdash_s A$ for all state constants $s$.
\end{definition}

For the soundness result we shall need this theorem of chapter \ref{chapter-learningbasedrealizability}.

\begin{theorem}\label{Realizability Theorem}
If $A$ is a closed formula of $\mathcal{L}$ provable in $\HA + \EM_1$, then there
exists $t\in \SystemTClass$ such that $t\Vvdash A$.
\end{theorem}

\section{Games, Learning and Realizability}\label{gamesrealizability}
\label{section-GamesLearningandRealizability}

In this section, we define the abstract notion of game, its 1-Backtracking version and Tarski games. We also prove our main theorem, connecting learning based realizability and 1-Backtracking Tarski games. 

\begin{definition}[Games]
\label{definition-Games} We define:
\begin{enumerate}

 \item
 A \emph{game} $G$ between two players is a quadruple \[(V,E_1,E_2, W)\]
where $V$ is a set, $E_1,E_2$ are  subsets of $V\times V$ such that
$Dom(E_1)\cap Dom(E_2)=\emptyset$, where $Dom(E_i)$ is the domain of $E_i$,  and $W$ is a set of sequences,
possibly infinite, of elements of $V$.  The elements of $V$ are
called \emph{positions} of the game; $E_1$, $E_2$ are the transition
relations respectively for player one and player two:
$(v_1,v_2)\in E_i$ means that player $i$ can legally move from the
position $v_1$ to the position $v_2$.\\

 \item We define a \emph{play} to be a walk, possibly infinite, in the
graph $(V,E_1\cup E_2)$, i.e. a sequence, possibly void,  $v_1::v_2::\ldots:: v_n::\ldots $ of elements of $V$  such that $(v_i, v_{i+1})\in E_1\cup E_2$ for every $i$. A play of the form $v_1:: v_2:: \ldots:: v_n::\ldots $ is said to \emph{start from}  $v_1$. A play is said to be
\emph{complete} if it is either infinite or is equal to $v_1::\ldots:: v_n$ and $v_n\notin Dom(E_1\cup E_2)$. $W$ is required to be a set of
complete plays. If $p$ is a complete play and $p\in  W$, 
%or if $p$is finite and its last element does belong to $Dom(E_2)$, 
we say that player one wins in $p$. If $p$ is a complete play and $p\notin
W$, 
%or if $p$ is finite and its last element does belong to
%$Dom(E_1)$, 
we say that player two wins in $p$.\\

 \item
 Let $P_G$ be the set of finite plays. Consider a function $f:
P_G\rightarrow V$. A play $v_1::\ldots:: v_n::\ldots$ is said to be
$f$-correct if $f(v_1::\ldots:: v_i)=v_{i+1}$ for every $i$ such that
$(v_i,v_{i+1})\in E_1$. $f$ is said to be a \emph{strategy} for player $i$ if for every play $p=v_1::\ldots :: v_n$ such that $v_n\in Dom(E_i)$, $v_1::\ldots ::v_n::f(p)$ is a play.\\

 \item
 A \emph{winning strategy} from position $v$ for player one is a strategy
$\omega: P_G\rightarrow V$ such that every complete
$\omega$-correct play $v::v_1::\ldots :: v_n::\ldots $  belongs to $W$.
 \end{enumerate}
 \end{definition}
 \textbf{Notation (Concatenation of Sequences).}  If for $i\in \NatSet, i=1,\ldots, n$  we have that $p_i=(p_i)_{0}:: \ldots :: (p_i)_{n_i}$ is a finite sequence of elements of length $n_i$,  with $p_1::\ldots ::p_n$ we denote the sequence \[(p_1)_0::\ldots ::(p_1)_{n_1}:: \ldots :: (p_k)_0::\ldots :: (p_k)_{n_k}\] where $(p_i)_j$ denotes the $j$-th element of the sequence $p_i$. \\\\
 Suppose that $a_1::a_2::\ldots :: a_n$ is a play of a game $G$,
representing, for some reason, a bad situation for player one (for
example, in the game of chess, $a_n$ might be a configuration of
the
chessboard in which player one has just lost his queen). Then,
learnt the lesson, player one might wish to erase some of his moves
and come back to the time the play was just, say, $a_1,a_2$ and
choose, say, $b_1$ in place of $a_3$; in other words, player one
might wish to \emph{backtrack}. Then, the game might go on as
$a_1 :: a_2 ::b_1::\ldots :: b_m$ and, once again, player one might want to
backtrack to, say, $a_1::a_2::b_1::\ldots :: b_i$, with $i< m$, and so
on... As there is no learning without remembering, player one
must keep in mind  the errors made during the play. This is the
idea
of 1-Backtracking games  (for more motivations, we refer the reader to \cite{BerCoq}) and here is our definition. %\\

\begin{definition}[1-Backtracking Games]
\label{definition-1BacktrackingGames} Let
$G=(V,E_1,E_2,W)$ be a game.

\begin{enumerate}

\item
We define $\oneback(G)$ as the game $(P_G, E_1',E_2', W')$, where:\\
\item $P_G$ is the set of finite plays of $G$\\

\item  \[E_2':=\{(p::a,\ p::a::b)\ |\   p\in P_G, p::a\in P_G,
(a,b)\in
E_2 \}\] and \[E_1':=\{(p::a,\ p::a::b)\ |\  p\in P_G, p::a\in P_G,  (a,b)\in
E_1\}\ \cup\] \[\{(p::a::q,\ p::a)\ |\ p, q\in P_G, p::a::q\in P_G, a\in Dom(E_1)\]\[
(q=q'::d\Rightarrow d\notin Dom(E_2)), p::a::q\notin W \};\]
\item $W'$ is  the set of finite complete plays $p_1::\ldots :: p_n$ of
$(P_G, E_1', E_2')$ such that $p_n\in W$.
\end{enumerate}
\end{definition}
\textbf{Note.} The pair $(p::a::q,\ p::a)$ in the definition above of $E_2'$ codifies a {\em backtracking move} by player one (and we point out that $q$ might be the empty sequence).\\

\textbf{Remark.} Differently from \cite{BerCoq}, in which both players are allowed to backtrack, we only consider the case in which only player one is supposed do that (as in \cite{Hayashi1}). It is not that our results would not hold: we claim that the proofs in this paper would work just as fine for the definition of 1-Backtracking Tarski games given in \cite{BerCoq}. However, as noted in \cite{BerCoq}, any player-one recursive winning strategy in our version of the game can be effectively transformed into a winning strategy for player one in the other version the game. Hence, adding backtracking for the second player does not increase the computational challenge for player one.
 Moreover, the notion of winner of the game given in \cite{BerCoq} is strictly non constructive and games played by player one with the correct winning strategy may even not terminate. Whereas, with our definition, we can formulate our main theorem as a program termination result: whatever the strategy chosen by player two, the game terminates with the win of player one. This is also the spirit of realizability and hence of this paper: the constructive information must be computed in a finite amount of time, not in the limit. \\

In the well known Tarski games, there are two players and a formula
on the board. The second player - usually called Abelard - tries to
show that the formula is false, while the first player - usually
called Eloise - tries to show that it is true. Let us see the
definition.%\\

\begin{definition}[Tarski Games]
\label{definition-TarskiGame} Let $A$ be a closed
implication and negation free arithmetical formula of $\mathcal{L}$. We define the
Tarski
game for $A$ as the game $T_A=(V, E_1, E_2, W)$, where:
\begin{enumerate}

\item
$V$ is the set of all subformula occurrences of $A$; that is, $V$ is the smallest set of formulas such that, if either $A\lor B$ or $A\land B$ belongs to $V$, then $A,B\in V$; if either $\forall x A(x)$ or $\exists x A(x)$ belongs to $V$, then $A(n)\in V$ for all numerals $n$. \\

\item
$E_1$ is the set of pairs $(A_1,A_2)\in V\times V$ such that   $A_1=\exists x
A(x)$  and $A_2=A(n)$,  or $A_1=A\lor B$ and either $A_2=A$ or
$A_2=B$;\\

\item
$ E_2$ is the set of pairs $(A_1,A_2)\in V\times V$ such that $A_1=\forall x
A(x)$  and $A_2=A(n)$,  or $A_1=A\land B$ and $A_2=A$ or
$A_2=B$;\\

\item
$W$ is the set of finite complete plays $A_1::\ldots ::A_n$ such that
$A_n=\True$.
\end{enumerate}
\end{definition}
\textbf{Note.} We stress that Tarski games are defined only for implication-and-negation-free arithmetical formulas. Indeed, $\oneback(T_A)$, when $A$ contains implications, would be much more involved and less intuitive (for a definition of Tarski games for every arithmetical formula see for example Lorenzen's \cite{Lorenzen}).\\

What we want to show is that if $t\Vvdash A$,
then $t$ gives to player one a recursive winning strategy in
$\oneback(T_A)$. The idea of the proof is the following. Suppose we play as player one. Our strategy is relativized to a knowledge state and we start
the game by fixing the actual state of knowledge as $\varnothing$.
Then we play in the same way as we would do in the Tarski game. For
example,  if there is $\forall x A(x)$ on the board and
$A(n)$ is chosen by player two, we recursively play the
strategy given by $tn$; if there is $\exists x A(x)$ on the
board, we calculate $\pi_0t[\varnothing]=n$ and play
$A(n)$ and recursively the strategy given by $\pi_1t$. If there is $A\lor B$ on the board, we calculate $\proj_0t[\varnothing]$, and according as to whether it equals $\True$ or $\False$, we play the strategy recursively given by $\proj_1t$ or $\proj_2t$.
If there is an atomic formula on the board, if it is true, we win; otherwise we extend the current state with the state $\varnothing \Cup t[\varnothing]$, we backtrack and play with respect to the new state of knowledge and trying to keep as close as possible to the previous game.
Eventually, we will reach a state large enough to enable our
realizer to give always correct answers and we will win. Let us consider first an example and then the formal definition of the winning strategy for Eloise.\\

\textbf{Example ($\EM_1$)}. Given a predicate $P$ of $\SystemT$, and its boolean negation predicate $\neg P$ (which is representable in $\SystemT$), the realizer $E_P$ of  \[\EM_1:=\forall x.\ \exists y\
P(x, y)\vee \forall y \neg P(x,y)\]
  is defined as \[\lambda \alpha^{\Nat}
\langle  \Chi_P\alpha,\ \langle
\Phic_P{\alpha},\
\varnothing \rangle ,\ \lambda m^{\Nat}\
\Add_P {\alpha}m\rangle   \]
We now compute a winning strategy for Eloise in the 1-Backtracking game associated to $\EM_1$. According to the rules of the game $\oneback(T_{\EM_1})$, Abelard is the first to move and, for some numeral $n$, chooses the formula

\[\exists y\ P(n, y)\vee \forall y \neg P(n,y)\]
Now is the turn of Eloise 
and she plays the strategy given by the term

\[\langle  \Chi_Pn,\ \langle
\Phic_Pn,\
\varnothing \rangle ,\ \lambda m^{\Nat}\
\Add_P nm\rangle   \]
Hence, she computes $\Chi_Pn[\varnothing]=\chi_P\varnothing n=\False$ (by definition \ref{definition-EquationalTheoryL1}), so she plays the formula 
\[\forall y \neg P(n,y)\]
and Abelard chooses $m$ and plays 

\[\neg P(n,m)\]
If $\neg P(n,m)=\True$, Eloise wins. Otherwise, she plays the strategy given by 
\[ (\lambda m^{\Nat}\ \Add_P nm )m[\varnothing]= \add_P \varnothing n m=\{\langle P,n,m\rangle\}\]
So, the new knowledge state is now $\{\langle P,n,m\rangle\}$ and she backtracks to the formula

\[\exists y\ P(n, y)\vee \forall y \neg P(n,y)\]
Now, by definition \ref{definition-EquationalTheoryL1}, $\Chi_Pn[\{\langle P,n,m\rangle\}]=\True$ and she plays the formula

\[\exists y\ P(n, y)\]
calculates the term \[\pi_0\langle
\Phic_Pn,\
{\varnothing} \rangle[\{\langle P,n,m\rangle\}]=\varphi_P\{\langle P,n,m\rangle\}n=m\]
plays $P(n,m)$ and wins.\\\\
\textbf{Notation.} In the following, we shall denote with upper case letters $A, B,C$
closed arithmetical formulas, with lower case letters $p,q,r$
plays of $T_A$ and with upper case letters $P,Q,R$  plays of
$\oneback(T_A)$ (and all those letters may be indexed by numbers). To avoid confusion with the plays of $T_A$, plays of 1Back($T_A$) will be denoted as $p_1,\ldots, p_n$ rather than $p_1::\ldots :: p_n$. Moreover, if $P=q_1,\ldots, q_m$, then $P, p_1,\ldots, p_n$ will denote the sequence $q_1, \ldots, q_m, p_1,\ldots p_n $.\\

We now define, given a play $p$ of $T_A$, a term $\rho(p)$, which we call ``the realizer associated to $p$'' and which represents the term that should be consulted by Eloise in a position $Q,p$ of the game $\oneback(T_A)$.

 \begin{definition}\label{adaptrealizer}
 Fix $u$ such that $u\Vvdash A$. Let  $p$ be a finite play of
$T_A$ starting with $A$. We define by induction on the length of
$p$ a term $\rho(p)\in \SystemTClass$ (read as `the realizer associated
to $p$')  in the following way: \begin{enumerate} 
\item If $p=A$, then
$\rho(p)=u$.\\

 \item If $p=(q:: \exists x B(x):: B(n))$ and
$\rho(q:: \exists x B(x))=t$, then $\rho(p)=\pi_1t$. \\
\item If $p=(q::
\forall x B(x):: B(n))$ and $\rho(q:: \forall x B(x))=t$,
then $\rho(p)=tn$. \\
\item  If $p=(q::  B_0\land B_1:: B_i)$
and $\rho(q:: B_0\land  B_1)=t$, then $\rho(p)=\pi_it$.\\ 
\item If
$p=(q::  B_1\lor B_2:: B_i)$ and $\rho(q:: B_1\lor  B_2)=t$, then
$\rho(p)=\proj_it$.\end{enumerate}
 Given a play $P=Q, q::B$ of $\oneback(T_A)$, we set
$\rho(P)=\rho(q::B)$.

\end{definition}

A play $P$ of $\oneback(T_A)$ may involve a number of backtracking moves by Eloise. In the winning strategy we are going to define, each time Eloise backtracks, she must extend the current state of knowledge by means of a realizer. In the above definition, we explain how Eloise calculates the state associated to $P$.

\begin{definition}\label{adaptstate}
 Fix $u$ such that $u\Vvdash A$. Let $\rho$ be as in definition \ref{adaptrealizer} and $P$ be a finite play of $\oneback(T_A)$ starting with $A$. We
define by induction on the length of $P$ a state $\Sigma(P)$ (read
as  `the state associated to $P$') in the following way:
\begin{enumerate} 
\item If $P=A$, then $\Sigma(P)=\varnothing$.\\
\item If $P=(Q, p::B, p::B::C)$ and $\Sigma(Q, p::B)=s$, then
$\Sigma(P)=s$.\\
 \item If $P=(Q, p::B::q, p::B)$ and $\Sigma(Q,
p::B::q)=s$ and $\rho(Q, p::B::q)=t$, then if $t:\State$, then
$\Sigma(P)=s\Cup t[s]$, else $\Sigma(P)=s$.\end{enumerate}

\end{definition}

We are now in a position to define a winning strategy for Eloise. Given a play $Q, p$ of $\oneback(T_A)$, she computes the state associated to $Q, p$ and then calls the realizer associated to $p$, which returns to her the next move to be performed.

\begin{definition}[Winning strategy for $\oneback(T_A)$]\label{definition-winningstrategy}
Fix $u$ such that $u\Vvdash A$. Let $\rho$ and $\Sigma$ be respectively as in definitions \ref{adaptrealizer} and \ref{adaptstate}. We define a function $\omega$ from the set of finite plays of
$\oneback(T_A)$ to set of finite plays of $T_A$; $\omega$ is intended
to be a recursive winning strategy from $A$ for player one in $\oneback(T_A)$.
 \begin{enumerate} 
 \item If $\rho(P,q::\exists x B(x))=t$,
$\Sigma(P,q::\exists x B(x))=s$ and $(\pi_0t)[s]={n}$,
then \[\omega(P,q::\exists x B(x))= q::\exists x B(x)::
B({n})\]

 \item If $\rho(P, q::B\lor C)=t$ and $\Sigma(P, q::B\lor C)=s$,
then if $(\proj_0t)[s]=\True$ then \[\omega(P, q::B\lor C)=q::B\lor
C::B\] else \[\omega(P, q::B\lor C)= q::B\lor C:: C\]

 \item If $A_n$ is atomic, $A_n=\False$, $\rho(P, A_1::\cdots::
A_n)=t$ and $\Sigma(P, A_1::\cdots ::A_n)=s$, then \[\omega(P,
A_1::\cdots:: A_n)= A_1::\cdots:: A_i\] where $i$ is equal to the smallest
$j< n$ such that $\rho( A_1::\cdots:: A_j)=w$ and either
\[A_j=\exists x C(x)\land  A_{j+1}= C({n}) \land
(\pi_0w)[s\Cup t[s]]\neq {n}\] 
or   \[A_j=B_1\lor B_2\land
A_{j+1}= B_1 \land (\proj_0w)[s\Cup t[s]]=\False\] 
or \[A_j=B_1\lor B_2\land A_{j+1}= B_2 \land
(\proj_0w)[s\Cup t[s]]=\True\] If such $j$ does not exist, we set $i=n$.\\
\item
In the other cases, $\omega(P,q)=q$.
\end{enumerate}

\end{definition}

 \begin{lem} \label{preservationlemma}
 \label{lemma-Completenessof1Backtracking}Suppose $u\Vvdash A$ and $\rho,\Sigma,\omega$ as in definition \ref{definition-winningstrategy}. Let $Q$ be a finite
$\omega$-correct play of $\oneback(T_A)$ starting with $A$, $\rho(Q)=t$, $\Sigma(Q)=s$. If
$Q=Q',q'::B$, then $t\Vvdash_s B$. \end{lem}

\proof 
By a straightforward induction on the length of
$Q$.
 \begin{enumerate}

  \item If $Q=A$, then $t=\rho(Q)=u\Vvdash_s A$.\\
  
  \item If $Q=P,q::\exists x B(x),q::\exists x B(x)::
B({n})$, then let $t'=\rho(P,q::\exists x B(x))$. By
definition of $\Sigma$, $s=\Sigma(P,q::\exists x B(x))$. Since $Q$
is $\omega$-correct and $(q::\exists x B(x),q::\exists x B(x)::
B({n}))\in E_1$, we have $\omega(P,q::\exists x
B(x))=q::\exists x B(x):: B({n})$ and so
${n}=(\pi_0t')[s]$. Moreover, by definition of $\rho$,
$t=\pi_1t'$; by induction hypothesis, $t'\Vvdash_s \exists x B(x)$;
so, $t=\pi_1 t' \Vvdash_s B({n})$.\\

   \item If $Q=P,q::B\lor C,q::B\lor C:: B$, then let
$t'=\rho(P,q::B\lor C)$. By definition of $\Sigma$,
$s=\Sigma(P,q::B\lor C)$. Since $Q$ is $\omega$-correct and
$(q::B\lor C,q::B\lor C:: B)\in E_1$, we have $\omega(P,q::B\lor
C)=q::B\lor C:: B$ and so $(\proj_0t')[s]=\True$. Moreover, by definition
of $\rho$, $t=\proj_1t'$; by induction hypothesis, $t'\Vvdash_s
B\lor C$; so, $t=\proj_1t' \Vvdash_s B$.
The other case is analogous.\\

  \item
 If $Q=P,q::\forall x B(x), q::\forall x B(x)::B({n})$,
then let $t'=\rho(P,q::\forall x B(x))$. By definition of $\Sigma$, $s=\Sigma(P, q::\forall x B(x))$. By definition of $\rho$,
$t=t'{n} $; by induction hypothesis, $t'\Vvdash_s \forall x
B(x)$; hence, $t=t'n \Vvdash_s B({n})$.\\

  \item
 If $Q=P,q::B\land C, q::B\land C::B$, then let
$t'=\rho(P,q::B\land
C)$. By definition of $\Sigma$, $s=\Sigma (P, q::B\land C)$. By definition of $\rho$, $t=\pi_0t'$; by induction hypothesis,
$t'\Vvdash_s B\land C$; hence, $t=\pi_0 t' \Vvdash_s B$. The other case is analogous.\\

 \item
 If $Q=P, A_1::\cdots:: A_n, A_1::\cdots:: A_i$, $i<n$, $A_n$ atomic, then $A_1=A$. Furthermore, if $\Sigma(P,A_1::\cdots:: A_n)=s'$
and  $t'=\rho(P,A_1::\cdots::A_n)$, then $s=s'\Cup t'[s']$. Let
$t_j=\rho(A_1::\cdots :: A_j)$, for $j=1,\ldots, i$. We prove by
induction on $j$ that $t_j\Vvdash_s A_j$, and hence the thesis.
If $j=1$, then $t_1=\rho(A_1)=\rho(A)=u\Vvdash_s A=A_1$.\\ If $j>1$,
by induction hypothesis $t_{k}\Vvdash_s A_{k}$, for every $k<j$. If either $A_{j-1}=\forall x
C(x)$ or $A_{j-1}=C_0\wedge C_1$, then
either $t_j=t_{j-1}{n}$ and $A_j=C({n})$, or
$t_j=\pi_mt_{j-1}$ and $A_j=C_m$: in both cases, we have
$t_j\Vvdash_s A_j$, since $t_{j-1}\Vvdash_s A_{j-1}$. Therefore, by definition of $\omega$ and $i$ and the $\omega$-\emph{correctness} of $Q$, the remaining possibilities are that either
$A_{j-1}=\exists x C(x)$, $A_j=
C({n})$, $t_j=\pi_1t_{j-1}$, with $(\pi_0t_{j-1})[s]= {n}$; or
$A_{j-1}=C_1\lor C_2$, $A_j=C_m$, $t_j=\proj_m t_{j-1}$ and
$(\proj_0t_{j-1})[s]=\True$ if and only if $m=1$; in both cases,  we
have $t_j\Vvdash_s A_j$. \end{enumerate}
 \qed

\begin{theorem}[Soundness Theorem]\label{theorem-Soundness} Let $A$ be a closed negation and implication free arithmetical formula. Suppose that $u\Vvdash A$ and consider the game
$\oneback(T_A)$. Let $\omega$ be as in definition \ref{definition-winningstrategy}. Then $\omega$ is a recursive winning strategy from $A$ for player one. \end{theorem}

\proof  We begin by showing that there is no infinite
$\omega$-correct play.\\ Let $P=p_1,\ldots, p_n,\ldots$ be, for the
sake of contradiction, an infinite $\omega$-correct play, with
$p_1=A$. Let $A_1::\cdots :: A_k$ be the \emph{longest} play of $T_A$ such that there exists $j$ such that for every
$n\geq j$, $p_n$ is of the form $A_1::\cdots ::A_k::q_n$. $A_1::\cdots::A_k$ is well defined, because:  $p_n$ is of the form $A::q'_n$ for every $n$;  the length of $p_n$ is at most the degree of the formula $A$; the sequence of maximum length is unique because any two such sequences are one the prefix of the other, and therefore are equal. Moreover,  let $\{n_i\}_{i\in\NatSet}$ be the infinite increasing sequence
of all indexes $n_i$ such that $p_{n_i}$ is of the form $A_1::\cdots:: A_k::q_{n_i}$ and
$p_{n_i+1}=A_1::\cdots :: A_k$ (indeed, $\{n_i\}_{i\in\NatSet}$
must be infinite: if it were not so, then there would be an index
$j'$ such that for every $n\geq j'$, $p_n=A_1::\cdots
::A_k::A_{k+1}::q$, violating the assumption on the maximal length of
$A_1::\cdots::A_k$). $A_k$, if not atomic, is a disjunction or an
existential statement.\\ Let now $s_i=\Sigma(p_1,\ldots, p_{i})$
and
$t=\rho(A_1::\cdots:: A_k)$. For every $i$, $s_i\leq s_{i+1}$, by definition of $\Sigma$.
There
are three cases:\\

 1) $A_k=\exists x B(x)$. Then, by the Stability Theorem (Theorem \ref{theorem-StabilityTheorem}), there exists
$m$ such that for every $a$, if $n_a\geq m$, then $(\pi_0t)[s_{n_a}]=(\pi_0t)[s_m]$. Let \[h:=(\pi_0t)[s_{n_a}\Cup t_1[s_{n_a}]]=(\pi_0t)[s_{n_a+1}]\]
 where $t_1=\rho(p_1,\ldots, p_{n_a})$. So let $a$ be such that $n_a\geq m$; then  \[p_{n_a+2}=\omega(p_1,\ldots,
p_{n_a+1})=\omega(p_1,\ldots, A_1::\cdots :: A_k)=A_1::\cdots ::A_k::B(h)\] Moreover, by hypothesis, and since $p_{n_a+1}=A_1::\cdots::A_k$, we have
 \[p_{n_{(a+1)}}=A_1::\cdots ::A_k::q_{n_{(a+1)}}= A_1::\cdots ::A_k::B(h)::q'\] 
 for some $q'$ and $p_{n_{(a+1)}+1}=A_1::\cdots ::A_k$: contradiction, since \[h=(\pi_0t)[s_{n_{a}+1}]=(\pi_0t)[s_{n_{(a+1)}+1}]=(\pi_0t)[s_{n_{(a+1)}}\Cup t_2[s_{n_{(a+1)}}]]\]
where  $t_2=\rho(p_1,\ldots, p_{n_{(a+1)}})$, whilst $h\neq (\pi_0t)[s_{n_{(a+1)}}\Cup t_2[s_{n_{(a+1)}}]]$ should hold, by definition of $\omega$ (point (3)).\\

2) $A_k=B\lor C$. This case is totally analogous to the preceding.\\

 3) $A_k$ is atomic. Then, for every $n\geq j$, $p_n=A_1::\cdots::
A_k$. So, for every $n\geq j$, $s_{n+1}=s_n\Cup t[s_n]$ and  hence, by Theorem \ref{Fixed Point Property} there
exists $m\geq j$ such that  $t[s_{m}]=\varnothing$. But $t
\Vvdash_{s_{m}} A_k$, by Lemma \ref{preservationlemma}; hence, $A_k$ must equal $\True$, and so it is
impossible that $(p_m, p_{m+1})=(A_1::\cdots ::A_k, A_1::\cdots
::A_k)\in E_1'$: contradiction.\\
Let now $p=p_1,\ldots, p_n$ be a
complete finite $\omega$-correct play. $p_n$ must equal 
$B_1::\cdots::B_k$, with $B_k$ atomic and $B_k=\True$: otherwise,
$p$ wouldn't be complete, since player one would lose the play $p_n$ in $T_A$ and hence would be allowed to backtrack by definition \ref{definition-1BacktrackingGames}.\\
\qed

\section{Examples}
\label{examples}
In this section we include two natural deduction classical proofs of two simple combinatorial statements, using only Excluded Middle for semi-decidable statements, then we extract a constructive content using our realizability semantics. For each of them, we interpret the program we shall extract using our game interpretation of the learning based realizability semantics.
\comment{
\textbf{Example ($\Sigma^0_1$ and $\Pi^0_2$ formulas).} Suppose
$u\Vdash \exists x P(x)$, with $P(x)$ atomic. Let $n$ be the
smallest natural number such that $s:=(\pi_1u)^n[\varnothing]$ is a
fixed point of $\pi_1u$. We have $u\Vdash_s \exists x P(x)$ and so
$\pi_1u\Vdash_s P(\overline{m})$, where $\overline{m}=(\pi_0u)[s]$.
Since $(\pi_1u)[s]=s$, $P(\overline{m})$ must be true. Hence, here
is the algorithm (in pseudo code) to find the witness for $\exists
x
P(x)$:\\ \\$s:=\varnothing$;\\ repeat $s:=(\pi_1u)[s]$ until
$(\pi_1u)[s]=s$;\\
return $(\pi_0u)[s];$\\ Suppose now $t\Vdash \forall x\exists y
P(x,y)$. Then, given $n\in\NatSet$, $t\overline{n}\Vdash\exists
yP(\overline{n},y)$. So we can apply the algorithm for $\Sigma^0_1$
formulas. Hence, the extracted algorithm is the following:\\
$u:=t\overline{n}$;\\ $s:=\varnothing$;\\ repeat $s:=(\pi_1u)[s]$
until
$(\pi_1u)[s]=s$;\\ return $(\pi_0u)[s];$\\
}

\subsection{Minimum Principle for Functions over Natural Numbers.} The
minimum principle states that every function $f$ over natural
numbers has a minimum value, i.e. there exists a $f(n)\in \NatSet$
such that for every $m\in\NatSet$ $ f(m)\geq f(n)$. We can prove
this principle in $\HA + \EM_1$, for any $f$ in the language. We assume $P(y,x)\equiv f(x)<y$, but, in order to
enhance readability, we will write $f(x)<y$ rather than the obscure
$P(y,x)$. We define:\\ $Lessef(n):= \exists \alpha f(\alpha)\leq
n$\\ $Lessf(n):=\exists \alpha f(\alpha)<n$\\ $Notlessf(n):=
\forall \alpha f(\alpha)\geq n$\\ Then we formulate - in equivalent form - the
minimum principle as: \[Hasminf:=\exists y.\ Notlessf(y)\wedge
Lessef(y)\] The informal argument goes as follows. We prove by induction on $n$ that for every $k$, if $f(k)\leq n$, then $f$ has minimum value. If $n=0$, we just observe that $f(k)\leq 0$, implies $f(k)$ is the minimum value of $f$. Suppose now $n>0$. If $Notlessf(f(k))$ holds true, we are done, $f(k)$ is the minimum of $f$. Otherwise, $Lessf(f(k))$ holds, and hence $f(\alpha)<f(k)\leq n$ for some $\alpha$ given by an oracle. Hence $f(\alpha)\leq f(k)-1\leq n-1$ and we conclude that $f$ has a minimum value by induction hypothesis. 

 Now we give the formal proofs, which are natural deduction trees, decorated with terms of $\SystemTClass$, as formalized in chapter \ref{chapter-learningbasedrealizability}. We first prove
that
$\forall n.\ (Lessef(n)\rightarrow Hasminf)\rightarrow
(Lessef(\suc(n))\rightarrow Hasminf)$ holds.

\def\proofSkipAmount{\vskip-2ex plus.1ex minus.1ex}
\begin{prooftree}
\small
\AxiomC{$E_P: \forall n.\ Notlessf(\suc(n))\lor Lessf(\suc(n))$}
\UnaryInfC{$E_Pn: Notlessf(\suc(n))\lor Lessf(\suc(n))$}

                   \AxiomC{$[Notlessf(\suc(n))]$}
                   \noLine
                   \UnaryInfC{$T_1$}
                   \noLine
                   \UnaryInfC{$Hasminf$}

                                        \AxiomC{$[Lessf(\suc(n))]$}
                                        \noLine
                                       \UnaryInfC{$T_2$}
                                       \noLine
                                        \UnaryInfC{$Hasminf$}

\TrinaryInfC{$D: Hasminf$}
\UnaryInfC{$\lambda w_2 D: Lessef(\suc(n))\rightarrow Hasminf$}
\UnaryInfC{$\lambda w_1\lambda w_2 D: (Lessef(n)\rightarrow
Hasminf)\rightarrow  (Lessef(\suc(n))\rightarrow Hasminf)$}
\UnaryInfC{$\lambda n \lambda w_1\lambda w_2 D: \forall n
(Lessef(n)\rightarrow Hasminf)\rightarrow  (Lessef(\suc(n)\rightarrow
Hasminf)$}
\end{prooftree}where for lack of space the term $D$ is defined later, $T_1$ is the tree

\def\proofSkipAmount{\vskip0ex plus.1ex minus.1ex}\begin{prooftree}
\small
\AxiomC{$v_1: Notlessf(\suc(n))$}
                   \AxiomC{$w_2: Lessef(\suc(n))$}
                   \BinaryInfC{$\langle v_1,w_2\rangle :
Notlessf(\suc(n))\land
Lessef(\suc(n))$}
                   \UnaryInfC{$\langle \suc(n),\langle v_1,w_2\rangle
\rangle : Hasminf$}
\end{prooftree}
and $T_2$ is the tree

\def\proofSkipAmount{\vskip-3ex plus.1ex minus.1ex}
\begin{prooftree}
\small
\AxiomC{$v_2: [Lessf(\suc(n))]$}
                                        \AxiomC{$w_1:
[Lessef(n)\rightarrow
Hasminf]$}

                                        \AxiomC{$[x_2: f(z)< \suc(n)]$}
                                        \UnaryInfC{$x_2: f(z)\leq
n$}
                                         \UnaryInfC{$\langle
z,x_2\rangle :
Lessef(n)$}
                                        \BinaryInfC{$w_1\langle
z,x_2\rangle :
Hasminf$}
                                        \BinaryInfC{$w_1\langle
\pi_0v_2,\pi_1v_2\rangle :
Hasminf$}
\end{prooftree}
We prove now that $Lessef(0)\rightarrow Hasminf$
\def\proofSkipAmount{\vskip-4ex plus.1ex minus.1ex}
 \begin{prooftree}
\small
 \AxiomC{$w: [Lessef(0)]$}

 \AxiomC{$x_1: [f(z)\leq 0]$}
 \UnaryInfC{$x_1: f(z)=0$}

 \UnaryInfC{$x_1: f(\alpha)\geq f(z)$}
 \UnaryInfC{$\lambda \alpha x_1: Notlessf(f(z))$}
                \AxiomC{$\varnothing: f(z)\leq f(z)$}
                \UnaryInfC{$\langle z,\varnothing\rangle : Lessef(f(z))$}
 \BinaryInfC{$\langle \lambda\alpha x_1,\langle z,\varnothing\rangle
\rangle : Notlessf(f(z))\land
Lessef(f(z))$}
 \UnaryInfC{$\langle f(z),\langle \lambda\alpha x_1,\langle
z,\varnothing\rangle \rangle \rangle : Hasminf$}
 \BinaryInfC{$\langle f(\pi_0w),\langle \lambda\alpha \pi_1
w,\langle \pi_0w,\varnothing\rangle \rangle \rangle :
Hasminf$}
 \UnaryInfC{$F:=\lambda w\langle f(\pi_0w),\langle \lambda\alpha
\pi_1
w,\langle \pi_0w,\varnothing\rangle \rangle \rangle : Lessef(0)\rightarrow
Hasminf$}
 \end{prooftree}
Therefore we can conclude with the induction rule that \[\lambda
\alpha^\Nat\ \rec F(\lambda n\lambda w_1\lambda w_2 D)\alpha: \forall x.
Lessef(x)\rightarrow Hasminf\] And now the thesis:

\def\proofSkipAmount{\vskip-1ex plus-1ex minus.1ex} \begin{prooftree}
\small
\AxiomC{$\varnothing: f(0)\leq f(0)$}
\UnaryInfC{$\langle 0, \varnothing\rangle : Lessef (f(0))$}
                  \AxiomC{$\lambda
\alpha^\Nat\ {\rec}F(\lambda n\lambda w_1\lambda w_2 D)\alpha: \forall x. Lessef(x)\rightarrow Hasminf$}
                  \UnaryInfC{$\rec F(\lambda n\lambda w_1\lambda w_2
D)f(0): Lessef(f(0))\rightarrow Hasminf$}
\BinaryInfC{$M:=\rec F(\lambda n\lambda w_1\lambda w_2 D)f(0)\langle
0,\varnothing\rangle :
Hasminf$}
\end{prooftree}

 Let us now  define \[D:={\tt if}\
\Chi_P \suc(n)\  {\tt then}\  w_1\langle \Phic_P \suc(n), \varnothing
\rangle \
{\tt else}\ \langle \suc(n),\langle \lambda \beta\ (\Add_P) 
\suc(n)\beta,w_2\rangle \rangle
\]
 Let $s$ be a state and let us consider $M$, the realizer of $Hasminf$, in the base case of the recursion and after in its general form during the computation: ${\rec}F(\lambda n\lambda
w_1\lambda w_2 D)f(0)\langle m,\varnothing\rangle
[s]$. If $f(0)=0$, \[M[s]={\rec}F(\lambda n\lambda w_1\lambda
w_2 D)f(0)\langle{0},\varnothing\rangle
[s]=\] \[=F\langle 0,\varnothing\rangle
=\langle f(0),\langle \lambda
\alpha \varnothing, \langle 0,\varnothing\rangle \rangle  \]
If $f(0)=\suc(n)$, we have two other cases. If $\chi_Ps\suc({n})=\True$, then \[{\rec}F(\lambda n\lambda w_1\lambda
w_2 D)\suc({n})\langle{m},\varnothing\rangle
[s]=\]\[=(\lambda n\lambda w_1\lambda w_2 D){n}({\rec}F(\lambda
n\lambda w_1\lambda w_2 D){n})\langle
{m},\varnothing\rangle [s]= \]\[={\rec}F(\lambda n\lambda
w_1\lambda w_2 D){n}\langle \Phic_P
(\suc({n})),\varnothing\rangle [s]\]
If $\chi_Ps\suc({n})=\False$, then \[{\rec}F(\lambda n\lambda
w_1\lambda w_2 D)\suc({n})\langle {m},\varnothing\rangle
[s]=\]\[=(\lambda n\lambda w_1\lambda w_2 D){n}({\rec}F(\lambda
n\lambda w_1\lambda w_2
D){n})\langle {m},\varnothing\rangle [s]= \]\[=
\langle \suc({n}),\langle \lambda \beta\ (\add_P)s \suc(n)\beta,\langle {m},\varnothing\rangle \rangle
\rangle \]
In the first case, the minimum value of $f$ has been found. In the second case, the operator ${\rec}$, starting from $\suc(n)$, recursively calls itself on $n$; in the third case, it reduces to its normal form. From these equations, we easily deduce the behavior of
the realizer of $Hasminf$.   In a pseudo imperative programming
language, for the witness of $Hasminf$ we would write:\\\\
$n:=f(0);$\\
{\tt while}  $(\chi_Psn=\True, i.e.\ \exists m \ such\ that\  f({m})<n\in s)$\\
{\tt do} $n:=n-1;$\\
{\tt return} $n;$ \\\\
Hence, when $f(0)>0$, we have, for some numeral $k$ 
\[M[s]=\langle k,\langle \lambda \beta\ (\add_P)s k\beta,\langle
\varphi_Psk,\varnothing\rangle \rangle \rangle \]
It is clear that $k$  is the
minimum value of $f$, according to the partial information provided by
$s$
about $f$, and that $f(\varphi_Psk)\leq k$. If $s$ is sufficiently complete, then $k$
is the true minimum of $f$.\\
The normal form of the realizer $M$ of $Hasminf$ is so simple that we can immediately extract the winning strategy $\omega$ for the 1-Backtraking version of the Tarski game for $Hasminf$.
Suppose the current state of the game is $s$. If $f(0)=0$, Eloise chooses the formula
\[Notlessf(0)\land Lessef (0)\]
and wins. If $f(0)>0$, she chooses, for $k$ defined as above, 

\[Notlessf(k)\land Lessef (k)=\forall \alpha\ f(\alpha) \geq k \land \exists \alpha\ f(\alpha)\leq k \]
If Abelard chooses $\exists \alpha\ f(\alpha)\leq k$, she wins, because she responds with $f(\varphi_Psk)\leq k$, which holds. Suppose hence Abelard chooses \[\forall \alpha\ f(\alpha) \geq k\] and then $f(\beta)\geq k$. If it holds, Eloise wins. Otherwise, she adds to the current state $s$
\[(\lambda \beta\ (\add_P)s k\beta)\beta= (\add_P)s k\beta=\{f(\beta)<k\}\]
 and backtracks to $Hasminf$ and then plays again. This time, she chooses 
 
 \[Notlessf(f(\beta))\land Lessef (f(\beta))\]
 (using $f(\beta)$, which was Abelard's counterexample to the minimality of $k$ and is smaller than her previous choice for the minimum value). After at most $f(0)$ backtrackings, she wins. \\

 \subsection{Coquand's Example.} We investigate now an example -  due to
Coquand - in our framework of realizability. We want to prove that
for every function over natural numbers and  for every
$a\in\NatSet$ there exists $x\in\NatSet$ such that $f(x)\leq
f(x+a)$. Thanks to the minimum principle, we can give a very easy
classical proof:
\comment{
\begin{prooftree}
\AxiomC{$h: Hasminf$}
            \AxiomC{$w: Notlessf(\mu)\land Lessef(\mu)$}
            \UnaryInfC{$\pi_1w: Lessef(\mu)$}
                  \AxiomC{$w: Notlessf(\mu)\land Lessef(\mu)$}
                  \UnaryInfC{$\pi_0 w: Notlessf (\mu)$}
                  \UnaryInfC{$\pi_0w(z+a): f(z+a)\geq \mu$}

                                   \AxiomC{$v: f(z)\leq \mu$}
                  \BinaryInfC{$\pi_0w(z+a)\Cup v: f(z)\leq f(z+a)$}
                  \UnaryInfC{$\langle z,\pi_0w(z+a)\Cup v\rangle :
\exists x
f(x)\leq f(x+a)$}
                  \UnaryInfC{$\lambda a \langle z,\pi_0w(z+a)\Cup
v\rangle :
\forall a \exists x f(x)\leq f(x+a)$}
                                    \BinaryInfC{$\lambda a
\langle \pi_0\pi_1w,\pi_0w(\pi_0\pi_1w+a)\Cup
\pi_1\pi_1w\rangle : \forall a \exists
x f(x)\leq f(x+a)$}
\BinaryInfC{$\lambda a
\langle \pi_0\pi_1\pi_1h,\pi_0\pi_1h(\pi_0\pi_1\pi_1h+a)\Cup
\pi_1\pi_1\pi_1h\rangle : \forall a \exists x f(x)\leq f(x+a)$}
\end{prooftree}}

\def\proofSkipAmount{\vskip-1ex plus.1ex minus.1ex}
%\hbox{ \kern-0.2cm
%\scriptsize
\begin{prooftree}
\scriptsize\AxiomC{$Hasminf$}
            \AxiomC{$ [Notlessf(\mu)\land Lessef(\mu)]$}
            \UnaryInfC{$ Lessef(\mu)$}
                  \AxiomC{$ [Notlessf(\mu)\land Lessef(\mu)]$}
                  \UnaryInfC{$ Notlessf (\mu)$}
                  \UnaryInfC{$ f(z+a)\geq \mu$}

                                   \AxiomC{$ [f(z)\leq \mu]$}
                  \BinaryInfC{$ f(z)\leq f(z+a)$}
                  \UnaryInfC{$\exists x
f(x)\leq f(x+a)$}
                  \UnaryInfC{$\forall a \exists x f(x)\leq f(x+a)$}
                                    \BinaryInfC{$ \forall a \exists
x f(x)\leq f(x+a)$}
\BinaryInfC{$\forall a \exists x f(x)\leq f(x+a)$}
\end{prooftree}
%\DisplayProof }
%\end{prooftree}
The extracted realizer is \[\lambda a
\langle \pi_0\pi_1\pi_1M,\pi_0\pi_1M(\pi_0\pi_1\pi_1M+a)\Cup
\pi_1\pi_1\pi_1h\rangle\] where $M$ is the realizer of $Hasminf$.
$m:=\pi_0\pi_1\pi_1M[s]$ is a point the purported minimum value $\mu:=\pi_0M$ of $f$ is attained at, accordingly to the information in the state $s$ (i.e. $f(m)\leq \mu$). So, if Abelard chooses

\[\exists x\ f(x)\leq f(x+a)\]
Eloise chooses

\[f(m)\leq f(m+a)\]
 We have to
consider the term  \[U[s]:=\pi_0\pi_1M(\pi_0\pi_1\pi_1M+a)\Cup
\pi_1\pi_1\pi_1M[s]\] which updates the current state $s$.  Surely,
$\pi_1\pi_1\pi_1M[s]=\varnothing$.  $\pi_0\pi_1M[s]$ is equal either
to $\lambda \beta\ (\add_P)s \mu\beta$
or to $\lambda \alpha \varnothing$.
 So,
what does $U[s]$
actually do? We have:
\[U[s]=\pi_0\pi_1M(\pi_0\pi_1\pi_1M+a)[s]
=\pi_0\pi_1M(m+a)[s]\]
with either $\pi_0\pi_1M(m+a)[s]=\varnothing$ or \[\pi_0\pi_1M(m+a)[s]= \{f(m+a)<f(m)\}\]
So $U[s]$  tests if $f(m+a)< f(m)$; if it is not the
case, Eloise wins, otherwise she enlarges the state $s$, including the information
$f(m+a)<f(m)$ and backtracks to $\exists x f(x)\leq f(x+a)$. Starting from the state $\varnothing$,
after $k+1$ backtrackings, it will be reached a state $s'$, which will
be of the form $\{ f((k+1)a)<f(ka),\ldots, f(2a)<f(a), f(a)<f(0)\}$
and Eloise will play $f((k+1)a)\leq f((k+1)a+a)$. Hence, the extracted
algorithm for Eloise's witness is the following:\\\\
$n:=0$;\\ {\tt while} $f(n)>f(n+a)$\\ {\tt do} $n:=n+a$;\\ {\tt return} $n$;

\section{Total Recursive Learning-Based Realizability}\label{completeness}

The realizability notion introduced in definition \ref{lemma-IndexedRealizabilityAndRealizability} is very interesting from the constructive point of view. But precisely for that reason, the system $\SystemTClass$  fails to realize every formula for which an Eloise recursive winning strategy exists in its associated 1-Backtracking Tarski game, as the following theorem implies:

\begin{theorem}[Incompleteness of $\SystemTClass$]
There is a $\Pi_2^0$ arithmetical sentence $A$ such that Eloise has recursive winning strategy in $\oneback(T_A)$, but no term of system $\SystemTClass$ realizes $A$.\end{theorem}

\proof Take any total recursive function $f:\Nat\rightarrow \Nat$ not representable by any term of system $\SystemTG$ of type $\Nat\rightarrow \Nat$.  Let $n$ be the code of $f$ in the enumeration of Turing machines assumed by Kleene's primitive recursive predicate $Txyz$. Then the formula $A:=\forall y\exists z Tnyz$ asserts the totality of $f$ and hence it is true. Clearly, Eloise has a winning recursive strategy in $\oneback(T_A)$: for any $y$, she may backtrack until she finds an $m$ such that $Tnym$ holds. Suppose, along the way of contradiction, that for some $t$ of system $\SystemTClass$, $t\Vvdash A$. Then, as proven in chapter \ref{chapter-constructiveanalysislearning}, there exists a term $u:\Nat\rightarrow \Nat$ of system $\SystemTG$ such that, for every numeral $l$, $Tnl(ul)$ holds. From this, it easily follows that $f$ can be coded by a term of system $\SystemTG$, which is a contradiction. 

\qed

The only way around this issue, which is the purpose of this section, is to extend our notion of realizability and increase the computational power of our realizers, in order to be able to represent any partial recursive function and in particular every recursive strategy of 1-Backtracking Tarski games. So, we choose to add to our calculus a fixed point combinator $\mathsf{Y}$, such that for every term $u:A\rightarrow A$, $\mathsf{Y}u=u(\mathsf{Yu})$, getting the full power of $\mathcal{PCF}$ (see for example Gunter \cite{Gunter}).

\begin{definition}[Systems $\PRclass$ and $\PRlearn$ ]
We define $\PRclass$ and $\PRlearn$ to be, respectively, the extensions of $\SystemTClass$ and $\SystemTLearn$ obtained by adding for every type $A$ a constant $\mathsf{Y}_A$ of type $(A\rightarrow A)\rightarrow A$ and a new equality axiom $\mathsf{Y_A}u=u(\mathsf{Y_A}u)$ for every term $u: A\rightarrow A$.

\end{definition}

Since in $\PRclass$ there is a schema for unbounded iteration, properties like convergence do not hold anymore, for terms may even not have a normal form. So we have to {\em ask} our realizers to be convergent. Hence, for each type $A$ of $\PRclass$ we define a set $\|A\|$  of
terms $u: A$ which we call the set of {\em stable terms} of
type $A$. We define stable terms by lifting the notion of
convergence from atomic types to arrow and product types.

\begin{definition}[Convergence for $\PCFclass$]
\label{definition-Convergence2} Assume
that $\{s_i\}_{i\in\NatSet} $ is a w.i. sequence of state constants,
and $u, v \in \PRclass$.
\begin{enumerate}

\item
  $u$ converges in $\{s_i\}_{i\in\NatSet}$ if there exists a normal form $v$ such that $\exists i
\forall j\geq i.u[s_j]=v$ in $\PRlearn$.

\item
$u$ converges if $u$ converges in every w.i. sequence of state constants.
\end{enumerate}
\end{definition}

\begin{definition}[Stable Terms]
\label{definition-StableTerms} Let
$\{s_i\}_{i\in\NatSet}$ be a w.i. chain of states and $s \in
\StateSet$. Assume $A$ is a
type. We define a set $\|A\|$ of terms
$t \in\PRclass$ of type $A$, by induction on $A$.
\begin{enumerate}

\item
$\|\State\|=\{t:\State \ |\ t\mbox{ converges}\}$\\

\item
$\|\Nat\|=\{t:\Nat \ |\ t\mbox{ converges}\}$\\

\item
$\|\Bool\|=\{t:\Bool \ |\ t\mbox{ converges}\}$\\

\item
$\|A\times B\|=\{t:A \times B \ |\
\pi_0t\in\|A\|,\pi_1t\in\|B\|\}$\\

\item
$\|A\rightarrow B\|=\{t:A\rightarrow B \ |\ \forall u\in
\|A\|, tu\in\|B\|\}$\\
\end{enumerate}

If $t\in\|A\|$, we say that $t$ is a {\em stable} term of type
$A$.
\end{definition}

Now we extend the notion of realizability with respect to $\PRclass$ and $\PRlearn$.

\begin{definition}[Total Recursive Learning-Based Realizability]
\label{lemma-IndexedRealizabilityAndRealizability2}
Assume $s$ is a state constant, $t\in \PRclass$ is a closed term of state $\varnothing$, $A \in \mathcal{L} $ is a closed formula, and $t\in \|[ A] \|$. Let $\vec{t} = t_1, \ldots, t_n : \Nat$.

\begin{enumerate}
\item
$t\Vvdash_s P(\vec{t})$ if and only if $t[s]  = \varnothing$ in $\PRlearn$ implies
$P(\vec{t})={\True}$\\

\item
$t\Vvdash_s{A\wedge B}$ if and only if $\pi_0t \Vvdash_s{A}$ and $\pi_1t\Vvdash_s{B}$\\

\item
$t\Vvdash_s {A\vee B}$  if and only if either $\proj_0t[{s}]={\True}$ in $\PRlearn$ and $\proj_1t\Vvdash_s A$, or $\proj_0t[{s}]={\False}$ in $\PRlearn$ and $\proj_2t\Vvdash_s B$\\

\item
$t\Vvdash_s {A\rightarrow B}$ if and only if for all $u$, if $u\Vvdash_s{A}$,
then $tu\Vvdash_s{B}$\\

\item
$t\Vvdash_s {\forall x A}$ if and only if for all numerals $n$,
$t{n}\Vvdash_s A[{n}/x]$\\
\item

$t\Vvdash_s \exists x A$ if and only for some numeral $n$, $\pi_0t[{s}]= {n}$ in $\PRlearn$ and $\pi_1t \Vvdash_s A[{n}/x]$\\
\end{enumerate}
We define $t \Vvdash A$ if and only if $t\Vvdash_s A$ for all state constants $s$.
\end{definition}

We observe that theorem \ref{Fixed Point Property} holds as well for the stable terms of $\PCFclass$, for it is a consequence of the Stability theorem \ref{theorem-StabilityTheorem}. Hence, the Soundness theorem \ref{theorem-Soundness}, which depends only on the definition of realizability, stability and  theorem \ref{Fixed Point Property}, also holds for realizers of $\PCFclass$. That is, we have 

\begin{theorem}[Soundness Theorem ($\PCFclass$)]\label{theorem-SoundnessPCF} Let $A$ be a closed negation-and-implication free arithmetical formula. Suppose that $u\in\PCFclass$ and $u\Vvdash A$ and consider the game $\oneback(T_A)$. Let $\omega$ be as in definition \ref{definition-winningstrategy}. Then $\omega$ is a recursive winning strategy from $A$ for player one. \end{theorem}

\section{Completeness}\label{section-completenessproof}

\subsection{Idea of the Proof}\label{subsection-ideaproof}

In this section, we prove our completeness theorem: if an implication-and-negation-free arithmetical formula has a winning recursive strategy in its associated 1-Backtracking Tarski game, then it is realizable by a term of $\PCFclass$.

The idea of the proof follows naturally from the very meaning of learning based realizability. In order to realize a formula, one has to provide in the first place a Kleene-Kreisel-style realizer of the formula, \emph{recursive in an oracle} for the Halting problem. This corresponds to the fact that the terms of $\SystemTClass$ contain symbols for non computable functions which are in the same Turing degree of the aforementioned oracle. That is why one can see learning based realizability as a way of ``programming with non computable functions''. Hence, one would like to apply directly Berardi et al. \cite{BerCoq} result: given an implication-and-negation-free arithmetical formula, if there exists a recursive winning strategy for Eloise in its associated 1-Backtracking Tarski game, then there also exists a winning strategy for Eloise in its associated Tarski game, \emph{recursive in an oracle} for the Halting problem. 

However, that result is not enough for our purposes. According to learning based realizability, together with an oracle-equipped Kleene-Kreisel-style realizer, one has also to provide an effective method for learning oracle values in a convergent way and show that the realizer is always defined, whatever oracle approximations are used. Hence, we have to refine Berardi et al. result and prove that oracle values can be learned by counterexamples and the program is not ``perturbed'' by the oracle approximations used. More precisely, we show that the ineffective oracle strategy given in \cite{BerCoq} can be made more effective by using the novel ideas of learning based realizability: we first approximate the strategy by allowing it to use in the computations only approximated oracles; then we show that good enough approximations can be attained by a process of intelligent learning by counterexamples. One this task will be accomplished, the completeness theorem will follow just by formalizing the argument.\\

We now give an informal overview of the construction to be carried out. This should serve the reader as a guide to the next technical sections. Suppose that $\omega$ is a recursive winning strategy for Eloise in $\oneback(T_A)$. We start by describing a winning strategy for Eloise in $T_A$, which is recursive in some oracle for the Halting problem. We begin with some terminology.

\begin{definition}[Improvable, Optimal Plays]\label{definition-optimal}
We say that a $\omega$-correct play \[Q, A_0::\ldots :: A_i\] of $\oneback(T_A)$, with $A_i$ of the form $\exists x B$ or $B\lor C$, is \emph{improvable} if there exists a $\omega$-correct play of $\oneback(T_A)$ of the form
\[Q, A_0::\ldots :: A_i, Q', A_0::\ldots :: A_i\]
and we call this latter play an \emph{improvement} of the former. Moreover, a play 
is said to be \emph{optimal} if it is not improvable.
\end{definition}
The reason why a play
\[Q, A_0::\ldots :: A_i, Q', A_0::\ldots :: A_i\]
is called an improvement of $Q, A_0::\ldots :: A_i$
is that the former gives more information to the strategy $\omega$ in order to choose the next move for Eloise. Moreover, if $Q, A_0::\ldots :: A_i$ is optimal, whatever $\omega$-correct continuation of the game we may consider, $\omega$ will not backtrack to $A_0::\ldots ::A_i$ anymore. Since any such continuation will extend the play 
\[Q, A_0::\ldots :: A_i, A_0::\ldots :: A_i:: A_{i+1}\]
where 
\[\omega(Q, A_0::\ldots :: A_i)=A_0::\ldots :: A_i:: A_{i+1}\]
the choice of $A_{i+1}$ operated by $\omega$ is the best possible. 

The oracle $\Chi_E$ that we consider answers to questions of the form: is the play $Q, A_0::\ldots :: A_i$ improvable? To facilitate computations we also consider an oracle $\Phic_E$ which given the code of a play $Q, A_0::\ldots :: A_i$ returns the code of a play
\[Q, A_0::\ldots :: A_i, Q', A_0::\ldots :: A_i\]
whenever $Q, A_0::\ldots :: A_i$ is improvable and returns a dummy code otherwise. Observe that the two oracles $\Chi_E$ and $\Phic_E$ are of the same Turing degree, so exactly one of them would suffice. Furthermore, observe that one can define a program taking as input a code of an improvable play $Q, A_0::\ldots :: A_i$ and returning the code of an optimal improvement of the form
\[Q, A_0::\ldots :: A_i, Q', A_0::\ldots :: A_i\]
(just iterate $\Phic_E$).

 Suppose now that $A_0::\ldots:: \exists x B$ is a position in the game $T_A$. How should Eloise  move? The idea, coming from \cite{BerCoq}, is to compute, using our oracles,  an optimal play of $\oneback(T_A)$ of the form \[Q, A_0::\ldots:: \exists x B\]  Then, Eloise should respond  by first computing \[\omega (Q, A_0::\ldots:: \exists x B)= A_0::\ldots:: \exists x B:: B(n)\] and then choosing the formula $B(n)$. The idea is that $B(n)$ is a good choice, since no backtracking to $A_0::\ldots::\exists x A$ is ever to be done, following the strategy $\omega$.
 
More precisely, Eloise, while playing, simultaneously constructs a sequence of plays $Q_0, \ldots, Q_k$ of $\oneback(T_A)$ in the following way. She first defines $Q_0$ to be an optimal improvement of $A$. Then, suppose $k>0$ and that she is in the position 
\[A_0::\ldots ::A_k\] of the game $T_A$ and has constructed a sequence of plays $Q_0,\ldots, Q_k$ of $\oneback(T_A)$ such that\\

i) each play $Q_{i}$ is of the form $Q_{i}', A_0::\ldots:: A_i$;\\

ii) For all $i<k$, $Q_{i+1}$ extends $Q_i$; \\

iii) For all $i$, $Q_i$ is optimal;\\

iv) For all $i$, $Q_i$ is $\omega$-correct.\\
 
Then if Abelard has to move and chooses $A_{k+1}$, Eloise defines $Q_{k+1}$ as an optimal improvement of 
\[Q_k, A_0::\ldots :: A_k::A_{k+1}\]
which she can compute using the oracles.
If Eloise has to move, she computes
\[\omega(Q_k)=A_0::\ldots::A_k::A_{k+1}\]
she chooses as next move $A_{k+1}$ and she defines  $Q_{k+1}$ as an optimal improvement of 
\[Q_k, A_0::\ldots :: A_k::A_{k+1}\]
which she can again compute using the oracles. It is clear that the sequence $Q_1, \ldots, Q_{k+1}$ still satisfies all properties i)-iv).

Suppose now that a complete play $A_0::\ldots:: A_n$ of $T_A$ has been played by Eloise following the above strategy and suppose by contradiction that $A_n=\False$. Then, since by iv) $Q_n$ is $\omega$-correct and $\omega$ is winning, we have \[\omega(Q_n)=A_0::\ldots :: A_i\]
for some $i\leq n$, which represents a backtracking move performed by $\omega$. Since by ii) $Q_n$ extends $Q_i$,  we have that $Q_n, A_0::\ldots:: A_i$ can be written as \[Q_i, Q, A_0::\ldots:: A_i\]
for some $Q$. As a consequence, $Q_i=Q_i', A_0::\ldots A_i$ can be improved, which contradicts the optimality of $Q_i$ stated at point iii).

Thus we have a winning strategy for Eloise, recursive in the oracles $\Chi_E, \Phic_E$. Now, however, we want Eloise to follow the very same strategy, but using \emph{approximations} of those oracles in the place of the original ones. Of course, responses of approximated oracles are not to be always trustable. However, we will prove that correct oracle values can be learned by counterexamples and therefore that the use of the oracles may be replaced by a learning mechanism. According to learning based realizability, we will have in particular to prove that whenever the ``approximated" strategy does not lead Eloise to win $T_A$, then some new value of the oracles can be learned: at least, from a failure Eloise corrects something old and gains something new, which is a perfect example of a \emph{self-correcting} strategy. 

Now, suppose that a complete play $A_0::\ldots:: A_n$ of $T_A$ has been played by Eloise following the new ``approximated" strategy.  We still may suppose that $Q_0, \ldots, Q_n$ satisfy i), ii), iv). However, they satisfy only the following weaker\\

iii)$'$ For all $i$, $Q_i$ is optimal, \emph{according} to the response of the current oracle \emph{approximations}.\\\\
In other words, it might happen that our approximated oracles believe $Q_i$ to be not improvable, whilst $Q_i$ actually \emph{is}. Since iii) does not hold any more, the previous argument - that has shown $A_n=\True$ - now fails and it might still occur that $A_n=\False$. How Eloise is to learn from this counterexample? She computes $\omega(Q_n)$, which is of the form $A_0::\ldots :: A_i$, since it represents a backtracking move performed by $\omega$. Then as before she writes $Q_n, A_0::\ldots:: A_i$ as \[Q_i', A_0:: \ldots :: A_i, Q, A_0::\ldots:: A_i\] As a consequence, she finds out that $Q_i', A_0::\ldots A_i$ can be extended as to contradict the approximated-oracle prediction of point iii)$'$ and hence collects a new value of the oracle.\\

In the next two sections, we spell out the details of the construction and prove that the above Eloise strategy  is sound and in fact convergent. In order to enhance readability and separate the important ideas from technicalities,  we split the construction into two parts. First, we define the concept of learning strategy, which represents a translation of our realizability notion into the language of game theory, and show that any winning strategy in $\oneback(T_A)$ can be transformed into a learning strategy in $T_A$.  Secondly, we show that in fact any learning strategy in $T_A$ can be translated into a learning based realizer of $A$.

\subsection{Winning 1-Backtracking Strategies into Learning Strategies}

For the rest of the paper, fix a closed implication-and-negation-free arithmetical formula $A$ and let $T_A$ be its associated Tarski game. Fix moreover a primitive recursive enumeration of the plays of $\oneback(T_A)$ and let $\omega$ be a winning recursive strategy from $A$ for player one in the game $\oneback(T_A)$. We assume, without loss of generality, that $\omega$ performs backtracking moves only in front of atomic formulas; that is, we assume that for every play $Q, A_0::\ldots:: A_n$, if \[\omega(Q, A_0::\ldots:: A_n)=A_0::\ldots :: A_i\] then $A_n$ is atomic. Clearly, any winning strategy can be transformed accordingly to this requirement:  any backtracking move can be delayed by dummy moves and be performed in front of an atomic formula.

First of all, we formalize the coding of plays into numbers which has to be represented in our calculus.

\begin{definition}[Abstract Plays of $T_A$, Coding terms]\label{definition-AbstractPlay}
A sequence of arithmetical formulas $A_0::\ldots:: A_n$ is said to be an \emph{abstract play} of $T_A$, if $A_0=A$ and for all $i$ \\

i) if $A_i=\forall x B$ or $A_i=\exists x B$, then $A_{i+1}=B$;\\

ii) if $A_i=B_0\land B_1$ or $A_i=B_0\lor B_1$, then $A_{i+1}=B_0$ or $A_{i+1}=B_1$.\\\\
Let moreover $p=A_0::\ldots:: A_k$ be any abstract play of $T_A$. By $|p|:\Nat$, we denote a term of $\PCFclass$ having as free variables precisely the variables occuring free in the formulas of $p$ and such that, for every sequence of numerals $\vec{n}$ and sequence of variables $\vec{x}$ comprising all the free variables of $p$, $|p|[\vec{n}/\vec{x}]$ is equal to the numeric code of the play $q=A_0[\vec{n}/\vec{x}]::\ldots:: A_k[\vec{n}/\vec{x}]$.

\qed
\end{definition}
We define now a translation of learning based realizability into
the language of Tarski games. The translation consists of a pair $(g_0,g_1)$ of terms of $\PCFclass$: the first one describes a strategy for Eloise in $T_A$, and the second one decides when Eloise should backtrack.

\begin{definition}[Learning Strategy]\label{definition-LearningStrategy}
Let $T_A=(V, E_1, E_2, W)$. Let $g=(g_0,g_1)$ be a pair of terms of $\PCFclass$ respectively of types $\Nat \rightarrow \Nat$ and $\Nat\rightarrow \State$.\\ For every state $s$, we say that a play $A_0::\ldots :: A_n$ of $T_A$ is $g[s]$-\emph{correct} if  for every $i<n$ such that $(A_i,A_{i+1})\in E_1$
\[g_0[s](|A_0::\ldots:: A_i|)=|A_{i+1}|\]
holds. \\$g$ is said to be a \emph{learning strategy} from $A$ if it satisfies the following conditions:

\begin{enumerate}
\item (Soundness) For every  state $s$ and play $p=A::A_0::\ldots :: A_n$ such that $A_n\in Dom(E_1)$, if $g_0[s](|p|)=|A_{n+1}|$, then $A::A_0::\ldots :: A_n:: A_{n+1}$ is a play.\\

\item (Convergence) For every play $p$ starting from $A$, $g_0(|p|)$ and $g_1(|p|)$ converge.\\

\item (Learning) For every state $s$ and complete $g[s]$-correct play $p$ starting from $A$, if \[g_1[s](|p|)=\varnothing\implies p\in W\]  

\end{enumerate}

\end{definition}

We observe that conditions (2) and (3) of definition \ref{definition-LearningStrategy} correspond respectively to the convergence property that $\PCFclass$ realizers must have and to the learning condition in realizability for atomic formulas.  

We now define a predicate $E: \Nat^2\rightarrow \Bool$ of $\SystemTG$, which codes the improvement relation between plays of $\oneback(T_A)$ we are interested in. In the following, our terms of $\PCFclass$ will make use only of the oracle constants $\Chi_E$ and $\Phic_E$, which are in the syntax of $\SystemTClass$ (recall chapter \ref{chapter-learningbasedrealizability}, definition \ref{definition-TermLanguageL1}) and hence of $\PCFclass$. Moreover, we define a term $\Psi:\Nat\rightarrow \Nat$, which will be our fundamental computational engine. Given a code of a $\omega$-correct play, $\Psi$ is intended to return the code itself or the code of an optimal improvement, according to the oracles $\Chi_E, \Phic_E$. $\Psi$ is in general non computable and by using $\Psi$ as it is, we could only write strategies recursive in the oracles $\Chi_E, \Phic_E$, as in \cite{BerCoq}. Therefore, we shall compute its approximations $\Psi[s]$, by which we will be able to write the ``approximated" strategy for Eloise we have discussed in section \ref {subsection-ideaproof}.

\begin{definition}[Improvement Relation, Optimality Operator $\Psi$]\label{definition-ExtensionTerm} Let $E: \Nat^2\rightarrow \Bool$ a predicate of G\"odel's $\SystemTG$ such that $Enm=\True$ iff $n$ and $m$ are numerals  coding respectively $\omega$-correct plays $P, p:: A$ and $P, p::A, Q, p::A$ of $\oneback(T_A)$, with $A=\exists x B(x)$ or $A=B_0\lor B_1$. \\
We want to define now a term $\Psi: \Nat\rightarrow \Nat$ of  $\PCFclass$ such that the following equation is provable in the equational theory of $\PCFclass$:
\[\Psi z=\mathsf{if }\ \Chi_Ez\ \mathsf{then}\ \Psi(\Phic_Ez)\ \mathsf{else}\ z\]
In order to do that, it is enough to let \[\alpha:=\lambda y^{\Nat\rightarrow \Nat} \lambda z^\Nat\mathsf{if }\ \Chi_Ez\ \mathsf{then}\ y(\Phic_Ez)\ \mathsf{else}\ z\]
and take $\Psi:=\mathsf{Y}(\alpha)$.

\end{definition}

The term $\Psi$, given a number $n$, checks whether $\Chi_En=\True$, i.e. whether there exists a number $m$ such that $Enm=\True$. In that case, it computes such an $m$ by calling $\Phic_En$, and continues the computation by calling itself on $m$; otherwise, it returns $n$. The termination of $\Psi n$ is guaranteed by the fact that there are no infinite $\omega$-correct plays and each recursive call made by $\Psi$ extends a current $\omega$-correct play (see \cite{BerCoq}). If $\Chi_E, \Phic_E$ are interpreted as oracles, $\Psi n$ returns an optimal improvement of the play coded by $n$. But when $\Chi_E, \Phic_E$ are approximated through a particular state $s$, in general $\Psi n[s]$ will return only an improvement of the play coded by $n$, or even $n$ itself. \\ We now have to prove  the crucial property that for any finite approximation of the oracles $\Chi_E, \Phic_E$ - that is, for any state $s$ - the term  $\Psi n[s]$ has a normal form and that $\Psi n$ converge: $\Psi$ is ``stable" with respect to oracle approximations. 

\begin{proposition}[Stability of $\Psi$]\label{proposition-StabilityPsi}
$\Psi \in \|\Nat\rightarrow \Nat\|$

\end{proposition}
\proof Let $\{s_i\}_{i\in\NatSet}$ be a w.i. chain of states. By definition \ref{definition-StableTerms}, we have to prove that, for every term $t\in \|\Nat\|$, $\Psi t$ converges. Since $t$ converges to a numeral, it is enough to show that  for every numeral $n$, $\Psi n$ converges. First of all, we observe that for any state $s$, $\chi_Esm$ is equal to $\True$ only for a finite number of arguments $m$. Moreover, by  definition \ref{definition-ExtensionTerm}
\[\Psi n[s]=\mathsf{if }\ \chi_Es n\ \mathsf{then}\ \Psi(\varphi_Esn)\ \mathsf{else}\ n\]
 Hence, by direct computation it can be seen that, for every $i\in \NatSet$, $\Psi n[s_i]$ has a normal form and it is equal, for some $k\in \NatSet$, to  $(\varphi_Es_i)^{k} n$, having defined by induction \[(\varphi_Es_i)^0 n:=n,\ (\varphi_Es_i)^{m+1} n:= \varphi_Es_i((\varphi_Es_i)^{m}n)\] Moreover, for every $m<k$ \begin{equation}\label{eq0}\chi_Es_i((\varphi_Es_i)^{m}n)=\True\end{equation} does hold and hence $(\varphi_Es_i)^{m+1} n$ codes a play properly extending the play coded by $(\varphi_Es_i)^{m}n$. \\
Now, if $\Psi n$ did not converge, then there would be two increasing infinite sequences of numbers $k_0, k_1, k_2,\ldots$ and $m_0, m_1, m_2,\ldots $ such that, for every $i\in\NatSet$, $\Psi n[s_{m_i}]=(\varphi_Es_{m_i})^{k_i}n$.  Furthermore, since $s_{m_i}\leq s_{m_{i+1}}$ and, by (\ref{eq0}), for every $m<k_i$ \[\langle E, (\varphi_Es_{m_{i}})^{m}n, \varphi_Es_{m_i}((\varphi_Es_{m_{i}})^{m}n) \rangle \in s_{m_i}\] we have that for every $m\leq k_i$   \[(\varphi_Es_{m_{i+1}})^{m}n=(\varphi_Es_{m_{i}})^{m}n\] as it can be seen by induction on $m$. Hence, for every $i$, letting $a=k_{i+1}-k_i$
\[(\varphi_Es_{m_{i+1}})^{k_{i+1}}n=(\varphi_Es_{m+1})^a((\varphi_Es_{m_{i+1}})^{k_{i}}n)=(\varphi_Es_{m+1})^a((\varphi_Es_{m_{i}})^{k_{i}}n)\]
 would be the code of a $\omega$-correct play of $\oneback(T_A)$ properly extending the play coded by $(\varphi_Es_{m_{i}})^{k_{i}}n$. Therefore, it would exist an infinite $\omega$-correct play of $\oneback(T_A)$, which is impossible since $\omega$ is a winning strategy by hypothesis.

\qed

We are going to define three terms $\Lambda, \Pi, \Omega$ that will implement the learning strategy for Eloise in $T_A$ sketched in section \ref{subsection-ideaproof}. For the sake of readability, we will describe only the properties that such terms must satisfy, without explicitly write down those terms. It is trivial, however, to actually code our definitions in $\PCFclass$. 

We start by defining a term $\Lambda:  \Nat\rightarrow \State$, which is supposed to code the function $g_1$ of definition \ref{definition-LearningStrategy}. $\Lambda[s]$ takes the code of a play of $\oneback(T_A)$ and builds a state containing values of the oracles $\Chi_E$ and $\Phic_E$ that can be drawn from the input play and are not already in $s$.

\begin{definition}[Learning Term]\label{definition-Lambda} Let  $\Lambda: \Nat\rightarrow \State$ be a term of $\PCFclass$ described as follows. $\Lambda$ takes as argument a numeral $m$. Then, it checks whether $m$ codes a $\omega$-correct play $Q=P, A_0::\ldots:: A_n$ of $\oneback(T_A)$, with $A_0::\ldots:: A_n$ complete play of $T_A$ and $A_n=\False$. If not, it returns a dummy state: $\varnothing$. Otherwise,  it computes $\omega(Q)=A_0::\ldots::A_i$, and, for any state $s$,  returns $\Lambda[s] m$, which equals the state containing all the triples $\langle E, n_0, n_1\rangle$ not belonging to $s$ and such that \\\\
i) $n_0$ codes a play $Q_0, A_0::\ldots:: A_i$;\\\\
ii) $n_1$ codes a play $Q_0, A_0::\ldots:: A_i, Q_1,A_0::\ldots:: A_i$; \\\\
iii) $Q_0, A_0::\ldots:: A_i, Q_1, A_0::\ldots:: A_i=Q, A_0::\ldots:: A_i$.
\end{definition}

\textbf{Note.} The term $\Lambda$ must make use of the constant $\Add_E$ to create its output $\Lambda[s]m$, which is a state, and that is why the triples $\langle E, n_0, n_1\rangle$ of the above definition \ref{definition-Lambda} are not in $s$.\\

Recall subsection \ref{subsection-ideaproof}.  We have explained that Eloise, while playing in $T_A$ and in position $A_0::\ldots :: A_i$, simultaneously constructs a sequence of  plays $Q_0, \ldots, Q_i$ satisfying some properties i)-iv). We now define a term  $\Pi:\Nat\rightarrow\Nat$ of $\PCFclass$, which will be used to construct that sequence. In particular,  $\Pi$ takes the code of a play $p=A_0::\ldots :: A_i$ and yields the code of a play $Q_i=Q, p$ of $\oneback(T_A)$. Again, if $\Pi$ is interpreted as a program recursive in the oracles $\Chi_E, \Phic_E$, it will yield an optimal play $Q, p$. But when $\Pi[s]$ is computed, the result $Q, p$ may not be optimal because $\Chi_E, \Phic_E$ have been approximated through the state $s$.

\begin{definition}[Sequence Constructor $\Pi$]\label{definition-Pii} We describe the behaviour of a term $\Pi: \Nat\rightarrow \Nat$ of $\PCFclass$, intended to take the code $|p|$ of a play $p$ of $T_A$ and return the code $|Q, p|$ of a play $Q, p$ of $\oneback(T_A)$. The definition of $\Pi[s](|p|)$ runs by induction over the length of $p$ and is distinguished by cases:
\begin{enumerate} 
\item If $p=A$, then
$\Pi[s](|p|)=\Psi[s](|A|)$.\\
 \item If $p=(q::B)$ and
\[\Pi[s](|q|)=|Q, q|\] then  \[\Pi[s](|p|)=\Psi[s](|Q, q, q::B|) \]

\comment{\item If $p=(q::
\forall x B(x)::B(n))$ and \[\Pi(|q:: \forall x B(x)|)=|Q, q::\forall x B(x)|\]
then \[\Pi(|p)|)=\Psi(|Q, q::\forall x B(x), q::\forall x B(x)::B(n)|)\] 
\item  If $p=(q::  B_0\land B_1:: B_i)$
and \[\Pi(|q:: B_0\land  B_1|)=|Q, q::B_0\land B_1|\] then $\Pi(|p|)=\Psi(|Q, q::B_0\land B_1, q::B_0\land B_1::B_i|)$.\\ 
\item If
$p=(q::  B_0\lor B_1:: B_i)$ and \[\Pi(|q:: B_0\lor  B_1|)=|Q, q:: B_0\lor B_1|\] then
$\Pi(|p|)=\Psi(|Q, q::B_0\lor B_1, q::B_0\lor B_1::B_i|)$.}
\end{enumerate}
\qed
\end{definition}

We now prove the convergence of $\Lambda$ and $\Pi$.

\begin{proposition}[Stability of $\Lambda$ and $\Pi$]\label{proposition-StabilityLambdaPi} $\Lambda \in \|\Nat\rightarrow \State\|$ and $\Pi\in \|\Nat\rightarrow \Nat\|$.

\end{proposition}

\proof Again, to prove that $\Lambda \in \|\Nat\rightarrow S\|$  it is enough to show that  for every numeral $n$, $\Lambda n$ converges.  Let then $\{s_i\}_{i\in\NatSet}$ be a w.i. chain of states. By definition \ref{definition-Lambda}, whatever $s_i$ is, $\Lambda[s_i]n$ construct a finite set of triples $\langle E, n_0, n_1\rangle$, which depends only on $n$, and then decides to output some of the triples, according as to whether they are in $s_i$ or not. This determination stabilizes for large enough $s_m$; that is, for all $m'\geq m$, $\Lambda[s_{m'}]n=\Lambda[s_{m}]$.

The convergence of $\Pi n$ follows by straightforward  induction on the length of the play coded by $n$ and by the convergence of $\Psi$.

\qed

We now put together the terms $\Psi, \Lambda, \Pi$ in order to define a learning strategy for Eloise in $T_A$.

\begin{definition}[Learning Strategy for $T_A$]\label{definition-Omega} We describe the behavior of a pair of terms $\Omega=(\Omega_0, \Omega_1)$ respectively of type $\Nat\rightarrow \Nat$ and $\Nat\rightarrow \State$ of $\PCFclass$, intended to represent a learning strategy for $T_A$.  The definition of $\Omega_i[s](|p|)$ is distinguished by cases:
\begin{enumerate} 

 \item If $p=q:: \exists x B(x)$ and
\[\Pi[s](|q:: \exists x B(x)|)=|Q, q::\exists x B(x)|\] and \[\omega(Q,q::\exists x B(x))=q::\exists x B(x)::B(n)\]  then $\Omega_0[s](|p|)=|B(n)|$.\\
 \item If
$p=q::  B_0\lor B_1$ and 
\[\Pi[s](|q:: B_0\lor B_1|)=|Q, q::B_0\lor B_1|\] and \[\omega(Q,q::B_0\lor B_1)=q::B_0\lor B_1::B_i\] then
$\Omega_0[s](|p|)=|B_i|$.\\
\item If $p=q::B$, with $B$ atomic, and \[\Pi[s](|q::B|)=|Q, q::B|\]
then $\Omega_1[s](|p|)= \Lambda[s](|Q,q::B|)$.\\
\item In all other (trivial) cases, $\Omega_i[s](|p|)$, for $i=0,1$, can be arbitrarily defined as $0$ or $\True$.

\end{enumerate}
\end{definition}

\begin{proposition}[Stability of $\Omega_0$ and $\Omega_1$]\label{proposition-StabilityOmega}  $\Omega_0\in \|\Nat\rightarrow \Nat\|$ and $\Omega_1\in \|\Nat\rightarrow \State\|$.

\end{proposition}

\proof Trivial, by proposition \ref{proposition-StabilityLambdaPi}.

\qed

We need first prove that the $\oneback(T_A)$ plays constructed by $\Pi$ are $\omega$-correct, when $\Pi$ starts from $\Omega$[s]-correct plays of $T_A$ built by $\Omega$ (this property correspond to property iv) of section \ref{subsection-ideaproof}).
\begin{lem}\label{lemma-omegacorrectness}
Suppose $p:=A_0::\ldots :: A_n$ is a $\Omega[s]$-correct play of $T_A$. Then
\[\Pi[s](|A_0::\ldots ::A_n|)=|Q, A_0::\ldots:: A_n|\]
 and  $Q, A_0::\ldots:: A_n$ is $\omega$-correct. 
\end{lem}

\proof Routine induction on $n$. If $n=0$, then $p=A_0$. By definition \label{definition-Pi} of $\Pi$ and $\Psi$
\[\Pi[s](|A_0|)=\Psi[s](|A_0|)=|Q, A_0|\]
with $Q, A_0$ $\omega$-correct by construction of $\Psi$.\\
Suppose now $n>0$. By induction hypothesis 
\begin{equation}\label{pieq}\Pi[s](|A_0::\ldots ::A_{n-1}|)=|Q, A_0::\ldots:: A_{n-1}|
\end{equation}
and $Q, A_0::\ldots:: A_{n-1}$ is $\omega$-correct. If $A_{n-1}=\exists x B$ or $A_{n-1}=B\lor C$, then by definition of $\Omega$, by equation (\ref{pieq}) and $\Omega$[s]-correctness of $A_0::\ldots ::A_n$, we have that 
\[\Omega_0[s](|A_0::\ldots :: A_{n-1}|)=|A_n|\]
with 
\[\omega(Q, A_0::\ldots :: A_{n-1})=A_0::\ldots:: A_n\]
By definition \ref{definition-Pii} and equation (\ref{pieq})
 \[\Pi[s](|A_0::\ldots ::A_{n}|)=\Psi[s](|Q,A_0::\ldots :: A_{n-1}, A_0::\ldots:: A_n|)\]
 which is $\omega$-correct. If $A_{n-1}=B\land C$ or $A_{n-1}=\forall x B$, then 
\[\Psi[s](|Q,A_0::\ldots :: A_{n-1}, A_0::\ldots:: A_n|)\]
is automatically $\omega$-correct.

\qed

We now prove the main theorem of this section: any recursive winning strategy $\omega$ for Eloise in $\oneback(T_A)$ can be translated into a learning strategy from $A$ for Eloise in $T_A$.
\begin{theorem}[1-Backtracking Strategies into Learning Strategies]
$\Omega$ is a learning strategy for $T_A$.
\end{theorem}
\proof The fact that $\Omega$ satisfies properties 1 and 2 of definition \ref{definition-LearningStrategy} is trivial and follows from proposition \ref{proposition-StabilityOmega}. So we prove property 3. Let $s$ be a state and assume $p=A_0::\ldots ::A_n$ is a complete $\Omega[s]$-correct play of $T_A$. Suppose that $A_0::\ldots:: A_n\notin W$ and hence $A_n=\False$. We have to prove that $\Omega_1[s]|p|\neq \varnothing$. By definition \ref{definition-Omega} of $\Omega$, we have \[\Omega_1[s]|p|=\Lambda [s](|Q, A_0::\ldots:: A_n|)\]
with \[\Pi[s](|A_0::\ldots ::A_n|)=|Q, A_0::\ldots:: A_n|\]
By lemma \ref{lemma-omegacorrectness} $|Q, A_0::\ldots:: A_n|$ is $\omega$-correct, since $A_0::\ldots:: A_n$ is $\Omega[s]$-correct.
By definition \ref{definition-Lambda} of $\Lambda$, $\Omega_1[s]|p|$ contains all triples 
\[\langle E, |Q_0, A_0::\ldots:: A_i|, |Q_0, A_0::\ldots:: A_i, Q_1, A_0::\ldots:: A_i| \rangle\]
not in $s$ such that  \[ \omega(Q, A_0::\ldots:: A_n)=A_0::\ldots:: A_i\] and \[Q, A_0::\ldots:: A_n,A_0::\ldots:: A_i=Q_0, A_0::\ldots:: A_i, Q_1,A_0::\ldots:: A_i\] As implied by very definition \ref{definition-Pii} of $\Pi$, for every $j<n$, $\Pi[s] (|A_0::\ldots:: A_{j+1}|)$ codes a play extending the play coded by $\Pi[s] (|A_0::\ldots:: A_{j}|)$. Furthermore, for some $Q', Q''$
\[\Pi[s](|A_0::\ldots:: A_i|)=\Psi [s]|Q', A_0::\ldots:: A_i|=|Q'', A_0::\ldots:: A_i|\]
So, there is some $Q'''$ such that
\[Q'', A_0::\ldots:: A_i, Q'''=Q, A_0::\ldots:: A_n \]
and most importantly, by definition of $\Psi [s]|Q', A_0::\ldots:: A_i|$
\[\chi_E s|Q'', A_0::\ldots:: A_i|=\False\]
Therefore the triple 
\[\langle E, |Q'', A_0::\ldots:: A_i|, |Q'', A_0::\ldots:: A_i, Q''', A_0::\ldots:: A_i|\rangle\]
belongs to $\Omega_1[s]|p|$, since it is not in $s$, by definition of $\Chi_E s$.

\qed

\subsection{Learning Strategies into Realizers}

In this section, we prove that the learning strategy $\Omega$ for $T_A$ can be translated into a learning based $\PCFclass$ realizer of $A$, thus proving our main completeness theorem.

We begin with a bit of coding.

\begin{definition}\label{definition-Code}
Let $\Omega_0^\Nat: \Nat\rightarrow \Nat$ and $\Omega_0^B: \Nat\rightarrow \Bool$ be  terms of $\PCFclass$  such that:

\begin{enumerate}

\item for every play $p::\exists x B$ of $T_A$ 
\[\Omega_0^\Nat|p::\exists x B|=n \iff \Omega_0|p::\exists x B|=|B(n)|\]
\item  for every play $p:: B_0\lor B_1$
 \[\Omega_0^B |p::B_0\lor B_1|=\True\iff\Omega_0 |p::B_0\lor B_1|=|B_0|\]

\end{enumerate}
\qed
\end{definition}

As a special case of the following  definition, we get a candidate realizer for $A$.

\begin{definition}[Realizer for $A$]\label{definition-Realizer}
Let $p$ be an abstract play of $T_A$. We define by induction and by cases a term $t_p$ of $\PCFclass$, with free variables among those occurring free in some formula of $p$, as follows:

\begin{enumerate}
\item $p=q::\forall x B$. \[t_{q::\forall x B}=\lambda x\ t_{q::\forall x B:: B}\]

\item $p=q::\exists x B$. \[t_{q::\exists x B}=\langle \Omega_0^\Nat|q::\exists x B|, t_{q::\exists x B:: B}[t_1/x]\rangle\]
where $t_1:= \Omega_0^\Nat|q::\exists x B|$.\\

\item $p= B_0\lor B_1$. \[t_{q:: B_0\lor B_1}=\langle \Omega_0^B|q::B_0\lor B_1|, t_{q:: B_0\lor B_1::B_0}, t_{q:: B_0\lor B_1::B_1}\rangle\]

\item $p=B_0\land B_1$. \[t_{q:: B_0\land B_1}=\langle  t_{q:: B_0\land B_1::B_0}, t_{q:: B_0\land B_1::B_1}\rangle\]

\item $p=q$, with $q$ complete. \[t_{q}=\Omega_1|q|\]

\end{enumerate}
\qed

\end{definition}

\begin{lem}[Completeness Lemma]\label{lemma-Completeness}\begin{enumerate}
\item Let $\vec{m}=m_0,\ldots, m_k$ a sequence of closed stable type-$\Nat$ terms of $\PCFclass$ and $\vec{x}=x_0,\ldots, x_k$ a sequence of variables containing all the free variables of $p::B$. Then  \[t_{p::B}[\vec{m}/\vec{x}]\in \|[ B] \|\]

\item Let $s$ be a state, $\vec{m}=m_0,\ldots, m_k$ a sequence of closed type-$\Nat$ terms of $\PCFclass$, $\vec{x}=x_0,\ldots, x_k$ a sequence of variables and  $(p::B)[\vec{m}[s]/\vec{x}]$ a $\Omega$[s]-correct play of $T_A$. Then
 \[t_{p::B}[\vec{m}/\vec{x}]\Vvdash_s B[\vec{m}[s]/\vec{x}]\]

\end{enumerate}
\end{lem}
\proof We prove (1) by induction on $p$ and by cases. We treat only three representative cases, those left out being obvious.

\begin{enumerate}

\item $B=\forall x C$. Let $n$ be a numeral. By inductive hypothesis, we have that  \[t_{p::\forall x C}[\vec{m}/\vec{x}]n=(\lambda x\ t_{p::\forall x C:: C}[\vec{m}/\vec{x}])n\]\[=t_{p::\forall x C:: C}[\vec{m},n/\vec{x},x]\in \|[ C] \|\]
Since $[ \forall x C] =\Nat\rightarrow [ C]$, we have that \[t_{p::\forall x C}[\vec{m}/\vec{x}]\in \|[ \forall x C] \|\]

\item $B=\exists x C$. Let \[t_1:=\Omega_0^\Nat|p::\exists x C|[\vec{m}/\vec{x}]\] Since the terms in $\vec{m}$ are stable by hypothesis and $\Omega_0$ is stable by proposition \ref{proposition-StabilityOmega} and by definition \ref{definition-Code} of $\Omega_0^\Nat$, we have that  $t_1$ is a stable term of type $\Nat$. 
By inductive hypothesis \[t_{p::\exists x C:: C}[\vec{m},t_1/\vec{x},x]\in \|[ C]\|\]
Since $[ \exists x C] =\Nat\times [ C]$ and  
 \[t_{p::\exists x C}[\vec{m}/\vec{x}]=\langle \Omega_0^\Nat|p::\exists x C|[\vec{m}/\vec{x}], t_{p::\exists x C:: C}[\vec{m},t_1/\vec{x},x]\rangle\]
we have \[t_{p::\exists x C}[\vec{m}/\vec{x}]\in \|[ \exists x C] \|\]

\item  $B$ atomic. Then \[t_{p::B}[\vec{m}/\vec{x}]= \Omega_1|p::B|[\vec{m}/\vec{x}]\] Since $[ B]=\State$ and $\Omega_1$ is stable by proposition \ref{proposition-StabilityOmega}, we have that $t_{p::B}[\vec{m}/\vec{x}]\in \|[ B] \|$.

\end{enumerate}

 We now prove (2) by induction on $p$ and by cases.
\begin{enumerate}
\item $B=\forall x C$. Let $n$ be a numeral. By inductive hypothesis, we have that  \[t_{p::\forall x C}[\vec{m}/\vec{x}]n=(\lambda x\ t_{p::\forall x C:: C}[\vec{m}/\vec{x}])n\]\[=t_{p::\forall x C:: C}[\vec{m},n/\vec{x},x]\Vvdash_s C[\vec{m},n[s]/\vec{x},x]\]
and hence \[t_{p::\forall x C}[\vec{m}/\vec{x}]\Vvdash_s \forall x C[\vec{m}[s]/\vec{x}]\]

\item $B=\exists x C$. Suppose \[t_1[s]:=\Omega_0^\Nat[s]|p::\exists x C|[\vec{m}[s]/\vec{x}]=n\] with $n$ numeral. By definition \ref{definition-Code} of $\Omega_0^\Nat$, we have that \[\Omega_0[s]|q::\exists x C|[\vec{m}[s]/\vec{x}]=|C[\vec{m}[s],n/\vec{x},x]|\]
and so $p::\exists x C::C[\vec{m}[s],n/\vec{x},x]$ is $\Omega[s]$-correct.
By inductive hypothesis \[t_{p::\exists x C:: C}[\vec{m},t_1/\vec{x},x]\Vvdash_s C[\vec{m}[s],n/\vec{x},x]\]
Since \[t_{p::\exists x C}[\vec{m}/\vec{x}]=\langle \Omega_0^\Nat|p::\exists x C|[\vec{m}/\vec{x}], t_{p::\exists x C:: C}[\vec{m},t_1/\vec{x},x]\rangle\]
we have \[t_{p::\exists x C}[\vec{m}/\vec{x}] \Vvdash_s \exists x C[\vec{m}[s]/\vec{x}]\]

\item $B=C_0\land C_1$. By inductive hypothesis
\[t_{p:: C_0\land C_1::C_0}[\vec{m}/\vec{x}] \Vvdash_s C_0[\vec{m}[s]/\vec{x}]\]
and 

\[t_{p:: C_0\land C_1::C_1}[\vec{m}/\vec{x}] \Vvdash_s C_1[\vec{m}[s]/\vec{x}]\]
Since \[t_{p:: C_0\land C_1}[\vec{m}/\vec{x}]=\langle  t_{p:: C_0\land C_1::C_0}[\vec{m}/\vec{x}], t_{p:: C_0\land C_1::C_1}[\vec{m}/\vec{x}]\rangle\]
we have  \[t_{p:: C_0\land C_1}[\vec{m}/\vec{x}]\Vvdash_s C_0\land C_1[\vec{m}[s]/\vec{x}]\]

\item $B=C_0\lor C_1$. Suppose \[t_1[s]:=\Omega_0^B[s]|p::C_0\lor C_1|[\vec{m}[s]/\vec{x}]=\True\] Then, by definition \ref{definition-Code} of $\Omega_0^B$, we have that \[\Omega_0[s]|q::C_0\lor C_1|[\vec{m}[s]/\vec{x}]=|C_0[\vec{m}[s]/\vec{x}]|\]
and hence $p::C_0\lor C_1::C_0[\vec{m}[s]/\vec{x}]$ is $\Omega[s]$-correct.
By inductive hypothesis \[t_{p:: C_0\lor C_1:: C_0}[\vec{m}/\vec{x}]\Vvdash_s C_0[\vec{m}[s]/\vec{x}]\]
Analogously, for $t_1[s]=\False$, we have
\[t_{p:: C_0\lor C_1:: C_1}[\vec{m}/\vec{x}]\Vvdash_s C_1[\vec{m}[s]/\vec{x}]\]
 Since \[t_{p:: C_0\lor C_1}[\vec{m}/\vec{x}]=\langle \Omega_0^B|p::C_0\lor C_1|[\vec{m}/\vec{x}], t_{p::C_0\lor C_1::C_0}[\vec{m}/\vec{x}], t_{p::C_0\lor C_1::C_1}[\vec{m}/\vec{x}]\rangle\]
we have \[t_{p:: C_0\lor C_1}[\vec{m}/\vec{x}] \Vvdash_s  C_0\lor C_1[\vec{m}[s]/\vec{x}]\]

\item $B$ atomic. Then \[t_{p::B}[\vec{m}/\vec{x}][s]= \Omega_1[s]|p::B|[\vec{m}[s]/\vec{x}]\] Since $p::B[\vec{m}[s]/\vec{x}]$ is $\Omega[s]$-correct and $\Omega$ is a learning strategy by lemma \ref{lemma-Completeness}, if $t_{p::B[\vec{m}/\vec{x}]}[s]=\varnothing$, then $B[\vec{m}[s]/\vec{x}]=\True$.\\
\end{enumerate}
\qed

\begin{theorem}[Completeness theorem]
Suppose there exists a recursive winning strategy for player one in $\oneback(T_A)$. Then there exists a term $t$ of $\PRclass$ such that $t\Vvdash A$.
\end{theorem}

\proof By Lemma \ref{lemma-Completeness}, point 1 and 2,  applied to $t_A$ and the empty sequence of terms.
\qed

\section{Conclusions}

We have proved a soundness and completeness result for total recursive learning based realizability with respect to 1-Backtracking game semantics, solving a conjecture left open in Aschieri \cite{Aschierigames}.

The contribution of the soundness theorem is semantical, rather than technical, and it should be useful to understand the significance and see possible uses of learning based realizability. We have shown how learning based realizers may be understood in terms of backtracking games and that this interpretation offers a way of eliciting constructive information from them. The idea is that playing games represents a way of challenging realizers; they react to the challenge by learning from failure and counterexamples. In the context of games, it is also possible to appreciate the notion of convergence, i.e. the fact that realizers stabilize their behaviour as they increase their knowledge. Indeed, it looks like similar ideas are useful to understand other classical realizabilities (see for example, Miquel \cite{Miq}). 

The proof of the completeness theorem has been definitely more technically challenging. In our view, moreover, it has two interesting features. In a sense, it is  the first application of the ideas of learning based realizability to a concrete non trivial classical proof, which is our version of the one given by Berardi et al. \cite{BerCoq}. This proof classically shows that if Eloise has recursive winning strategy in the 1-Backtracking Tarski game associated to a formula $A$, then she also has a winning strategy in the Tarski game associated to $A$ (but a strategy only recursive in an oracle for the Halting problem). Since the existence of this latter strategy implies the truth of $A$, the argument can be seen as a proof of $A$ in some version of intuitionistic Arithmetic with $\EM_1$. We managed to associate a constructive content to this seemingly ineffective proof and found out that it hides a learning mechanism to gain correct oracle values from failures and counterexamples. We have then transformed this learning mechanism into a learning based realizer of $A$.   Secondly, we have shown the interesting theoretical result that backtracking strategies in 1-Backtracking games can interpreted as learning realizers. We have thus successfully established a close non trivial relationship between two interpretations of classical proofs: game semantics and learning based realizability.

\chapter{Constructive Analysis of Learning in Peano Arithmetic}\label{chapter-constructiveanalysislearning}

\begin{abstract}
In this chaper we give a constructive analysis of learning as it arises in various computational interpretations of classical Peano Arithmetic, such as our learning based realizability, Avigad's  update procedures and epsilon substitution method. In particular, we show how to compute in G\"odel's system $\SystemTG$ upper bounds on the length of learning processes, which are themselves represented in $\SystemTG$ through learning based realizability. The result is achieved by the introduction of a new non standard model of G\"odel's $\SystemTG$, whose new basic objects are pairs of non standard natural numbers (convergent sequences of natural numbers) and moduli of convergence, where the latter are objects giving constructive information about the former. As foundational corollary, we obtain that learning based realizability is a constructive interpretation of Heyting Arithmetic plus excluded middle over $\Sigma_1^0$ formulas (for which it was designed) and of all Peano Arithmetic when combined with G\"odel's double negation translation. As byproduct of our approach, we also obtain a new proof of Avigad's theorem for update procedures and thus of termination of epsilon substitution method for $\PA$.

\end{abstract}

\maketitle

\section{Introduction}\label{section-intro}

The aim of this chapter is to carry out a detailed and complete constructive analysis of learning, as it arises in learning based realizability for $\HA+ \EM_1$ and in Avigad's \cite{Avigad} axiomatization of the epsilon substitution method for Peano Arithmetic through the concept of update procedure. The importance of this analysis is both practical and foundational. In the first place, we explicitly show how to compute upper bounds to the length of learning processes, thus providing the technology needed to analyze their computational complexity. Secondly, we answer positively to the foundational question of whether learning based realizability can be seen as an interpretation of classical Arithmetic into intuitionistic Arithmetic.

Our constructive framework is G\"odel's system $\SystemTG$ and our metatheory will be purely intuitionistic. Our analysis will be accompIished by  restating and then reproving constructively the following convergence theorem.

\begin{theorem}[Convergence]\label{theorem-con}
Let $t:\funnat$ be a closed term of $\SystemTG$. Let $s: \seqfun$ be any closed  term of $\SystemTG$ representing a weakly increasing chain of functions: that is, assume that for every numerals $n\leq m$, $s_n\leq s_m$\footnote{Define $s_n\leq s_m$ iff for all numerals $l$, $s_n(l)\neq 0$ implies $s_n(l)=s_m(l)$. See the premise to definition \ref{definition-ordering} for intuitive meaning.} holds. Then, there exists an $n$ such that for all $m\geq n$, $t(s_n)=t(s_m)$.

\end{theorem}

 The intuitive meaning of the convergence theorem is the following. It is intended to be an analysis of oracle computations. That is, given a non computable function $f:\funn$ one would like to ``compute" $t(f)$. Since this is not effectively possible, in order to obtain significant results one may try nevertheless to define a weakly increasing chain $s$ of functions with the property that for all numerals $n$, $s_n\leq f$. Such a chain can be seen as a sequence of more and more refined approximations of $f$ and can for example be constructed by means of \emph{learning processes} as they arise in learning based realizability or epsilon substitution method (see Mints \cite{Mints}). The theorem says that if $t$ is computed with respect to such a sequence of approximations, then a stable answer about the value of $t(f)$ is eventually obtained. 

The  convergence theorem is already interesting in itself,  but its special significance lies in its consequences,  which we now describe and shall prove in the final part of the chapter. Since most of them cannot be proved if the convergence theorem is not first restated and then proven constructively, they provide an important motivation for working in this direction.  

A first consequence of the convergence theorem is that any learning process represented by a learning based realizer always terminates. Formally:

\begin{theorem}[Zero Theorem]  Let $t$ be a type $\State$ term of $\SystemTClass$. Define $s_0:=\varnothing$ and, for every natural number $n$, $s_{n+1}:=s_n\Cup t[s_n]$. Then, there is an $n$ such that $t[s_n]=\varnothing$. 
\end{theorem}

If the convergence theorem is proven constructively, also the above Zero theorem can be and so one obtains a constructive analysis of the numbers of learning steps required to complete the learning process. It has as a constructive consequence the following theorem:

\begin{theorem}[Program Extraction via Learning Based Realizability]\label{theorem-ext} Let $t$ be a term of $\SystemTClass$ and suppose that $t\Vvdash \forall x^\Nat \exists y^\Nat Pxy$, with $Pxy$ atomic. Then, from $t$ one can define a term $u$ of G\"odel's system $\SystemTG$ such that for every numeral $n$, $Pn(un)=\True$.
\end{theorem}
The above theorem sharpens the result obtained in chapter \ref{chapter-learningbasedrealizability} and in Aschieri and Berardi \cite{AB}. There, we have proved as well that from any $t$ such that $t\Vvdash\forall x^\Nat \exists y^\Nat Pxy$ one can extract a computable function $v$ such that for every numeral $n$, $Pn(vn)=\True$. However, the extracted $v$ made use of unbounded iteration, while the $u$ of theorem \ref{theorem-ext} is a ``bounded" algorithm, that is, a program not explicitly using any kind of unbounded iteration. This is an important point from a foundational point of view: the algorithms extracted via learning based realizability \emph{construct} witnesses, rather than \emph{searching} for them. Let us make clear that, however, $u$ - from the computational point of view - is equal to $v$. In fact, $u$ results from $v$ just by replacing its only unbounded iteration with a primitive recursive one (of an appropriate type). Thus, $u$ just adds to $v$ information about the computational complexity of the learning process generated by $v$. For practical purposes, therefore, $v$ is as efficient as $u$. 

As corollary, one obtains the important result that from classical proofs in Peano Arithmetic $\PA$ of $\forall\exists$-formulas one can extract bounded algorithms via learning based realizability $\Vvdash$. This is done by, first, extracting a realizer from any given proof and, then, by applying theorem \ref{theorem-ext}. In other words, one is able to give a novel proof of the following theorem due to G\"odel (through its Dialectica interpretation, see e.g. \cite{Kohl}):

\begin{theorem}[Provably Total Functions of $\PA$] If  $\PA\vdash \forall x^\Nat \exists y^\Nat Pxy$, then there exists a term $u$ of G\"odel's system $\SystemTG$ such that for every numeral $n$, $Pn(un)=\True$. 
\end{theorem}
The novelty, here, is the technique employed to prove the theorem and the new understanding of extracted programs as realizers able to learn in a constructive way.

From a constructive proof of the convergence theorem one can also provide new constructive proofs of Avigad's \cite{Avigad} \emph{fixed point theorem} for $n$-ary update procedures and hence of the \emph{termination of the epsilon substitution method} for $\PA$. Hence, one also obtains a \emph{constructive analysis of learning} in Peano Arithmetic. The novelty, here, is the use of type theory to reason about the learning processes generated by update procedures and hence epsilon substitution method. 

Theorem \ref{theorem-con} can be proven easily, but ineffectively, in second order logic:\\

\textit{Proof of theorem \ref{theorem-con} (Ineffective)}. The informal idea of the proof is the following. Terms of system $\SystemTG$ use only a finite number of  values of their function arguments. If we ``apply" $t$ to the least upper bound $f_s$ of the sequence $s$ (w.r.t the relation $\leq$ of definition \ref{definition-ordering}), we find that the finite part of $f_s$ effectively used in the computation of $t(f_s)$ is already contained in some $s_k$. So,  for every $h\geq k$, $t(s_h)=t(s_k)$.\\  Let us see the details. As proven by Kreisel (for a proof see Schwichtenberg \cite{Schwichtenberg}), $t$ has a modulus of continuity $\mathcal{C}$, which is a term of system $\SystemTG$ of type $\funnat$ such that the following statement is provable in extensional $\HA^\omega$:

\begin{equation}\label{s1}\forall f^\funn, g^\funn. (\forall x^\Nat\leq (\mathcal{C}f)\ f(x)=g(x))\rightarrow t(f)=t(g)\end{equation}
By using the comprehension axiom, we can define the least upper bound $f_s$ of the sequence $s$ as follows

\[f_s(n)= \begin{cases}
m   &\text{if $\exists i$ such that $s_i(n)=m\neq 0$}\\

0  &\text{otherwise}

\end{cases}
\]
Let $\mathcal{C}^M$ be the denotation of $\mathcal{C}$ in the full set theoretic model $M$ of extensional $\HA^\omega$ (see Kohlenbach \cite{Kohl}). Then there exists an $n$ such that for all $m\geq n$

\[\forall x^\Nat\leq (\mathcal{C}^Mf_s)\ s_n(x)=s_m(x)\]
By \ref{s1}, we get that for all $m\geq n$, $t(s_n)^M=t(s_m)^M$. Hence by soundness of the model with respect to formal equality  of extensional $\HA^\omega$, $t(s_n)$ and $t(s_m)$ normalize to the same numeral, since $t(s_n)=a$ and $t(s_m)=b$, with $a,b$ numerals, implies $a^M=t(s_n)^M=t(s_m)^M=b^M$ and then $a=b$.

\qed

The convergence theorem is therefore true, but one cannot hope to prove it constructively as it is stated. In fact, it is a formula of the form $\forall\exists \forall$ and the simplest incompleteness of intuitionistic reasoning as compared to classical reasoning arises precisely for that kind of formulas. It is known, for example, that classical finite type Peano Arithmetic $\PA^\omega$ proves the formula $\forall f^\funn\exists x^\Nat \forall y^\Nat f(x)\leq f(y)$, while intuitionistic Heyting Arithmetic $\HA^\omega$ does not. In our case, one could associate to any Turing machine a weakly increasing sequence $s:\Nat\rightarrow (\funn)$ such that for all $m$, $s_m(n)=0$ if $n\neq 0$, and $s_m(0)=1$ if the machine terminates on input $n$ in less than $m$ steps, $s_m(0)=0$ otherwise. A constructive proof of the convergence theorem relatively to the term $\lambda f^\funn f(0)$ would compute the limit of the sequence $\lambda m^\Nat s_ m(0)$, thus determining whether the Turing machine terminates on input $n$. By producing such a sequence $s$ for every Turing machine, we would have a solution for the Halting problem.\\ \\

\textit{Synopsis of the chapter}. In the rest of the chapter, we develop a technology for constructively reasoning about convergence in G\"odel's system $\SystemTG$ and proving a classically equivalent form of the convergence theorem. All proofs will be constructive and all their constructive content will be made explicit. This constructive effort results in a longer and a bit more complex presentation than it could be. However, if one is not interested in full explicit details,  the techniques used may be simplified in order to yield quite short and powerful constructive proofs of the main results of the chapter.

Our approach has a semantical content. In fact, we starts from considering a kind of constructive non standard model for Peano Arithmetic and then we reinterpret G\"odel's system $\SystemTG$ constants in order to manipulate the new individuals of the model. The reinterpretation of system $\SystemTG$ will turn out to be particularly suited to perform the computations we need to do for constructively reasoning about convergence.  From the high level point of view, the proof techiques used amount to a combination of Kreisel's no-counterexample interpretation and Tait's reducibility/logical-relations method. With the first one, we can constructively reason about convergence. With the second, we prove the soundness of the model with respect to our purposes. 
   
 In detail, the plan of the chapter is the following.
 
  In section \S \ref{section-TheTermCalculus2}, we recall details of G\"odel system $\SystemTG$. 
  
  In section \S \ref{section-notionconvergence} we define the first ingredient of our approach, which is a constructive notion of convergence for sequences of objects, due to Berardi \cite{BerPer}. It is a no-counterexample interpretation of the classical notion of convergence, but it is different from the usual interpretation. Its main advantage is that it is very efficient from the computational point of view, since it enables \emph{programming with continuations} and hence the writing of powerful and elegant realizers of its constructive content, which we will call \emph{moduli of convergence}. Intuitively, a modulus of convergence for a convergent function $f:\Nat\rightarrow A$ will be a term able to find suitable intervals in which $f$ is constant; moreover,  the length of those intervals   will depend on a continuation. At the end of the section we use Berardi's notion of convergence to reformulate the convergence theorem (see theorem \ref{theorem-weakconvergence}).
  
 In section \S \ref{section-nonstandardmodel}, we introduce the second ingredient of our approach: a model that extends the usual full finite type structure generated over natural numbers by replacing naturals by pairs $\langle \mathcal{N}, f\rangle$ of a non standard natural number $f$ (which is a function $\funn$ as in ultrapower models of Peano Arithmetic) and its modulus of convergence $\mathcal{N}$. We also syntactically define a semantics $\model{\_}$ (where $s$ is a weakly increasing chain of functions) mapping terms of $\SystemTG$ in to elements of the model and in section \S \ref{section-adeguacy} we show that, thanks to $\model{\_}$, we can evaluate every term $t: (\funn)\rightarrow \Nat$, into a pair $\langle \mathcal{N}, f\rangle$ such that $\mathcal{N}$ is a modulus of convergence for the function $f=\lambda n^\Nat t(s_n)$. 
 
 In section \S \ref{section-consequences}, we prove all the corollaries of the convergence theorem that we have discussed before.  

\section{Term Calculus}\label{section-TheTermCalculus2}

In this chapter we will prove results that hold for any ``simple" extension of G\"odel's system $\SystemTG$ (see chapter \ref{chapter-technicalpreliminaries}). In this section, we recall the definition and results we shall need and introduce some useful notation.\\

\textbf{Notation}. For notational convenience and to define in a more readable way terms of type $A\times B\rightarrow C$, for any variables $x_0: A$ and $x_1: B$ we define 
\[\lambda \langle x_0, x_1\rangle^{A\times B} u:= \lambda x^{A\times B} u[\pi_0x/x_0\ \pi_1x/x_1]\]
where $x$ is a fresh variable not appearing in $u$. We observe that for any terms $t_0, t_1$
\[(\lambda \langle x_0, x_1\rangle^{A\times B} u)\langle t_0, t_1\rangle= u[t_0/x_0\ t_1/x_0]\]

Often, as in chapter \ref{chapter-learningbasedrealizability}, it is useful to add to system $\SystemTG$ new constants and atomic types, together with a set of algebraic reduction rules we call ``functional''.

\begin{definition}[Functional set of rules]\label{definition-functional}
Let $C$ be any set of constants, each one of some type $A_1\rightarrow \ldots \rightarrow A_n\rightarrow A$, for some atomic types $A_1,\ldots,A_n, A$. We say that $\mathcal{R}$ is a {\em functional set of reduction rules} for $C$ if $\mathcal{R}$ consists, for all $c\in C$ and all closed normal terms ${a_1}:A_1,\ldots, {a_n}:A_n$ of $\SystemT$, of one and exactly one rule $c {a_1}\ldots {a_n}\mapsto {a}$, where ${a}:A$ is a closed normal term of $\SystemT$.
\end{definition}

If a system $\SystemT$ is obtained from G\"odel's $\SystemTG$ by adding a recursive set $C$ of constants and a recursive functional set of rules for $C$, we we call $\SystemT$ a \emph{simple extension of} $\SystemTG$.  As in chapter \ref{chapter-learningbasedrealizability}, by a standard reducibility argument it can be proved that $\SystemT$ is strongly normalizing and has Church-Rosser property. Moreover, any atomic-type term of any simple extension $\SystemT$ of $\SystemTG$ is equal either to a numeral, if it is of type $\Nat$, or to a boolean, if it is of type $\Bool$, or to a constant of type $A$, if it is of type $A$. All results of this paper hold whatever simple extension of $\SystemTG$ is chosen. Let us fix one.

\begin{definition}[System $\SystemT$]
From now on, we denote with $\SystemT$  an arbitrarily chosen simple extension of G\"odel's system $\SystemTG$. We also assume that $\SystemT$ contains constants for deciding equality of constants of atomic type.
\end{definition}

Throughout the paper, the intended interpretation of the natural number $0$ will be as a ``default'' value. That is,  when we do not have any information about what value a function has on argument $n$, we assume that it has value $0$. That being said, it is natural to consider a function $f_1: \NatSet\rightarrow \NatSet$ to be \emph{extending} another function $f_2:\NatSet\rightarrow \NatSet$, whenever it holds that for every $n$ such that $f_1(n)$ is a non default value (and hence different from $0$), then $f_1(n)=f_2(n)$. $f_2$ may hence have a non default value at some argument where $f_1$ has a default value, but it agrees with $f_1$ at the arguments where $f_1$ has not default value. So, $f_2$ carries more information than $f_1$.

\begin{definition}[Ordering Between Functions and Terms]\label{definition-ordering}
Let $f_1, f_2$ be functions $\NatSet\rightarrow \NatSet$. We define \[f_1\leq f_2 \iff \forall n\in\NatSet\ f_1(n)\neq 0\Rightarrow f_1(n)=f_2(n)\]
Moreover, if $t_1, t_2$ are closed terms of $\SystemT$ of type $\Nat\rightarrow \Nat$ representing respectively functions $g_1, g_2:\NatSet\rightarrow \NatSet$, we will write $t_1\leq t_2$ if and only if $g_1\leq g_2$.
\end{definition}

In the following, we will write ``$s\in$ w.i." if is $s:\seqfun$ is a closed term representing a weakly increasing sequence of functions, that is, if for all numerals $n,m$,  $n\leq m$ implies $s_n\leq s_m$.

\section[No-Counterexample Interpretation and Berardi's Notion of Convergence]{The No-Counterexample Interpretation and Berardi's Notion of Convergence}\label{section-notionconvergence}

In this paper, we are interested in arithmetical formulas stating convergence of natural number sequences. Classically, we consider a sequence of natural numbers to be convergent if it is definitely constant, that is, if the there is an element of the sequence which is equal to all successive elements of the sequence. Hence, we will consider formulas of the form  \begin{equation}\label{ef}(\forall z^A)\ \exists x^\Nat \forall y^\Nat P(z,x,y)\end{equation}
Since that kind of formulas cannot generally be proven constructively, a common standpoint is to consider classically equivalent but constructively weak enough statements, as in Kreisel \emph{no-counterexample interpretation}:

\[(\forall z^A)\ \forall f^{\funn}\exists x^\Nat P(z,x,f(x))\]
If the statement \ref{ef} (with $A=\Nat$) is provable in $\PA$, then one can constructively extract from any proof a term $t:\funnat$ of system $\SystemT$ such that 

\[(\forall z^\Nat)\ \forall f^{\funn} P(z,t(f),f(t(f)))\]
holds (see for example Kohlenbach \cite{Kohl}).  In our cases, we have to deal with formulas of the form

\[\exists x^\Nat \forall y^\Nat\geq x f(x)=f(y)\]
where $f$ is a term of type $\funn$, and hence we may be tempted to consider their no-counterexample interpretation

\begin{equation}\label{e1}\forall h^{\funn} \exists x^\Nat\ h(x)\geq x \rightarrow f(x)=f(h(x))\end{equation}
If one introduces the notation \[f\downarrow [n,m] \overset{\text{def}}{\equiv} \forall x^\Nat.\ n\leq x\leq m\rightarrow f(x)=f(n)\]
one often finds in literature the following equivalent version of \ref{e1}:

\begin{equation}\label{con0}\forall h^{\funn} \exists x^\Nat\ f\downarrow [x, h(x)]\end{equation}
which is the no-counterexample interpretation of
\[\exists x^\Nat \forall y^\Nat\geq x f\downarrow [x, y]\]

While the above notion of convergence \ref{con0} would be enough for our purposes, it seems not to allow straightforward compositional reasoning when one has to deal with non trivial \emph{interaction} of convergent functions. Even when there is no complex interaction, the needed reasoning is not direct. For example, one may want to prove that if two functions $f,g$ converge in the sense of \ref{con0}, one can systematically find intervals in which they are \emph{both} constant. That is, if

\[\forall h^{\funn} \exists x^\Nat f\downarrow [x, h(x)]\land \forall h^{\funn} \exists x^\Nat g\downarrow [x, h(x)]\]
 then one may want to prove that
 \[\forall h^{\funn} \exists x^\Nat f\downarrow [x, h(x)]\land g\downarrow [x, h(x)]\]
 The above implication is provable in a non overly complicated way, but when interaction increases (as we shall see in proposition \ref{proposition-merge} below), one begins to feel the need for a more suitable formulation of convergence.
 
Berardi \cite{BerPer} introduced a notion of convergence  especially suited for managing interaction of convergent functions. If one consider the formula

\[\forall z^\Nat\exists x^\Nat\geq z \forall y^\Nat\geq x f\downarrow [x, y]\]
(with the intent of expressing very redundantly the fact that there are infinite points of convergence for $f$) one obtains a very strong notion of constructive convergence by taking its no-counterexample interpretation

\[\forall z^\Nat \forall h^{\funn} \exists x^\Nat\geq z\ h(x)\geq x \rightarrow f\downarrow [x, h(x)]\]
which after skolemization becomes

\[\forall h^{\funn} \exists \alpha^\funn\geq \id \forall z^\Nat\ h(z)\geq z\rightarrow   f\downarrow [\alpha(z), h(\alpha(z))]\]  
which is equivalent to

\begin{equation}\label{con1}\forall h^{\funn}\geq\id \exists \alpha^\funn\geq \id\ \forall z^\Nat f\downarrow [\alpha(z), h(\alpha(z))]
\end{equation}
 where we have used  the notation
 
 \[\alpha^\funn \geq \id \overset{\text{def}}{\equiv} \forall x^\Nat \alpha(x)\geq x\]
 We observe that \ref{con0} and \ref{con1} are constructively equivalent. However, from a computational point of view, their realizers are quite different: the realizers of \ref{con1} are able to interact directly with each other, as we will see.
 
 We are now ready to formally define a constructive notion of convergence for sequences of numbers: a sequence of objects $f: \Nat\rightarrow A$  is convergent if for any $h^\funn\geq \id$ there are infinitely many intervals $[n,h(n)]$ in which $f$ is constant.
\begin{definition}[Convergence (Berardi \cite{BerPer})]
\label{definition-Convergence}
 Let $f: \seqin{A}$ be a closed term of $\SystemT$, with $A$ atomic type. We say that $f$ \emph{converges} if

\[\forall h^{\funn}\geq\id\ \exists \alpha^\funn\geq \id\ \forall z^\Nat f\downarrow [\alpha(z), h(\alpha(z))]\]  

\end{definition}
\textbf{Notation.} If $t: A\rightarrow B$ and $u: A$ we shall often write $t_u$ in place of  $tu$, for notational convenience or for highlighting that $t_u$ is an element of a collection of type-$B$ terms parametrized by terms of type $A$. \\

We now make explicit the constructive information associated to the above notion of convergence, through the concept of \emph{modulus of convergence}. A modulus of convergence takes an $h: \funn$ and returns an enumeration of intervals $[n,h(n)]$ in which $f$ is constant. It is a intuitionistic realizer of the notion of convergence.
\begin{definition}[Modulus of Convergence]\label{definition-mc}
 Let $f: \seqin{A}$ be a closed term of $\SystemT$, with $A$ atomic type. A term $\mathcal{M}:(\funn)\rightarrow (\funn)$ of $\SystemT$ is a \emph{modulus of convergence} for $f$ if 
 \begin{enumerate}
 \item $\forall h^{\funn}\geq\id\ \mathcal{M}_h\geq \id$\\
 \item  $\forall h^\funn\geq\id\  \forall z^\Nat f \downarrow [\mathcal{M}_h(z), h(\mathcal{M}_h(z))]$\\
 \end{enumerate}
 If $h:\funn\geq\id$ and $\forall z^\Nat f \downarrow [\mathcal{N}(z), h(\mathcal{N}(z))]$, $\mathcal{N}$ is said to be an $h$-\emph{modulus of convergence} for $f$.
 
 \end{definition}
 We observe that by definition, if one has a modulus of convergence $\mathcal{M}$ for a function $f$, one can find an infinite number of intervals of any desired length in which $f$ is constant. For example, if one wants to find an interval of length $5$, one just defines the function $h(x)=x+5$ and compute $n:=\mathcal{M}_h(0)$. Then,  $f$ is constant in $[n,n+5]$. Clearly, a modulus of convergence carries a lot of constructive information about $f$.
 
 \subsection{Intuitive Significance of the Concept of Modulus of Convergence and Restatement of the Convergence Theorem} 
 As we said, Berardi's notion of convergence works remarkably well when convergent functions interact together, for instance, in the definition of a new function. The fact that Berardi's notion is a no-counterexample interpretation of the classical notion of convergence, explains why it \emph{works}. We can intuitively describe the reasons why it does it \emph{well} as follows. 
 
 A first reason is purely computational. Given a function $h^\funn\geq \id$ and a modulus of convergence $\mathcal{M}$, we can interpret the role of $h$ in the computation of $\mathcal{M}_h$ as that of a \emph{continuation}.  Constructively, when a new convergent function is defined from other convergent functions, one will need to produce intervals in which the new function is constant. One will try to achieve the goal by finding intervals in which the functions involved in the definition are \emph{all} constant. The problem is that one may be able to find such intervals for every single function, but not for them all together. For example, if one defines the function
\[\beta:=\lambda x^\Nat f(g(x), x)\] 
then $\beta$ is convergent if $g$ and $\lambda x^\Nat f(n,x)$ are such for every choice of $n$. But an interval in which $g$ is constant need not be an interval in which $\beta$ too is constant, because we have to find some interval in which \emph{both} $g$ is constantly equal to some $m$ and $\lambda x^\Nat f(m,x)$ is constant. We solve the problem through the use of continuations.  

We start by observing that it seems there is a strict sequence of tasks to be performed. First, one tries to find an $m_1$ such that $g$ is constant in, say, $[m_1,l_1]$ with $m_1<l_1$. Then, he computes $g(m_1)=n$ and pass $n$ to a ``continuation" $h: \funn$ which returns an $h(n)=m_2<l_2$ such that $\lambda x^\Nat f(n, x)$ is constant in $[m_2, l_2]$.  If $m_1<m_2<l_1$, a non trivial interval in which $g$ is constant has been found. But if $m_2>l_1$? Then, $g$ may assume different values in all points of the interval $[m_2, l_2]$ and one cannot hope that $\beta$ is going to be convergent in $[m_2, l_2]$. We anticipate the solution contained in the proof of proposition \ref{proposition-merge}, by letting $m_1=\mathcal{M}_k(0)$, where $\mathcal{M}$ is a modulus of convergence for $g$ and, for example, $k(x)=h(g(x))+1$. Then, by definition of modulus of convergence, $g$ is constant in $[m_1, k(m_1)]$ and letting $l_1=k(m_1)$ we obtain that  \[m_2=h(n)=h(g(m_1))<h(g(m_1))+1=l_1\] as required. In other words, we use $k$ and hence $h$ as \emph{continuations}, thanks to $\mathcal{M}$. 

 The issue we are facing may be further exemplified by the following sequential game between $k$ players. Suppose there are convergent functions $f_1, f_2, \ldots, f_k$ ot type $\funn$ on the board and an arbitrarily chosen number $m$. Players make their moves in order, starting from player one and finishing with player $k$. A play of the game, is an increasing sequence of numbers $m, m_1, m_2, \ldots, m_k$, with $m_i$ the move of player $i$. Player $i$ wins if $f_i$ is constant in an interval $[m_k, l_k]$, for some $l_k>m_k$. A strategy for player $i$ is just a function $h$ over natural numbers, taking the move of the player $i-1$ (or the integer $m$ if $i=1$) and returning the move of player $i$. The fact that the winning condition depends on the move of player $k$ makes very difficult for players $1,\ldots k-1$ to win. In this game, each player hopes that  in the resulting final interval its own function will be constant but his hope is frustrated by the following ones, which are trying to accomplish the same task but with respect to their own functions.  However, player $i$ has a winning strategy effectively computable if he knows the strategies of all subsequent players $i+1, \ldots, k$ (we cannot assume the trivial winning strategy returning the point of stabilization of $f_i$ to be effectively computable, since $f_i$ is arbitrary)
 
 We are now in a position to tell another reason why moduli of convergence are so useful. A winning strategy for player $i$ can be computed by a convergence module. More precisely, it can be proved, as consequence of proposition \ref{proposition-pair}, that if players $i+1, \ldots, k$ play strategies $h_{i+1},\ldots, h_k$, then $h_i:=\mathcal{M}_{h_{k}\circ\cdots\circ h_{i+1}}$ is a winning strategy for player $i$ against $h_{i+1},\ldots, h_k$, whenever $\mathcal{M}$ is a modulus of convergence for $f_i$. Therefore, if a modulus of convergence for each function $f_1, \ldots, f_k$ is given, one can compute a particularly desirable  instance of \emph{Nash equilibrium}, that is, a sequence of functions $h_1, h_2, \ldots, h_k$ such that, if every player $i$ plays according to the strategy $h_i$, every play will be won by every player. Therefore, at the end of the interaction, every participant will have accomplished its own task.
 
 We now formulate the promised restatement of theorem \ref{theorem-con} that we shall be able to prove.
 
 \begin{theorem}[Weak Convergence]\label{theorem-weakconvergence}
 Let $t: (\funn)\rightarrow \State$ be a closed term of $\SystemT$, with $\State$ atomic type. Then we can effectively define a closed term $\mathcal{M}: ( \Nat\rightarrow (\funn))\rightarrow (\funn)\rightarrow (\funn)$ of $\SystemT$, such that the following holds: for all  $s:  \Nat\rightarrow (\funn)$, $s\in$ w.i. and numerals $n$, $\mathcal{M}s$ is a modulus of convergence for $\lambda m^\Nat t_n(s_m)$.
\end{theorem} 
  
 \subsection{Basic Operations with Moduli of Convergence}
 We now prove a couple of propositions, both to illustrate the use of moduli of convergence and to provide lemmas we will need in the following. First, we show that given two terms $f_1$ and $f_2$, if each one of them has a modulus of convergence, then there is a modulus of convergence that works simultaneously for both of them. In particular, we can define a binary operation $\sqcup$ between moduli of convergence such that, for every pair of moduli $\mathcal{M}, \mathcal{N}$, $\mathcal{M}\sqcup \mathcal{N}$ is ``more general" than both $\mathcal{M}$ and  $\mathcal{N}$. Here, for every $\mathcal{M}_1, \mathcal{M}_2$, we call $\mathcal{M}_2$ more general than $\mathcal{M}_1$, if for every term $f$, if $\mathcal{M}_1$ is a modulus of convergence for $f$ then also $\mathcal{M}_2$ is a modulus of convergence for $f$.  We this terminology, we may see $\mathcal{M}\sqcup \mathcal{N}$ as an upper bound of the set $\{\mathcal{M}, \mathcal{N}\}$, with respect to the partial order induced by the relation ``to be more general than". The construction of the pair ${\mathcal{M}}_{h\circ {\mathcal{N}}_h},\mathcal{N}_h$ below may also be seen as a Nash equilibrium for the two player version of the game we have discussed above.
 
 \begin{proposition}[Joint Convergence]\label{proposition-pair}
 Let $\mathcal{M}$ and $ \mathcal{N}$ be moduli of convergence respectively for $f_1$ and $f_2$. Define \[\mathcal{M}\sqcup \mathcal{N}:=\lambda h^\funn\lambda z^\Nat \mathcal{N}_h({\mathcal{M}}_{h\circ {\mathcal{N}}_h}(z)) \]Then  $\mathcal{M}\sqcup\mathcal{N}$ is a modulus of convergence for both $f_1$ and $f_2$.
  \end{proposition}
  \proof 
  Set \[\mathcal{L}:=\mathcal{M}\sqcup\mathcal{N}\]
  First, we check property 1 of definition \ref{definition-mc}. For all $h^\funn\geq\id$, $\mathcal{N}_h\geq \id$ by definition \ref{definition-mc} point (1) and so $h\circ \mathcal{N}_h\geq \id$. Thus, for all $h^\funn\geq\id$ and $z^\Nat$  \[\mathcal{L}_h(z)= \mathcal{N }_h({\mathcal{M}}_{h\circ {\mathcal{N}}_h}(z)) \geq  z\]
  since $\mathcal{M}$ has property (1) of definition \ref{definition-mc} and hence $  {\mathcal{M}}_{h\circ {\mathcal{N}}_h}\geq \id$. Therefore, for all $h^\funn\geq \id$, $\mathcal{L}_h\geq \id$ and we are done.\\
  Secondly, we check property 2 of definition \ref{definition-mc}. Fix a term $h^\funn\geq\id$ and a numeral $z$. We have that  \begin{equation}\label{i}f_1\downarrow [\mathcal{M}_{h\circ \mathcal{N}_h}(z), h\circ \mathcal{N}_h(\mathcal{M}_{h\circ \mathcal{N}_h}(z))]\end{equation}
  since $\mathcal{M}$ is a module of convergence for $f_1$.  Moreover, 
  
  \begin{equation}\label{i0}f_2 \downarrow [\mathcal{N}_h(\mathcal{M}_{h\circ \mathcal{N}_h}(z)), h( \mathcal{N}_h(\mathcal{M}_{h\circ \mathcal{N}_h}(z)))]\end{equation}
  since $\mathcal{N}$ is a modulus of convergence for $f_2$.  But the starting point of the interval in \ref{i0} is greater or equal to the starting point of the interval in \ref{i}, for $\mathcal{N}_h\geq \id$, while their ending points are equal. Hence also
\[f_1 \downarrow [\mathcal{N}_h(\mathcal{M}_{h\circ \mathcal{N}_h}(z)), h( \mathcal{N}_h(\mathcal{M}_{h\circ \mathcal{N}_h}(z)))] \]
and hence both $f_1$ and $f_2$ are constant in the interval $[\mathcal{L}_h(z), h(\mathcal{L}_h(z))] $ by definition of $\mathcal{L}$.
  
  \qed

We now consider a situation in which a family $\{f_n\}_{n\in\Nat}$ of convergent terms interacts with a convergent term $g$ and we show the result of the interaction is still a convergent term. In the following, we call ``object of type $A$" any closed normal term of type $A$.

\begin{proposition}[Merging of Functions]\label{proposition-merge}
Let $f:A\rightarrow (\Nat\rightarrow A)$ be a closed term, with $A$ atomic, and $\mathcal{N}: A \rightarrow (\funn)\rightarrow (\funn)$ be such that for every object $a$ of type $A$, $\mathcal{N}_a$ is a modulus of convergence for $f_a$. Let moreover $g: \Nat\rightarrow A$ and let $\mathcal{M}$ be a modulus of convergence for $g$. Define 

\[\mathcal{H}_1(\mathcal{M},\mathcal{N},g):=\lambda h^\funn\lambda z^\Nat \mathcal{N}_h'({\mathcal{M}}_{h\circ {\mathcal{N}}_{h}'}(z)) \]
with 

\[\mathcal{N}_h':=\lambda n^\Nat (\mathcal{N}_{g(n)}h)n\]
Then $\mathcal{H}_1(\mathcal{M},\mathcal{N},g)$ is a modulus of convergence for \[\lambda n^\Nat f_{g(n)}(n)\]
\end{proposition}

\proof Property (1) of definition \ref{definition-mc} follows by the same reasoning used in proposition \ref{proposition-pair}. We check property (2) of definition \ref{definition-mc}. Set $\mathcal{L}:=\mathcal{H}_1(\mathcal{M},\mathcal{N},g)$. The idea is that $\mathcal{L}$ has to produce an interval $i$ in which $g$  is constant and equal to $a$, while the interval produced by $\mathcal{N}_a$ in which $f_a$ is constant will be contained in $i$. $\mathcal{L}$ does the job by using $\mathcal{N}_h'$ as a \emph{continuation}.\\
Fix a term closed $h^\funn\geq\id$ and $z$ a numeral. We have that  \begin{equation}\label{i1}g\downarrow [\mathcal{M}_{h\circ \mathcal{N}_h'}(z), h\circ \mathcal{N}_h'(\mathcal{M}_{h\circ \mathcal{N}_h'}(z))]\end{equation}
since $\mathcal{M}$ is a module of convergence for $g$.  In particular, 
  
  \begin{equation}\label{i2}g \downarrow [\mathcal{N}_h'(\mathcal{M}_{h\circ \mathcal{N}_h'}(z)), h( \mathcal{N}_h'(\mathcal{M}_{h\circ \mathcal{N}_h'}(z)))]\end{equation}
since $\mathcal{N}_h'\geq \id$. Say that for all $n$ in the intervals in \ref{i1} and \ref{i2}, $g(n)=a$. By definition  of $\mathcal{N}_h'$
\[ [\mathcal{N}_h'(\mathcal{M}_{h\circ \mathcal{N}_h'}(z)), h( \mathcal{N}_h'(\mathcal{M}_{h\circ \mathcal{N}_h'}(z)))]\]

\begin{equation}\label{i3}= [\mathcal{N}_ah(\mathcal{M}_{h\circ \mathcal{N}_h'}(z)), h( \mathcal{N}_ah(\mathcal{M}_{h\circ \mathcal{N}_h'}(z)))]\end{equation}
Since $\mathcal{N}_a$ is a modulus of convergence for $f_a$, we have 

\[f_a \downarrow  [\mathcal{N}_ah(\mathcal{M}_{h\circ \mathcal{N}_h'}(z)), h( \mathcal{N}_ah(\mathcal{M}_{h\circ \mathcal{N}_h'}(z)))]\]
But for all $x$ in the interval \ref{i3}, \[(\lambda n^\Nat f_{g(n)}(n))x=f_a(x)\]
Hence 
\[\lambda n^\Nat f_{g(n)}(n) \downarrow[\mathcal{N}_h'(\mathcal{M}_{h\circ \mathcal{N}_h'}(z)), h( \mathcal{N}_h'(\mathcal{M}_{h\circ \mathcal{N}_h'}(z)))] \]
and so $\lambda n^\Nat f_{g(n)}(n)$ is constant in the interval $[\mathcal{L}_h(z), h(\mathcal{L}_h(z))] $ by definition of $\mathcal{L}$.

\qed

\section{Computations with non Standard Natural Numbers}\label{section-nonstandardmodel}

For technical convenience we add now to system $\SystemT$ a constant $\Phic:\funn$ with no associated reduction rules. In this way, each term $t: A$ can be viewed as functionally depending on $\mathsf{\Phi}$, but it is still considered as having type $A$, instead the more complicated $(\funn)\rightarrow A$. Of course, terms of atomic type are not in general equal to a constant or a numeral, if they contain $\Phic$.   
\begin{definition}[Evaluation at $u$]
Let $t$ be a term. For any term $u: \funn$, we denote with $t[u]$ the term $t[u/\Phic]$.

\end{definition}

Adopting this notation, what we want prove is that if $t: A$, with $A$ atomic, and $s\in$ w.i., then the function $\lambda m^\Nat t[s_m]$ constructively converges, that is, it has a modulus of convergence. A natural attempt for achieving the goal is to recursively decompose the problem. For example, suppose we want to study the convergence of the function
\[^*(+t_1t_2):=\lambda m^\Nat +t_1t_2[s_m]: \funn\]
where $s\in$ w.i.. and $+:\Nat\rightarrow \funn$ represents a constant of $\SystemT$ encoding the operation of addition of natural numbers. Since $t_1:\Nat$ and $t_2:\Nat$ may have complex structure, it is natural to recursively study the functions

\[^*t_1:=\lambda m^\Nat t_1[s_m]: \funn\]
 and 
\[^*t_2:=\lambda m^\Nat t_2[s_m]: \funn\]
 But if we want to study the function $^*(+t_1t_2)$ as a combination of $^*t_1$ and $^*t_2$, it is clear that $+$ cannot be interpreted as itself, but as a function $^*+$ of $^*t_1$ and $^*t_2$. We would like the following equation to hold
\[ ^*(+t_1t_2)=  {^*+} ^*t_1^*t_2\]
 As a consequence of our notation, also the following equation must be true for all numerals $n$
 \[ ^*(+t_1t_2)(n)=+(^*t_1(n))(^*t_2(n))\]
These considerations impose us to define 
\[^*+:=\lambda g_1^\funn \lambda g_2^\funn \lambda m^\Nat  +g_1(n)g_2(n)\]
At a first look, this may seem a rather strange way of doing computations. But it turns out that it is \emph{strongly} not the case. $^*t_1$ and $^*t_2$ may be interpreted as \emph{hypernatural} numbers and $^*+$ as the operation of addition of hypernaturals as they are defined in ultrapower non standard models of Peano Arithmetic. 

\subsection{Non Standard Models of Arithmetic}\label{subsection-nonstandardmodels}

The first non standard model of Arithmetic is due to Skolem \cite{Skolem}. The universe of that model is indeed made of functions $\NatSet\rightarrow \NatSet$, but we instead describe a variant of the Skolem construction, which is the ultrapower construction (see for example Goldblatt \cite{Goldblatt}).

 Fix a non principal ultrafilter $\mathcal{F}$ over $\NatSet$. First, define an equivalence relation $\simeq$ between functions $\NatSet \rightarrow \NatSet$ as follows:
\[f_1\simeq f_2 \iff \{x\in\NatSet\ |\ f_1(x)=f_2(x) \}\in \mathcal{F}\]
(The intuition here is that an ultrafilter collects the ``big" subsets of $\NatSet$ and hence two functions are to be considered equal if they have equal values for ``great many" arguments. For example, two functions which, as sequences, converge to the same natural number are considered equal, for they agree on a cofinite set of $\NatSet$, which must belong to every non principal ultrafilter). Secondly, define
\[^*\NatSet := (\NatSet\rightarrow \NatSet)_\simeq\]
that is, $^*\NatSet$ is the set of all natural number functions partitioned under the equivalence relation $\simeq$. Finally, set
\[\begin{aligned}^*0&:= \lambda n^\NatSet 0\\
^*S &:=\lambda n^\NatSet S(n)\\
^*+ &:=\lambda f_1^{\NatSet\rightarrow \NatSet}\lambda f_2^{\NatSet\rightarrow \NatSet}\lambda n^\NatSet f_1(n)+f_2(n)\\
^*\cdot &:=\lambda f_1^{\NatSet\rightarrow \NatSet}\lambda f_2^{\NatSet\rightarrow \NatSet}\lambda n^\NatSet f_1(n)\cdot f_2(n)
\end{aligned}\]
where $S, +, \cdot$ are the usual operations over natural numbers. In general, if one wants to define  the non standard version of a standard function $f:\NatSet^k\rightarrow \NatSet$, he simply lets

\[^*f:=\lambda f_1^{\NatSet\rightarrow \NatSet}\ldots \lambda f_k^{\NatSet\rightarrow \NatSet} \lambda n^\NatSet f(f_1(n),\ldots, f_k(n))\]

 It can be proved that the structure
\[( ^*\NatSet, {^*0}, {^*S}, {^*+}, {^*\cdot})\]
is a model of Peano Arithmetic as similar to the usual structure of natural numbers as to satisfy \emph{precisely} the same sentences which are true under the usual interpretation. Formally, it is \emph{elementarily equivalent} to the structure of natural numbers. 

Elements of $^*\NatSet$ are usually called hypernatural numbers. Since they are so similar to natural numbers, it perfectly makes sense to think about defining a model of system $\SystemT$ over hypernaturals. Indeed, Berardi \cite{Ber2005} used hypernaturals, under a weaker equivalence relation, to construct an intuitionistic model for $\Delta_0^2$ maps and Berardi and de' Liguoro \cite{BerardiLiguoroMonadi} used them to interpret a fragment of classical primitive recursive Arithmetic.

\subsection{A non Standard Model for the System $\SystemTG_0$}\label{subsection-modelforT0}

In order to approach gradually our final construction, we first give a definition of a non standard model for $\SystemTG_0$, which is G\"odel's $\SystemTG$ restricted to having only a recursion operator $\rec$ of type $\Nat\rightarrow (\Nat\rightarrow \funn)\rightarrow \funn$ and choice operator $\ifn:\Bool\rightarrow\funn\rightarrow\Nat$. Hence, $\SystemTG_0$ represents the primitive recursive functions. 

The definition of the model is purely syntactical and this is the key for our approach to go through. In fact, we are defining an internal model, that is a representation of $\SystemTG_0$ into $\SystemTG$ itself. First, define the new type structure as: 

\[\begin{aligned}^*\Nat &:= \Nat\rightarrow \Nat\\
^*\Bool &:=\Nat\rightarrow \Bool\\
^*(A\rightarrow B) &:= {^*A}\rightarrow {^*B}\\
^*(A\times B) &:= {^*A}\times {^*B}
\end{aligned}\]
From a semantical point of view, we interpret natural numbers as functions. Since the construction is syntactical, there is no need to describe an equivalence relation between those functions. But, accordingly to which equivalence relation one has in mind, the definition we are going to give will make sense or not from the \emph{semantical} point of view. For the results of this chaper, we have no utility in putting extra effort to define a model for $\SystemTG_0$, which is also a model for Peano Arithmetic. Hence, we may assume that $^*\Nat$ represents just all functions over $\NatSet$ without any partition.

Now, for every term $u: T$,  define a term $^*u$ of type $^*T$ by induction as follows

\[\begin{aligned}
{^*0}&:= \lambda m^\Nat\ 0\\
{^*\True}&:= \lambda m^\Nat\ \True\\
{^*\False}&:= \lambda m^\Nat\ \False\\
{^*\mathsf{S}} &:=\lambda f^{^*\Nat}\lambda m^\Nat\ \mathsf{S}(f(m))\\
{^*\ifn}&:=\lambda g^{^*\Bool } \lambda f_1^{^*\Nat} \lambda f_2^{^*\Nat} \lambda m^\Nat\ \ifn{f(m)}{f_1(m)}{f_2(m)}\\
{^*\rec}&:=\lambda f_1^{^*\Nat } \lambda f_2^{{^*\Nat}\rightarrow {^*\Nat}\rightarrow {^*\Nat}} \lambda g^{^*\Nat} \lambda m^\Nat\ (\rec_B f_1(\lambda n^\Nat f_2(\lambda x^\Nat n))g(m))(m)\\
{^*(x^A)}&:=x^{{^*A}}\\
{^*(ut)}&:={^*u}{^*t}\\
{^*(\lambda x^A u)}&:=\lambda x^{^*A} {^*u}\\
{^*\langle u, t\rangle}&:=\langle {^*u}, {^*t}\rangle\\
{^*(\pi_i u)}&:=\pi_i {^*u}
\end{aligned}
\]
with the type $B$ of $\rec_B$ equal to ${^*\Nat}\rightarrow (\Nat\rightarrow {^*\Nat}\rightarrow {^*\Nat})\rightarrow {\Nat}\rightarrow {^*\Nat}$.\\

The definition of the constants and the functions $^*\mathsf{S}$ and $^*\ifn$ is exactly the one used in the construction of ultrapower models of natural numbers. The definition of $^*\rec$ is different because involves higher type arguments, but it is a straightforward generalization of the ultrapower construction. Intuitively, $^*\rec f_1f_2g$ has to iterate $f_2$ a number of times given by $g$. But since $g$ is now an hypernatural number, the concept ``$g$ times'' makes no direct sense. Hence, $^*\rec$ also picks as input a number $m$, transform $g$ into $g(m)$ and iterates $f_2$ a number of times given by $g(m)$. But since $f_2$ is of type $^*\Nat\rightarrow {^*\Nat}\rightarrow {^*\Nat}$, the function given to $\rec_B$ is not directly $f_2$, but a term $\lambda n^\Nat f_2(\lambda x^\Nat n)$ that transforms $n$ into its hypernatural counterpart and gives it as the first argument of $f_2$. After all this work is done, one obtains a hypernatural
\[h:=\rec_B f_1(\lambda n^\Nat f_2(\lambda x^\Nat n))g(m)\]
So if $^*\rec$ stopped here, it would not return the right type of object. Hence, it returns $h(m)$, consistently to the fact that $g$ has been instantiated to $m$ previously.

The above construction can be generalized to G\"odel's $\SystemTG$, with a little more effort to be put in the generalization of $^*\rec$ and $^*\ifn$ to all types. A version of $\SystemT$ just manipulating hypernaturals is not enough for our purposes, and will be included in our final construction, so details are postponed to the next sections.

\subsection{System $\SystemT$ over Hypernaturals with Moduli of Convergence}

In the context of this work, we are not interested into the whole collection of hypernatural numbers, but only in those who are \emph{convergent}. Moreover, we want also to produce, for each one of these convergent hypernaturals, a modulus of convergence. The idea therefore is to put more constructive information into the model of hypernatural numbers and to define operations that preserve this information. The new objects we are going to consider are \emph{hypernatural numbers with moduli of convergence}. They can be represented as pairs 
\[\langle \mathcal{N}, f\rangle\]
where $f:\NatSet\rightarrow \NatSet$ is an hypernatural number and $\mathcal{N}: (\NatSet\rightarrow \NatSet)\rightarrow (\NatSet\rightarrow \NatSet)$ is modulus of convergence for $f$, as in definition \ref{definition-mc}. The resulting model will be a full type structure generated over this basic objects and their equivalent in the other atomic types. We will call it the \emph{model of hypernaturals with moduli}.
In the following, for any $s\in$w.i., $\model{u}$ will be the denotation of a term $u$ of $\SystemT$ in this new model and the aim of this sections is to syntactically define the interpretation function $[\![ \_ ]\!]_s$.

In order to construct such a model, we will have to define new operations that, first, generalize the ones over hypernaturals  we have previously studied and, secondly, are also able to combine moduli of convergence.

 For example, how to define the non standard version $\model{+}$ of addition? The summands are two objects of the form $\langle \mathcal{N}_1, f_1\rangle$ and $\langle \mathcal{N}_2, f_2\rangle$.  The second component of the sum will be the non standard sum 
 \[f_1{^*+}f_2:=\lambda m^\Nat f_1(m)+f_2(m)\] of $f_1$ and $f_2$. The first component will be a modulus of convergence for $f_1{^*+}f_2$, and so a simultaneous modulus of convergence for both $f_1$ and $f_2$ is enough. From proposition \ref{proposition-pair}, we know how to compute it with $\sqcup$ from $\mathcal{N}_1$ and $\mathcal{N}_2$. We can thus define

\[\model{+}\langle \mathcal{N}_1, f_1\rangle\langle \mathcal{N}_2, f_2\rangle:=\langle \mathcal{N}_1\sqcup \mathcal{N}_2,  f_1{^*+}f_2\rangle\]

We now launch into the definition of our syntactically described model for the whole system $\SystemT$. First we define the intended interpretation $\M_T$ of every type $T$.

\begin{definition}[Interpretation of Types]
For every type $T$ of system $\SystemT$, we define a type $\M_T$ by induction on $T$ as follows.

\begin{enumerate}
\item $T=A$, with $A$ atomic. Then \[\M_A:=  ((\funn)\rightarrow(\funn))\times (\Nat\rightarrow A)\]
\item $T=A\rightarrow B$. Then 
 \[\M_{A\rightarrow B}:= \M_A\rightarrow  \M_B\]
 \item $T=A\times B$. Then \[ \M_{A\times B}:= \M_A \times \M_B\]
\end{enumerate}

\end{definition}

If $A$ is atomic,  the interpretation $\M_A$ of $A$ is the set of pairs formed by a function $\NatSet \rightarrow A$ and its modulus of convergence (if it happens to have one). This is an accord with our view that whenever a $s\in$w.i. is fixed, a term $t$ of atomic type can be interpreted as a function $\lambda m^\Nat t[s_m]$ paired with a modulus of convergence. The model of hypernaturals with moduli can be seen as the collection of sets denoted by types $\M_T$, for $T$ varying on all types of $\SystemT$.

We now define a logical relation between the terms of our intended model of hypernaturals with moduli and the terms of system $\SystemT$. It formally states what properties any denotation of any term of $\SystemT$ should have. It  formalizes of our previous description of what the model should contain. 
\begin{definition} [Generalized Modulus of Convergence]
Let $t$ and $\mathcal{M}$ be  closed terms of $\SystemT$ and $s\in$w.i.. We define the relation $\mathcal{M} \gmc t$ -  representing the notion ``$\mathcal{M}$ is a \emph{generalized modulus of convergence} for $t$" -  by induction on the type $T$ of $t$ as follows:\\

\begin{enumerate}
\item $T=A$, with $A$ atomic. Let $\mathcal{M}: \M_A$.  Then
 \[\mathcal{M} \gmc t \iff \text{$\mathcal{M}=\langle\mathcal{L}, g\rangle$, $\mathcal{L}$ is a modulus of convergence for $g$ and $g\overset{\text{ext}}{=}\lambda n^\Nat t[s_n]$}\]
where we have defined $(g\overset{\text{ext}}{=}\lambda n^\Nat t[s_n]) \equiv \text{for all numerals $m$, $g(m)=t[s_m]$}$.\\
\item $T=A\rightarrow B$. Let $\mathcal{M}: \M_{A\rightarrow B}$. Then 
\[\mathcal{M} \gmc t \iff (\forall u^A.\ \mathcal{N} \gmc u \implies \mathcal{MN} \gmc tu)\]

\item $T=A\times B$. Let $\mathcal{M}: \M_{A\times B}$. Then 
\[\mathcal{M} \gmc t \iff ( \pi_0\mathcal{M} \gmc \pi_0t \land \pi_1\mathcal{M} \gmc \pi_1t)\]
\end{enumerate}

\end{definition}

 The aim of the rest of this section is to syntactically define a semantic interpretation $\model{\_}$ of the terms of $\SystemT$ into the model of hypernaturals with moduli, such that for every term $u: A$  and $s\in$w.i., $\model{u}\gmc u$. This means that, if $A$ is atomic, $u$ is evaluated in a pair $\langle\mathcal{L}, g\rangle$ such that $g\overset{\text{ext}}{=}\lambda n^\Nat u[s_n]$ and $\mathcal{L}$ is a modulus of convergence of $g$. Then, if given any term $t: (\funn)\rightarrow A$, we set $u:=t\Phic$ and consider $\model{u}$, we automatically obtain a constructive proof of theorem \ref{theorem-weakconvergence}.  
 
 In the following, we will make repeated use of the fact that the notion of generalized modulus of convergence is consistent with respect to equality.

\begin{lem}[Equality Soundness]\label{lemma-equalitysoundess}
Suppose $\mathcal{M}_1\gmc t_1$, $\mathcal{M}_1=\mathcal{M}_2$ and $t_1=t_2$. Then $\mathcal{M}_2\gmc t_2$.
\end{lem}
\proof Trivial induction on the type of $T$.

\qed

We now define a fundamental operation on moduli of convergence. The construction is a generalization of the one in proposition \ref{proposition-merge}.

\begin{definition}[Collection of Moduli Turned into a Single Modulus]\label{definition-gmc} Let $\mathcal{N}:A\rightarrow \M_T$ and $\langle \mathcal{M}, g\rangle: \M_{A}$, with $A$ atomic. We define by induction on $T$ and by cases a term $\mathcal{H}(\langle \mathcal{M}, g\rangle, \mathcal{N})$ of type $ \M_T$.
\begin{enumerate}

\item $T$ atomic. Then
\[\mathcal{H}( \langle \mathcal{M}, g\rangle, \mathcal{N}):= \langle\mathcal{H}_1(\mathcal{M}, \lambda a^A \pi_0(\mathcal{N}_a), g),\ \lambda n^\Nat f_{g(n)}(n)\rangle\]
with $f:=\lambda a^A\pi_1\mathcal{N}_a$ and $\mathcal{H}_1$ as in proposition \ref{proposition-merge}.\\
\item $T=C\rightarrow B$. Then

\[\mathcal{H}( \langle \mathcal{M}, g\rangle, \mathcal{N}):=\lambda \mathcal{L}^{\M_C}\mathcal{H}( \langle \mathcal{M}, g\rangle, \lambda a^A \mathcal{N}_a\mathcal{L}) \]

\item $T=C\times B$. Then

\[\mathcal{H}( \langle \mathcal{M}, g\rangle, \mathcal{N}):=  \langle \mathcal{H}( \langle \mathcal{M}, g\rangle, \lambda a^A\pi_0\mathcal{N}_a), \mathcal{H}( \langle \mathcal{M}, g\rangle, \lambda a^A \pi_1\mathcal{N}_a)\rangle  \]

\end{enumerate}
\end{definition}
If we call ``object of type $A$" any closed normal term of type $A$, then the role of the term $\mathcal{H}$ is to satisfy the following lemma, which is one the most important pieces of our construction. It provides a way of constructing the semantics of a term $ut$, with $t$ of atomic type $A$, if one is able to define a semantics for $t$ and for $ua$ for every object $a$ of type $A$.

\begin{lem}\label{lemma-merge}  Let $u$ and $t$ be  closed terms respectively of types $A\rightarrow T$ and $A$, with $A$ atomic. Suppose that for every object $a$ of type $A$, $\mathcal{N}_a\gmc ua$ and $\langle \mathcal{M}, g\rangle \gmc t$. Then $\mathcal{H}(\langle \mathcal{M}, g\rangle,  \mathcal{N})\gmc ut$.

\end{lem}
\proof By induction on $T$ and by cases.
\begin{enumerate}
\item $T$ atomic. We have
\[\mathcal{H}(\langle \mathcal{M}, g\rangle, \mathcal{N}):= \langle\mathcal{H}_1(\mathcal{M}, \lambda a^A \pi_0(\mathcal{N}_a), g), \lambda n^\Nat f_{g(n)}(n)\rangle\]
with \[f:=\lambda a^A\pi_1\mathcal{N}_a\] and \[g\overset{\text{ext}}{=}\lambda n^\Nat t[s_n]\]
for by hypothesis  $\langle \mathcal{M}, g\rangle\gmc t$. Moreover, for every object $a$ of type $A$ \[f_a\overset{\text{ext}}{=}\lambda n^\Nat ua[s_n]\] since by hypothesis  $\mathcal{N}_a\gmc ua$.  We must show that  \[\mathcal{H}_1(\mathcal{M}, \lambda a^A \pi_0\mathcal{N}_a, g)\] is a modulus of convergence for the function $\lambda n^\Nat f_{g(n)}(n)$ and that $\lambda n^\Nat f_{g(n)}\overset{\text{ext}}{=} \lambda n^\Nat ut[s_n]$. For this last part, indeed, for every numeral $m$, there is an object $a=g(m)$ such that \[\begin{aligned}
(\lambda n^\Nat f_{g(n)}(n))m &=f_a(m)\\
                                     &\overset{\text{ext}}{=}u[s_m](a)\\
                                     &=u[s_m](g(m))\\
                                     &\overset{\text{ext}}{=}u[s_m]((\lambda n^\Nat t[s_n])m)\\
                                     &= u[s_m](t[s_m])\\
                                     &=(\lambda n^\Nat ut[s_n])m
 \end{aligned}\]
Now, since $\langle \mathcal{M}, g\rangle \gmc t$,  $\mathcal{M}$ is a modulus of convergence for $g$. Moreover, for every object $a$ of type $A$, $\mathcal{N}_a\gmc ua$ by hypothesis, and therefore $\pi_0\mathcal{N}_a$ is a modulus of convergence for $\lambda n^\Nat\ ua[s_n]\overset{\text{ext}}{=}f_a$.  By proposition \ref{proposition-merge}, we obtain that \[\mathcal{H}_1(\mathcal{M},\lambda a^A \pi_0\mathcal{N}_a, g)\] is modulus of convergence for $\lambda n^\Nat f_{g(n)}(n)$, and we are done. \\

\item $T=C\rightarrow B$. Let $v: C$ and suppose $\mathcal{L}\gmc v$. We have to show that

\[\mathcal{H}( \langle \mathcal{M}, g\rangle, \mathcal{N})\mathcal{L}= \mathcal{H}( \langle \mathcal{M}, g\rangle, \lambda a^A\mathcal{N}_a\mathcal{L})\gmc utv \]
 But for every object $a$ of type $A$, $\mathcal{N}_a\gmc ua$. Therefore, for every object $a$ of type $A$ \[\mathcal{N}_a\mathcal{L}\gmc uav=(\lambda m^A umv)a\]
By induction hypothesis
\[\mathcal{H}( \langle \mathcal{M}, g\rangle, \lambda a^A \mathcal{N}_a\mathcal{L})\gmc (\lambda m^Aumv)t=utv\]
which is the thesis.\\

\item $T=C\times B$.  We have to show that, for $i=0,1$, 

\[\begin{aligned}\pi_i\mathcal{H}( \langle \mathcal{M}, g\rangle, \mathcal{N})
 &=\pi_i \langle \mathcal{H}( \langle \mathcal{M}, g\rangle, \lambda a^A \pi_0\mathcal{N}_a), \mathcal{H}( \langle \mathcal{M}, g\rangle, \lambda a^A \pi_1\mathcal{N}_a)\rangle\\
 &= \mathcal{H}( \langle \mathcal{M}, g\rangle, \lambda a^A \pi_i\mathcal{N}_a) \rangle\gmc \pi_i(ut)
\end{aligned} \]
 Now, for every object $a$ of type $A$, $\mathcal{N}_a\gmc ua$. Therefore, for every object $a$  
 of type $A$ \[\pi_i\mathcal{N}_a\gmc \pi_i(ua)=(\lambda m^A \pi_i(um))a\]
By induction hypothesis
\[\mathcal{H}( \langle \mathcal{M}, g\rangle, \lambda a^A \pi_i\mathcal{N}_a)\gmc (\lambda m^A \pi_i(um))t=\pi_i(ut)\]
which is the thesis.

\end{enumerate}
\qed

We are now in a position to define for each constant $c$ of $\SystemT$ a term $\model{c}$, which is intended to satisfy the relation $\model{c} \gmc c$. $\model{c}$ can be seen as the non standard version of the operation denoted by $c$.

 \begin{definition}[Generalized Moduli of Convergence for Constants]\label{definition-gmcc}
We define for every constant $c:T$ a closed term $\model{c}:\M_T$, accordingly to the form of $c$.
\begin{enumerate}
\item $c:A$, $A$ atomic. For any closed term $u$ of atomic type, define
\[\mathcal{M}_{\id,u}:=\langle\lambda h^\funn\lambda m^\Nat m, \lambda n^\Nat u\rangle\] Then

\[\model{c}:=\mathcal{M}_{\id,c}\]

\item $c=\Phic:\funn$. Let \[\mathcal{N}:=\lambda n^\Nat\langle\lambda h^\funn\lambda m^\Nat \ifthen{s_m(n)=s_{h(m)}(n)}{m}{h(m)},\lambda m^\Nat\ s_m(n)\rangle \]
Then 
\[\model{\Phic}:=\lambda  \langle \mathcal{M}, g\rangle^{\M_\Nat}\mathcal{H}( \langle \mathcal{M}, g\rangle, \mathcal{N})\]

\item $c\neq \Phic, c\neq \ifn$, $c:A_0\rightarrow \cdots \rightarrow A_m\rightarrow A$, with $A, A_i$ atomic for $i=0,\ldots,m$.\\ If $m>0$, then  define 

\[\begin{aligned}
\model{c}:=\lambda  \langle \mathcal{L}_0, g_0\rangle^{\M_{A_0}} \ldots \lambda  \langle \mathcal{L}_m, g_m\rangle^{\M_{A_m}} \langle&\mathcal{L}_0\sqcup\mathcal{L}_1\sqcup \ldots \sqcup \mathcal{L}_m,\\ &\lambda n^\Nat c(g_0(n))\ldots (g_m(n))\rangle
\end{aligned} \] 
assuming left association for $\sqcup$.
\\If $m=0$ ($c: A_0\rightarrow A$), define

\[\model{c}:= \lambda  \langle \mathcal{M}, g\rangle^{\M_{A_0}}  \langle \mathcal{M},\lambda n^\Nat c(g(n))\rangle\]

\item $c=\mathsf{R}_T$, $\mathsf{R}_T$ recursor constant with $T=A\rightarrow (\Nat\rightarrow A\rightarrow A )\rightarrow \Nat\rightarrow A$. Define
\[\mathcal{N}:=\lambda n^{\Nat}\ \mathsf{R}_U\mathcal{I}(\lambda n^\Nat  \mathcal{L}\mathcal{M}_{\id,n})n \]
with 
\[U:=\M_A\rightarrow (\Nat\rightarrow \M_A\rightarrow \M_A )\rightarrow \Nat\rightarrow M_A\] Then

\[\model{\rec_T}:=\lambda \mathcal{I}^{\M_A}\lambda \mathcal{L}^{\M_{\Nat\rightarrow A\rightarrow A }}\lambda \langle \mathcal{M}, g\rangle^{\M_\Nat}\mathcal{H}(\langle \mathcal{M}, g\rangle, \mathcal{N})\]

\item $c=\ifn_T$ with $T: \Bool\rightarrow A\rightarrow A\rightarrow A$. Define 

\[\mathcal{N}:=\lambda b^\Bool \ifthen{b}{\mathcal{L}_1}{\mathcal{L}_2}, \]
Then
\[\model{\ifn_T}:=\lambda \langle \mathcal{M}, g\rangle^{\M_\Bool}\lambda \mathcal{L}_1^{\M_A}\lambda \mathcal{L}_2^{\M_A} \mathcal{H}(\langle \mathcal{M}, g\rangle, \mathcal{N})\]
\end{enumerate}
\end{definition}

The definition of $\model{c}$ is a generalization of the operations done with hypernaturals. In case (1), we transform basic objects into their hypernatural, hyperboolean and hyperconstant counterparts (we call them hyperobjects) all paired with their trivial moduli of convergence.

 In case (2), the interpretation $\model{\Phic}$ of $\Phic$ is obtained by first defining uniformly on the numeral parameter $n$  a collection of interpretations  $\mathcal{N}_n=\model{\Phic n}$ of $\Phic n$, and then using the term $\mathcal{H}$ to put together the interpretations in such a way that $\model{\Phic}\langle \mathcal{M}, g\rangle$ interprets $\model{\Phic t}$ whenever $\langle \mathcal{M}, g\rangle=\model{t}$. 
 
  In case (3), we provide the non standard version of the function represented by $c$, which is a function $\model{c}$ which combines both hyperobjects and their moduli of convergence.
  
   In case (4) and (5) we have generalized the ideas of subsection \ref{subsection-modelforT0}. In particular, for $T=A\rightarrow (\Nat\rightarrow A\rightarrow A )\rightarrow \Nat\rightarrow A$ and $A$ atomic, the definition of $\model{\rec_T}$ is exactly the same of $^*\rec$ in subsection \ref{subsection-modelforT0}, enriched with the information of how to combine moduli of convergence. In fact, if we consider the term 
\[\begin{aligned}&\model{\rec_T}\mathcal{I} \mathcal{L}\langle \mathcal{M}, g\rangle\\
=&\langle\mathcal{H}_1(\mathcal{M}, \lambda n^\Nat \pi_0(\mathcal{N}_n), g),\ \lambda n^\Nat f_{g(n)}(n)\rangle
\end{aligned}\]
with $f:=\lambda n^\Nat \pi_1(\mathcal{N}_n)$, 
its right projection is equal to
\[\lambda n^\Nat f_{g(n)}(n)\]
which is equal to 
\[\lambda n^\Nat (\pi_1\mathsf{R}_U\mathcal{I}(\lambda n^\Nat  \mathcal{L}\mathcal{M}_{\id,n})g(n))(n)\]
which correspond exactly to the term 
\[^*\rec f_1f_2g= \lambda n^\Nat\ (\rec_B f_1(\lambda n^\Nat f_2(\lambda x^\Nat n))g(n))(n)\]
of subsection  \ref{subsection-modelforT0}.

We now prove that for any constant $c$, $\model{c}$ is a generalized modulus of convergence for $c$.

\begin{proposition}\label{proposition-constants}
For every constant $c$, $\model{c} \gmc c$.
\end{proposition}

\proof We proceed by cases, accordingly to the form of $c$.
\begin{enumerate}
\item $c=\Phic$. Let $t:\Nat$ and suppose $\langle \mathcal{M}, g\rangle\gmc t$. We have to prove that \[\model{\Phic}\langle \mathcal{M}, g\rangle \gmc \Phic t\]  By definition \ref{definition-gmcc} of $\model{\Phic}$
\[\model{\Phic}\langle \mathcal{M}, g\rangle= \mathcal{H}(\langle \mathcal{M}, g\rangle, \mathcal{N}) \]
with  

\[\mathcal{N}:=\lambda n^\Nat \langle\lambda h^\funn\lambda m^\Nat \ifthen{s_m(n)=s_{h(m)}(n)}{m}{h(m)}, \lambda m^\Nat s_m(n) \rangle\]
Since $\langle \mathcal{M}, g\rangle \gmc t$, if we prove that  for every numeral $n$, $\mathcal{N}_n \gmc \Phic n$, we obtain by lemma \ref{lemma-merge} that $\mathcal{H}(\langle \mathcal{M}, g\rangle, \mathcal{N})\gmc \Phi t$  and we are done. So let us show that, given a numeral $n$, $\pi_0\mathcal{N}_n$ is a modulus of convergence for the function \[\pi_1\mathcal{N}_n=\lambda m^\Nat s_m(n)\overset{\text{ext}}{=}\lambda m^\Nat\ \Phic(n)[s_m]\]
We have to prove that given any closed term $h^\funn\geq\id$ and numeral $n$,

\[\lambda m^\Nat s_m(n) \downarrow [(\pi_0\mathcal{N}_n)_h(z), h((\pi_0\mathcal{N}_n)_h(z))]\]
We have two possibilities:\\ i) $s_z(n)=s_{h(z)}(n)$. Since $s\in$ w.i., we have either $s_{h(z)}(n)=0$ and so
\[\forall y^\Nat.\  z\leq y\leq h(z)\implies  s_y(n)=0\]  
or $s_z(n)=s_{h(z)}(n)\neq 0$ and so

\[\forall y^\Nat.\  z\leq y\leq h(z)\implies  s_y(n)=s_z(n)\]  
Therefore

\[\begin{aligned}\lambda m^\Nat s_m(n) \downarrow &[z, h(z)]\\
                                                         =&[(\pi_0\mathcal{N}_n)_h(z), h((\pi_0(\mathcal{N}_n)_h(z))]
\end{aligned}
\] 
by definition of $\mathcal{N}$.\\
ii) $s_z(n)\neq s_{h(z)}(n)$.  Since $s\in$ w.i., we have $s_z\leq s_{h(z)}$ and hence $0=s_z(n)$. So $s_{h(z)}(n)=s_{h(h(z))}(n)$ and as above  

\[\begin{aligned}\lambda m^\Nat s_m(n) \downarrow &[h(z), h(h(z))]\\
                                                         =&[(\pi_0\mathcal{N}_n)_h(z), h((\pi_0\mathcal{N}_n)_h(z))]\end{aligned}\]

\item $c\neq \Phic $, $c\neq \ifn$, $c:A_0\rightarrow \cdots \rightarrow A_m\rightarrow A$.\\
i) $m>0$. Suppose $t_i: A_i$  and $\langle \mathcal{L}_i, g_i\rangle\gmc t_i$ for all $i=0,\ldots, m$. We have to prove that 
\[\model{c}\langle \mathcal{L}_0, g_0\rangle\ldots \langle \mathcal{L}_m, g_m\rangle\gmc ct_1\ldots t_m\]
We have that $g_i\overset{\text{ext}}{=}\lambda n^\Nat t_i[s_n]$ for $i=0,\ldots, m$. Moreover, since  by definition \ref{definition-gmcc} of $\model{c}$  
\[\begin{aligned}&\model{c}\langle \mathcal{L}_0, g_0\rangle\ldots \langle \mathcal{L}_m, g_m\rangle\\
=&\langle\mathcal{L}_0\sqcup\mathcal{L}_1\sqcup \ldots\sqcup\mathcal{L}_m,\lambda n^\Nat c(g_0(n))\ldots (g_m(n))\rangle\end{aligned}\]
 we must show that 
\[\begin{aligned} &\pi_0(\model{c}\langle \mathcal{L}_0, g_0\rangle\ldots \langle \mathcal{L}_m, g_m\rangle)\\
=&\mathcal{L}_0\sqcup \mathcal{L}_1\sqcup\ldots \sqcup\mathcal{L}_m
\end{aligned}\] is a modulus of convergence for  
\[\begin{aligned}&\lambda n^\Nat c(g_0(n))\ldots (g_m(n))\\
                          \overset{\text{ext}}{=}&\lambda n^\Nat c(t_1[s_n])\ldots (t_m[s_n])\\
                          =&\lambda n^\Nat ct_1\ldots t_m[s_n]
\end{aligned}\]
Since for $i=0,\ldots, m$, $\mathcal{L}_i$ is a modulus of convergence for $g_i$,  by repeated application of proposition \ref{proposition-pair} we deduce that $\mathcal{L}_0\sqcup\mathcal{L}_1\sqcup \ldots \sqcup\mathcal{L}_m$ is modulus of convergence for all $g_1,\ldots, g_m$ simultaneously. Hence for all closed terms $h:\funn\geq \id$ and numerals $z$, and for $i=0,\ldots, m$

\[ g_i \downarrow [(\mathcal{L}_0\sqcup \mathcal{L}_1\sqcup\ldots \sqcup\mathcal{L}_m)_h(z), h((\mathcal{L}_0\sqcup \mathcal{L}_1\sqcup\ldots \sqcup\mathcal{L}_m)_h(z))]\]
and therefore

\[\lambda m^\Nat\ c(g_1(m))\ldots (g_n(m))\downarrow [(\mathcal{L}_0\sqcup \mathcal{L}_1\sqcup\ldots \sqcup\mathcal{L}_m)_h(z), h((\mathcal{L}_0\sqcup \mathcal{L}_1\sqcup\ldots \sqcup\mathcal{L}_m)_h(z))]\]
which is the thesis.\\
ii) $m=0$. Straightforward simplification of the argument for i).\\

\item $c: A$, $A$ atomic. By definition \ref{definition-gmcc} 
\[\model{c}=\langle\lambda h^\funn\lambda m^\Nat m, \lambda n^\Nat c\rangle\]
We have therefore to prove that  $\lambda h^\funn\lambda m^\Nat m$ is a modulus of convergence for $\lambda n^\Nat c$, which is trivially true, and that  $\lambda n^\Nat c[s_n]=\lambda n^\Nat c$, which is also trivial. We conclude $\model{c}\gmc c$.\\

\item $c=\mathsf{R}_T$, $\mathsf{R}_T$ recursor constant with $T=A\rightarrow (\Nat\rightarrow A\rightarrow A )\rightarrow \Nat\rightarrow A$. Suppose $\mathcal{I}\gmc u: A$, $\mathcal{L}\gmc v: \Nat\rightarrow A\rightarrow A $ and  $\langle \mathcal{M}, g\rangle\gmc t:\Nat$.  We have to prove that 
\[\model{\rec_T}\mathcal{I}\mathcal{L}\langle \mathcal{M}, g\rangle=\mathcal{H}(\langle \mathcal{M}, g\rangle, \mathcal{N})\gmc \rec_Tuvt
\]
where
\[\mathcal{N}:=\lambda n^{\Nat}\ \rec_U\mathcal{I}(\lambda n^\Nat  \mathcal{L}\mathcal{M}_{\id,n})n \]
If we show that for all numerals $n$, $\mathcal{N}_n\gmc  \rec_Tuvn$, by lemma \ref{lemma-merge} we obtain that  \[\mathcal{H}(\langle \mathcal{M}, g\rangle, \mathcal{N})\gmc \rec_Tuvt\] We prove that by induction on $n$. \\
If $n=0$, then \[\mathcal{N}_0= \rec_U\mathcal{I}(\lambda n^\Nat  \mathcal{L}\mathcal{M}_{\id,n})0=\mathcal{I} \gmc u= \rec_Tuv0\]
If $n=\suc(m)$, then

\[\begin{aligned}\mathcal{N}_{\suc(m)}
&=\mathsf{R}_U\mathcal{I}(\lambda n^\Nat  \mathcal{L}\mathcal{M}_{\id,n})\suc(m)\\
&=(\lambda n^\Nat  \mathcal{L}\mathcal{M}_{\id,n})m(\mathsf{R}_U\mathcal{I}(\lambda n^\Nat  \mathcal{L}\mathcal{M}_{\id,n})m)\\
&= \mathcal{L}\mathcal{M}_{\id,m}(\mathsf{R}_U\mathcal{I}(\lambda n^\Nat  \mathcal{L}\mathcal{M}_{\id,n})m)\\
&=\mathcal{L}\mathcal{M}_{\id,m}\mathcal{N}_m
\end{aligned}\]
By induction hypothesis, $\mathcal{N}_m\gmc \mathsf{R}_Tuvm$. Moreover, $\mathcal{M}_{\id,m}\gmc m $ and by hypothesis $\mathcal{L}\gmc v$. Hence 

\[\mathcal{L}\mathcal{M}_{\id,m}\mathcal{N}_m\gmc vm(\mathsf{R}_Tuvm)=\mathcal{R}_Tuv\mathsf{S}(m)\]
which is the thesis.\\

\item $c=\mathsf{if}_T$, with $T: \Bool\rightarrow A\rightarrow A\rightarrow A$. Suppose $\mathcal{L }_1\gmc u_1: A$, $\mathcal{L }_2\gmc u_2: A$ and $\langle \mathcal{M}, g\rangle\gmc t: \Bool$. We have to prove that 
\[\model{\ifn_T}\langle \mathcal{M}, g\rangle\mathcal{L}_1\mathcal{L}_2=\mathcal{H}(\langle \mathcal{M}, g\rangle, \mathcal{N})\gmc \mathsf{if}_Ttu_1u_2\]
where
\[\mathcal{N}:=\lambda b^\Bool \ifthen{b}{\mathcal{L}_1}{\mathcal{L}_2}\]
If we show that for all $a\in \{\True, \False\}$, $\mathcal{N}_a\gmc (\lambda b^\Bool \mathsf{if}_Tbu_1u_2)a$, by lemma \ref{lemma-merge} we obtain that  \[\mathcal{H}(\langle \mathcal{M}, g\rangle, \mathcal{N})\gmc (\lambda b^A \mathsf{if}_Tbu_1u_2)t= \mathsf{if}_Ttu_1u_2\]
We prove that by cases.\\
If $a=\True$, then
\[\mathcal{N}_a=\mathcal{L}_1\gmc u_1= (\lambda b^A \mathsf{if}_Tbu_1u_2)a\]
If $a=\False$, then
\[\mathcal{N}_a=\mathcal{L}_2\gmc u_2= (\lambda b^A \mathsf{if}_Tbu_1u_2)a\]
Hence, we have the thesis.

\end{enumerate}
\qed

We are finally ready to define the interpretation of every term of $\SystemT$ in our model of hypernaturals with moduli.

\begin{definition}[Generalized Moduli of Convergence for Terms of $\SystemT$] For every term $v: T$ of  system $\SystemT$ and $s\in$w.i., we define a term $\model{v}: \M_T$ by induction on $v$ and by cases as follows:
\begin{enumerate}

\item  $v=c$, with $c$ constant. We define $\model{c}$ as in definition \ref{definition-gmcc}.\\

\item $v=x^A$, $x$ variable. Then
\[\model{x^A}:=x^{\M_A}\]

\item $v=ut$. Then
\[\model{ut}:=\model{u} \model{t}\]

\item $v=\lambda x^A u$. Then
\[\model{\lambda x^A u}:=\lambda x^{\M_A}\ \model{u}\]

\item $v=\langle u, t\rangle$. Then
\[\model{\langle u, t\rangle}:=\langle \model{u}, \model{t}\rangle\]

\item $v=\pi_iu$. Then 
\[\model{\pi_i u}:=\pi_i \model{u}\]

\end{enumerate}

\end{definition}

\section{Adequacy Theorem}\label{section-adeguacy}
We are now able to prove our main theorem. For every closed term $u$, $\model{u}$ is an inhabitant of the model of hypernaturals with moduli of convergence.
\begin{theorem}[Adequacy Theorem]\label{theorem-adeguacy}
 Let $w: A$ be a term of $\SystemT$ and
let $x_1^{A_1},\ldots,x_n^{A_n}$ contain
all the free variables
of $w$. Then, for all $s\in$w.i.
\[\lambda x_1^{\M_{A_1}}\ldots \lambda x_n^{\M_{A_n}} \model{w}\gmc \lambda x_1^{A_1}\ldots \lambda x_n^{A_n} w\]
\end{theorem}
\proof
Let  $t_1: A_1,\ldots,t_n: A_n$ be arbitrary terms. We have to prove that
\[\mathcal{M}_1\gmc t_1,\ldots, \mathcal{M}_n\gmc t_n \implies \model{w}[\mathcal{M}_1/x_1^{\M_{A_1}}\ldots \mathcal{M}_n/x_n^{\M_{A_n}}]\gmc w[t_1/x_1^{A_1}\ldots t_n/x_n^{A_n}]\]For any term $v$, we set
\[\overline{v}:=v[t_1/x_1^{A_1}\cdots t_n/x_n^{A_n}]\] and \[\overline{\model{v}}:=\model{v}[\mathcal{M}_1/x_1^{\M_{A_1}}\ldots \mathcal{M}_n/x_n^{\M_{A_n}}]\] With that notation, we have to prove that $\overline{\model{w}}\gmc \overline{w}$. The proof is by
induction on $w$ and proceeds by cases, accordingly to the form of $w$.\\
\begin{enumerate}
\item $w=c$, with $c$ constant. Since $\model{c}$ is closed and $c$ does not have free variables, by proposition \ref{proposition-constants}
\[\overline{\model{w}}=\overline{\model{c}}=\model{c}\gmc c=\overline{c}=\overline{w}\]
which is the thesis.\\
\item

$w=x_i^{A_i}$, for some $1\leq i \leq n$. Then
\[\overline{\model{w}}=x_i^{\M_{A_i}}[\mathcal{M}_1/x_1^{\M_{A_1}}\ldots \mathcal{M}_n/x_n^{\M_{A_n}}]=\mathcal{M}_i\gmc t_i=x_i^{A_i}[t_1/x_1^{A_1}\cdots t_n/x_n^{A_n}]=\overline{w}\]
which is the thesis.\\
\item
 $w=ut$. By induction hypothesis, $\overline{\model{u}}\gmc\overline{u}$ and $\overline{\model{t}}\gmc\overline{t}$. So
\[\overline{\model{ut}}=\overline{\model{u}}\overline{\model{t}} \gmc\overline{u}\overline{t}=\overline{w}\]
which is the thesis.\\

\item
  $w=\lambda x^{A}u$. Let $t: A$ and suppose $\mathcal{M}\gmc t$. We have to prove that
$\overline{\model{w}}\mathcal{M}\gmc \overline{w}t$. By induction hypothesis \[\overline{\model{\lambda x^A u}}\mathcal{M}=(\lambda x^{\M_A}\overline{\model{u}})\mathcal{M}=\overline{\model{u}}[\mathcal{M}/x^{\M_A}]\gmc \overline{u}[t/x^{A}]=\overline{w}t \] 
which is the thesis.\\
\item  $w=\langle u_0,u_1\rangle$.  By induction hypothesis, $\overline{\model{u_0}}\gmc \overline u_0$ and $\overline{\model{u_1}}\gmc \overline u_1$. Therefore, for $i=0,1$
\[\pi_i\overline{\model{\langle u_0, u_1\rangle}}=\pi_i\langle\overline{\model{u_0}}\overline{\model{u_1}}\rangle=\overline{\model{u_i}}\gmc \overline{u}_i=\pi_i\overline{w}\]
which is the thesis.\\
\item $w=\pi_iu$, with $i\in\{0,1\}$. By induction hypothesis, $\overline{\model{u}}\gmc\overline{u}$. Therefore,
\[\overline{\model{\pi_i u}}=\pi_i\overline{\model{u}}\gmc \pi_i\overline{u}=\overline{w}\]
which is the thesis.
\end{enumerate}
\qed

\section{Consequences of the Adequacy Theorem}\label{section-consequences} 
In this section, we spell out most interesting consequences of adequacy theorem. 
\subsection{Weak Convergence Theorem}

We can finally prove the constructive version of theorem \ref{theorem-con}, our main goal. The following theorem is even stronger of the previously enunciated  theorem \ref{theorem-weakconvergence}, because it states that one can find moduli of convergence for any uniformly defined collection of terms.

\begin{theorem}[Weak Convergence Theorem for Collection of Terms]\label{theorem-c}
Let $t: \Nat\rightarrow (\funn)\rightarrow \State$ be a closed term of $\SystemT$ not containing $\mathsf{\Phi}$, with $\State$ atomic type. Then we can effectively define a closed term $\mathcal{M}: \Nat\rightarrow  ( \Nat\rightarrow (\funn))\rightarrow (\funn)\rightarrow (\funn)$ of $\SystemT$, such that for all  $s: \Nat\rightarrow (\funn)$, $s\in$ w.i. and numerals $n$, $\mathcal{M}_ns$ is a modulus of convergence for $\lambda m^\Nat t_n(s_m)$.
\end{theorem}

\proof Let  \[\mathcal{M}:=\lambda y^\Nat \lambda s^{\seqfun}\pi_0((\lambda x^{\M_\Nat}\model{t_{x^\Nat}\mathsf{\Phi}})\mathcal{M}_{\id, y})\] By the adequacy theorem \ref{theorem-adeguacy},
\[\lambda x^{\M_\Nat}\model{t_{x^\Nat}\mathsf{\Phi}}\gmc \lambda x^\Nat t_{x^\Nat}\mathsf{\Phi}\]
Since for every numeral $n$, $\mathcal{M}_{\id, n}\gmc n$, we have 
\[(\lambda x^{\M_\Nat}\model{t_{x^\Nat}\mathsf{\Phi}})\mathcal{M}_{\id, n}\gmc  t_n\mathsf{\Phi}\]
By definition of generalized modulus of convergence and of $\mathcal{M}$, $\mathcal{M}_ns$ is a modulus of convergence for $\lambda m^\Nat t_n\mathsf{\Phi}[s_m]=\lambda m^\Nat t_n(s_m)$.

\qed

\subsection{Learning Based Realizability and Provably Total Functions of $\PA$ }

If $\sigma$ is a state of knowledge, we mantain the notation \[u[\sigma]:=[\varphi_P\sigma/\Phic_P\ \chi_P\sigma/\Chi_P\ \add_P\sigma/\Add_P]\] of chapter \ref{chapter-learningbasedrealizability}, since there is no confusion with the notation $u[v]$ of the previous section, which was defined for $v$ of type $\funn$.
\begin{theorem}[Zero Theorem for Collections of Terms of $\SystemTClass$]\label{theorem-zerotclass}
Let $t:\Nat \rightarrow \State$ be a closed term of $\SystemTClass$, where $\State$ is the atomic type of knowledge states. Then, there exists a term  $\mathsf{Zero}: \Nat\rightarrow \State$ of $\SystemTLearn$ such that for every numeral $n$, $t_n[\mathsf{Zero}_n]=\varnothing$.
\end{theorem}

\proof 
We may assume $t$ contains as oracles only constants $\Chi_P, \Phic_P, \Add_P$ for some fixed predicate $P$. The general case is analogous and involves only a little bit more of trivial coding. \\
In the first place, we have to carry out some simple coding. We have to represent states of knowledge by terms of type $\funn$. This is straightforward since a state represents a function over $\NatSet$.  Define
\[f :=\lambda \sigma^\State\lambda n^\Nat \ifthen{\chi_P\sigma n}{\varphi_P\sigma n+1}{0}\]
$f$ takes a state $\sigma$ and returns the function coding it; if $n$ is a numeral, when $\chi_P\sigma n=\True$, $f_\sigma n$ returns $\varphi_P\sigma n+1$ and not just $\varphi_P\sigma n$, in order to ensure that $0$ is returned only when it is just a trivial value, i.e. when $\chi_P\sigma=\False$.  \\
Given a term $g:\funn$, intended to represent a state $\sigma$, define terms $\varphi_P^g, \chi_P^g,\add_P^g$, intended to code respectively  $\varphi_P\sigma, \chi_P\sigma, \add_P\sigma$, as follows
\[\begin{aligned}
\varphi_P^{g}&:=\lambda n^\Nat \ifthen{g(n)=0}{0}{g(n)-1}\\
\chi_P^g &:= \lambda n^\Nat \ifthen{g(n)=0}{\False}{\True}\\
\add_P^g &:= \lambda n^\Nat \lambda m^\Nat \ifthen{g(n)=0}{\add_P\varnothing nm}{\varnothing}
\end{aligned}\]
Moreover, for every term $u$, define
\[u^g:=u[\varphi_P^g/\Phic_P\ \chi_P^g/\Chi_P\ \add_P^g/\Add_P]\]
It is easy to see that for all terms $u$ and all states $\sigma$, $u[\sigma]=u^{f_\sigma}$. \\
We can now give the important part of the argument. Fix a numeral $n$. Let $\sigma_m:=\varnothing$ and $\sigma_{m+1}:=\sigma_m\Cup t_n[\sigma_m]$ be a recursive definition of a sequence of states (which can be coded in $\SystemT$).  Our goal is to write down a term of system $\SystemTLearn$ which is able to find a state $\sigma_{k+1}$ such that $t_n[\sigma_k]=t_n[\sigma_{k+1}]$: as we will see, this condition implies $t[\sigma_{k+1}]=\varnothing$. Let
\[\update:=\lambda m^\Nat \lambda g^{\funn} t_m^{g}\]
By the weak convergence theorem  \ref{theorem-c} applied to $\update$,  there exists a term $\mathcal{M}$ of $\SystemTLearn$ such that, for every numeral $l$ and $s\in $ w.i., $\mathcal{M}_ls$ is a modulus of convergence for $\lambda m^\Nat \update_l(s_m)$. Set $s:=(\lambda m^\Nat f_{\sigma_m})$. Then $s\in$ w.i., since $\sigma_0\leq \sigma_1\leq \sigma_2\cdots$. Therefore $\mathcal{M}_ns$ is a modulus of convergence for $\lambda m^\Nat\ \update_n(s_m)$. If we choose $h:=\lambda m^\Nat m+1$ and set $k:=\mathcal{M}_nsh$, we have that $\update_n(s_k)=\update_n(s_{k+1})$ by definition of modulus of convergence. We thus obtain, by definition of $\update$ ad $s$, that $t_n^{f_{\sigma_k}}=t_n^{f_{\sigma_{k+1}}}$ and so $t_n[\sigma_k]=t_n[\sigma_{k+1}]$. Let
  \[\mathsf{Zero}:= \lambda m^\Nat \sigma_{(\mathcal{M}_msh)+1}\] 
  We have 

\[\begin{aligned}\mathsf{Zero}_n\Cup t_n[\mathsf{Zero}_n]=&\sigma_{k+1}\Cup t_n[\sigma_{k+1}]\\
&=(\sigma_k\Cup t_n[\sigma_k])\Cup t_n[\sigma_{k+1}]\\
&=(\sigma_k\Cup t_n[\sigma_k])\Cup t_n[\sigma_{k}]\\
 &=\sigma_k\Cup t_n[\sigma_k]\\
 &=\sigma_{k+1}\\
 &=\mathsf{Zero}_n
  \end{aligned}\]
  Since $t_n[\mathsf{Zero}_n]$ is consistent and disjoint with $\mathsf{Zero}_n$, we conclude $t_n[\mathsf{Zero}_n]=\varnothing$ and obtain the thesis.

\qed

\comment{We used a version of G\"odel's system $\SystemT$ over hypernaturals with moduli as a convenient computational device. We have interpreted the calculations we needed to perform as operations over non standard integers, and with this interpretation in mind we have been able to explain  every detail of our syntactic construction. }
As anticipated in the introduction, from learning based realizers we can extract algorithms of system $\SystemTG$.

\begin{thm}[Program Extraction via Learning Based Realizability]\label{theorem-extraction} Let $t$ be a term of $\SystemTClass$ and suppose that $t\Vvdash \forall x^\Nat \exists y^\Nat Pxy$, with $P$ atomic. Then, from $t$ one can effectively define a term $u$ of G\"odel's system $\SystemTG$ such that for every numeral $n$, $Pn(un)=\True$.
\end{thm}

\proof Let \[v:=\lambda m^\Nat \pi_1(tm)\]
$v$ is of type $\Nat\rightarrow \State$. By theorem \ref{theorem-zerotclass}, there exists a term $\mathsf{Zero}: \Nat\rightarrow \State$ of $\SystemTLearn$ such that $v_n[\mathsf{Zero}_n]=\varnothing$ for every numeral $n$. Define 
\[w:=\lambda m^\Nat (\pi_0(tm)[\mathsf{Zero}_m])\]
and fix a numeral $n$. By unfolding the definition of realizability with respect to the state $\mathsf{Zero}_n$, we have that 

\[tn\Vvdash_{\mathsf{Zero}_n} \exists y^\Nat Pny\] 
and hence
\[\pi_1(tn)\Vvdash_{\mathsf{Zero}_n}\ Pn(wn)\]
that is to say
\[v_n[\mathsf{Zero}_n]=\varnothing\implies Pn(wn)=\True\]
and therefore \[Pn(wn)=\True\] We observe that $w:\funn$ and $w$ is a term of $\SystemTLearn$. By standard coding of states into natural numbers and of all other constants of $\SystemTLearn$ into terms of G\"odel's $\SystemTG$, one can code every term of $\SystemTLearn$ into G\"odel's $\SystemTG$. Hence, there exists a term $u$ of G\"odel's $\SystemTG$ such that for all numerals $n$, $u(n)=w(n)$, which is the thesis. 

\qed

\begin{remark}
We point out that the term $u$ of theorem \ref{theorem-extraction}, modulo some trivial coding of states into numbers, bears a strong resemblance to the term $t$ from which it is defined. In particular, $u$ is  straightforwardly obtained from a modulus of convergence $\mathcal{M}$ that carries the constructive information associated to the convergence of $t$. In turn, $\mathcal{M}$ is obtained via the translation $\model{\_}$, which just replaces  type-$A$ constants and variables of $t$ with new terms of type $\M_A$. As an instance, recursion constants $\rec_A$ are replaced with recursion constants  $\rec_{\M_A}$. Therefore, the type of recursion goes through a constant increase of $2$, since $\M_A$ is obtained by changing the basic types $C$ with $\M_C$.\\ We conjecture it is not possible to amend our translation as to preserve the types of recursion constants. This should be due to Avigad's theorem: if one is able to find zeros of finite update procedures, than one can compute every provably total function of $\PA$. But in the next subsections, we show how to compute finite zeros of update procedures in G\"odel's system $\SystemTG$, thanks to moduli of convergence. If the translation $\model{\_}$ did not increase the recursion type, then the term computing the zero of an update procedure would have the same recursion level of the latter. But since primitive recursive functionals are enough as update procedures, one would get a contradiction to Avigad's theorem, because one could interpret all provably total functions of Arithmetic with primitive recursion.

\end{remark}
We are now able to prove a version of the classic theorem of G\"odel, characterizing the class of functions provably total in $\PA$ as the class of functions representable in system $\SystemTG$. 

\begin{theorem}[Provably Total Functions of $\PA$] If $\PA\vdash \forall x^\Nat \exists y^\Nat Pxy$, then there exists a term $u$ of G\"odel's system $\SystemTG$ such that for every numeral $n$, $Pn(un)=\True$. 
\end{theorem}
\proof Starting from the assumption that  
\[\PA\vdash \forall x^\Nat \exists y^\Nat Pxy\]
by Kolmogorov double negation translation (see for instance \cite{Sorensen}), we have that 
\[\HA\vdash \forall x^\Nat \lnot\lnot \exists y^\Nat Pxy\]
Therefore
\[\HA+\EM_1\vdash \forall x^\Nat \exists y^\Nat Pxy\]
and so there is a term $t$ of $\SystemTClass$ such that
\[t\Vvdash\forall x^\Nat \exists y^\Nat Pxy \]
By theorem \ref{theorem-extraction}, there exists a term $u$ of G\"odel's system $\SystemTG$ such that for all numerals $n$
\[Pnu(n)=\True\]

\qed

\subsection{Zeros for Update Procedures}

Thanks to the adeguacy theorem, we are able to give a new constructive proof of Avigad's theorem for update procedures. Here, we give a slightly different definition of update procedure. This is not a limitation, since the update procedures which are \emph{actually} used by Avigad \cite{Avigad} in proving 1-consistency of $\PA$ still fall under the following definition.

\begin{definition}[Update Operator, Typed Update Procedure]\label{definition-finiteupdateproc}
Fix a primitive recursive bijective  coding $| \_ |: (\NatSet^3\cup \{\varnothing\}) \rightarrow \NatSet$ of  $\varnothing$ and of triples of natural numbers into natural numbers. Define a binary operation $\oplus$ which combine functions $f: \funn$ and codes of the form $|(1, n ,m)|$ of pairs of natural numbers and returns a function $\funn$ as follows
\[f \oplus |(1,m,n)|:=\lambda x^\Nat \ifthen{x=m}{n}{f(x)}\]
For convenience, define also $f\oplus |\varnothing|=f$.

 A \emph{typed update procedure of ordinal} $k\in\NatSet$ (also said a $k$-ary typed update procedure) is a term $\update: (\funn)^k\rightarrow \Nat$ of G\"odel's $\SystemTG$ such that the following holds: 
\begin{enumerate}
\item for all sequences ${f}={f}_1,\ldots, f_k$ of closed type-$\funn$ terms of $\SystemTG$, $\update {f}=|(i,n,m)| \implies 1\leq i\leq k$.\\
\item for all sequences ${f}={f}_1,\ldots, {f}_k$ and ${g}={g}_1,\ldots, {g}_k$ of closed type-$\funn$ terms of $\SystemTG$ and for all $1\leq i<k$, if \\

i) for all $j<i$, ${f}_j={g}_j$; \\

ii) $\update {f}=|(i,n,m)|$, ${g}_i(n)=m$ and $\update {g}=|(i,h,l)|$\\\\
then $h\neq n$. 
\end{enumerate}
If $\update$ is a $k$-ary typed update procedure, a \emph{zero} for $\update$ is a sequence ${f}={f}_1,\ldots, f_k$ of closed type-$\funn$ terms of $\SystemTG$ such that $\update f=|\varnothing|$.

\end{definition}
Every unary update procedure gives rise to a learning process, i.e. a weakly increasing  chain of functions.
\begin{proposition}[Learning Processes from Unary Update Procedures]\label{proposition-updateincreasing}Let $\update$ be a unary update procedure and define by recursion
$s_0:=0^\funn:=\lambda x^\Nat 0$ and $s_{k+1}:=s_k\oplus \update s_k$. Then $s\in$w.i..
\end{proposition}
\proof Suppose $s_i(n)=m\neq 0$. We have to prove that for all $j$,  $ s_{i+j}(n)=m$. We proceed by induction on $j$. Suppose $j>0$. Since $s_0=0^\funn$ and $s_i(n)\neq 0$,  it must be that for some $i_0<i$, $\update s_{i_0}=|(1, n, m)|$. By induction hypothesis, $s_{i+j-1}(n)=m$. By definition  \ref{definition-finiteupdateproc}, point  2-ii), $\update s_{i+j-1}\neq |(1, n, l)|$ for all $l$. Since \[s_{i+j}=s_{i+j-1}\oplus \update s_{i+j-1}\] it must be that $s_{i+j}(n)=m$. 

\qed

We first prove that unary typed update procedures have zeros. 

\begin{theorem}[Zero Theorem for Unary Typed Update Procedures]\label{theorem-fixupdproc}
Let  $\update :  (\funn)^{k+1}\rightarrow \Nat$ be a term of $\SystemTG$  such that for all closed type-$\funn$ terms $f_1,\ldots, f_k$ of $\SystemTG$, $\update f_1\ldots f_k$ is a typed unary update procedure. Then one can constructively define a closed term $\varepsilon: (\funn)^k\rightarrow (\funn)$ of $\SystemTG$ such that for all closed type-$\funn$ terms $f_1,\ldots, f_k$ of $\SystemTG$ \[\update f_1\ldots f_k(\varepsilon f_1\ldots f_k)= |\varnothing|\] 
\end{theorem}

\proof First, for any term $h:\funn$, define 
\[\mathcal{L}^h:= \lambda \langle \mathcal{M}, g\rangle^{\M_{\Nat}}  \langle \mathcal{M},\lambda n^\Nat h(g(n))\rangle\]  
The same proof of proposition \ref{proposition-constants} (in the case of constants of type $\funn$) shows that for all closed terms $h$ of $\SystemTG$ and $s\in$w.i., $\mathcal{L}^h\gmc h$. Define

\[\mathcal{N}_s :=\lambda h_1^\funn\ldots\lambda h_k^\funn \model{\update}\mathcal{L}^{h_1}\ldots\mathcal{L}^{h_k}\model{\Phic}\]
and fix closed type-$\funn$ terms $f_1,\ldots, f_k$ of $\SystemTG$.
 By the adeguacy theorem  \ref{theorem-adeguacy}, for all $s\in$w.i., $\model{\update}\gmc \update$ and hence
 \[\mathcal{M}_s:=\mathcal{N}_sf_1\ldots f_k\gmc \update f_1\ldots f_k{\Phic}\]
 So, for all $s\in $ w.i., $\pi_0(\mathcal{M}_s)$ is a modulus of convergence for $\lambda m^\Nat \update f_1\ldots f_k(s_m)$. Define by recursion a term $s$ such that $s_0:=0^\funn$ and $s_{n+1}:=s_n\oplus \update f_1\ldots f_ks_n$. Then $s\in$w.i. by proposition \ref{proposition-updateincreasing}, since $\update f_1\ldots f_k$ is a typed unary update procedure.
So $\pi_0(\mathcal{M}_s)$ is a modulus of convergence for $\lambda m^\Nat\ \update f_1\ldots f_k(s_m)$. If we choose $h:=\lambda m^\Nat m+1$ and set $j:=\pi_0(\mathcal{M}_s)h$, we have that 
\[\update f_1\ldots f_k(s_j)=\update f_1\ldots f_k(s_{j+1})\]
 by definition of modulus of convergence. So let 

\[\varepsilon f_1\ldots f_k:=  s_{(\pi_0(\mathcal{M}_s)h)+1}\] 
Then
 \[\begin{aligned}s_{j+2}&=s_{j+1}\oplus \update f_1\ldots f_k(s_{j+1})\\
  &=(s_j\oplus \update f_1\ldots f_k(s_j))\oplus \update f_1\ldots f_k(s_{j+1})\\
 &=(s_j\oplus \update f_1\ldots f_k(s_j))\oplus \update f_1\ldots f_k(s_{j})\\
 &=s_j\oplus \update f_1\ldots f_k(s_j)\\
  &=s_{j+1}
  \end{aligned}\]
 and hence it must be that \[\update f_1\ldots f_k(\varepsilon f_1\ldots f_k)=\update f_1\ldots f_k(s_{j+1})=|\varnothing|\] which is the thesis.
 
  \qed

We are now able to prove the zero theorem for n-ary typed update procedures, following the idea of Avigad's original construction.

\begin{theorem}[Zero Theorem for n-ary Typed Update Procedures]\label{theorem-zeronupdate}
Let $\update :  (\funn)^{k}\rightarrow \Nat$ be a typed update procedure of ordinal $k\in\NatSet$. Then one can constructively define terms $\varepsilon_1,\ldots, \varepsilon_k$ of G\"odel's $\SystemTG$ such that 

\[\update \varepsilon_1\ldots \varepsilon_k=|\varnothing|\]
\end{theorem}
\proof By induction on $k$. The case $k=1$ has been treated in theorem \ref{theorem-fixupdproc}. Therefore, suppose $k>1$. Define 
\[\update_k:=\lambda g_1^\funn\ldots \lambda g_{k-1}^\funn \ifthen{\update g_1\ldots g_{k-1}=|(k,n,m)|}{|(1,n,m)|}{|\varnothing|}\]
 Since for all closed type-$\funn$ terms $f_1,\ldots, f_{k-1}$ of $\SystemTG$, $\update_kf_1\ldots f_{k-1}$ is a unary typed update procedure, by theorem \ref{theorem-fixupdproc}, we can constructively define a term $\varepsilon_k$ of $\SystemTG$ such that for all closed type-$\funn$ terms $f_1,\ldots, f_{k-1}$ of $\SystemTG$ \[\update_k f_1\ldots f_{k-1}(\varepsilon_k f_1\ldots f_{k-1})= |\varnothing|\] 
and hence 
 \[\update f_1\ldots f_{k-1}(\varepsilon_k f_1\ldots f_{k-1})\neq |(k,n,m)|\] 
This implies that 
\[\lambda g_1^\funn \ldots \lambda g_{k-1}^\funn \update g_1\ldots g_{k-1}(\varepsilon g_1\ldots g_{k-1})\] 
is a typed update procedure of ordinal $k-1$. By induction hypothesis, we can constructively define terms $\varepsilon_1,\ldots, \varepsilon_{k-1}$ of G\"odel's $\SystemTG$ such that 
 \[\update \varepsilon_1\ldots \varepsilon_{k-1}(\varepsilon_k \varepsilon_1\ldots \varepsilon_{k-1})= |\varnothing|\] 
 which is the thesis.

\qed

An important corollary of theorem \ref{theorem-zeronupdate}, is the termination of the epsilon substitution method for first order Peano Arithmetic.

\begin{theorem}[Termination of Epsilon Substitution Method for $\PA$]
The $H$-process (as defined in Mints \cite{Mints}) of the epsilon substitution method for $\PA$ always terminates.
\end{theorem}
\proof See Avigad \cite{Avigad}
\qed

\chapter{Learning in Predicative Analysis}\label{chapter-learninganalysis}

\begin{abstract} We give an axiomatization of the concept of learning, as it implicitly appears in various computational interpretations of Predicative classical second order Arithmetic. We achieve our result by extending Avigad's notion of update procedure to the transfinite case.

\end{abstract}

\section{Introduction}
The aim of this chapter is to provide an abstract description of learning as it is appears in various computational interpretations of predicative fragments of classical second order Arithmetic. Our account has a twofold motivation and interest.

 Its first purpose is to provide a foundation that will serve to extend learning based realizability to predicative fragments of Analysis: a possible path to follow is the one suggested here. In particular, we describe the learning processes that arise when extending the approach of learning based realizability to predicative Arithmetic and prove their termination. This is achieved by introducing the notion of transfinite update procedure.
 
 Secondly, we continue the work of Avigad on update procedures \cite{Avigad} and - as anticipated - extend them to the transfinite case. The concept of transfinite update procedure may be seen as an axiomatization of learning as  implicitly used in the epsilon substitution method formulated in the work of Mints et al (\cite{Mints}, \cite{Mints2}). The notion is useful to understand the core of the epsilon method and its fundamental ideas without having to deal with the complicated formalism and non relevant details. Moreover, as interesting byproduct of the conceptual analysis of the epsilon method through transfinite update procedures, one can provide a combinatorial statement that is equivalent to the 1-consistency of various fragments of predicative Arithmetic in a very scalable way. In particular, the informal statement ``$\mathsf{U}(\alpha)$: every transfinite update procedure of ordinal less than $\alpha$ has a finite zero" very rapidly acquires logical complexity even at small ordinals. For example:
 \begin{enumerate}
 \item $\mathsf{U}(2)$ corresponds to the 1-consistency of $\HA+\EM_1$.
 \item $\mathsf{U}(n+1)$, with $n\in \omega$,  corresponds to the 1-consistency of $\HA+\EM_n$ (excluded middle over $\Sigma_n^0$ formulas).
 \item $\mathsf{U}(\omega)$ corresponds to the 1-consistency of $\PA$
 \item  $\mathsf{U}(\omega\cdot 2)$ corresponds to the 1-consistency of $\EA$ (Elementary Analysis)
 \item $\mathsf{U}(\omega^\omega)$ corresponds to the 1-consistency of $\PA^2$ plus $\Delta_1^1$-comprehension axiom.
   \end{enumerate}
  (1), (4) are treated in this thesis, (3) has been proved by Avigad \cite{Avigad}, (2) will follow by extending learning based realizability. (5) should follow by straightforward extension of the methods we will use to prove (4). In order to make precise the statement $\mathsf{U}(\alpha)$ and prove refined versions of (1)-(5), one has to choose a formal system in to which represent update procedures. All the update procedures used in this thesis may be assumed to be represented in system $\SystemTG$.  \\
 
 \emph{Plan of the Chapter}. In section \S \ref{section-transfiniteupdateprocedures} we introduce and motivate the concept of transfinite update procedure and give a very short non constructive proof of the existence of finite zeros. 
 
 In section \S \ref{section-learningprocesses} we explain the notion of ``learning process generated by an update procedure" and prove that every learning process terminates with a zero for the associated update procedure. The result represents a more constructive proof of the existence of finite zeros and the learning processes are ``optimal", in the sense that one could provide constructively the expected ordinal bounds to their length  and of the size of finite zeros (by applying techniques of Mints \cite{Mints})
 
 In sections \S  \ref{section-typedupdateprocedures} and \S\ref{section-barrecursivezeros} we formalize the notion of update procedure in typed lambda calculus plus bar recursion and prove the existence of zeros for update procedures of ordinal less than $\omega^\omega$ by writing down simple bar recursive terms.
 
 In section \S \ref{section-casestudyea} we devote ourselves to a case study: we show that $\mathsf{U}(\omega\cdot 2)$ implies the 1-consistency of $\EA$ by proving that it implies the termination of $H$-processes of the epsilon substitution method for $\EA$.

\section{Transfinite Update Procedures for Predicative Analysis}\label{section-transfiniteupdateprocedures}
From the computational point of view, classical predicative second order Arithmetic poses very difficult problems. Axioms of comprehension ask for functions $g$ able to decide truth of formulas:
\[\exists g^\funn \forall x^\Nat.\ g(x)=0\leftrightarrow \phi(x)\] 
 Axioms of (countable) choice ask for functions $g$ computing existential witnesses of truth of formulas:
 \[(\forall x^\Nat \exists y^A\ \phi(x,y))\rightarrow \exists g^{\Nat\rightarrow A} \forall x^{\Nat} \phi(x, g(x))\] 
 A Kleene-style realizability interpretation for even the most simple form of the excluded middle 
 \[\EM_1:\forall n^\Nat. \forall y^\Nat \lnot Pny \lor \exists x^\Nat Pnx\] asks for deciding the truth of semidecidable formulas of the form $\exists x^\Nat\ Pnx$, with $P$ decidable. 

In general, learning based computational interpretations of predicative fragments of classical analysis (such as our learning based realizability, epsilon substitution method, update procedures, Herbrand analysis) provide answers to the above computational challenges by the following three-stage pattern:
\begin{enumerate}
\item They identify a sequence $F$ - possibly transfinite - of non computable functions $\fun$.\\

\item They define classical witnesses for provable formulas by using programs recursive in $F$.\\

\item They define learning procedures through which it is possible to find, for every particular computation, a suitable finite approximation of the functions of $F$ such that one can effectively compute the witnesses defined at stage two.
\end{enumerate} 
The functions in the $F$ of stage (1) are the computational engine of the interpretation. Given the difficulty of computing witnesses in classical Arithmetic, they are always non computable. It is  therefore not surprising that given this additional computational power, one  is able to define  at stage (2) witnesses for classical formulas. If we picture the sequence $F$ as a sequence of infinite stacks of numbers, the learning process of point (3) finds a ``vertical" approximation of $F$:  functions of $F$ are infinite stacks of numbers whereas their finite approximations are finite stacks. Moreover, a crucial point is that the sequence $F$ is not an arbitrary sequence. In a sense, $F$ is also ``horizontally" approximated:  for every ordinal $\alpha$, the recursion theoretic Turing degree of $F_\alpha$ is approximated by the degrees of $F_\beta$, for $\beta<\alpha$. This property is very important: in this way, the values of the functions in $F$ can be gradually approximated and learned.

 More precisely, $F$ can be seen a sequence of functions obtained by transfinite iteration of recursion theoretic \emph{jump} operator (see for example Odifreddi \cite{Odifreddi}). That is, for every $\beta$, if $\beta$ is a successor, $F_\beta$ has the same Turing degree of an oracle for the halting problem for the class of functions recursive in $F_{\beta-1}$ (jump); if $\beta$ is limit, $F_\beta$ has the same Turing degre of the function mapping the code of a pair $(\alpha, n)$, with $\alpha<\beta$,  into $F_\alpha(n)$ (join or $\beta$-jump).  
A fundamental property of such a sequence is that the assertion $F_\beta(n)=m$ depends only on the values of the functions $F_\alpha$, for $\alpha<\beta$, and the values of $F_\beta$ are learnable in the limit\footnote{In the sense of Gold \cite{Gold}: $F_\beta(n)=m \iff \lim_{k\to\infty}g(n,k)=m$} by a program $g$ recursive in the join of $F_\alpha$ for $\alpha<\beta$, which is a guarantee that learning processes will terminate. 

We now give an informal example of the kind of analysis which is needed to carry out the first stage of a learning based interpretation, in the case of $\EA$. A complete treatment is postponed to section \ref{section-casestudyea}.
\begin{example}[Elementary Analysis]
Consider a subsystem of second order Peano Arithmetic in which second order quantification is intended to range over arithmetical sets and hence over arithmetical formulas (formulas with only numerical quantifiers and possibly set parameters). Since one has to interpret excluded middle over arbitrary formulas, it is necessary to provide \emph{at least} programs that can decide truth of formulas. Numerical quantifiers correspond to Turing jumps. That is, if we have a program $t$ (with the same function parameters of $\phi$) such that for every pair of naturals $n, m$
\[t(n,m)=\True \iff \phi(n,m)\]
then the truth of 
\[\exists x^\Nat \phi(n,x)\]
is equivalent to the termination of a program  $Q(n)$ exhaustively checking
\[t(n,0), t(n,1), t(n,2), \ldots\]
until it finds - if there exists - an $m$ such that $t(n,m)=\True$. Applying the jump operator to the Turing degree $t$ belongs to, one can write down a program $\chi_t$ which is able to determine whether $Q(n)$ terminates. That is
\[\chi_t(n)=\True \iff \exists x^\Nat \phi(n,x)\]
Similarly, one eliminates universal numerical quantifiers, thanks to the fact that $\forall \equiv \lnot \exists \lnot $. Iterating these reasoning and applying $2k$ times the jump operator - and given a recursive enumeration $\phi_0, \phi_1, \ldots,$ of arithmetical formulas - one can obtain for every $\Sigma_{2k}^0$ formula 
\[\phi_n(m):=\exists x_1^\Nat \forall y_1^\Nat \ldots \exists x_k^\Nat \forall y_k^\Nat\ P(m, x_1,y_1, \ldots, x_k, y_k)\]
a program $t_n$ such that
\[t_n(m)=\True \iff \phi_n(m)\]
Using the $\omega$-jump operator, one can write down a program $u$ such that
\[u(n,m)=\True \iff t_n(m)=\True\]
and hence
\[u(n,m)=\True \iff \phi_n(m)\]
Now a $\Sigma_1^1$ formula
\[\exists f^{\Nat\rightarrow \Bool} \phi_i \]
 - provided we assume that $f^{\Nat\rightarrow \Bool}$ ranges over arithmetical predicates - can also be expressed as
\[\exists n^\Nat t_i[\lambda m^\Nat u(n, m)/f]\]
So applying again a jump operator to the recursive degree of $t_i[\lambda m^\Nat u(n, m)/f]$, one is able to write a program determining the truth value of $\exists f^{\Nat\rightarrow \Bool} \phi_i$. Iterating this reasoning, one can decide the truth of arbitrary $\Sigma_n^1$ formulas. \\
Summing up, in order to decide  truth in Elementary Analysis, one needs to apply the jump operator $\omega+\omega$ times and thus produces a sequence $F$ of non computable functions $F$ of length $\omega+\omega$. All the programs that we have described are recursive in some initial segment of $F$.

\end{example}

We are now in a position to understand the following axiomatization of the learning procedures cited in point (3) above.  

\begin{defi}[Transfinite Update Procedures]\label{transfiniteupdateprocedures}
Let $\alpha\geq 1$ be a numerable ordinal. An {\em update procedure of ordinal $\alpha$} is a function $\update:( \alpha\rightarrow (\NatSet\rightarrow \NatSet))\rightarrow (\alpha \times \NatSet \times \NatSet)\cup \{\emptyset\}$ such that:

\begin{enumerate}
\item $\update$ is continuous. i.e. for any $f: \alpha\rightarrow (\NatSet\rightarrow \NatSet)$ there is a finite subset $A$ of $\alpha \times \NatSet$ such that for every $g: \alpha\rightarrow (\NatSet\rightarrow \NatSet)$ if $f_\gamma (n) =g_\gamma (n)$ for every $(\gamma, n)\in A$, then $\update f=\update g$.\\

\item For all functions $f,g:  \alpha\rightarrow (\NatSet\rightarrow \NatSet)$ and every ordinal $\beta\in \alpha$, if \\

i) for every $\gamma < \beta$, $f_{\gamma}=g_{\gamma}$;\\

ii) $\update f = (\beta, n, m)$, $g_\beta (n)=m$ and $ \update g= (\beta, i, j)$;\\\\
 then $i\neq n$.
\end{enumerate}
\end{defi}

The concept of transfinite update procedure is a generalization of Avigad's notion of update procedure \cite{Avigad}. A transfinite update procedure, instead of taking just a finite number of function arguments, may get as input an arbitrary transfinite sequence of functions, which are intended to approximate a target sequence $F$; as output, it may return an update $(\beta, n, m)$, which means that the $\beta$-th function taken as argument is an inadequate approximation of $F_\beta$ and must be updated as to output $m$ on input $n$. Condition (2) in definition \ref{transfiniteupdateprocedures} is a little bit stronger than Avigad's requirement, which would be:\\\\
(2)' For all functions $f,g:  \alpha\rightarrow (\NatSet\rightarrow \NatSet)$ and every ordinal $\beta\in \alpha$, if \\

i) for every $\gamma < \beta$, $f_{\gamma}=g_{\gamma}$;\\

ii) $\update f = (\beta, n, m)$, $g_\beta n=m$, $f_\beta \leq g_{\beta}$ and $ \update g= (\beta, i, j)$ (where $f_\beta \leq g_{\beta}$ is defined as $f_\beta(x)\neq 0$ implies $f_\beta(x)=g_\beta(x)$);\\\\
 then $i\neq n$.\\\\
But in fact, the update procedures which are \emph{actually} used by Avigad \cite{Avigad} in proving 1-consistency of $\PA$ still fall under our definition. \\Condition (2) means that the values for the $\beta$-th function depend only on the values of functions of ordinal less than $\beta$ in the input sequence and an update procedure returns only updates which are \emph{relatively verified} and hence need not to be changed. In this sense, if $\update f = (\beta, n, m)$, one has \emph{learned} that $F_\beta(n)=m$; so if $g_\beta$ is a candidate approximation of $F_\beta$ and $g_\beta(n)=m$, then $\update g$ does not represent a request to modify the value of $g_\beta$ at point $n$, whenever $f$ and $g$ agree on all ordinals less than $\beta$. \\
 We remark that the choice of the type for an update procedure is somewhat arbitrary: we could have chosen it to be
 \[( \alpha\rightarrow (X\rightarrow Y))\rightarrow (\alpha \times X \times Y)\cup \{\emptyset\}\] 
 \emph{as long as} the elements of the sets $X$ and $Y$ can be coded by finite objects. Since such coding may always be performed by using natural numbers, we choose to consider $X=Y=\NatSet$.  
  
The use of transfinite update procedures made by learning based computational interpretations of classical arithmetic can be described as follows. Suppose those interpretations are given a provable formula with an attainable computational meaning, for example one of the form $\forall x^\Nat \exists y^\Nat\ Pxy$, with $P$ decidable. Then, for every numeral $n$, they manage to  define a term $t_n: (\alpha\rightarrow (\NatSet\rightarrow \NatSet))\rightarrow \NatSet$ and an update procedure $\update_n$ of ordinal $\alpha$ such that 
\[\update_n (f) =\emptyset \implies Pn(t_n(f))\]
for all $f: \alpha\rightarrow (\NatSet\rightarrow \NatSet) $. The idea is that a witness for the formula $\exists y^\Nat Pny$ is calculated by $t_n$ with respect to a particular approximation $f$ of the sequence $F$ we have previously described. If the formula $Pn(t_n(f))$ is true, there is nothing to be done. If it is false, then $\update_n (f)= (\beta, n, m)$ for some $\beta, n, m$: a new value for $F_\beta$ is learned.  This is what we call ``learning by counterexamples": from every failure a new positive fact is acquired. We have studied an instance of this kind of learning in the chapter on learning based realizability for $\HA+\EM_1$ (for the case of $\alpha=1$), when we defined realizability for atomic formulas: in that case the pair $(n,m)$ was produced by the realizer of the excluded middle. We will see another instance in section \ref{section-casestudyea} in the case $\alpha=\omega+k$, with $k\in \NatSet$: the triple $(\beta, n ,m)$ will be produced through the evaluation of axioms for epsilon terms. 

The effectiveness of the above approach depends on the fact that every update procedure has a finite zero, as defined below.

\begin{defi}[Finite Functions, Finite Zeros, Truncation and Concatenation of Function Sequences]\label{definition-truncconc}
Let $\update$ be an update procedure of ordinal $\alpha$. 
\begin{enumerate}
\item
$f: \alpha\rightarrow (\NatSet\rightarrow \NatSet)$ is said to be a {\em finite function} if the set of $(\gamma, n)$ such that $f_\gamma n\neq 0$ is finite.\\

\item A {\em finite zero} for $\update$ is a finite function $f: \alpha\rightarrow (\NatSet\rightarrow \NatSet)$ such that $\update f=\emptyset$.\\
\item Let  $f: \alpha \rightarrow (\fun)$ and $\beta<\alpha$. Let  $f_{<\beta}: \beta\rightarrow(\fun)$  be the \emph{truncation} of $f$ at $\beta$:
\[f_{<\beta}:=\gamma\in\beta \mapsto f_\gamma  \]

\item Let $\alpha_1, \alpha_2$ be two ordinals, $f_1: \alpha_1\rightarrow (\fun)$ and $f_2:\alpha_2\rightarrow (\fun)$. Then the \emph{concatenation} $f_1*f_2: (\alpha_1+\alpha_2)\rightarrow (\fun)$ of $f$ and $g$ is defined as:
\[(f*g)_\gamma(n):=
\begin{cases}
f_\gamma(n) &\text{$\mathsf{if}$ } \gamma <\alpha_1 \\
g_\beta(n)  &\text{$\mathsf{if}$ $\gamma=\alpha_1+\beta<\alpha_1+\alpha_2$ } 
\end{cases}
\]
\item With a slight abuse of notation, a function $f:\fun$ will be sometimes identified with the corresponding length-one sequence of functions $0\mapsto (n\in \NatSet \mapsto f(n))$.

\end{enumerate}
\end{defi}

We now prove that every update procedure has a finite zero. We will give other more and more constructive proofs of this theorem, that will allow to compute finite zeros for update procedures and thus witnesses for classically provable formulas, thanks to learning based interpretations. But for now we are only interested into understanding the reason of the theorem's \emph{truth} and give a very short non constructive proof. All the subsequent proofs can be seen as more and more sophisticated and refined constructivizations of  the following argument.
\begin{thm}[Zero Theorem for Update Procedures of Ordinal $\alpha$]
Let $\update$ be an update procedure of ordinal $\alpha$. Then $\update$ has a finite zero. 
\end{thm}
\proof

 We define, by transfinite induction,  a function $f: \alpha\rightarrow (\fun)$ as follows. Suppose we have defined $f_\gamma: \fun$, for every $\gamma<\beta$. Define the sequence $f_{<\beta} : \beta \rightarrow (\fun)$ of them all 
 \[f_{<\beta}:=\gamma \in \beta \mapsto f_\gamma\]
 Then define
 \[f_\beta(x)=
\begin{cases}
0 &\text{$\mathsf{if }\ \forall g^{(\alpha-\beta)\rightarrow (\fun)}\ \forall z^\NatSet\ \update({f_{<\beta}}*g)\neq \langle \beta, x, z\rangle$ }\\
y &\text{$\mathsf{otherwise }$}, \text{ for some $y$ such that  $\exists g^{(\alpha-\beta)\rightarrow (\fun)}\ \update(f_{<\beta}*g)= \langle \beta, x,y \rangle$}

\end{cases}
\]
 By axiom of choice and classical logic, for every $\beta$, $f_{<\beta}$ and $f_\beta$ are well defined. So we can let
  \[f:=f_{<\alpha}\]
 Suppose $\update (f)=\langle \beta, x, z\rangle$, for some $\beta<\alpha$: we show that it is impossible. For some $h:{(\alpha-(\beta+1))\rightarrow (\fun)}$, $f=f_{<\beta}*f_\beta*h$. Hence, for some $g:{(\alpha-\beta)\rightarrow (\fun)}$ \[\update(f_{<\beta}*g)=\langle \beta, x, y\rangle\land f_\beta(x)=y\]
 by definition of $f$. But $\update$ is an update procedure and so
 \[(\update(f_{<\beta}*g)=\langle \beta, x, y\rangle\land f_\beta(x)=y\land \update (f_{<\beta}*f_\beta*h)=\langle \beta, x, z\rangle)\implies x\neq x\]
 which is impossible. We conclude that $\update(f)=\emptyset$ and, by continuity, that $\update$ has a finite zero.

\qed

\section{Learning Processes Generated by Transfinite Update Procedures}\label{section-learningprocesses}
In this section we show that every update procedure $\update$ generates a learning process and this learning process always terminates with a finite zero of $\update$. This result is an abstract version of the termination of the $H$-process as defined in the various versions of epsilon substitution method (see Mints et al. \cite{Mints}). The proof of termination is non-constructive and is similar to the one in  Mints et al. \cite{Mints} (which however is by contradiction while ours is not).\\

If $\update$ is an update procedure and $\update(f)=\langle \gamma, n, m\rangle$, then the value of $f_\gamma$ at argument $n$ must be updated as to be equal to $m$. But as explained in the introduction, we may imagine that all the values of all the functions $f_\beta$, with $\beta>\gamma$, depend on the values of the \emph{current} $f_\gamma$. Therefore, if we change some of the values of $f_\gamma$, we must erase all the values of all the functions $f_\beta$, for $\beta>\gamma$, because they may be inconsistent with the new values of $f_\beta$. In a sense, $f$ is a fragile structure, that may be likened to an \emph{house of cards}: if we change some layer all the higher ones collapse. We define an update operator $\oplus$ that performs those operations.

\begin{defi}[Controlled Update of Functions]

Let $f:\alpha\rightarrow (\fun)$ and $\langle \gamma, n, m\rangle\in \alpha\times\NatSet\times \NatSet$. We define a function $f\oplus\langle \gamma, n, m\rangle: \alpha\rightarrow (\fun)$ such that

\[(f\oplus\langle \gamma, n, m\rangle)_\beta (x):=
\begin{cases}
 f_\beta (x) &\text{$\mathsf{if }$ $\beta< \gamma$ or ($\beta=\gamma$  and $x\neq n$)}\\
 m  &\mathsf{if }\text{ $\gamma=\beta$ and $x=n$}\\ 
 0 & \text{$\mathsf{otherwise}$}\\

\end{cases}
\]
We also define $f\oplus \emptyset:=f$.

\end{defi}

We now define the concept of ``learning process generated by an update procedure $\update$". It may be thought as a process of updating and learning new values of functions, which is \emph{guided} by $\update$. It corresponds to the step three of the learning based computational interpretations of classical arithmetic we have described in the  introduction.  Intuitively, such a learning process starts from the always zero function $\zerof^{\alpha}$. If $\update$ says that some value of $\zerof^{\alpha}$ must be updated - i.e. $\update(\zerof^{\alpha})=\langle \gamma, n ,m\rangle$ -  then the learning process generates the function $\update^{(1)}:=\zerof^{\alpha}\oplus \langle\gamma, n ,m\rangle$. Similarly, if $\update$ says that some value of $\update^{(1)}$ must be updated - i.e. $\update(\update^{(1)})=\langle \gamma', n' ,m'\rangle$ - then the learning process generates the function $\update^{(2)}:=\update^{(1)}\oplus \langle\gamma', n' ,m'\rangle$. The process goes on indefinitely in the same fashion.
\begin{defi}[Learning Processes Generated by $\update$]\label{definition-learningprocess}
Let $\update$ be an update procedure of ordinal $\alpha$. For every $n\in\NatSet$, we define a function $\update^{(n)}: \alpha\rightarrow (\fun)$ by induction as follows:

\[\begin{aligned}\update^{(0)}&:= \zerof^{\alpha}:=\gamma \in\alpha \mapsto (n\in \NatSet\mapsto 0)\\
\update^{(n+1)}&:=\update^{(n)}\oplus \update(\update^{(n)} )
\end{aligned}\] 
Moreover, a function $f: \alpha\rightarrow (\fun)$ is said to be $\update$-{\em generated} if there exists an $n$ such that $f=\update^{(n)}$.

\end{defi}
The aim of the learning process generated by $\update$ is to find a finite zero for $\update$. Indeed, if for some $n$, $\update(\update^{(n)})=\emptyset$, then for all $m\geq n$, $\update^{(m)}=\update^{(n)}$ and we thus say that the learning process \emph{terminates}. We now devote ourselves to the proof that learning processes always terminate. In other words, every update procedure $\update$ has a $\update$-generated finite zero.

Given an update procedure $\update$, its useful to define a new ``simpler" update procedure, obtained from $\update$ by fixing some initial segment of its input, ignoring all updates relative to this fixed part of the input and adjusting their indexes. 

\begin{defi}\label{defiinstance}
Let $\update$ be an update procedure of ordinal $\alpha$. Then, for any function $g:\beta\rightarrow (\fun)$, with $\beta<\alpha$,  define a function    
\[\update_{g}: ((\alpha-\beta)\rightarrow (\fun))\rightarrow (\alpha-\beta)\times \NatSet\times\NatSet \cup \{ \emptyset \}
\]
as follow
\[\update_g(f)=
\begin{cases} 
\langle \gamma, n,m\rangle &\text{$\mathsf{if }$ } \update(g*f)=\langle \beta + \gamma, n,m\rangle \\
\emptyset &\mathsf{otherwise}

\end{cases}
\]
(We point out that if $\beta=0=\emptyset$, $\update_g=\update$ as it should be) 
\end{defi}

Indeed $\update_g$ as defined above is an update procedure.
\begin{fact}\label{fact-ovvio}
 Let $\update$ be an update procedure of ordinal $\alpha$. Then, for any function $g:\beta\rightarrow (\fun)$, with $\beta<\alpha$:
\begin{enumerate} \item $\update_g$ is an update procedure of ordinal $\alpha-\beta$.\\

\item For every $h:\fun$, $\update_{g*h}=(\update_{g})_h$.
\end{enumerate}
\end{fact}
\proof Immediate.

\qed

The strategy of our termination proof can be described as follows. Given an update procedure $\update$ of ordinal $\alpha$, we shall define a sequence of functions $g: \alpha\rightarrow (\fun)$ such that a ``reduction lemma'' can be proved: if, for some $\beta<\alpha$, $\update_{g_{<\beta}}$ has a $\update_{g_{<\beta}}$-generated finite zero, then for some $\gamma<\beta$ also $\update_{g_{<\gamma}}$ has a $\update_{g_{<\gamma}}$-generated finite zero (for definition of $g_{<\beta}$, recall definition \ref{definition-truncconc} ). But the greater the ordinal $\beta$ the easier is to compute with a learning process a finite zero for $\update_{g_{<\beta}}$ because the sequence $g_{<\beta}$ becomes so long that the input for $\update_{<\beta}$ becomes short. So we shall be able to show that for some large enough $\beta$, $\beta<\alpha$, $\update_{g_{<\beta}}$ has a $\update_{g_{<\beta}}$-generated finite zero, which proves the theorem in combination with the reduction lemma. This technique can be seen as a generalization of Avigad's \cite{Avigad} to the transfinite case.

We now prove the reduction lemma in the limit case.

\begin{lem}[Reduction Lemma, Limit Case]\label{lemma-transfiniteinstance}
Let $\update$ be an update procedure of ordinal $\alpha$ and $g: \beta \rightarrow (\fun)$, with $\beta$ limit ordinal and $\beta<\alpha$. Then
\begin{enumerate} \item 
If $f:(\alpha-\beta)\rightarrow (\fun)$ is $\update_g$-generated, then there exists $\gamma<\beta$ such that $\mathsf{0}^{\beta-\gamma}*f$ is $\update_{g_{<\gamma}}$-generated.\\
\item If $\update_g$ has a $\update_g$-generated finite zero, then there exists $\gamma<\beta$ such that $\update_{g_{<\gamma}}$ has a $\update_{g_{<\gamma}}$-generated finite zero.
\end{enumerate}
\end{lem}
\proof (1) Let $n$ be the smallest among the $i$ such that $f=\update_g^{(i)}$. If $k<n$ and $\update_g(\update_g^{(k)})=\emptyset$, then 
\[\update_g^{(k+1)}=\update_g^{(k)}\oplus \update_g(\update_g^{(k)})=\update_g^{(k)}\]
So we have that for all $k<n$, $\update_g(\update_g^{(k)})\neq \emptyset$. Since $\beta$ is a limit ordinal and $\update$ is continuous, there exists a $\gamma<\beta$ such that for every $k\leq n$
\[\update(g*\update_g^{(k)})=\update(g_{<\gamma}*\mathsf{0}^{\beta-\gamma}*\update_g^{(k)})\]
and
\begin{equation}\label{delta}\update(g_{<\gamma}*\mathsf{0}^{\beta-\gamma}*f)=\langle \delta, n , m\rangle \land \delta<\beta \implies \delta<\gamma
\end{equation}
We prove by induction on $k\leq n$ that 
\[\update_{g_{<\gamma}}^{(k)}=\mathsf{0}^{\beta-\gamma}*\update_g^{(k)}\]
which is the thesis.
If $k=0$, \[\update_{g_{<\gamma}}^{(0)}=\mathsf{0}^{\alpha-\gamma}=\mathsf{0}^{\beta-\gamma}*\mathsf{0}^{\alpha- \beta}=\mathsf{0}^{\beta-\gamma}*\update_g^{(0)}\]
If $k+1\leq n$, then for some $\delta, l, m$ 
\[\update_g(\update_g^{(k)})=\langle \delta, l, m\rangle\]
Then, by definition of $\update_g$, 
\[\update(g*\update_g^{(k)})=\langle \beta+\delta, l, m\rangle\]
and hence 
\[\update(g_{<\gamma}*\mathsf{0}^{\beta-\gamma}*\update_g^{(k)})=\langle \beta+\delta, l, m\rangle=\langle \gamma+ (\beta-\gamma)+\delta, l, m\rangle\]
Therefore, by definition of $\update_{g_{<\gamma}}$
\[ \update_{g_{<\gamma}}(\mathsf{0}^{\beta-\gamma}*\update_g^{(k)})=\langle (\beta-\gamma)+\delta, l, m\rangle\]
By the help of induction hypothesis, we conclude that
\[\begin{aligned}\update_{g_{<\gamma}}^{(k+1)}&=\update_{g_{<\gamma}}^{(k)}\oplus \update_{g_{<\gamma}}(\update_{g_{<\gamma}}^{(k)})\\
&=\update_{g_{<\gamma}}^{(k)}\oplus \update_{g_{<\gamma}}(\mathsf{0}^{\beta-\gamma}*\update_g^{(k)})\\
&=\update_{g_{<\gamma}}^{(k)}\oplus\langle (\beta-\gamma)+\delta, l, m\rangle\\
&=(\mathsf{0}^{\beta-\gamma}*\update_g^{(k)})\oplus\langle (\beta-\gamma)+\delta, l, m\rangle\\
&=\mathsf{0}^{\beta-\gamma}*(\update_g^{(k)}\oplus\langle \delta, l, m\rangle)\\
&=\mathsf{0}^{\beta-\gamma}*(\update_g^{(k)}\oplus \update_g(\update_g^{(k)}))\\
&=\mathsf{0}^{\beta-\gamma}*\update_g^{(k+1)}
\end{aligned}\]
(2) We continue the previous proof. Suppose that $f$  is also a finite zero of $\update_g$.  If
\[\emptyset=\update(g*f)=\update(g_{<\gamma}*\mathsf{0}^{\beta-\gamma}*f)\]
then by definition of $\update_{g_{<\gamma}}$
\[\update_{g_{<\gamma}}(\mathsf{0}^{\beta-\gamma}*f)=\emptyset\]
Therefore, suppose
\[\langle \delta, l ,m\rangle=\update(g*f)=\update(g_{<\gamma}*\mathsf{0}^{\beta-\gamma}*f)\]
Since $\update_g(f)=\emptyset$, by definition of $\update_g$ we have that $\delta <\beta$.
By \ref{delta}, we have also that $\delta<\gamma$. So, by definition of $\update_{g_{<\gamma}}$
\[\update_{g_{<\gamma}}(\mathsf{0}^{\beta-\gamma}*f)=\emptyset\]
which is the thesis.

\qed

We now prove the reduction lemma in the successor case.

\begin{lem}[Reduction Lemma, Successor Case]\label{lemma-updateinstance}
Let $\update$ be an update procedure of ordinal $\alpha$. Define $g: \fun$ as follows:

\[g(x):=
\begin{cases}
y &\text{$\mathsf{if }$ } \exists i.\ \update(\update^{(i)})= \langle 0, x, y \rangle\land i=\mathsf{min}\{n\ |\ \exists z\ \update(\update^{(n)})=\langle 0, x, z\rangle\}\\
0 &\text{$\mathsf{otherwise }$ }
\end{cases}
\]
Then:
\begin{enumerate}
\item For every finite function $g_0 \leq\footnote{As in chapter 4, $g_0\leq g$ iff for all $x$ $g_0(x)\neq 0\implies g_0(x)=g(x)$} g$, if $g_0*\zerof^{\alpha-1}$ is $\update$-generated and  $f$ is  $\update_{g_0}$-generated, then $g_0*f$ is $\update$-generated.\\

\item If $\update_{g}$ has a $\update_{g}$-generated finite zero, then $\update$ has a $\update$-generated finite zero.
\end{enumerate}
\end{lem}
\proof

(1) By induction on the number $m$ such that $\update_{g_0}^{(m)}=f$. If $f=\update_{g_0}^{(0)}=\zerof^{\alpha-1}$, then $g_0*f=g_0*\zerof^{\alpha-1}$ is $\update$-generated by hypothesis.\\
If \[f=\update_{g_0}^{(k+1)}=\update_{g_0}^{(k)}\oplus \update_{g_0}(\update_{g_0}^{(k)})\] then $g_0* \update_{g_0}^{(k)}$ is $\update$-generated by inductive hypothesis, i.e. for some $n$, $g_0* \update_{g_0}^{(k)}=\update^{(n)}$ . We have two cases:\\

i) $\update_{g_0}(\update_{g_0}^{(k)})=\langle \gamma, x,z \rangle$. Then, by definition of $\update_{g_0}$ \[\update (g_0*\update_{g_0}^{(k)})=\langle \gamma+1, x,z \rangle\] 

Therefore, 
\[\begin{aligned}\update^{(n+1)}&=\update^{(n)}\oplus\update(\update^{(n)})\\
&=(g_0* \update_{g_0}^{(k)})\oplus \update (g_0*\update_{g_0}^{(k)})\\
&=(g_0* \update_{g_0}^{(k)})\oplus\langle \gamma+1, x,z \rangle\\
&=g_0*( \update_{g_0}^{(k)}\oplus\langle \gamma, x,z \rangle)\\
&= g_0*f
\end{aligned}\]

ii)  $\update_{g_0}(\update_{g_0}^{(k+1)})=\emptyset$. Then 
\[f= \update_{g_0}^{(k)}\oplus \emptyset= \update_{g_0}^{(k)}\] 

Hence $g_0*f$ is $\update$-generated by induction hypothesis. \\

 (2) Let $f$ be a $\update_{g}$-generated finite zero  of $\update_{g}$. By definition of $g$, for every $m$ such that, for some numbers $x,y$,
  \[\update(\update^{(m)}) = \langle 0,x,y\rangle\] we have that $\update^{(m+1)}=h_m*\zerof^{(\alpha-1)}$ for some finite function $h_m\leq g$. $h_m*\zerof^{(\alpha-1)}$ is $\update$-generated by definition. The sequence of all $h_m$ is increasing, and, by definition of $g$, the limit of all $h_m$ is indeed $g$. Thus, for every finite function $h\leq g$ we have $h\leq h_m\leq g$ for some finite $h_m$ such that $h_m*\zerof^{(\alpha-1)}$ is $\update$-generated.

By assumption $f$ is $\update_g$-generated, that is, $f = \update_g^{(n)}$ for some $n$. By continuity of $\update(g*f)$ in $g$, and by continuity of $\update_g^{(m)}$ in $g$ for any $m\leq n$, we deduce that there is some finite function $h \leq g$ such that for all functions $h\leq g_0\leq g$ the two conditions below hold: \[\update(g_0 * f) = \update(g *f)\]  \[f = \update_{g_0}^{(n)}\] that is, $f$ is an $\update_{g_0}$-generated finite zero of $\update_{g_0}$. By the discussion above, we may choose some finite $g_0$ such that $h\leq g_0\leq g$ and $g_0*\zerof^{(\alpha-1)}$ is $\update$-generated.  By point (1), $g_0*f$ is $\update$-generated: 
 \begin{equation}\label{ee0}g_0*f=\update^{(n)}
 \end{equation}
  for some $n$. Suppose
  \begin{equation}\label{ee1}\update(g_0*f)=\langle 0, x, z\rangle
  \end{equation}
  for some $x, z$: we show it is impossible and hence obtain that $\update(g_0*f)=\emptyset$, by the fact that $\update_{g_0}(f)=\emptyset$ and definition of $\update_{g_0}$. Combining (\ref{ee0}) and (\ref{ee1}), we obtain
 \[\update(\update^{(n)})=\langle 0, x, z\rangle\] 
 By definition of $g$, for some $m\leq n$
 \[\update(\update^{(m)})=\langle 0, x, y\rangle\land g(x)=y\]
This last fact plus (\ref{ee1}) imply that $x\neq x$, since by definition $\update$ is an update procedure: impossible.
 
 \qed

We are now able to prove the main theorem: update procedures generate terminating learning processes.

\begin{thm}[Termination of Learning Processes]
Let $\update$ be an update procedure of ordinal $\alpha$. Then, $\update$ has a finite zero. In particular, there exists $k\in\NatSet$ such that $\update(\update^{(k)} )=\emptyset$.

\end{thm}
\proof
We define, by transfinite induction,  a function $g: \alpha\rightarrow (\fun)$ as follows. Suppose we have defined $g_\gamma: \fun$, for every $\gamma<\beta<\alpha$. Define the sequence $g_{<\beta} : \beta \rightarrow (\fun)$ of them all 
 \[g_{<\beta}:=\gamma \in \beta \mapsto g_\gamma\]
 Then define
 \[g_\beta(x):=
\begin{cases}
y &\text{$\mathsf{if }$ } \exists i.\ \update_{g_{<\beta}}(\update_{g_{<\beta}}^{(i)})= \langle 0, x, y \rangle\land i=\mathsf{min}\{n\ |\ \exists z\ \update_{g_{<\beta}}(\update_{g_{<\beta}}^{(n)})=\langle 0, x, z\rangle\}\\
0 &\text{$\mathsf{otherwise }$ }
\end{cases}
\]
 By axiom of choice and classical logic, for every $\beta$, $g_{<\beta}$ and $g_\beta$ are well defined. So we can let
  \[g:=g_{<\alpha}\]
 We first want to show that there exists a $\beta$ such that $\update_{g_{<\beta}}$ has a $\update_{g_{<\beta}}$-generated finite zero. We have two cases:\begin{enumerate}
 \item $\alpha$ is a successor ordinal. Then, by fact \ref{fact-ovvio}, $\update_{g<\alpha}$ is an update procedure of ordinal $1$, which has a $\update_{g<\alpha}$-generated finite zero (see chapter 4).\\
 
 \item $\alpha$ is a limit ordinal. Then, by continuity of $\update$, there is some $\beta<\alpha$ such that for all $\beta\leq \delta < \alpha$
 \[\update(g)=\update(g_{<\delta}*\zerof^{(\alpha-\delta)})\]
If $\update(g)=\emptyset$, then by definition \ref{defiinstance}
\[\update_{g_{<\beta}}(\zerof^{(\alpha-\beta)})=\emptyset\]
and we are done. If $\update(g)=\langle \gamma, n, m\rangle$, without loss of generality we can assume we have chosen $\beta$ such that $\gamma<\beta$. Again by definition \ref{defiinstance} of $\update_{g_{<\beta}}$
\[\update_{g_{<\beta}}(\zerof^{(\alpha-\beta)})=\emptyset\]
and we are done.
 \end{enumerate}
 Let now \[\beta_0:=\mathsf{min}\{\beta\ |\ \update_{g_{<\beta}}\text{ has a $\update_{g_{<\beta}}$-generated finite zero} \}
 \]
 $\beta_0$ cannot be a successor, otherwise if we let $\beta_0=\beta_1+1$\[g_{<\beta_0}=g_{<(\beta_1)}*g_{\beta_1}\]
 and hence by fact \ref{fact-ovvio} point (2)
 \[\update_{g_{<\beta_0}}=(\update_{g_{<\beta_1}})_{g_{\beta_1}}\]
 and by reduction lemma \ref{lemma-updateinstance} $\update_{g_{<\beta_1}}$ would have a $\update_{g_{<\beta_1}}$-generated finite zero. But $\beta_0$ also cannot be a limit ordinal, otherwise by reduction lemma \ref{lemma-transfiniteinstance}, for some $\gamma<\beta_0$, $\update_{g_{<\gamma}}$ would have a $\update_{g_{<\gamma}}$-generated finite zero. We conclude that $\beta_0=0$. Since $\update_{g_{<0}}=\update$, we obtain the thesis.
 \qed

\section{Spector's System $\SystemB$ and Typed Update Procedures of Ordinal $\omega^k$}\label{section-typedupdateprocedures}

Zeros of transfinite update procedures cannot in general be computed in G\"odel's system $\SystemTG$: as we will show, already update procedures of ordinal $\omega+k$, with $k\in\omega$, can be used to give computational interpretation to Elementary Analysis and hence their zeros can be used to compute the functions provably total in Elementary Analysis.  We will show however that Spector's system $\SystemB$ is enough to compute zeros.

\begin{defi}[Bar Recursion Operator, Spector's System $\SystemB$, Type Level of Bar Recursion]
In the following, we will work with Spector's system $\SystemB$ which is G\"odel's $\SystemTG$ augmented with constants $\BR_{\tau,\sigma}, \Psi_{\tau, \sigma}$ respectively of type \[ T_1\rightarrow T_2\rightarrow T_3 \rightarrow T_4 \rightarrow \tau\] 
and
\[ T_1\rightarrow T_2\rightarrow T_3 \rightarrow T_4 \rightarrow \Bool\rightarrow \tau\]with 
\[T_1=(\Nat \rightarrow \sigma)\rightarrow \Nat\] 
\[T_2=\sigma^*\rightarrow\tau\]
\[T_3=\sigma^*\rightarrow (\sigma\rightarrow \tau)\rightarrow \tau\]
\[T_4=\sigma^*\]
where $\sigma^*$ is a type representing finite sequences of objects of type $\sigma$. The meaning of $\BR_{\tau,\sigma}$ is defined by the equation

\begin{equation}\label{barrec}\mathsf{BR_{\tau,\sigma}}YGHs\overset{\text{$\tau$}}{=}
\begin{cases}
Gs &\text{if $ Y\hat{s}<|s|$}\\
Hs(\lambda x^\sigma \mathsf{BR}_{\tau,\sigma}YGH(s*x)) &\text{$\mathsf{otherwise}$}
\end{cases}
\end{equation}
where $s*x$ denotes the finite sequence $s$ followed by $x$, $\hat s$ denote the function mapping $n$ to $s_n$, if $n < |s|$, to $0^\sigma$ otherwise, where $s_n$ is the $n$-th element of $s$ and $|s|$ is the number of elements in $s$. If $\sigma,\tau, Y,G,H$ are determined by the context, we we will just write $\BR(s)$ in place of $\BR_{\tau,\sigma}YGHs$.\\
$\BR_{\tau,\sigma}$ is said to be {\em bar recursion of type} $\sigma$. The {\em type level} of bar recursion $\BR_{\tau, \Nat }$ of type $\Nat$ (said also type $0$), is the  type level  of the constant $\BR_{\tau, \Nat }$, that is, assuming $\Nat^*=\Nat$, $\mathsf{max}(1, \mathsf{typelevel}(\tau))+2$.\\ 
In order to obtain a strongly normalizing system such that equation \ref{barrec} holds, we have to add to system $\SystemB$ the following reduction rules (see Berger \cite{Berger}):
\[\begin{aligned}\BR_{\tau,\sigma}YGHs &\mapsto\Psi_{\tau,\sigma} YGHs(Y\hat{s}<|s|)\\
\Psi_{\tau,\sigma} YGHs(\True)&\mapsto Gs\\
 \Psi_{\tau,\sigma} YGHs(\False)&\mapsto Hs(\lambda x^\sigma \mathsf{BR}_{\tau,\sigma}YGH(s*x))\\
\end{aligned}\]
where $<$ is a term coding the correspondent relation on natural numbers.
\end{defi}

Since we are interested only in computable update procedures, we now fix a system for representing them. For the aim of computationally interpreting Elementary Analysis, update procedures can be assumed to belong to system $\SystemTG$. However, for more powerful systems one may need more capable update procedures, so we define them to belong to $\SystemB$. Here, we limit ourselves to the ordinal $\omega^k$, for $k\in \omega$, since this ordinal is enough to interpret Elementary Analysis and even fragments of Ramified Analysis (see for example, Mints et al. \cite{Mints2})

\begin{defi}[Representation of Ordinals and Typed Update Procedures of Ordinal $\omega^k$] We will represent ordinal numbers of the form $\omega^k$, with $k\in\omega$, by exploiting  the order isomorphism between $\omega^k$ and $\NatSet^k$ lexicographically ordered. So, for $k\in\omega$, $k>0$, we set \[[\omega^0]:=\nu, [\omega^k]:=\Nat^k\]  where $\nu$ is the empty string and

 \[[\omega^0 \rightarrow (\NatSet\rightarrow \NatSet)]:= \Nat\rightarrow \Nat\] and, if $k\in\omega$ \[[\omega^{k+1} \rightarrow (\NatSet\rightarrow \NatSet)]:= \Nat \rightarrow [\omega^k\rightarrow( \NatSet\rightarrow \NatSet)]\]
 where $\Nat$ is the type representing $\NatSet$ in typed lambda calculus. Define moreover
 \[[(\omega^k\times \NatSet\times \NatSet)\cup \{\emptyset\}]:=[\omega^k]\times \Nat\times\Nat\]
 Unfortunately, $\emptyset$ does not have a code. So we have to use an injective coding $|\_|$ of the set $(\omega^k\times \NatSet\times \NatSet)\cup \{\emptyset\}$ into the set of closed normal terms of type $[(\omega^k\times \NatSet\times \NatSet)\cup \{\emptyset\}]$. To fix ideas, we define $|(\beta,n,m)|=\langle\beta', n+1,m+1\rangle$, with $\beta':\Nat^k$ the code of $\beta$, and $|\emptyset|= \langle 0,\ldots, 0\rangle$.\\
A {\em typed update procedure of ordinal $\omega^k$} is a term of Spector's system $\SystemB$ of type:\[[\omega^{k} \rightarrow (\NatSet\rightarrow \NatSet)]\rightarrow  [(\omega^k\times \NatSet\times \NatSet)\cup \{\emptyset\}]\]
 satisfying point (2) of definition \ref{transfiniteupdateprocedures}, where for simplicity function quantification is assumed to range over functions definable in system $\SystemB$. Equality as it appears in the definition is supposed to be extensional.
\end{defi}

\section[Bar Recursion Proof of Zero Theorem for Typed Upd. Proc. of Ordinal $\omega^k$]{Bar Recursion Proof of the Zero Theorem for Typed Update Procedures of Ordinal $\omega^k$}\label{section-barrecursivezeros}
In this section we give a constructive proof of the Zero theorem for typed update procedures of ordinal less than $\omega^k$. In particular we show that finite zeros can be computed with bar recursion of type $1$. We start with the base case.

\begin{thm}[Zero Theorem for Update Procedures of Ordinal 1=$\omega^0$]\label{fixedpointtheorem} Let $\update$ be a typed update procedure of ordinal $1$. Then $\update$ has a finite zero $\sigma$. Moreover, $\sigma$ can be calculated as the normal form of a bar recursive term $\mathsf{Zero}(\update)$ (defined uniformly on the parameter $\update$) of system $\SystemTG$ plus bar recursion of type $0$.
\end{thm}

The result follows by Oliva \cite{Oliva}. We give below another proof, which is a  simplification of Oliva's one, made possible by the slightly stronger condition we have imposed on the notion of update procedure.\\  The informal idea of the construction - but with some missing justifications - is the following. We reason over the well-founded tree of finite sequences of numbers $s$ such that $\update(\hat{s})=|(n,m)|$ and $ n\geq |s|$.  We want to construct a function $\sigma: 
\Nat\rightarrow \Nat$ which is a zero of $\update$. Suppose that we have constructed a ``good" initial approximation $\sigma(0)*\cdots *\sigma(i)$ of $\sigma$; we want to prove that it can be extended to a long enough approximation of $\sigma$.  Our first step is to continue with $\sigma(0)*\cdots *\sigma(i)*0$. If this is a good guess, by well-founded induction hypothesis, we can extend  $\sigma(0)*\cdots *\sigma(i)*0$ to a complete approximation $\sigma(0)*\cdots*\sigma(n)$ of $\sigma$, with $n>i$. Since we are not sure that our previous guess was lucky, we compute $\update(\sigma(0)*\cdots*\sigma(n))$. If for all $m$ \[\update(\sigma(0)*\cdots*\sigma(n))\neq|(i+1, m)|\] then our approximation for $\sigma(i+1)$ is adequate, and we claim that $\sigma(0)*\cdots*\sigma(n)$ is the approximation of $\sigma$ we were seeking.  Otherwise \[\update(\sigma(0)*\cdots*\sigma(n))=|(i+1, m)|\] for some $m$: $\update$ tells us that our guess for the value of $\sigma(i+1)$ is wrong. But now we know that $\sigma(0)*\cdots *\sigma(i)*m$ is a good initial approximation of $\sigma$ and we have made  progress. Again by well-founded induction hypothesis, we conclude that we can extend $\sigma(0)*\cdots *\sigma(i)*m$ to a good approximation of $\sigma$.\\\\
\emph{Proof of Theorem} \ref{fixedpointtheorem}. We formalize and complete the previous informal argument. In the following $s$ will be a variable for finite sequences of numbers.  Using bar recursion of type $0$, we can define a term which builds directly the finite zero we are looking for and is such that:
 
\[\mathsf{BR}(s)=
\begin{cases}
\hat{s} &\text{if $\update \hat{s}=|(n, m )|$ and $n<|s|$}\\
\hat{s}  & \text{if }\update \hat{s}=|\emptyset|\\
\mathsf{BR}(s*m) &\text{if $\update ({\mathsf{BR}(s*0)})=|(|s|,m)|$ }\\
\mathsf{BR}(s*0) &\text{if $\update ({\mathsf{BR}(s*0)})\neq |(|s|,m)|$ for all $m$}

\end{cases}
\]
(we assume that $\mathsf{BR}(s)$ checks in order every condition in its definition and executes the action corresponding to the first satisfied condition). We let $\sigma$ as the normal form of \[\mathsf{Zero(\update)}:={\mathsf{BR}(\langle\rangle)}\] 
where $\langle \rangle$ is the empty sequence. Let us prove that $\sigma$ is a finite zero of $\update$. Suppose $\update\sigma=|(n,m)|$: by showing that this is impossible, we obtain that $\update\sigma=|\emptyset|$. The normalization of $\BR(\langle \rangle)$ leads to the following chain of equations:
\[\begin{aligned}\BR(\langle \rangle)&=\BR(\sigma(0))\\
&=\BR(\sigma(0)*\sigma(1))\\
&\ldots\\
&\ldots\\
&\ldots\\
&=\BR(\sigma(0)*\cdots*\sigma(i))\\
&=\widehat{\sigma(0)*\cdots*\sigma(i)}\\
&=\sigma
\end{aligned}\]
  with \[n<|\sigma(0)*\cdots * \sigma(i)|=i+1\] In particular 
\[\BR(\langle\rangle)=\BR(\sigma(0)*\cdots * \sigma(n-1))\]
Now, we have two cases:
\begin{enumerate}
\item
 $\update ({\BR(\sigma(0)*\cdots * \sigma(n-1)*0)})=|(n,l)|$. Then 

\[\BR(\langle\rangle)=\BR(\sigma(0)*\cdots * \sigma(n-1)*l)\]
and so $\sigma(n)=l$, which is impossible, by definition \ref{transfiniteupdateprocedures} of update procedure, point (2), for $\update\sigma=|(n,m)|$.\\
\item for all $l$,  $\update ({\BR(\sigma(0)*\cdots * \sigma(n-1)*0))}\neq |(n,l|)$. Then by definition  \[\BR(\sigma(0)*\cdots * \sigma(n-1))=\BR(\sigma(0)*\cdots * \sigma(n-1)*0)\] Therefore \[|(n,m)|=\update \sigma=\update({\mathsf{BR}(\langle\rangle)})=\update ({\BR(\sigma(0)*\cdots * \sigma(n-1)*0))}\]
again impossible, by assumption of this case.
\end{enumerate}
 We have then proved that $\sigma$ is the sought finite zero. \comment{Finally, as proven in Schwichtenberg (\cite{Schwichtenberg}), system $\SystemTG$ is closed under the rule of bar recursion of type $0$; hence, it follows that $\mathsf{Fix}(\update)$ can be taken as a term of system $\SystemTG$.}

\qed 

We now prove that every typed update procedure of ordinal $\omega$ has a finite zero.

\begin{thm}[Zero Theorem for Typed Update Procedures of Ordinal $\omega$]\label{fixedpointomega} Let $\update$ be a typed update procedure of ordinal $\omega$. Then $\update$ has a finite zero $\sigma$. Moreover, $\sigma$ can be calculated as the normal form of a bar recursive term $\mathsf{Zero}_\omega(\update)$  (defined uniformly on the parameter $\update$) of system $\SystemTG$ plus bar recursion of type $1:=\Nat\rightarrow \Nat$.
\end{thm}
\proof The finite function $\sigma: [\omega\rightarrow (\NatSet\rightarrow \NatSet)]$ we are going to construct can be represented as a finite function sequence $\sigma(0)*\sigma(1)*\cdots *\sigma(n)$, for a large enough $n$. In the following $s$ is a variable ranging over finite sequences of natural number functions. Using bar recursion of type $1$, we can define in a most simple way a term which builds directly the finite zero we are looking for. We present the construction gradually. To begin with, suppose we are able to define - uniformly on $s$ - terms $\BR(s)$ and $g_s: (\funn)$ satisfying the following equation for every $s$: 
\[\mathsf{BR}(s)=
\begin{cases}
\hat{s} &\text{if $\update \hat{s}=|(\gamma, n, m)|$ and $\gamma<|s|$}\\
\hat{s}  & \text{if }\update \hat{s}=|\emptyset|\\
\mathsf{BR}(s*g_s) &\text{$\mathsf{otherwise}$, where $\forall n,m\ $ $\update(\BR(s*g_s))\neq (|s|, n, m)$}\\

\end{cases}
\]
Let \[\sigma:=\mathsf{Zero_\omega(\update)}:={\mathsf{BR}(\langle\rangle)}\] We prove that $\sigma$ is a finite zero of $\update$. We show this by proving that $\update\sigma=(\gamma,n,m)$ is impossible. As in the proof of theorem \ref{fixedpointtheorem} \[\mathsf{BR}(\langle\rangle)=\BR(\sigma(0)*\cdots*\sigma(i))=\widehat{\sigma(0)*\cdots*\sigma(i)}\] with $\gamma<i+1$. Let \[r:={\sigma(0)*\cdots* \sigma (\gamma-1)}\] By some computation 
\[\begin{aligned}\update\sigma&=\update(\BR(\langle\rangle))\\
&=\update(\BR(\sigma(0)*\cdots * \sigma(\gamma-1)))\\
&=\update(\BR(r))\\
&=\update(\BR(r*g_r))
\end{aligned}\]
Since by construction for all $n,m$ \[\update(\BR(r*g_r))\neq |(|r|, n,m)|=|(\gamma, n, m)|\] we obtain that $\update\sigma\neq (\gamma,n,m)$: impossible.\\
It remains to show that a $g_s$ such that appears in the definition of $\BR(s)$ exists. Indeed, it is enough to set
\[g_s:=\mathsf{Zero}(\lambda f^{\Nat\rightarrow \Nat} \update_{|s|}({\BR(s*f)}))\]
where, for $i\in\NatSet$, we have defined   \[\update_i:= \lambda f^{\Nat\rightarrow(\Nat\rightarrow \Nat)}.\ \mathsf{if }\ \update(f)=|(i,n,m)|\ \mathsf{then}\ |(n,m)|\ \mathsf{else}\ |\emptyset|\] We prove now that in fact $ \update(\BR(s*g_s))\neq |(|s|, n ,m)|$ for all $n,m$. First, observe again that for every $s$ \[\BR(s)=\widehat{s*h_1*\cdots*h_n}\] for some terms $h_1,\ldots, h_n$ of type $\funn$.  Now, fix any finite sequence $s$ of type-$\funn$ terms. We want to show that \[F_{s}:= \lambda f^{\Nat\rightarrow \Nat} \update_{|s|}({\BR(s*f)})\] is an update procedure of ordinal $1$.   Suppose $F_sg_1=|(n,m)|$, $g_2(n)=m$ and $F_sg_2=|(h,l)|$. Then, by definition of $F_s$, it must be that \[\update ({\BR(s*g_1)})=|(|s|,n,m)|\] and \[\update ({\BR(s*g_2)})=|(|s|,h,l)|\]
Moreover, 
\[{\BR(s*g_2)}_{|s|}(n)=g_{2}(n)=m\]
Since $\update$ is an update procedure, $h\neq n$ must hold; therefore $F_s$ is an update procedure of ordinal $1$. But by definition of $g_s$, $\mathsf{Zero}$ and theorem \ref{fixedpointtheorem}, this means that 
\[|\emptyset|=F_s(\mathsf{Zero}(F_s))=\update_{|s|}(\BR(s*g_s))\]
By definition of $\update_{|s|}$ it must be true that $ \update(\BR(s*g_s))\neq |(|s|, n ,m)|$ for all $n,m$.

\qed

The previous argument can be generalized in order to prove the Zero theorem for typed update procedures of ordinal $\omega^k$.

\begin{thm}[Zero Theorem for Typed Update Procedures of Ordinal $\omega^k$, with $k\in\omega$] Let $\update$ be a typed update procedure of ordinal $\omega^k$. Then $\update$ has a finite zero $\sigma$. Moreover, $\sigma$ can be calculated as the normal form of a bar recursive term $\mathsf{Zero}_{\omega^k}(\update)$  (defined uniformly on the parameter $\update$) of system $\SystemTG$ plus bar recursion of some type $A$, where $\mathsf{typelevel}(A)=1$.
\end{thm}
\proof By induction on $k$. The cases $k=0,1$ have already been taken care. Now, we want to prove the thesis for $k+1$, with $k>0$. The finite function $\sigma: [ \omega^{k+1}\rightarrow (\NatSet\rightarrow \NatSet)]$ we are going to construct can be represented as a finite function sequence $\sigma(0)*\sigma(1)*\cdots *\sigma(n)$, for a large enough $n$, with each $\sigma(i)$ of type $[\omega^{k}\rightarrow (\NatSet\rightarrow \NatSet)]$.  In the following $s$ represents a sequence of functions of type $[\omega^{k}\rightarrow (\NatSet\rightarrow \NatSet)]$. We recall that an ordinal less than $\omega^{k+1}$ is coded as a pair $(\gamma, \beta)$, with $\gamma\in \NatSet$ and $\beta \in \NatSet^k$. 
We again present the construction gradually. To begin with, suppose we are able to define - uniformly on $s$ - terms $\BR(s)$ and $g_s: [\omega^{k}\rightarrow (\NatSet\rightarrow \NatSet)]$ satisfying the following equation for every $s$: 
\[\mathsf{BR}(s)=
\begin{cases}
\hat{s} &\text{if $\update \hat{s}=|((\gamma, \beta), n, m)|$ and $\gamma<|s|$}\\
\hat{s}  & \text{if }\update \hat{s}=\emptyset\\
\mathsf{BR}(s*g_s) &\text{$\mathsf{otherwise}$, where $\forall \beta, n,m\ $ $\update(\BR(s*g_s))\neq ((|s|, \beta), n, m)$}\\

\end{cases}
\]
Let \[\sigma:=\mathsf{Zero_{\omega^{k+1}}(\update)}:={\mathsf{BR}(\langle\rangle)}\] We prove that $\sigma$ is a finite zero of $\update$. We show this by proving that $\update\sigma=((\gamma, \beta),n,m)$ is impossible. As in the proof of theorem \ref{fixedpointtheorem} \[\mathsf{BR}(\langle\rangle)=\BR(\sigma(0)*\cdots*\sigma(i))=\widehat{\sigma(0)*\cdots*\sigma(i)}\] with $\gamma<i+1$. Let \[r:={\sigma(0)*\cdots* \sigma (\gamma-1)}\] By some computation 
\[\begin{aligned}\update\sigma&=\update(\BR(\langle\rangle))\\
&=\update(\BR(\sigma(0)*\cdots * \sigma(\gamma-1)))\\
&=\update(\BR(r))\\
&=\update(\BR(r*g_r))
\end{aligned}\]
Since by construction 
\[\update(\BR(r*g_r))\neq |((|r|,\beta), n,m)|=|((\gamma, \beta), n, m)|\] we obtain that $\update\sigma\neq ((\gamma, \beta),n,m)$: impossible.\\
It remains to show that a $g_s$ such that appears in the definition of $\BR(s)$ exists. Indeed, it is enough to set
\[g_s:=\mathsf{Zero}_{\omega^k}(\lambda f^{[\omega^k\rightarrow (\NatSet\rightarrow \NatSet)]} \update_{|s|}({\BR(s*f)}))\]
where, for $i\in\NatSet$, we have defined   \[\update_i:= \lambda f^{[\omega^{k+1}\rightarrow (\NatSet\rightarrow \NatSet)]}.\ \mathsf{if }\ \update(f)=((i,\delta),n,m)\ \mathsf{then}\ (\delta,n,m)\ \mathsf{else}\ \emptyset\] We prove now that in fact $ \update(\BR(s*g_s))\neq |((|s|, \beta), n ,m)|$ for all $\beta, n,m$. First, observe that for every $s$ \[\BR(s)=\widehat{s*h_1*\cdots*h_n}\] for some terms $h_1,\ldots, h_n$ of type $\funn$.  Now, fix any finite sequence $s$ of type-$\funn$ terms. First, we want to show that \[F_{s}:= \lambda f^{[\omega^k\rightarrow (\NatSet\rightarrow \NatSet)]} \update_{|s|}({\BR(s*f)})\] is an update procedure of ordinal $\omega^k$.  Suppose for some $\delta$ of type $[\omega^k]$: $F_sg_1=|(\delta,n,m)|$, $\forall \delta_0<\delta.\ (g_1)_{\delta_0}=(g_2)_{\delta_0}$, $(g_2)_\delta(n)=m$ and $F_sg_2=|(\delta,h,l)|$. Then, by definition of $F_s$, it must be that \[\update ({\BR(s*g_1)})=|((|s|,\delta),n,m)|\] and \[\update ({\BR(s*g_2)})=|((|s|,\delta),h,l)|\]
Since $\update$ is an update procedure and \[{\BR(s*g_2)}|s|\delta( n)=(g_2)_\delta(n)=m\] then $h\neq n$ must hold; therefore $F_s$ is an update procedure of ordinal $\omega^k$. But by definition of $g_s$ and induction hypothesis, this means that 
\[\update_{|s|}(\BR(s*g_s))=F_s(\mathsf{Zero}_{\omega^k}(F_s))=|\emptyset|\]
By definition of $\update_{|s|}$ it must be true that $ \update(\BR(s*g_s))\neq |((|s|, \beta), n ,m)|$ for all $\beta, n,m$.

\qed

\section{Case Study: Elementary Analysis}\label{section-casestudyea}

In this section, we give a three-step description of the epsilon substitution method for Elementary Analysis. Every step corresponds to one of the three stages in which - according to section \ref{section-transfiniteupdateprocedures} - learning based computational interpretations of predicative classical Arithmetic can be decomposed. As main foundational result one obtains a constructive proof that the zero theorem for update procedures of ordinal less that $\omega\cdot 2$ implies the 1-consistency of Elementary Analysis. More precisely, any zero of such an update procedure can be used to compute witnesses for $\Pi_0^2$ formulas.

The content of this section is based on Mints et al. \cite{Mints} and may be considered as an informal survey and a general guide to the reading of the epsilon substitution method in the light of our ideas on learning. Neither full details nor full proofs will be provided, but our description should be clear enough for the reader to gain an understanding of the basic ideas underpinning the epsilon method and its learning based interpretation.\\

 We first define the language of Elementary Analysis, which is  a fragment of second order Arithmetic in which second order quantification ranges over arithmetical formulas (possibly with free set variables).

\begin{definition}[Language $\mathcal{L}_{\EA}$of $\EA$]
The \emph{terms} of $\mathcal{L}_\EA$ are inductively defined as follow:
\begin{enumerate}
\item Numerical variables $x, y, z,\ldots$ are terms of type $0$.\\

\item  Set variables $X, Y, Z, \ldots$ are terms of type $1$.\\

\item $\zerof$ is a term of type $0$.\\

\item If $t$ is a term of type $0$, $\suc(t)$ is term of type $0$.

\end{enumerate}
The \emph{formulas} of $\mathcal{L}_\EA$ are inductively defined as follows:
\begin{enumerate}
\item For every natural number $n$, there is a denumerable set of $n$-ary predicate constants, one for every computable predicate over $n$-uples of natural numbers. If $P$ is a $n$-ary predicate constant and $t_1, \ldots, t_n$ are terms of type $0$, then $Pt_1\ldots t_n$ is an atomic formula.\\
\item If $t$ is a term of type $0$ and $X$ a variable of type $1$, then $t\in X$ is an atomic formula.\\
\item If $A$ and $B$ are formulas, then $A\land B$, $A\rightarrow  B$, $\lnot A$ are formulas.\\
\item If $A$ is a formula and $v$ is a variable, $\exists v A$ is a formula and $\forall v A$ is defined as $\lnot \exists v\lnot A$.\\

\end{enumerate}
If $A$ is a formula and $z$ a variable of type $0$ free in $A$, then $\lambda z A$ is a \emph{lambda set}; $\lambda z A$ is said to be \emph{arithmetical} if it contains no bound set variables. The formula $B(\lambda z A/X)$ is defined as the formula obtained from $B$ by substituting each atomic formula $t\in X$ of $B$ with $A(t/z)$, as usual without capture of variables.

\end{definition}

We now define the axioms and inference rules of $\EA$.

\begin{definition} [Axioms and Inference Rules of $\EA$]The \emph{axioms} of $\EA$ are formulas of $\mathcal{L}_{\EA}$ defined as follows:
\begin{enumerate}
\item Propositional tautologies are axioms.\\

\item Definitions of predicate constants are axioms, e.g. for $add$ predicate
\[add(x, 0, x) \text{ and } add(x, y, z)\rightarrow add(x, \suc(y), \suc(z))\]
\item $x=x$ and $x=y\rightarrow A(x)\rightarrow A(y)$ are equality axioms.\\
\item $\lnot \suc(x)=0$ and $\suc(x)=\suc(y)\rightarrow x=y$ are axioms.\\
\item $A(0)\rightarrow  (\forall x. A(x)\rightarrow A(\suc(x)))\rightarrow \forall x A(x)$ is the induction axiom scheme.\\
\item $A(t/x)\rightarrow \exists x A$ is an axiom for every term $t$ of type $0$. \\
\item $A(T/X)\rightarrow \exists X A$ is an axiom if $T$ is a set variable or an arithmetical lambda set. 

\end{enumerate}
The \emph{inference rules} of $\EA$ are modus ponens 
\[\begin{array}{c}   A\rightarrow B\ \ \  A \\ \hline
 B
\end{array}\ \ \ \ \]
and 
\[\begin{array}{c}      A \rightarrow C\\ \hline
 \exists v A \rightarrow C
\end{array}\ \ \ \ \]
with the standard proviso that $v$ does not occur free in $C$.
\end{definition}
We are now ready to take the first step of a learning based interpretation.
\subsection{First Stage: Identification of a Sequence of non Computable Functions $F$}
We now define the sequence of non computable functions needed to give a computational interpretation of $\EA$. We do that by first introducing the concept of epsilon term. 

\begin{definition}[Language $\mathcal{L}_{\EA\epsilon}$] We define by simultaneous induction the \emph{terms} and the \emph{formulas}  of $\mathcal{L}_{\EA\epsilon}$:
\begin{enumerate}
\item Numerical variables $x, y, z,\ldots$ are terms of type $0$.\\

\item  Set variables $X, Y, Z, \ldots$ are terms of type $1$.\\

\item $\zerof$ is a term of type $0$.\\

\item If $t$ is a term of type $0$, $\suc(t)$ is term of type $0$.\\

\item  For every $n$-ary predicate constant $P$ of $\mathcal{L}_{\EA}$, if $t_1, \ldots, t_n$ are terms of type $0$, then $Pt_1\ldots t_n$ is an atomic formula.\\
\item If $t$ is a term of type $0$ and $T$ a term of type $1$, then $t\in T$ is an atomic formula.\\
\item If $A$ and $B$ are formulas, then $A\land B$, $A\rightarrow  B$, $\lnot A$ are formulas.\\
\item If $A$ is a formula  and $v$ is a variable, then $\epsilon v A$ is an \emph{epsilon term} (of type equal to the type of $v$) and $v$ is considered bound in $\epsilon v A$. \\
\end{enumerate}
If $A$ is a formula and $z$ a variable of type $0$ free in $A$, then $\lambda z A$ is a \emph{lambda term} (but we do not ask it is in $\mathcal{L}_{\EA\epsilon}$). A formula or a lambda term or a term is said to be an \emph{expression of} $\EA\epsilon$ and is \emph{arithmetical} if it contains no bound set variables and \emph{canonical} if it is closed (i.e. no free variables occurs in it) and does not have closed epsilon terms as subterms. $\lambda z A$ is \emph{regular} if it is of the form $(\lambda z B)[t_1/v_1\ldots t_n/v_n]$ with $\lambda z B$ arithmetical and $t_1, \ldots, t_n$ any terms.
\end{definition}

Canonical epsilon terms are the ones that are assigned a meaning. The intended denotation of a canonical epsilon term $\epsilon x A$ is the least number $n$ such that $A(n)$ is true while the denotation of $\epsilon X A$ is an arithmetical canonical lambda term $\lambda z G$ (which represents an arithmetical set) such that $A(\lambda z G/X)$ is true. In other words, the following \emph{critical} formulas should be true:
\[A(t/x)\rightarrow A(\epsilon x A/x)\]
\[A(T/X)\rightarrow A(\epsilon X A/X)\]
where $t$ is any term of type $0$ and $T$ is an epsilon term  of type $1$ or a regular lambda term. The notion of truth for formulas of $\EA\epsilon$ requires further explanation: in order to evaluate their truth, first epsilon terms must be evaluated and hence eliminated.

\begin{definition}[Substitutions, Evaluations of Epsilon Terms] We define:
\begin{enumerate}
\item An \emph {epsilon substitution} $S$ is a function from the set of canonical epsilon terms to the set of numerals and arithmetical canonical lambda terms such that $S(\epsilon x A)$ is always a numeral and  $S(\epsilon X A)$ is always an arithmetical canonical lambda term.\\
\item Let $t_1, t_2$ be  expressions of $\EA\epsilon$ and $S$ an epsilon substitution. We write $t_1\mapsto_S t_2$ if $t_2$ is obtained from $t_1$ either by substituting one of its canonical epsilon subterms $\epsilon x A$ with $S(\epsilon x A)$ or replacing one of its subformulas $t\in \epsilon X A$ with $G(t)$, where $S(\epsilon X A)=\lambda z G$.\\
\item An expression $t$ of $\EA\epsilon$ is said to be in $S$-normal form if there is no $t_1$ such that $t\mapsto_S t_1$. We indicate with $|t|_S$ the unique $S$-normal form of $t$, which exists by theorem \ref{theorem-normalization} below.

\end{enumerate}
\end{definition}
The truth value of a closed formula $A$ of $\EA\epsilon$ is the truth value of $|A|_S$ (which does not contain epsilon terms): so it is always relative to a epsilon substitution $S$. 

The relation $\mapsto_S$ is well founded and has Church-Rosser property (see Mints et al. \cite{Mints}).
\begin{theorem}[Normalization and Church-Rosser]\label{theorem-normalization}
For every $S$, the relation $\mapsto_S$ is well founded and every expression $t$ of $\EA\epsilon$ has an $S$-normal form.
\end{theorem}

We now need to measure the ``computational strength" of an expression $t$ of $\EA\epsilon$. Intuitively, from the computational point of view, an epsilon term $\epsilon x A$ represent a recursion theoretic jump, because in general one has to enumerate all natural numbers in order to decide if an $n$ exists such that $|A(n)|_S$ is true. So, closed arithmetical expressions will have a computational strength below $\omega$, because no more than a finite number of jumps is done inside them by their epsilon subterms. Instead, an epsilon term $\epsilon X A$ must have a computational strength of at least $\omega$ because one has to know all the values of arithmetical canonical epsilon terms in the first place if he wants to determine whether  there exists a canonical arithmetical lambda term $\lambda z G$ such that $|A(\lambda z G/X)|_S$ is true. Indeed one can assign to expressions a computational strength which is always less that $\omega\cdot 2$. This is done through the so called \emph{rank function}, which we introduce only by exposing the properties that it must have (for the actual definition and details, see Mints et al. \cite{Mints}).

\begin{theorem}[Rank Function]\label{theorem-rank}
There exists a function $\rank$ from the set of expressions of $\EA\epsilon$ to $\omega\cdot 2$ such that the following holds. For every epsilon substitution $S$ and ordinal $\alpha$, denote with $S_{\leq \alpha}$ the function mapping $e$ to $S(e)$ if $\rank(e)\leq \alpha$, to $0$ or $0^1:=(\lambda z. z=z)$ otherwise (according to the type of $e$); then
\begin{enumerate}
\item For every canonical epsilon term $\epsilon v A$  \[\rank(\epsilon v A)>\rank (A(e/v))\]
whenever $e$ is a numeral or an arithmetical canonical lambda term.\\

\item For every expression $e$ of $\EA\epsilon$ and epsilon substitutions $S_1, S_2$, if $\rank(e)=\alpha$ and $(S_1)_{\leq \alpha}=(S_2)_{\leq \alpha}$, then $|e|_{S_1}=|e|_{S_2}$. 
\end{enumerate}
\end{theorem}

The above theorem \ref{theorem-rank} is crucial and its meaning is the following. Given any canonical epsilon term $\epsilon x A$, any substitution $S$ and natural number $n=S(\epsilon x A)$, in order to check whether $n$ is a correct denotation for $\epsilon x A$, we have to determine the truth value of 
\[|A(n/x)|_S\]
Since $\rank(\epsilon x A)$ is strictly greater than $\rank(A(n/x))$, the truth of $|A(n/x)|_S$ depends only on the values that $S$ assigns to epsilon terms of rank strictly \emph{less} than that of $\epsilon x A$. So the meaning of $\epsilon x A$ is predicatively determined by the meaning of epsilon terms of lower rank. The same holds for canonical epsilon terms of type $1$. 

We are now able to define the sequence of functions that will enable us to define classical witness for formulas in $\EA$ and that we shall try to approximate. Suppose that for every ordinal $\alpha<2\cdot \omega$ we have a primitive recursive enumeration $\epsilon_0^\alpha, \epsilon_1^\alpha, \ldots$ of canonical epsilon terms of rank equal to  $\alpha$ and a primitive recursive enumeration  $\lambda_0, \lambda_1, \ldots, $ of canonical arithmetical lambda terms. We associate to  any function $f: \omega\cdot 2\rightarrow (\NatSet\rightarrow \NatSet)$ the  epsilon substitution $S_f$ such that:  
\[f_\alpha (n)=
\begin{cases}
m & \text{if }S_f(\epsilon_n^\alpha)=m\land \epsilon_n^\alpha= \epsilon x A\\ 
l & \text{if } S_f(\epsilon_n^\alpha)=\lambda_l \land \epsilon_n^\alpha= \epsilon X A
\end{cases}
\]
It is easy to see - using classical logic - that there exists an epsilon substitution $\mathsf{S}$ which makes true every critical formula $C$, i.e. $|C|_S$ is true: just start by assigning values to canonical epsilon terms of rank $1$, then to those of rank $2$ and so on. Then our target collection $F:\omega\cdot 2\rightarrow (\fun)$ of functions  is the one such that $\mathsf{S}=S_F$.

\subsection{Second Stage: Definition of Classical Witnesses by Programs Recursive in $F$}
Using epsilon terms one can define classical witness for any provable formula of $\EA$ by first translating formulas of $\EA$ into formulas of $\EA\epsilon$, using the equivalence 
\[\exists v A\equiv A(\epsilon v A/v)\]

\begin{definition}[Translation of Formulas of $\EA$ into Formulas $\EA\epsilon$]We define a translation of Formulas of $\EA$ into Formulas $\EA\epsilon$ by induction as follows:
\begin{enumerate}
\item If $P$ is atomic, $P^*:=P$.\\
\item $(\lnot A)^*:=\lnot A^*$.\\
\item $( A\land B)^*:= A^*\land B^*$.\\
\item $( A\rightarrow B)^*:= A^*\rightarrow B^*$.\\
\item $(\exists v A)^*:=A^*(\epsilon v A^*/v)$
\end{enumerate}
\end{definition}

We now define the axioms and inference rules for $\EA\epsilon$.

\begin{definition} [Axioms and Inference Rules of $\EA\epsilon$]The \emph{axioms} of $\EA\epsilon$ are formulas of $\mathcal{L}_{\EA\epsilon}$ defined as follows:
\begin{enumerate}
\item Propositional tautologies are axioms.\\

\item Definitions of predicate constants are axioms, e.g. for $add$ predicate
\[add(x, 0, x) \text{ and } add(x, y, z)\rightarrow add(x, \suc(y), \suc(z))\]
\item $x=x$ and $x=y\rightarrow A(x)\rightarrow A(y)$ are equality axioms.\\
\item $\lnot \suc(x)=0$ and $\suc(x)=\suc(y)\rightarrow x=y$ are axioms.\\
\item Minimality axioms: $\epsilon x A=\suc(t)\rightarrow \lnot A(t)$ \\
\item Critical formulas:
\[\lnot s=0\rightarrow s=\suc(\epsilon x s=\suc(x))\]
\[A(t/x)\rightarrow A(\epsilon x A/x)\]
\[A(T/X)\rightarrow A(\epsilon X A/X)\]
where $t$ is any term of type $0$ and $T$ is either a term of type $1$ or a regular lambda term.
\end{enumerate}
 The only \emph{inference rule} of $\EA\epsilon$ is modus ponens.

\end{definition}
The following theorem shows that $\EA$ can be embedded in the quantifier free system $\EA\epsilon$. This allows one to extract witnesses for existential statements provable in $\EA$. As usual $\vdash$ denotes provability.

\begin{theorem}[Classical Witnesses for Provable Formulas of $\EA$]\label{theorem-witnessesea} The following holds:
\begin{enumerate}
\item Suppose that  $\EA\vdash A$. Then $\EA\epsilon \vdash A^*$. \\

\item Suppose that  $\EA\vdash \exists x A$, with $A$ atomic. Then there exists a finite sequence $C_1, \ldots, C_n$ of closed critical formulas of $\EA\epsilon$ and a closed term $t$ of $\EA\epsilon$ such that if for all $i$, $|C_i|_S$ is true, then $|t|_S=n$ and  $A(n)$ is true.

\end{enumerate}
\end{theorem}

By the above theorem \ref{theorem-witnessesea}, it is now clear that witnesses for $\EA$ can be computed by programs recursive in $F$, because $F$ represent an epsilon substitution $S_F$ which makes every critical formula true. 

\subsection{Stage Three: Learning Processes Approximating $F$}

In order to compute witnesses for $\EA$ is necessary to find a good finite approximation of $F$, in order to satisfy some finite set of critical formulas. This is the point when update procedures come into the scene. The fundamental property of a closed critical formula $C$, for example of the form
\[A(t/x)\rightarrow A(\epsilon x A/x)\]
is that from the fact $|C|_S$ is false one can always learn something. Suppose $|C|_S$ is false. In this case, if $S(\epsilon x |A|_S)=m$, one has that $|A(m/x)|_S$ is false. Fortunately, if $|t|_S=n$, the formula $|A(n/x)|_S$ is true since $|C|_S$ is false and so $|A(n_0/x)|_S$ is true for some minimal $n_0\leq n$. So one learns a new value $n_0$ that can be assigned to $\epsilon x |A|_S$. Analogously, if $C$ is a closed critical formula of the form
\[A(T/X)\rightarrow A(\epsilon X A/X)\]
and $|C|_S$ is false, if we suppose $S(\epsilon X |A|_S)=\lambda z G$, one has that $|A(\lambda z G/X)|_S$ is false. Fortunately, if $|T|_S=\lambda z H$, the formula $|A(\lambda z H/X)|_S$ happens to be true. So one learns a new value $\lambda z H$ that can be assigned to $\epsilon X |A|_S$. Observe that is not a priori obvious that $\lambda z H$ is an arithmetical lambda term, but in fact  this can be proved given the assumptions we have made on $T$ (again, for details see Mints et al. \cite{Mints}).

From now on, fix a finite sequence of critical formulas $C_0, \ldots, C_N$ (for brevity only of the two forms considered above). We want to define an update procedure of ordinal $\omega+k$ out of it (with $k\in\omega$), any of whose finite zeros will represent an epsilon substitution making all the critical formulas true. 

\begin{definition}[Update Procedure for $C_0, \ldots, C_N$]\label{definition-updatecritical}
We define an update procedure $\update$ of ordinal $\omega+k$. Let $f: \omega+k\rightarrow (\fun)$. If for every $i$, $|C_i|_{S_f}$ is true, we set 
\[\update(f)=\emptyset\]
Otherwise, consider the first $i$ such that $|C_i|_{S_f}$ is false. If 
\[C_i=A(t/x)\rightarrow A(\epsilon x A/x)\]
and $|t|_{S_f}=m$ and $\epsilon x |A|_{S_f}=\epsilon_n^\alpha$, we set 
\[\update(f)=\langle \alpha, n, m_0\rangle\]
where $m_0\leq m$ is the smallest among the $i$ such that $|A(i/x)|_{S_f}$ is true (which exists since $|A(m/x)|_{S_f}$ is true).
If 
\[C_i=A(T/X)\rightarrow A(\epsilon X A/X)\]
and $|T|_{S_f}=\lambda z H=\lambda_m$ and $\epsilon X |A|_{S_f}=\epsilon_n^\alpha$, we set
\[\update(f)=\langle \alpha, n, m\rangle\]
\end{definition}
We now show that $\update$ is a well defined update procedure.

\begin{theorem}[Adequacy of $\update$]\label{theorem-adeqaucyofu}$\update$ is an update procedure of ordinal $\omega+k$, for some $k\in \omega$.
\end{theorem}
\proof We skip the proof that $\update$ is well defined, which amounts to show that for every substitution $S$, if $\rank(e)\geq \omega$, then $\rank(|e|_S)\leq \rank(e)$: this ensures an upper bound on the rank of the terms that are evaluated by $\update$ in its computations and so $\update$ never updates values for epsilon terms of rank greater than $\omega+k$, if $k$ is chosen large enough. \\
We prove instead that $\update$ is an update procedure. $\update$ is continuous, since for every $f$ only a finite number of values of $S_f$ and hence of $f$ are used to compute $\update(f)$. Suppose now that
$\update(f)=\langle \beta, n, m\rangle$, for all $\gamma<\beta$ $f_\gamma=g_\gamma$, $g_\beta(n)=m$ and $\update(g)=\langle \beta, h, l\rangle$. Suppose $h= n$: we have to prove it is impossible. 
Consider the first $i$ such that $|C_i|_{S_f}$ is false. Suppose $C_i$ is of the form
\[C_i=A(t/x)\rightarrow A(\epsilon x A/x)\]
Then by definition of $\update(f)$, $\epsilon x |A|_{S_f}=\epsilon_n^\beta$ and $m$ is the smallest among the $i$ such that $|A(i/x)|_{S_f}$ is true. Furthermore, consider the first $j$ such that $|C_j|_{S_g}$ is false. By definition of $\update(g)$ and since $h=n$ we have
\[C_j=B(t/v)\rightarrow B(\epsilon v B/v)\]
with  $\epsilon v |B|_{S_g}=\epsilon_n^\beta$. Therefore, $|A|_{S_f}=|B|_{S_g}$. Moreover, let \[\delta:=\rank(|B|_{S_g}(m/x))\] By theorem \ref{theorem-rank}, point (1), $\delta<\beta=\rank(\epsilon x |B|_{S_g})$. By hypothesis \[(S_g)_{\leq \delta}=(S_f)_{\leq \delta}\] So by theorem \ref{theorem-rank}, point (2)\[|A(m/x)|_{S_f}=||A|_{S_f}(m/x)|_{S_f} =||B|_{S_g}(m/x)|_{S_f}=||B|_{S_g}(m/x)|_{S_g}=|B(m/x)|_{S_g}\]
Since $g_\beta(n)=m$, we have
\[S_g(\epsilon v |B|_{S_g})=S_g(\epsilon_n^\beta)=m\]
Thus
\[|B(\epsilon v B/v)|_{S_g}=|B(m/x)|_{S_g}\]
But $|A(m/x)|_{S_f}$ is true by construction and so $ |B(\epsilon v B/v)|_{S_g}$ itself must be true, which contradicts the assumption that $|C_j|_{S_g}$ is false.\\
An analogous reasoning yields a contradiction when $C_i$ is of the form
\[A(T/x)\rightarrow A(\epsilon X A/X)\]

\qed

\begin{thm}[1-Consistency of Elementary Analysis]
If for all $k\in\omega$ every update procedure of ordinal $\omega +k $ has a finite zero, then Elementary Analysis is 1-consistent.
\end{thm}

\proof By theorem \ref{theorem-witnessesea}, it is enough to show that, given any finite sequence  of critical formulas ($C_1, \ldots, C_N$ without loss of generality), there exists a finite epsilon substitution that makes true every formula. This amounts to show that  $\update$ has a finite zero, which is true by theorem \ref{theorem-adeqaucyofu} and hypothesis.
\qed

\section{Further Work}

Much remains to be done and we plan to address the following issues in the future. 

Our constructive proof of the zero theorem for typed update procedures of ordinal $\omega^k$ is not optimal, in the sense that, by Howard \cite{Howard}, it should be possible to use only system $\SystemTG$ plus bar recursion of type $0$.

Moreover, a more self contained proof that the zero theorem for update procedures of ordinal less that $\omega\cdot 2$ implies the consistency of $\EA$ is a major aim.

%\appendix
%    Include appendix "chapters" here.
%\include{}

\backmatter
%    Bibliography styles amsplain or harvard are also acceptable.

%    See note above about multiple indexes.
\printindex

\end{document}

%-----------------------------------------------------------------------
% End of amsbook.template
%-----------------------------------------------------------------------